\documentclass[11pt,reqno]{amsart}

\usepackage{caption,amsmath,amssymb,amsfonts,amsthm,latexsym,graphicx,multirow,color,enumerate}
\usepackage[all]{xy}
\usepackage{hyperref}

\numberwithin{table}{section}

\newcommand{\A}{\mathrm{A}} \newcommand{\AGL}{\mathrm{AGL}} \newcommand{\AGaL}{\mathrm{A\Gamma L}}  \newcommand{\ASL}{\mathrm{ASL}} \newcommand{\Aut}{\mathrm{Aut}}
 \newcommand{\bbF}{\mathbb{F}} \newcommand{\bfO}{\mathbf{O}}
\newcommand{\C}{\mathrm{C}} \newcommand{\Cen}{\mathbf{C}}  \newcommand{\calC}{\mathcal{C}} \newcommand{\calN}{\mathcal{N}} \newcommand{\calP}{\mathcal{P}} \newcommand{\calO}{\mathcal{O}}   \newcommand{\Co}{\mathrm{Co}} 
\newcommand{\D}{\mathrm{D}} 

\newcommand{\F}{\mathrm{F}} 
\newcommand{\G}{\mathrm{G}} \newcommand{\GaG}{\mathrm{\Gamma G}} \newcommand{\GaL}{\mathrm{\Gamma L}} \newcommand{\GaSp}{\mathrm{\Gamma Sp}} \newcommand{\GaO}{\mathrm{\Gamma O}} \newcommand{\GaU}{\mathrm{\Gamma U}}  \newcommand{\GL}{\mathrm{GL}} \newcommand{\GO}{\mathrm{O}}  \newcommand{\GU}{\mathrm{GU}}
  
\newcommand{\J}{\mathrm{J}}

\newcommand{\M}{\mathrm{M}} \newcommand{\magma}{\textsc{Magma}}
\newcommand{\N}{\mathrm{N}}
\newcommand{\Nor}{\mathbf{N}}
\newcommand{\Out}{\mathrm{Out}}
\newcommand{\Pa}{\mathrm{P}}  \newcommand{\PGL}{\mathrm{PGL}} \newcommand{\PGaL}{\mathrm{P\Gamma L}}    \newcommand{\ppd}{\mathrm{ppd}} \newcommand{\POm}{\mathrm{P\Omega}} \newcommand{\PSL}{\mathrm{PSL}}    \newcommand{\PSp}{\mathrm{PSp}} \newcommand{\PSU}{\mathrm{PSU}}
\newcommand{\Q}{\mathrm{Q}}
 \newcommand{\Rad}{\mathbf{R}}  
\newcommand{\Setminus}{-} \newcommand{\SiL}{\mathrm{\Sigma L}}  \newcommand{\SL}{\mathrm{SL}}  \newcommand{\Soc}{\mathrm{Soc}} \newcommand{\Sp}{\mathrm{Sp}}  \newcommand{\SU}{\mathrm{SU}} \newcommand{\Suz}{\mathrm{Suz}} \newcommand{\Sy}{\mathrm{S}}  \newcommand{\Sz}{\mathrm{Sz}}
\newcommand{\Tr}{\mathrm{Tr}}
\newcommand{\Z}{\mathbf{Z}} 

\newtheorem{theorem}{Theorem}[section]
\newtheorem{lemma}[theorem]{Lemma}
\newtheorem{proposition}[theorem]{Proposition}
\newtheorem{corollary}[theorem]{Corollary}
\newtheorem{problem}[theorem]{Problem}

\theoremstyle{definition}
\newtheorem{definition}[theorem]{Definition}

\newtheorem*{remark}{Remark}
\newtheorem{hypothesis}[theorem]{Hypothesis}

\newtheorem{MyExample}[theorem]{Example}
\newenvironment{example}
  {\pushQED{\qed}\begin{MyExample}}
  {\popQED\end{MyExample}}

\begin{document}

\title[Factorizations]{Factorizations of Almost Simple Groups with Applications}

\author{Robert M.~Guralnick, Cai Heng Li, Lei Wang, Binzhou Xia}
\address[Guralnick]
{Department of Mathematics\\University of Southern California\\Los Angeles, CA 90089-2532\\USA
\newline Email: {\tt guralnic@usc.edu}
}
\address[Li]
{SUSTech International Center for Mathematics, and Department of Mathematics, Southern University of Science and Technology\\Shenzhen 518055, Guangdong\\P. R. China
\newline Email: {\tt lich@sustech.edu.cn}
}
\address[Wang]
{School of Mathematics and Statistics\\Yunnan University\\Kunming 650091, Yunnan\\P. R. China
\newline Email: {\tt wanglei@ynu.edu.cn}
}
\address[Xia]
{School of Mathematics and Statistics\\The University of Melbourne\\Parkville, VIC 3010\\Australia
\newline Email: {\tt binzhoux@unimelb.edu.au}
}

\begin{abstract}
The classification of factorizations $G=HK$ of finite almost simple groups, proposed by Wielandt in 1979 and pursued through several partial classifications, has remained open in one major case. We settle that case: for every finite almost simple classical group $G$, we determine all factorizations $G=HK$ in which both $H$ and $K$ have a unique nonsolvable composition factor, completing the classification of factorizations of finite almost simple groups. Among other consequences of this classification, we prove that the smallest dimension of a biperfect bicrossproduct Hopf algebra is $41287680$.


\end{abstract}

\maketitle

%%%%%%%%%%%%%%%%%%%%%%%%%%%%%%%%%%%%%%%%%%%%%%%%%%%%%%%%%%%%%%%%%%%
%%%%%%%%%%%%%%%%%%%%%%%%%%%%%%%%%%%%%%%%%%%%%%%%%%%%%%%%%%%%%%%%%%%

\section{Introduction}

%%%%%%%%%%%%%%%%%%%%%%%%%%%%%%%%%%%%%%%%%%%%%%%%%%%%%%%%%%%%%%%%%%%
%%%%%%%%%%%%%%%%%%%%%%%%%%%%%%%%%%%%%%%%%%%%%%%%%%%%%%%%%%%%%%%%%%%

An expression $G=HK$ of a group $G$ as the product of subgroups $H$ and $K$ is called a \emph{factorization} of $G$, where $H$ and $K$ are called \emph{factors}. A group $G$ is said to be \emph{almost simple} if $S\leqslant G\leqslant\Aut(S)$ for some nonabelian simple group $S$, where $S=\Soc(G)$ is the \emph{socle} of $G$. The aim of this paper is to solve the long-standing open problem:

\begin{problem}\label{PrbXia1}
Classify factorizations of finite almost simple groups.
\end{problem}

Determining all factorizations of almost simple groups stands as a fundamental problem in the study of simple groups, which was formally proposed by Wielandt~\cite[6(e)]{Wielandt1979} in 1979 on simple groups but had received attention much earlier in the literature (see for instance~\cite{Ito1953}).
Besides the importance in group theory, it also has numerous applications in other branches of mathematics, such as combinatorics~\cite{CJT2007,GGP,GLX2017,Li2006,LP2008,LX}, Hopf algebra~\cite{EGGS2000,Kac1968,LWX2023,Takeuchi1981} and number theory~\cite{FGS1993,GS1995}.

In what follows, all groups are assumed to be finite if there is no special instruction.
Problem~\ref{PrbXia1} for exceptional groups of Lie type was solved by Hering, Liebeck and Saxl~\cite{HLS1987} in 1987.
For the other families of almost simple groups, a landmark was achieved later on by Liebeck, Praeger and Saxl~\cite{LPS1990}, classifying the \emph{maximal factorizations}, namely, those with both factors maximal.
%
%\begin{theorem}[Liebeck-Praeger-Saxl]\label{ThmMaximal}
%Let $G$ be an almost simple group and let $A$ and $B$ be maximal subgroups of $G$. Then the triples $(G,A,B)$ with $G=AB$ are precisely known.
%\end{theorem}
Since then, this work has served as a cornerstone for subsequent investigations into factorizations of almost simple groups. For example, the factorizations of alternating and symmetric groups were also classified in~\cite{LPS1990}, the factorizations of sporadic almost simple groups were classified by Giudici in~\cite{Giudici2006}, and the factorizations of almost simple groups with a trivial intersection of the two factors were determined in~\cite{BL2021} and~\cite{LWX2023}, complemented by prior advancements~\cite{Baumeister2005,Baumeister2007,FGS1993,Li2003,Li2006,LPS2010,WW1980}.

Recently, factorizations of almost simple groups with a factor having at least two nonsolvable composition factors were classified in~\cite{LX2019}, and those with a factor being solvable were described in~\cite{LX} and~\cite{BL2021}.
Thus Problem~\ref{PrbXia1} is reduced to the classification of factorizations of almost simple classical groups such that both factors have a unique nonsolvable composition factor.

The main result of this paper, encapsulated in Theorem~\ref{ThmMain}, attains this objective, which classifies the factorizations $G=HK$ of almost simple classical groups $G$ such that both $H$ and $K$ have a unique nonsolvable composition factor.
For the convenience in stating the main result, let us establish some terminology. The undefined notation in this paper follows~\cite[\S2.1]{LX} or otherwise is standard.

The \emph{solvable residual} of a group $G$ is denoted by $G^{(\infty)}$; that is, $G^{(\infty)}$ denotes the smallest normal subgroup of $G$ such that $G/G^{(\infty)}$ is solvable.
For an almost simple classical group $G$, the following three operations on subgroups $H$ and $K$ of $G$ do not change the property that $G=HK$:
\begin{itemize}
\item swap $H$ and $K$;
\item replace $H$ and $K$ by $H^\alpha$ and $K^\alpha$ for some $\alpha\in\Aut(G)$ (in particular, $\alpha\in\Aut(G^{(\infty)})$);
\item replace $H$ and $K$ by $H^x$ and $K^y$ for some $x,y\in G$.
\end{itemize}
Starting from an almost simple classical group $G$ with subgroups $H_0$ and $K_0$, one may apply a finite number of the above operations and obtain subgroups $H_1$ and $K_1$ of $G$. In this case, we say the triples $(G,H_0,K_0)$ and $(G,H_1,K_1)$ are \emph{equivalent} (see Definition~\ref{DefEquiv} for a more formal definition).

\begin{definition}\label{DefTight}
A group $X$ is said to \emph{tightly contain} a group $Y$ if $X\geqslant Y$ and $X^{(\infty)}=Y^{(\infty)}$.
For a group $G$ with subgroups $H$ and $K$, we say that $(G,H,K)$ \emph{tightly contains} $(G_0,H_0,K_0)$ if
\begin{equation}\label{EqnTight}
G/G^{(\infty)}=\big(HG^{(\infty)}/G^{(\infty)}\big)\big(KG^{(\infty)}/G^{(\infty)}\big)
\end{equation}
and $G$, $H$ and $K$ tightly contain $G_0$, $H_0$ and $K_0$ respectively.
\end{definition}

%For groups $X$ and $Y$, we say that $X$ \emph{tightly contains} $Y$ if $X\geqslant Y$ and $X^{(\infty)}=Y^{(\infty)}$.

Note that~\eqref{EqnTight} is obviously a necessary condition for $G=HK$, and the determination of the factorizations of $G/G^{(\infty)}$ for a classical group $G$ is straightforward, as $G/G^{(\infty)}$ has a rather simple structure (in particular solvable).

\begin{theorem}\label{ThmMain'}
Let $G$ be an almost simple classical group, and let $H$ and $K$ be core-free subgroups of $G$ such that both $H$ and $K$ have a unique nonsolvable composition factor. Then $G=HK$ if and only if $(G,H,K)$ tightly contains some triple $(\overline{G_0},\overline{H_0},\overline{K_0})$ up to equivalence, where $\,\overline{\phantom{\varphi}}\,$ is the quotient modulo scalars and $(G_0,H_0,K_0)$ lies in Tables~$\ref{TabLinear}$, $\ref{TabUnitary}$,~$\ref{TabUnitary-2}$, $\ref{TabOmega}$, $\ref{TabOmegaMinus}$, $\ref{TabOmegaPlus}$,~~$\ref{TabOmegaPlus2}$, $\ref{TabSymplectic1}$ and $\ref{TabSymplectic2}$.
\end{theorem}

\begin{remark}
The triples $(G_0,H_0,K_0)$ in Tables~\ref{TabLinear}, \ref{TabUnitary}, \ref{TabUnitary-2}, \ref{TabOmega}, \ref{TabOmegaMinus}, \ref{TabOmegaPlus}, \ref{TabOmegaPlus2}, \ref{TabSymplectic1} and \ref{TabSymplectic2} are those satisfying the description of the corresponding examples, lemmas or propositions (with $G_0$, $H_0$ and $K_0$ sometimes denoted as $G$, $H$ and $K$) whose labels are displayed in the last column. The intersection $H_0\cap K_0$ is also provided in these tables, where the appearance of $\Sp_m(q)$ indicates that $m$ is even, and we set $\Sp_0(q)=1$ for convenience.
Although, roughly speaking, Theorem~\ref{ThmMain'} lists the \emph{minimal} (with respect to tight containment) factorizations $(G_0,H_0,K_0)$ to keep the result concise, all the factorizations $G=HK$ can be obtained following the last paragraph of Subsection~\ref{Sec2}.
\end{remark}

More remarks on Theorem~\ref{ThmMain'} can be found in Section~\ref{SecXia3}, including an illustration of the approaches to prove the theorem. The proof itself is carried out in Sections~\ref{SecXia1}--\ref{SecSymplectic02}, one section for each family of classical groups.
%The proof strategy begins with a reduction to potential candidates for $H^{(\infty)}$ and $K^{(\infty)}$ based on the maximal factorizations given in~\cite{LPS1990}. However, the core and substantial component of this work, presented in Sections~\ref{SecXia1} to~\ref{SecSymplectic02}, is the explicit determination of which corresponding pairs $(H, K)$ actually yield a factorization. For this, geometric approaches serve as the primary technique.
To address the unipotent radicals of parabolic subgroups---the most challenging cases in the classification---we consider representations of classical groups over defining characteristics. The resolution of these cases relies on the novel linearized polynomial model introduced in~\cite{FLWXZ} to describe the submodules of the unipotent radical.

A somewhat surprising consequence of Theorem~\ref{ThmMain'} is that at least one of the two factors in a factorization of an almost simple classical group is always large, remaining close to being maximal. This phenomenon is made explicit in the next two corollaries, which provide quantitative and qualitative versions of this property.

\begin{corollary}\label{CorOrder}
Let $G$ be an almost simple group whose socle is a group of Lie type of rank $k$. Then for every factorization of $G$, there is a factor of order at least $|G|^{c(k)}$, where $c(k)$ is a function of $k$ such that $\lim_{k\to\infty}c(k)=1$.
\end{corollary}

It is worth remarking on Corollary~\ref{CorOrder} that, in a factorization $G=HK$ of a finite group, one of the factors has order at least $|G|^{1/2}$, namely, $\max\{|H|,|K|\}\geqslant|G|^{1/2}$, and this bound is tight in general.

\begin{corollary}\label{CorTranGps}
Let $G$ be an almost simple classical group of Lie type with socle $L$, and let $G=HK$ be a factorization. Then up to conjugation in $\Aut(L)$, interchanging $H$ and $K$ if necessary, one of the following holds:
\begin{enumerate}[{\rm(a)}]
\item $\Nor_L(K^{(\infty)})$ is the stabilizer in $L$ of a $1$-space;
%and $(L,\Nor_L(K^{(\infty)}))$ lies in the table:
%\[
%\begin{array}{|c|ccccccc|}
%\hline
%L & \PSL_n(q) & \PSU_n(q) & \PSp_{2m}(q) & \Omega_{2m+1}(q) & \Omega_{2m+1}(q) & \Omega_7(q) & \POm_{2m}^\varepsilon(q) \\
%\hline
%\Nor_L(K^{(\infty)}) & \Pa_1 & \N_1 & \Pa_1 & \N_1^-,\, \N_1^+ & \N_1^- & \N_1^+,\, \Pa_1 & \N_1,\, \Pa_1 \\
%\hline
%\textup{Conditions} & & & & q\textup{ even} & q\textup{ odd} & q\textup{ odd} & \varepsilon\in\{+,-\} \\
%\hline
%\end{array}
%\]
\item $\Nor_L(K^{(\infty)})$ is the stabilizer in $L$ of a nondegenerate $2$-space;
%and $(L,H^{(\infty)},\Nor_L(K^{(\infty)}))$ lies in the table:
%\[
%\begin{array}{|c|cccccc|}
%\hline
%L & \Sp_{2m}(q) & \Omega_{2m+1}(q) & \POm_{2m}^-(q) & \POm_{2m}^+(q) & \Omega_{2m}^+(q) & \Omega_{2m}^+(q) \\
%\hline
%H^{(\infty)} & \Sp_m(q^2) & \Sp_m(q^2) & \SU_m(q) & (\Pa_m)^{(\infty)} & \SL_m(q) & \SU_m(q) \\
%\hline
%\Nor_L(K^{(\infty)}) & \N_2 & \N_2 & \N_2^+ & \N_2^- & \N_2^- & \N_2^+\\
%\hline
%\textup{Conditions} & q\in\{2,4\} & q\in\{2,4\} & m\textup{ odd} & & q\in\{2,4\} & q\in\{2,4\},\,m\textup{ even} \\
%\hline
%\end{array}
%\]
\item $\Nor_L(K^{(\infty)})$ is the stabilizer in $L$ of an antiflag;
%and $(L,H^{(\infty)},K^{(\infty)})$ lies in the table:
%\[
%\begin{array}{|c|ccc|}
%\hline
%L & \PSL_n(q) & \PSL_{2m}(q) & \PSL_6(q) \\
%\hline
%H^{(\infty)} & \lefthat\Sp_n(q)' & \lefthat\SL_{m}(q^2),\, \Sp_{m}(q^2) & \G_6(q)' \\
%\hline
%K^{(\infty)} & \SL_{n-1}(q) & \SL_{2m-1}(q) & \SL_5(q) \\
%\hline
%\textup{Conditions} & & q\in\{2,4\} & q\textup{ even} \\
%\hline
%\end{array}
%\]
\item $L=\PSL_4(q)$, or $\PSp_4(q)$ with odd $q$, or $\Sp_6(q)$ with even $q$, or $\Omega_7(q)$ with odd $q$;
\item $L$ is in a known list of finitely many groups.
\end{enumerate}
\end{corollary}

Corollary~\ref{CorTranGps} motivates the study of subgroups of classical groups that act transitively on certain families of subspaces. In fact, this topic has garnered significant interest, with numerous partial results established over several decades~\cite{GGP, Hering1974,Hering1985,Kantor1973,KantorLiebler1982,LPS2010,Perin1972}. The classification of such subgroups is now ultimately completed by Theorem~\ref{ThmMain'}, as each of them is precisely a factor in the factorization of the classical group with the other factor being the stabilizer of the subspace.

Beyond its role in classifying transitive subgroups on subspaces, the classification presented in this paper has served as a principal tool in several recent works: it is utilised in~\cite{CGP2025,YFX2023} to study $2$-arc-transitive digraphs, in~\cite{LP2025} to explore arc-transitive $2$-cell embeddings of graphs, and most recently in~\cite{BGKLPW} to investigate regular subgroups of primitive groups. In this paper, we give a further application to Hopf algebras.

A group factorization $G=HK$ is called \emph{exact} if $H\cap K=1$. Such factorizations arise naturally in the theory of Hopf algebras: the construction of Kac~\cite{Kac1968} and Takeuchi~\cite{Takeuchi1981} associates with an exact factorization $G=HK$ a semisimple Hopf algebra $\mathcal H(G,H,K)$ of dimension $|G|$, called a \emph{bicrossproduct} Hopf algebra. In 2000, Etingof, Gelaki, Guralnick and Saxl~\cite{EGGS2000} proved that $\mathcal H(G,H,K)$ is biperfect if and only if both $H$ and $K$ are perfect and self-normalizing. Giving a negative answer to~\cite[Question~7.5]{EG2000}, the authors then used the exact factorization
\[
\M_{24}=\PSL_2(23)\,(2^4{:}\A_7)
\]
to construct a biperfect Hopf algebra of dimension $|\M_{24}|=244823040$, and asked whether a biperfect Hopf
algebra of smaller dimension exists~\cite[Question~3.4(2)]{EGGS2000}.
Moreover, it is proved in~\cite{LWX2023} that the above factorization of $\M_{24}$ is the unique exact factorization of a finite simple group giving rise to a biperfect Hopf algebra.
In Section~\ref{sec:Application}, we answer the question of Etingof et al.~in the affirmative and, more precisely, use the classification of factorizations of almost simple groups to determine the smallest possible dimension within the bicrossproduct construction.

\begin{theorem}
The smallest dimension of a biperfect bicrossproduct Hopf algebra is $41287680$.
\end{theorem}

%%%%%%%%%%%%%%%%%%%%%%%%%%%%%%%%%%%%%%%%%%%%%%%%%%%%%%%%%%%%%%%%%%%
%%%%%%%%%%%%%%%%%%%%%%%%%%%%%%%%%%%%%%%%%%%%%%%%%%%%%%%%%%%%%%%%%%%

\section{Preliminaries, approaches, and remarks}\label{SecXia3}

%%%%%%%%%%%%%%%%%%%%%%%%%%%%%%%%%%%%%%%%%%%%%%%%%%%%%%%%%%%%%%%%%%%
%%%%%%%%%%%%%%%%%%%%%%%%%%%%%%%%%%%%%%%%%%%%%%%%%%%%%%%%%%%%%%%%%%%

%In this section we collect notation and preliminary results needed throughout the paper.
Most preliminary results of this section will be used repeatedly but usually without mentioning in this paper.

\subsection{Remarks on the main result}\label{Sec2}
\ \vspace{1mm}

To be precise, a \emph{classical group} in this paper refers to a group $G$ that satisfies
\begin{equation}\label{EqnXia18}
\Omega\leqslant G\leqslant\Gamma{:}\langle\iota\rangle\quad\text{or}\quad\Omega/\Z(\Omega)\leqslant G\leqslant\Aut(\Omega/\Z(\Omega)),
\end{equation}
where $\Omega$ is a group in the following table, $\Gamma$ is the corresponding conformal semilinear group defined in~\cite[\S1.6.2]{BHR2013}, and $\iota$ is the transpose-inverse if $\Omega=\SL_n(q)$ with $n\geqslant3$ while $\iota=1$ otherwise.
We call a group $G$ with $\Omega\leqslant G\leqslant\Gamma{:}\langle\iota\rangle$ an \emph{almost quasisimple classical group}.
\[
\begin{array}{c|cccccc}
\hline
\Omega & \SL_n(q) & \SU_n(q) & \Sp_{2m}(q) & \Omega_{2m+1}(q) & \Omega_{2m}^+(q) & \Omega_{2m}^-(q) \\
\hline
\text{Conditions} & n\geqslant2 & n\geqslant3 & m\geqslant2 & m\geqslant3 & m\geqslant4 & m\geqslant4 \\
 & q\geqslant4\text{ if }n=2 & (n,q)\neq(3,2) & (m,q)\neq(2,2) & q\text{ odd} &  & \\
\hline
\end{array}
\]

In a factorization $G=HK$, the factors $H$ and $K$ are called a \emph{supplement} of $K$ and $H$, respectively, in $G$.
%Group factoriztions are naturally related to transitive actions.
Let $K$ be a subgroup of a group $G$, and let
\[
\Delta=K\backslash G=\{Kg\mid g\in G\}
\]
be the set of right cosets of $K$ in $G$. Then $G$ acts transitively by right multiplication on $\Delta$, and $K=G_\delta$ is the stabilizer of a point $\delta\in\Delta$. Clearly, $H$ is transitive on $\Delta$ if and only if $G=HK$. Indeed, we have more equivalent conditions for $G=HK$ in the following well-known (and straightforward) lemma.

\begin{lemma}\label{LemXia2}
Let $H$ and $K$ be subgroups of $G$. Then the following are equivalent:
\begin{enumerate}[{\rm (a)}]
\item $G=HK$;
\item $G=H^\alpha K^\alpha$ for any $\alpha\in\Aut(G)$;
\item $G=H^xK^y$ for any $x,y\in G$;
\item $|G||H\cap K|=|H||K|$;
\item $|G|\leqslant|H||K|/|H\cap K|$;
\item $H$ acts transitively by right multiplication on $K\backslash G$;
\item $K$ acts transitively by right multiplication on $H\backslash G$.
\end{enumerate}
\end{lemma}

%\begin{remark}\label{RmkXia1}
%Lemma~\ref{LemXia2} plays a fundamental role in the study of group factorizations. For example,
%\begin{itemize}
%\item due to part~(c) we may consider specific representatives of a conjugacy class of subgroups according to our requirements when studying factorizations of a group;
%\item given a group $G$ and its subgroups $H$ and $K$, in order to verify whether $G=HK$, we only need to compute the orders of $G$, $H$, $K$ and $H\cap K$ by part~(d) or~(e), which enables us to search group factorizations efficiently in \magma~\cite{BCP1997}.
%\end{itemize}
%\end{remark}

Lemma~\ref{LemXia2} motivates the following definition, which generalizes equivalence for almost simple classical groups introduced before Definition~\ref{DefTight}.

\begin{definition}\label{DefEquiv}
Let $G$ be a classical group with subgroups $H_1$, $H_2$, $K_1$ and $K_2$, we say that $(G,H_1,K_1)$ is \emph{equivalent} to $(G,H_2,K_2)$ if there exist $\alpha\in\Aut(\overline{G}^{(\infty)})$ and $x,y\in G$ such that $\overline{G}=\overline{G}^\alpha$ and either $\overline{H_1^x}=\overline{H_2}^\alpha$ and $\overline{K_1^y}=\overline{K_2}^\alpha$ or $\overline{H_1^x}=\overline{K_2}^\alpha$ and $\overline{K_1^y}=\overline{H_2}^\alpha$, where $\,\overline{\phantom{\varphi}}\,$ is the quotient modulo scalars.
\end{definition}

\begin{remark}
When considering possible factorizations for classical groups, equivalent triples can be viewed as the same.
Recall the concept of tight containment given in Definition~\ref{DefTight}.
From now on in this paper, whenever we say $(G,H,K)$ tightly contains $(G_0,H_0,K_0)$, we mean $(G,H,K)$ tightly contains some triple that is equivalent to $(G_0,H_0,K_0)$.
\end{remark}

To classify the factorizations of almost simple classical groups, we first establish the following theorem that classifies the factorizations of almost quasisimple classical groups.

\begin{theorem}\label{ThmMain}
Let $G$ be an almost quasisimple classical group, and let $H$ and $K$ be subgroups of $G$ not containing $G^{(\infty)}$ such that both $H$ and $K$ have a unique nonsolvable composition factor. Then $G=HK$ if and only if either $(G,H,K)$ tightly contains some $(G_0,H_0,K_0)$ in Tables~$\ref{TabLinear}$, $\ref{TabUnitary}$, $\ref{TabOmega}$, $\ref{TabOmegaMinus}$, $\ref{TabOmegaPlus}$, $\ref{TabSymplectic1}$ and $\ref{TabSymplectic2}$ or $(\overline{G},\overline{H},\overline{K})$ tightly contains some $(G_0,H_0,K_0)$ in Tables~$\ref{TabUnitary-2}$ and~$\ref{TabOmegaPlus2}$, where $\,\overline{\phantom{\varphi}}\,$ is the quotient modulo scalars.
\end{theorem}

%The statement of Theorem~\ref{ThmMain} is basically a presentation of the minimal factorizations with respect to tight containment. {\color{red}Xiao Xia: Single out those few that are not???}

Theorem~\ref{ThmMain} follows from a combination of the main theorems in Sections~\ref{SecXia1} to~\ref{SecSymplectic02}. For a group homomorphism $\varphi\colon G\to G^\varphi$, the factorizations of $G^\varphi$ are precisely $G^\varphi=H^\varphi K^\varphi$ from the factorizations $G=HK$ of $G$. In this way, Theorem~\ref{ThmMain} gives the factorizations of $\overline{G}$ such that both factors have a unique nonsolvable composition factor. Note that, according to our definition of classical groups in~\eqref{EqnXia18}, an almost simple classical group is not necessarily the quotient of an almost quasisimple classical group modulo scalars. In fact, the only exceptions are the $4$-dimensional symplectic groups in even characteristic and the $8$-dimensional orthogonal groups of plus type, in which cases, $\overline{\Gamma}$ is a subgroup of $\Aut(\Omega/\Z(\Omega))$ of index $2$ or $3$, respectively.
Nevertheless, for these exceptions, our argument in Sections~\ref{SecOmegaPlus05} and~\ref{SecSymplectic02}, with minimal adjustments to accommodate quotients modulo scalars, shows that there are no such factorizations of $\Aut(\Omega/\Z(\Omega))$. Thus, Theorem~\ref{ThmMain'} follows.

Finally, we remark on how to determine all the factorizations from the minimal (with respect to tight containment) ones in Theorem~\ref{ThmMain'} or Theorem~\ref{ThmMain}.
Observe that, given a factorization $G_0=H_0K_0$ and an almost quasisimple classical group $G$ that tightly contains $G_0$, for any subgroups $H$ and $K$ of $G$ containing $H_0$ and $K_0$ respectively, $G=HK$ if and only if~\eqref{EqnTight} holds. Thus, all the factorizations $G=HK$ corresponding to a triple $(G_0,H_0,K_0)$ in Theorem~\ref{ThmMain} are precisely given by the triples $(G,H,K)$ with $G$ tightly containing $G_0$,
\[
H_0\leqslant H\leqslant\Nor_G(H_0^{(\infty)})\ \text{ and }\ K_0\leqslant K\leqslant\Nor_G(K_0^{(\infty)})
\]
such that~\eqref{EqnTight} holds. Similarly, one can obtain all the factorizations of almost simple groups corresponding to the triples $(\overline{G_0},\overline{H_0},\overline{K_0})$ in Theorem~\ref{ThmMain'}, or equivalently, such factorizations are precisely given by $(\overline{G},\overline{H},\overline{K})$ with $(G,H,K)$ as in Theorem~\ref{ThmMain}.

\subsection{Remarks on the literature}
\ \vspace{1mm}

A subgroup of a classical group $G$ is called $max^+$ if it is maximal and does not contain $G^{(\infty)}$, and is called $max^-$ if it is a maximal one among the subgroups not containing $G^{(\infty)}$.
Note that $\max^+$ subgroups are necessarily $\max^-$, but the converse is not true.
For example, the stabilizer of a $1$-space in $\PSL_3(3)$ is a $\max^-$ subgroup of $\Aut(\PSL_3(3))=\PSL_3(3).2$ but not a $\max^+$ subgroup of $\Aut(\PSL_3(3))$.
A factorization of $G$ with $\max^\varepsilon$ factors is said to be a $max^\varepsilon$ factorization, where $\varepsilon\in\{+,-\}$.
The $\max^+$ factorizations of almost simple groups are classified in~\cite{LPS1990}, and the $\max^-$ factorizations of almost simple groups that are not $\max^+$ are determined in~\cite{LPS1996}.

For the reader's benefit, we collect some corrigenda of the published work in the literature on factorizations of almost simple groups.
\begin{enumerate}[{\rm(I)}]
\item In part~(b) of Theorem~2 in~\cite{HLS1987}, $A_0$ can also be $\G_2(2)$, $\SU_3(3)\times2$, $\SL_3(4).2$ or $\SL_3(4).2^2$ besides $\G_2(2)\times2$.
\item It is pointed out in~\cite[Remark~1.4]{GGS2024} that Table~4 of~\cite{LPS1990} should contain the $\max^+$ factorizations $G=\N_2^-[G](\mathrm{O}_8^-(2^f).f)$ (see~\cite[\S2.1]{LX} for the notation $\N_2^-[G]$) for $G=\Omega_8^+(4^f).2f$ with $f\in\{1,2\}$.
\item It is pointed out in~\cite{GGP} that a regular subgroup of order $2^9\cdot3\cdot7^2$ in $\PGaL_3(8)$ is missing from Table~16.1 of~\cite{LPS2010}. In fact, this is covered by the classification~\cite[Proposition~4.1]{LX2019} of factorizations $G=HK$ of almost simple groups with solvable factors: the triples $(G,H,K)$ here are $(\PGaL_3(8).\calO,73{:}(9\times\calO_1),2^{3+6}{:}7^2{:}(3\times\calO_2))$ with $\calO\leqslant\C_2$ and $\calO_1,\calO_2\leqslant\calO$ such that $\calO=\calO_1\calO_2$. Besides the case $\calO=\calO_1=\calO_2=1$, the cases $(\calO,\calO_1,\calO_2)=(2,1,2)$ and $(\calO,\calO_1,\calO_2)=(2,2,1)$ also give regular subgroups of primitive groups and hence are also missing from Table~16.1 of~\cite{LPS2010}. The case $\calO=\calO_1=\calO_2=2$ give a $\max^+$ factorization that is overlooked in~\cite{LPS1990}.
\item In Table~1 of~\cite{LX2019}, the triple $(L,H\cap L,K\cap L)=(\Sp_{6}(4),(\Sp_2(4)\times\Sp_{2}(16)).2,\G_2(4))$ is missing, and for the first two rows $R.2$ should be $R.P$ with $P\leqslant2$. Also, see the remark after our Lemma~\ref{LemXia9} for a correction of~\cite[Lemma~5.1]{LX2019}.
\item The second and fourth authors have recognized, through communication with Mikko Korhonen, that the following triples for $(G,H,K)$ should have been in~\cite[Proposition~4.3]{LX}:
    \[
    (\Sy_9,\Sy_4\times\AGL_1(5),\PGaL_2(8)),\ \ (\A_9,(\A_4\times\D_{10})^{\boldsymbol{\cdot}}2,\PGaL_2(8)).
    \]
    This oversight occurred due to a careless reading of the list of $k$-homogeneous permutation groups for $k\geqslant4$ in the second last paragraph of the proof.
\end{enumerate}

\subsection{Arithmetic results}
\ \vspace{1mm}

For positive integers $a$ and $b$, denote the greatest common divisor of $a$ and $b$ by $(a,b)$ when there is no confusion with the pair $(a,b)$. For a positive integer $n$ and a prime $p$, denote by $n_p$ the largest $p$-power dividing $n$, and denote $n_{p'}=n/n_p$.
%For convenience, we say that $n$ is divisible by $a/b$ if and only if $nb$ is divisible by $a$. For example, $5$ is divisible by $5/2$ but not by $3/2$.

Let $a>1$ and $k>1$ be integers. A prime number $r$ is called a \emph{primitive prime divisor} of the pair $(a,k)$ if $r$ divides $a^k-1$ but does not divide $a^j-1$ for any positive integer $j<k$. We record the celebrated Zsigmondy's theorem here on the existence of primitive prime divisors.

\begin{theorem}[Zsigmondy~\cite{Zsigmondy1892}]
If $k\geqslant2$, then $(a,k)$ has a primitive prime divisor except for either $k=2$ and $a=2^e-1$ is a Mersenne prime, or $(a,k)=(2,6)$.
%Moreover, if $r$ is a primitive prime divisor of $(a,k)$ then $r>k$.
\end{theorem}

If $a$ is prime and $(a,k)\neq(2,6)$, then a primitive prime divisor of $(a,k)$ is simply called a \emph{primitive prime divisor} of $a^k-1$, and we denote the set of such primitive prime divisors by
\[
\ppd(a^k-1).
\]
For example, $\ppd(15)=\ppd(2^4-1)=\{5\}$. Moreover, set $\ppd(63)=\{7\}$.
We say that an integer is \emph{divisible by} $\ppd(p^b-1)$ if it is divisible by all the primes in $\ppd(p^b-1)$, and that an integer is \emph{coprime to} $\ppd(p^b-1)$ if it is coprime to all the primes in $\ppd(p^b-1)$.

Checking the orders of outer automorphism groups of finite simple groups leads to the following result (see~\cite[Page~38, Proposition~B]{LPS1990}).

\begin{lemma}\label{LemXia6}
Let $L$ be a simple group of Lie type over $\bbF_q$. If $k\geqslant3$ and $(q,k)\neq(2,6)$, then $\ppd(q^k-1)$ is coprime to $|\Out(L)|$.
\end{lemma}

For a group $X$, the \emph{solvable radical} of $X$ is denoted by $\Rad(X)$, which is the largest solvable normal subgroup of $X$.
If $G$ is an almost simple group with socle $L$ and $X$ is a subgroup of $G$, then $\Rad(X\cap L)=\Rad(X)\cap L$, and so $|\Rad(X)|/|\Rad(X\cap L)|=|\Rad(X)L/L|$ divides $|\Out(L)|$.
Then, by Lemma~\ref{LemXia6}, we obtain the following lemma by checking $|\Rad(A\cap L)|$ and $|\Rad(B\cap L)|$ for the triples $(L,A\cap L,B\cap L)$, where $G=AB$ is a $\max^-$ factorization (see~\cite{LPS1990} and~\cite{LPS1996}).

\begin{lemma}\label{LemXia22}
Let $G$ be an almost simple group with socle $L$, where $L$ is a simple group of Lie type over $\bbF_q$, and let $G=AB$ be a $\max^-$ factorization with nonsolvable factors $A$ and $B$. If $k\geqslant3$ and $(q,k)\neq(2,6)$, then $\ppd(q^k-1)$ is coprime to both $|\Rad(A)|$ and $|\Rad(B)|$.
\end{lemma}

\subsection{Some actions of classical groups}\label{SecXia2}
\ \vspace{1mm}

For classical groups $G$, the groups $\Pa_i[G]$, $\N_i[G]$, $\N_i^+[G]$ and $\N_i^-[G]$ are defined in~\cite[\S2.1]{LX}.
Let $\mathcal{P}_i[G]$, $\mathcal{N}_i[G]$, $\mathcal{N}_i^+[G]$ and $\mathcal{N}_i^-[G]$ be the set of right cosets of $\Pa_i[G]$, $\N_i[G]$, $\N_i^+[G]$ and $\N_i^-[G]$ in $G$, and let $\mathcal{P}_i^{(\infty)}[G]$, $\mathcal{N}_i^{(\infty)}[G]$, $\mathcal{N}_i^{+(\infty)}[G]$ and $\mathcal{N}_i^{-(\infty)}[G]$ be the set of right cosets of $\Pa_i[G]^{(\infty)}$, $\N_i[G]^{(\infty)}$, $\N_i^+[G]^{(\infty)}$ and $\N_i^-[G]^{(\infty)}$ in $G$, respectively. When the group $G$ is clear from the context, we simply write them as $\calP_i$, $\calN_i$, $\calN_i^+$, $\calN_i^-$, $\calP_i^{(\infty)}$, $\calN_i^{(\infty)}$, $\calN_i^{+(\infty)}$ and $\calN_i^{-(\infty)}$.
It is worth noting that these sets of right cosets can be identified with the sets of certain subspaces. For example, $\mathcal{P}_1[\SL_n(q)]$ is the set of $1$-spaces in a vector space of dimension $n$ over $\bbF_q$, and $\mathcal{N}_1[\SU_n(q)]$ is the set of nonsingular $1$-spaces in a unitary space of dimension $n$ over $\bbF_{q^2}$, while $\mathcal{N}_1^{(\infty)}[\SU_n(q)]$ is the set of nonsingular vectors therein.
%Determining transitive subgroups of $G$ on these sets is an important special case for determining factorizations of $G$, which we explain below by a couple of examples.

%Let $G=\PGaL_n(q)$. Then $G$ is transitive on $\calP_1[G]$, the set of 1-subspaces of $\mathbb{F}_q^n$. The subgroups $H$ of $G$ that are transitive on $\calP_1[G]$ are classified by Hering (see~\cite[Lemma~3.1]{LPS2010}), and each $H$ gives rise to a factorization $G=H\Pa_1[G]$, see Table~\ref{TabLinear}.

An \emph{antiflag} of a vector space $V$ is an unordered pair $\{\langle v\rangle,W\}$, where $v$ is a nonzero vector in $V$ and $W$ is a hyperplane in $V$ not containing $v$.
%The antiflag-transitive subgroups of $\PGaL_n(q)$ are classified by Cameron and Kantor~\cite{CK1979} (see~\cite{Kantor} for a correction), each of which gives rise to a factorization of $\PGaL_n(q)$. In this paper we will adopt a slightly different version of the classification from~\cite[Theorem~3.2]{LPS2010}, see Table~\ref{TabLinear}.
For such a pair $v$ and $W$, we call $\{v,W\}$ a \emph{refined antiflag}. In the same spirit as above, we use $\calP_{1,n-1}^{(\infty)}$ to denote the set of refined antiflags in an $n$-dimensional vector space. Transitive actions on $\calP_{1,n-1}^{(\infty)}$ are said to be \emph{refined-antiflag-transitive}.

%For some other classical groups $G$, subgroups of $G$ that are transitive on a $G$-orbit of certain subspaces are described in~\cite[Chapter~4]{LPS2010}. Recently, the work in~\cite[Chapter~4]{LPS2010} has been extended by Giudici, Glasby and Praeger~\cite{GGP}.

It is worth remarking that it is elementary to write down the group structures of $\Pa_i$, $\N_i$, $\N_i^+$ and $\N_i^-$. Let us take $\Pa_1[\Sp_{2m}(q)]$, the stabilizer in $\Sp_{2m}(q)$ of a $1$-space, as an example. Write
\[
\Pa_1[\Sp_{2m}(q)]=R{:}(S\times T)=(Q.q^{2m-2}){:}(S\times T)
\]
in the canonical way, where $Q=q$, $R=Q.q^{2m-2}$, $S=\Sp_{2m-2}(q)$ and $T=\GL_1(q)$. Here $\langle Q,S\rangle=Q\times S$, while $R$ is abelian if $q$ is even, and nonabelian with center $Q$ if $q$ is odd. Moreover, the following property is needed in Proposition~\ref{LemSymplectic12} and fairly easy to prove.

\begin{lemma}\label{LemXia21}
Let $(m,q)\neq(2,2)$ or $(2,3)$. Then every proper subgroup of $R$ normalized by $S$ is contained in $Q$. In particular, $\Pa_1[\Sp_{2m}(q)]^{(\infty)}=R{:}S$.
\end{lemma}
%
%\begin{proof}
%Let $X$ be a subgroup of $R$ normalized by $S$, and let
%\[
%R(b,u)=
%\begin{pmatrix}
%1&0&0\\
%-B_{m-1}u^\mathsf{T}&I_{2m-2}&0\\
%b&u&1
%\end{pmatrix}
%\,\text{ for }b\in\bbF_q\text{ and }u\in\bbF_q^{2m-2}.
%\]
%Suppose for a contradiction that $X<R$ and $X\nleqslant Q$. Then $X$ contains some $R(b,u)$ with $u\neq0$. Since
%\[
%\begin{pmatrix}
%1&&\\
%&A&\\
%&&1
%\end{pmatrix}^{-1}
%R(b,u)
%\begin{pmatrix}
%1&&\\
%&A&\\
%&&1
%\end{pmatrix}
%=R(b,uA)
%\ \text{ for }
%\begin{pmatrix}
%1&&\\
%&A&\\
%&&1
%\end{pmatrix}
%\in S
%\]
%and $X$ is normalized by $S=\Sp_{2m-2}(q)$, we have
%\[
%X\supseteq R(b,u)^S=\{R(b,v)\mid v\in\bbF_q^{2m-2},\,v\neq0\}.
%\]
%This implies that
%\begin{align*}
%X&\supseteq\{R(b,u+w)R(b,w)^{-1}\mid w\in\bbF_q^{2m-2},\,u+w\neq0,\,w\neq0\}\\
%&=\{R(uB_{m-1}w^\mathsf{T},u)\mid w\in\bbF_q^{2m-2},\,w\neq-u,\,w\neq0\}=\{R(c,u)\mid c\in\bbF_q\},
%\end{align*}
%and so
%\[
%X\supseteq\{R(c,u)\mid c\in\bbF_q\}^S=\{R(c,v)\mid c\in\bbF_q,\,v\in\bbF_q^{2m-2},\,v\neq0\}.
%\]
%As a consequence, $|X|$ is at least $q(q^{2m-2}-1)$ and hence $|X|$ does not properly divide $|R|=q^{2m-1}$, which contradicts our assumption that $X<R$.
%\end{proof}

We conclude this subsection with some representation-theoretic notions in Subsections~\ref{SecUnitaryPm} and~\ref{SecOmegaPlus04}. For a vector space $V$, let $\bigwedge^2(V)$ denote the alternating square of $V$; that is,
\[
\mbox{$\bigwedge^2(V)=(V\otimes V)/\langle v\otimes v\mid v\in V\rangle$}.
\]
For positive integers $i$ and $j$, if $M$ is an $\bbF_{p^i}\GL_\ell(p^i)$-module with underlying vector space $\bbF_{p^i}^\ell$, then let $M^{(p^j)}$ denote the $\bbF_{p^i}\GL_\ell(p^i)$-module $\bbF_{p^i}^\ell$ with the $\GL_\ell(p^i)$ action defined by $(v,g)\mapsto vg^{\phi^j}$, where $\phi$ is the Frobenius map taking $p$th power.
If $M$ is an $\bbF_{q^k}\GL_\ell(q^k)$-module for some prime power $q$ and positive integers $\ell$ and $k$ such that $M^{(q)}$ and $M$ are identical $\GL_\ell(q^k)$-modules, then $M$ can be realized over $\bbF_q$.
For example, the modules $U_{(b)}^\sharp(i)$ defined by~\eqref{EqnUnitary1} and~\eqref{EqnOrthonal1}, respectively, can be realized over $\bbF_q$.

\subsection{Illustration of approaches}
\ \vspace{1mm}

Before working on different families of classical groups from next section, we briefly sketch our approach to classify their factorizations.

\begin{lemma}\label{LemXia10}
Let $G=MK$ be a factorization, and let $H$ be a subgroup of $M$. Then $G=HK$ if and only if $M=H(M\cap K)$.
\end{lemma}

\begin{remark}
Lemma~\ref{LemXia10} shows how we embed the factorization $M=H(M\cap K)$ of a smaller group $M$ into the larger factorization $G=MK$ with a factor $M$ to formulate the ``refined'' factorization $G=HK$.
\end{remark}

If $G=HK$ with $H$ and $K$ contained in subgroups $A$ and $B$ of $G$, respectively, then Lemma~\ref{LemXia10} gives $A=H(A\cap B)$ and $B=(A\cap B)K$. Further, if $M$ is a normal subgroup of $A$ and $B$ is a normal subgroup of $B$, then by taking quotients we have
\begin{equation}\label{EqnXia17}
A/M=(HM/M)((A\cap B)M/M)\quad\text{and}\quad B/N=((A\cap B)N/N)(KN/N).
\end{equation}
Starting from the $\max^-$ factorizations $G=AB$ of classical groups $G$, our analysis of all possible factorizations $G=HK$ is to determine which subgroups $H$ of $A$ and subgroups $K$ of $B$ satisfy $G=HK$. For this purpose, we will often apply~\eqref{EqnXia17} with $M=\Rad(A)$ and $N=\Rad(B)$, in which case the ensuing result is of fundamental importance. The result itself is a corollary of the classification of $\max^-$ factorizations in~\cite[Theorem~A]{LPS1990} and~\cite{LPS1996} (see also~\cite[Proposition~2.17]{LX}).

\begin{lemma}\label{LemXia7}
Let $G$ be a classical group, and let $G=AB$ with $\max^-$ subgroups $A$ and $B$ of $G$. If $A$ has a unique nonsolvable composition factor, then $A/\Rad(A)$ is almost simple and $(A\cap B)\Rad(A)/\Rad(A)$ is core-free in $A/\Rad(A)$. If $B$ has a unique nonsolvable composition factor, then $B/\Rad(B)$ is almost simple and $(A\cap B)\Rad(B)/\Rad(B)$ is core-free in $B/\Rad(B)$.
\end{lemma}

This lemma tells us that, if $A$ has exactly one nonsolvable composition factor, then taking $M=\Rad(A)$, the group $HM/M$ either contains $\Soc(A/M)$ or is a supplement of $(A\cap B)M/M$ in a core-free factorization of $A/M$. In either case we are led to a description for $HM/M$, and in the same vein we get a description for $KN/N$ with $N=\Rad(B)$.
%This will be the foundation of our analysis for the potential factorizations $G=HK$.
This provides us with a recursive approach to facilitate the delineation of factorizations $G=HK$, which turns out to produce the candidates for $H^{(\infty)}$ and $K^{(\infty)}$ rather quickly.
For example, starting from the $\max^-$ factorization $G=AB$ with $G^{(\infty)}=\Omega_{2m}^-(q)$, $A^{(\infty)}=\SU_m(q)$ ($m$ odd) and $B=\N_1[G]$ (see~Proposition~\ref{Prop:MaxOmegaminu}), we will show in Lemma~\ref{LemOmegaMinusRow2,3} that the candidates for $(H^{(\infty)},K^{(\infty)})$ are
\begin{equation}\label{EqnXia19}
(\SU_m(q),\Omega_{2m-1}(q))\ \text{ and }\ (\SU_m(q),\Omega_{2m-2}^-(q)).
\end{equation}

After identifying the candidates for $H^{(\infty)}$ and $K^{(\infty)}$, the next step is to determine which subgroups $H$ and $K$ of $\Nor_G(H^{(\infty)})$ and $\Nor_G(K^{(\infty)})$ indeed satisfy $G=HK$. This will be accomplished by essentially calculating $H\cap K$ (so one can apply part~(d) of Lemma~\ref{LemXia2}). Given the need to know $H\cap K$ (or at least its order) explicitly, this aspect of the process is subtle and constitutes the main focus of the paper. The argument in this part predominantly takes a geometric approach, whereas  leveraging representations of classical groups over defining characteristic will be the primary tool for handling factors from certain parabolic subgroups (see Subsections~\ref{SecUnitaryPm} and~\ref{SecOmegaPlus04}).
As for our example in the previous paragraph, the geometric argument in Example~\ref{ex:OmegaMinus01} and the majority of Subsection~\ref{nonsingular 2-spaces} leads to the characterization that, for $(H^{(\infty)},K^{(\infty)})$ in~\eqref{EqnXia19}, $G=HK$ if and only if $(G,H,K)$ tightly contains the triple $(G_0,H_0,K_0)$ in the first row of the table in Proposition~\ref{prop:OmegaMinus-02} or any row of the table in Proposition~\ref{LemOmegaMinus09}.

%%%%%%%%%%%%%%%%%%%%%%%%%%%%%%%%%%%%%%%%%%%%%%%%%%%%%%%%%%%%%%%%%%%
%%%%%%%%%%%%%%%%%%%%%%%%%%%%%%%%%%%%%%%%%%%%%%%%%%%%%%%%%%%%%%%%%%%

\section{Linear groups}\label{SecXia1}

%%%%%%%%%%%%%%%%%%%%%%%%%%%%%%%%%%%%%%%%%%%%%%%%%%%%%%%%%%%%%%%%%%%
%%%%%%%%%%%%%%%%%%%%%%%%%%%%%%%%%%%%%%%%%%%%%%%%%%%%%%%%%%%%%%%%%%%

This section is devoted to factorizations of linear groups. Throughout this section, let $q=p^f$ with $p$ prime, let $V=\bbF_q^n$, and let $G$ be a classical group with $L=G^{(\infty)}=\SL(V)=\SL_n(q)$, where $n\geqslant2$ and $(n,q)\neq(2,2)$ or $(2,3)$.

\begin{theorem}\label{ThmLinear}
Let $H$ and $K$ be nonsolvable subgroups of $G$ not containing $G^{(\infty)}$. Then $G=HK$ if and only if $(G,H,K)$ tightly contains some triple $(G_0,H_0,K_0)$ in Table~$\ref{TabLinear}$.
\end{theorem}

We shall describe examples of the factorizations in Subsection~\ref{sec:Linear-ex}, and then prove Theorem~\ref{ThmLinear} in Subsection~\ref{SecLinearProof}.

\subsection{Examples}\label{sec:Linear-ex}
%\ \vspace{1mm}
%
%The examples of factorizations of $G$ are divided into two types:
%\begin{itemize}
%\item one factor stabilizes a $1$-subspace or an antiflag;
%\item two small exceptions with $G^{(\infty)}=\SL_2(9)$ or $\SL_3(4)$.
%\end{itemize}
%We describe them separately.

\subsubsection{Transitive groups on nonzero vectors}

The subgroups of $\GaL_n(q)$ that are transitive on the set of nonzero vectors are classified by Hering~\cite{Hering1985} (see also \cite{Liebeck1987}). The cases needed here give the following factorizations.

\begin{example}\label{ex:Linear}
For each subgroup $H$ of $G$ listed in the table below, $H$ is transitive on $V\Setminus\{0\}$ with the described stabilizer $H_v$, which gives the factorization $G=HG_v$ with $H\cap G_v=H_v$.
\[
\begin{array}{llll}
\hline
G & H & G_v & H_v \\
\hline
\SL_{ab}(q) & \SL_a(q^b) & q^{ab-1}{:}\SL_{ab-1}(q) & q^{ab-b}{:}\SL_{a-1}(q^b) \\
 & \Sp_a(q^b)' &  & [q^c]{:}\Sp_{a-2}(q^b)\\
 & \G_2(q^b)'\text{ ($a=6$, $q$ even)} &  & [(q^{5b},q^{6b}/4)]{:}\SL_2(q^b) \\
\SL_4(2) & \A_7 & 2^3{:}\SL_3(2) & \PSL_2(7)\\
\SL_4(3) & 2^{\boldsymbol{\cdot}}\Sy_5 & 3^3{:}\SL_3(3) & 3\\
 & 8\circ\SL_2(5)=8^{\boldsymbol{\cdot}}\A_5 & & \Sy_3\\
 & 2^{1+4}{}^{\boldsymbol{\cdot}}\A_5 & & \SL_2(3)\\
\SL_6(3) & \SL_2(13) & 3^5{:}\SL_5(3)& 3\\
\hline
\end{array}
\]
The point stabilizer in $\G_2(q^b)'$ is described in~\cite[Lemma~5.1]{Cooperstein1981}, the factorizations of $\SL_4(2)$, $\SL_4(3)$ and $\SL_6(3)$ are verified in \magma~\cite{BCP1997}, and $c$ is defined as
\begin{equation}\label{EqnLinear01}
c=
\begin{cases}
ab-b&\text{ if }(q,a,b)\neq(2,4,1)\\
2&\text{ if }(q,a,b)=(2,4,1).
\end{cases}
\end{equation}
\end{example}

\subsubsection{Refined-antiflag-transitive groups}

Recall the definition of antiflags and refined antiflags in Section~\ref{SecXia3}. Let $V=\bbF_q^{2m}$, and let
\[
\calP_{1,2m-1}^{(\infty)}=\{\{v,W\}\mid V=\langle v\rangle\oplus W\}
\]
be the set of refined antiflags of $V$. Then $L=\SL_{2m}(q)$ is transitive on $\calP_{1,2m-1}^{(\infty)}$ with
\[
L_{\{v,W\}}=\SL_{2m-1}(q)\ \text{ and }\ |\calP_{1,2m-1}^{(\infty)}|=|L|/|L_{\{v,W\}}|=q^{2m-1}(q^{2m}-1),
\]
and we may identify $\calP_{1,2m-1}^{(\infty)}$ with the right coset set $L_{\{v,W\}}\backslash L=\SL_{2m-1}(q)\backslash\SL_{2m}(q)$. For a group $G$, each refined-antiflag-transitive subgroup $H$ of $G$ gives rise to a factorization $G=HG_{\{v,W\}}$. We next describe them.
%Antiflag transitive groups have been classified, see \cite[Theorem II]{Antiflag-trans-gps}.

Let $V$ be associated with a nondegenerate alternating form $\beta$, and let $H=\Sp(V)$. Pick a standard basis $e_1,f_1,e_2,f_2,\dots,e_m,f_m$ for $(V,\beta)$. Let $W_1=\langle f_1,e_2,f_2,\dots,e_m,f_m\rangle$, so that $(e_1,W_1)$ is a refined antiflag. Then $H_{\{e_1,W_1\}}=H_{e_1}\cap H_{W_1}$ stabilizes $e_1^\perp\cap W_1=\langle e_2,f_2,\dots,e_m,f_m\rangle$, and is faithful on $\langle e_2,f_2,\dots,e_m,f_m\rangle$. It follows that
\[
H_{\{e_1,W_1\}}=\Sp(e_1^\perp\cap W_1)=\Sp(\langle e_2,f_2,\dots,e_m,f_m\rangle)=\Sp_{2m-2}(q),
\]
and so $|\{\{e_1,W_1\}^g\mid g\in\Sp_{2m}(q)\}|=|\Sp_{2m}(q)|/|\Sp_{2m-2}(q)|=|\calP_{1,2m-1}^{(\infty)}|$. Thus, $\calP_{1,2m-1}^{(\infty)}=\{\{e_1,W_1\}^g\mid g\in\Sp_{2m}(q)\}$, and $H=\Sp_{2m}(q)$ is transitive on $\calP_{1,2m-1}^{(\infty)}$.

\begin{example}\label{ex:Sp<SL}
With the above notation, the refined-antiflag-transitivity of $\Sp_{2m}(q)$ gives
\[
\SL_{2m}(q)=L=HL_{\{e_1,W_1\}}=\Sp_{2m}(q)\,\SL_{2m-1}(q),
\]
where the intersection of the two factors is $H_{\{e_1,W_1\}}=\Sp_{2m-2}(q)$. For $(m,q)=(2,2)$, the group $\Sp_4(2)'\cong\A_6$ is also refined-antiflag-transitive. Consequently, the overgroup $\A_7$ of $\Sp_4(2)'\cong\A_6$ in $\SL_4(2)\cong\A_8$ is refined-antiflag-transitive as well. Thus we obtain the factorizations
\begin{align*}
\SL_4(2)&=\Sp_4(2)'\,\SL_3(2),\\
\SL_4(2)&=\A_7\,\SL_3(2).
\end{align*}
The intersection of the two factors is $3$ or $7{:}3$, respectively.
\end{example}

For $H=\G_2(q)<\Sp_6(q)$ with $q$ even, $H_{\{v,W\}}=\SL_2(q)$, see~\cite[4.3.6]{Wilson2009}. This gives rise to the factorizations in the next example.

\begin{example}\label{ex:G_2<Sp}
For even $q$, we have the factorization
\[
\SL_6(q)=\G_2(q)\,\SL_5(q),
\]
with the intersection $\G_2(q)\cap\SL_5(q)=\SL_2(q)$, as in row~6 of Table~\ref{TabLinear}.
\end{example}

\begin{remark}
We note that the commutator subgroup $\G_2(2)'=\PSU_3(3)$ is not refined-antiflag-transitive, verified by computation in \magma~\cite{BCP1997}.
\end{remark}

Identify the underlying set of the vector space $V=\bbF_q^{2m}$ with that of the vector space $V_\sharp=\bbF_{q^2}^m$. Let $v_1,\dots,v_m$ be a basis of $V_\sharp$, and $\lambda$ a generator of $\bbF_{q^2}^\times$. Then $v_1,\lambda v_1,\dots,v_m,\lambda v_m$ is a basis of $V$.
Let $W_1=\langle\lambda v_1,v_2,\lambda v_2,\dots,v_m,\lambda v_m\rangle$, a hyperplane of $V$. Then $(v_1,W_1)$ is a refined antiflag of $V$.
Let $L=\SL(V)=\SL_{2m}(q)$, and $S=\SL(V_\sharp)=\SL_m(q^2)$. Since $\langle v_2,\dots,v_m\rangle_{\bbF_{q^2}}$ is the only hyperplane of $V_\sharp$ that is contained in $W_1$, the stabilizer of $\{v_1,W_1\}$ in $S$ is
\begin{equation}\label{eq:SL-stab}
S_{\{v_1,W_1\}}=\SL(V_\sharp)_{v_1,W_1}=\SL(\langle v_2,\dots,v_m\rangle_{\bbF_{q^2}})=\SL_{m-1}(q^2).
\end{equation}
It follows that
\begin{equation}\label{eq:SL-stab-1}
|S|/|S_{(v_1,W_1)}|=q^{2m-2}(q^{2m}-1)=|\calP_{1,2m-1}^{(\infty)}|/q,
\end{equation}
and so $S$ is intransitive on $\calP_{1,2m-1}^{(\infty)}$. Moreover, if a group $H$ with $H^{(\infty)}=S$ is transitive on $\calP_{1,2m-1}^{(\infty)}$, then the normal subgroup $S$ has $q$ orbits on $\calP_{1,2m-1}^{(\infty)}$. We will show that, in the case where $q=2$ or $4$, some overgroups of $S$ are indeed transitive on $\calP_{1,2m-1}^{(\infty)}$.

Let $\psi\in\SiL(V_\sharp)$ be such that
\[
\psi\colon \ a_1v_1+\dots+a_mv_m\mapsto a_1^pv_1+\dots+a_m^pv_m
\]
for $a_1,\dots,a_m\in\bbF_{q^2}$, where $q=p^f$ with prime $p$. Then $|\psi|=2f$. In the following, we write $\omega=\{v_1,W_1\}$ for short.

\begin{lemma}\label{LemLinear6}
Let $V_\sharp$ and $\omega$ be as above, and let $q\in\{2,4\}$. Then
\[
\SiL(V_\sharp)_{\omega}=\SiL(V_\sharp)\cap\SiL(V)_{\omega}=\SL(V_\sharp)\cap\SL(V)_{\omega}=\SL(V_\sharp)_{\omega}=\SL_{m-1}(q^2).
\]
\end{lemma}

\begin{proof}
Let $H=\SiL(V_\sharp)=\SiL_m(q^2)$, so that $H=S{:}\langle\psi\rangle$ with $S=\SL(V_\sharp)$. Since $S<\SL(V)\leqslant\SiL(V)$, we have $S\cap \SiL(V)_\omega=S_\omega$. Thus $H_\omega=H\cap\SiL(V)_{\omega}\trianglerighteq S\cap\SiL(V)_\omega=S_\omega$. Suppose that $H_\omega>S_\omega$. Then $\psi^f s\in H_\omega$ for some $s\in S_{\omega}$. It follows that $v_1^s=v_1^{\psi^fs}=v_1$ and
\[
(\lambda v_1)^{\psi^fs}=(\lambda^{2^f}v_1)^s=\lambda^{2^f}v_1^s=\lambda^{2^f}v_1=\lambda^qv_1.
\]
Since $\bbF_{q^2}$ is a $2$-dimensional vector space over $\bbF_q$ with a basis $1,\lambda$, there exist $a,b\in\bbF_q$ such that  $\lambda^q=a\lambda+b$. As $\lambda$ is a generator of $\bbF_{q^2}$, neither $a$ nor $b$ is zero. Hence $(\lambda v_1)^{\psi^fs}=\lambda^qv_1=a\lambda v_1+bv_1\notin W_1$, contradicting the condition that $\psi^fs$ stabilizes $W_1$.
Thus $\SiL(V_\sharp)\cap\SiL(V)_{\omega}=H\cap\SiL(V)_{\omega}=H_\omega=S_{\omega}=\SL_{m-1}(q^2)$ by~\eqref{eq:SL-stab}.
\end{proof}

Let $\phi$ be the field automorphism of $L=\SL(V)=\SL_{2m}(q)$ of order $f$ that fixes each of $v_1,\lambda v_1,\dots,v_m,\lambda v_m$, and let $\gamma$ be the graph automorphism of $L=\SL(V)$ induced by the transpose of inverse with respect to the basis $v_1,\lambda v_1,\dots,v_m,\lambda v_m$. Then both $\phi$ and $\gamma$ fix the refined antiflag $\{v_1,W_1\}$, and $\gamma$ normalizes $\SL(V_\sharp)$.
Let $G=L{:}\langle\phi,\gamma\rangle$. Since $L$ is transitive on $\calP_{1,2m-1}^{(\infty)}$, it follows that $G$ is transitive on $\calP_{1,2m-1}^{(\infty)}$, with the stabilizer
\[
G_\omega=G_{\{v_1,W_1\}}=L_{\{v_1,W_1\}}{:}\langle\phi,\gamma\rangle.
\]

\begin{lemma}\label{LemLinear16}
Let $V_\sharp$, $\omega$, $\phi$ and $\gamma$ be as above, and let $\SL(V_\sharp)\trianglelefteq H<G=\SL(V){:}\langle\phi,\gamma\rangle$ with $q\in\{2,4\}$. Then $H$ is transitive on $\calP_{1,2m-1}^{(\infty)}$ if and only if $H$ contains $\psi$ or $\psi\gamma$. Moreover, if $H=\SL(V_\sharp){:}\langle\psi\rangle$ or $\SL(V_\sharp){:}\langle\psi\gamma\rangle$, then $H_\omega=\SL_{m-1}(q^2)$.
\end{lemma}

\begin{proof}
Let $H=\SL(V_\sharp){:}\langle\rho\rangle$, where $\rho=\psi$ or $\psi\gamma$. Then an element of $H$ has the form $s\rho^i=s\psi^i\gamma^{i'}$, where $s\in\SL(V_\sharp)$, and $i'=0$ or $i$. Moreover, an element of $G_\omega=L_{\{v_1,W_1\}}{:}\langle\phi,\gamma\rangle$ has the form $t\phi^j\gamma^k$ with $t\in L_\omega$. Thus an element $g$ of $H_\omega=H\cap G_\omega$ has the form
\begin{equation}\label{EqnLinear02}
g=s\psi^i\gamma^{i'}=t\phi^j\gamma^k,
\end{equation}
where $s\in \SL(V_\sharp)$, $t\in L_\omega$, and $i,j,k$ are suitable integers with $i'=0$ or $i$. It follows that
\[
\GaL(V)\gamma^{i'}=\GaL(V)s\psi^i\gamma^{i'}=\GaL(V)g=\GaL(V)t\phi^j\gamma^k=\GaL(V)\gamma^k.
\]
Hence $\gamma^{i'}=\gamma^k$, and so $s\psi^i=t\phi^j$, which lies in $\SiL(V_\sharp)\cap\SiL(V)_\omega=\SL(V_\sharp)_\omega=\SL_{m-1}(q^2)$ by Lemma~\ref{LemLinear6}. This implies that $i=j=0$, which then yields $\gamma^{i'}=1$ and so $\gamma^k=1$. Hence
\[
g=s=t\in\SL(V_\sharp)_\omega=\SL_{m-1}(q^2).
\]
Consequently, $H_\omega=\SL_{m-1}(q^2)$. Now
\[
|H|/|H_\omega|=|\SL_m(q^2){:}(2f)|/|\SL_{m-1}(q^2)|=2fq^{2m-2}(q^{2m}-1)=q^{2m-1}(q^{2m}-1)=|\calP_{1,2m-1}^{(\infty)}|.
\]
We conclude that $H$ is transitive on $\calP_{1,2m-1}^{(\infty)}$.

The above paragraph shows that a group $H$ with $\SL(V_\sharp)\trianglelefteq H<G=\SL(V){:}\langle\phi,\gamma\rangle$ is transitive on $\calP_{1,2m-1}^{(\infty)}$ if it contains $\psi$ or $\psi\gamma$. To prove the converse, we show that the group $H=\SL(V_\sharp){:}\langle\psi^2,\gamma\rangle$ is not transitive on $\calP_{1,2m-1}^{(\infty)}$. In fact, since $\gamma$ fixes the refined antiflag $\omega=\{v_1,W_1\}$, we have $H_\omega\geqslant\SL(V_\sharp)_\omega{:}\langle\gamma\rangle$. It follows that
\begin{align*}
|H|/|H_\omega|&\leqslant|\SL_m(q^2){:}(f\times2)|/|\SL_{m-1}(q^2){:}2|\\
&=fq^{2m-2}(q^{2m}-1)<q^{2m-1}(q^{2m}-1)=|\calP_{1,2m-1}^{(\infty)}|,
\end{align*}
and so $H$ is not transitive on $\calP_{1,2m-1}^{(\infty)}$.
\end{proof}

Lemma~\ref{LemLinear16} gives rise to factorizations in the next example.

\begin{example}\label{ex:SL_m<SL_2m}
Let $\psi$, $\phi$ and $\gamma$ be as above, and let $q\in\{2,4\}$. Then we have the factorizations
\begin{align*}
\SL_{2m}(q){:}\langle\phi\rangle&=(\SL_m(q^2){:}\langle\psi\rangle)\,(\SL_{2m-1}(q){:}\langle\phi\rangle),\\
\SL_{2m}(q){:}\langle\phi\gamma\rangle&=(\SL_m(q^2){:}\langle\psi\gamma\rangle)\,(\SL_{2m-1}(q){:}\langle\phi\gamma\rangle),
\end{align*}
where the intersection of the two factors is $\SL_{m-1}(q^2)$, as in rows~3--5 of Table~\ref{TabLinear}. In particular, both $\SL_m(q^2){:}\langle\psi\rangle$ and $\SL_m(q^2){:}\langle\psi\gamma\rangle$ are refined-antiflag-transitive.
\end{example}

By Lemma~\ref{LemXia10}, we obtain the subsequent example by combining the factorizations in Example~\ref{ex:SL_m<SL_2m} with those in Examples~\ref{ex:Sp<SL} and~\ref{ex:G_2<Sp}.

\begin{example}\label{ex:Sp<SL-02}
Let $\psi$, $\phi$ and $\gamma$ be as above, and let $q\in\{2,4\}$.
\begin{enumerate}[{\rm(a)}]
\item If $m$ is even, then for $\Sp_m(q^2)<\Sp_{2m}(q)$ we have the factorizations
\begin{align*}
\SL_{2m}(q){:}\langle\phi\rangle&=(\Sp_m(q^2){:}\langle\psi\rangle)\,(\SL_{2m-1}(q){:}\langle\phi\rangle),\\
\SL_{2m}(q){:}\langle\phi\gamma\rangle&=(\Sp_m(q^2){:}\langle\psi\gamma\rangle)\,(\SL_{2m-1}(q){:}\langle\phi\gamma\rangle),
\end{align*}
where the intersection of the two factors is $\Sp_{m-2}(q^2)$, as in rows~3--5 of Table~\ref{TabLinear}; in particular, both $\Sp_m(q^2){:}\langle\psi\rangle$ and $\Sp_m(q^2){:}\langle\psi\gamma\rangle$ are refined-antiflag-transitive.
\item If $m=6$, then for $\G_2(q^2)<\Sp_6(q^2)<\Sp_{12}(q)$ we have the factorizations
\begin{align*}
\SL_{12}(q){:}\langle\phi\rangle&=(\G_2(q^2){:}\langle\psi\rangle)\,(\SL_{11}(q){:}\langle\phi\rangle),\\
\SL_{12}(q){:}\langle\phi\gamma\rangle&=(\G_2(q^2){:}\langle\psi\gamma\rangle)\,(\SL_{11}(q){:}\langle\phi\gamma\rangle),
\end{align*}
where the intersection of the two factors is $\SL_2(q^2)$, as in rows~7--9 of Table~\ref{TabLinear}; in particular, both $\G_2(q^2){:}\langle\psi\rangle$ and $\G_2(q^2){:}\langle\psi\gamma\rangle$ are refined-antiflag-transitive.\qedhere
\end{enumerate}
\end{example}

Subgroups of $\GaL_n(q){:}\langle\gamma\rangle$ that are transitive on the set of antiflags of $\bbF_q^n$ are classified by Cameron and Kantor (see~\cite[Theorem~3.2]{LPS2010}). Now we state a classification of nonsolvable transitive groups on the set of refined antiflags, which turns out to coincide with that of the antiflag-transitive groups.

\begin{theorem}\label{refined-antiflags}
Let $H$ be a nonsolvable subgroup of $\GaL_n(q){:}\langle\gamma\rangle$. Then $H$ is refined-antiflag-transitive if and only if $H$ tightly contains $H_0^\alpha$ for one of the following $H_0$ as described in Examples~$\ref{ex:Sp<SL}$,~$\ref{ex:G_2<Sp}$,~$\ref{ex:SL_m<SL_2m}$,~$\ref{ex:Sp<SL-02}$ and some $\alpha\in\GaL_n(q){:}\langle\gamma\rangle$.
\begin{enumerate}[{\rm(a)}]
\item $H_0=\SL_n(q)$, $\Sp_n(q)'$, or $\G_2(q)$ with $q$ even and $n=6$;
\item $q=2^f\in\{2,4\}$, and $H_0=\SL_{n/2}(q^2){:}(2f)$, $\Sp_{n/2}(q^2){:}(2f)$, or $\G_2(q^2){:}(2f)$ with $n=12$;
\item $n=4$, $q=2$, and $H_0=\A_7$.
\end{enumerate}
Moreover, $H$ is refined-antiflag-transitive if and only if it is antiflag-transitive.
\end{theorem}

\begin{proof}
Let $\calP_{1,2m-1}$ be the set of antiflags and $\calP_{1,2m-1}^{(\infty)}$ the set of refined antiflags of $\bbF_q^n$. Suppose that $H$ is transitive on $\calP_{1,2m-1}^{(\infty)}$.
Since $\calP_{1,2m-1}$ is a $(\GaL_n(q){:}\langle\gamma\rangle)$-invariant partition of $\calP_{1,2m-1}^{(\infty)}$, it follows that $H$ is transitive on $\calP_{1,2m-1}$.
By~\cite[Theorem~3.2]{LPS2010}, $H$ satisfies one of the following:
\begin{enumerate}[{\rm(i)}]
\item $H^{(\infty)}=\SL_n(q)$, $\Sp_n(q)'$, or $\G_2(q)'$ with $n=6$ and $q$ even;
\item $q\in\{2,4\}$, and $H^{(\infty)}=\SL_{n/2}(q^2)$, $\Sp_{n/2}(q^2)$, or $\G_2(q^2)$ with $n=12$; moreover $H$ contains a full group of field automorphisms in each case;
\item $n=4$, $q=2$, and $H=\A_7$.
\end{enumerate}
For~(i) and~(iii), either $H$ tightly contains $H_0=H^{(\infty)}$, or $(n,q)=(6,2)$ and $H$ tightly contains $H_0=\G_2(2)$, as in part~(a) or~(c) of the lemma.
For~(ii), since $H$ contains a full group of field automorphisms, $H$ tightly contains the group $H_0$ in part~(b) of the lemma,

From the above paragraph we also see that every antiflag-transitive group $H$ satisfies one of~(i)--(iii) and so tightly contains $H_0$ in one of~(a)--(c). By Examples~\ref{ex:Sp<SL},~\ref{ex:G_2<Sp},~\ref{ex:SL_m<SL_2m} and~\ref{ex:Sp<SL-02}, if $H$ tightly contains such an $H_0$ then $H$ is antiflag-transitive. Thus a subgroup of $\GaL_n(q){:}\langle\gamma\rangle$ is refined-antiflag-transitive if and only if it is antiflag-transitive, and the theorem holds.
\end{proof}

\subsubsection{Other examples}

Now we describe the examples in the last two rows of Table~\ref{TabLinear}.

\begin{example}\label{ex:SL-exception1}
The group $G=\SL_2(9)$ has precisely two conjugacy classes of subgroups isomorphic to $\SL_2(5)=2^{\boldsymbol{\cdot}}\A_5$. Let $H$ and $K$ be two subgroups from these two classes respectively. Then $G=HK$ with $H\cap K=2^{\boldsymbol{\cdot}}\D_{10}$.
\end{example}

The following example is known in~\cite[Theorem~3.3]{LX} and verified in \magma~\cite{BCP1997}.

\begin{example}\label{ex:SL-exception2}
The group $G=\SL_3(4){:}\langle\phi\gamma\rangle$ has three conjugacy classes of subgroups isomorphic to $\SL_3(2){:}2$, and let $K$ be such a subgroup. Then there are precisely two conjugacy classes of maximal subgroups $H$ of $G$ isomorphic to $3^{\boldsymbol{\cdot}}\M_{10}$ such that $G=HK$. For each such pair $(H,K)$ we have $H\cap K=\Sy_3$.
\end{example}

\subsection{Classifying the factorizations}\label{SecLinearProof}

Let $G$ be a group with $G^{(\infty)}=L=\SL_n(q)$, where $n\geqslant2$ and $(n,q)\neq(2,2)$ or $(2,3)$.
Suppose $G=HK$ such that $H$ and $K$ are nonsolvable and neither $H$ nor $K$ contains $L$.
We first treat a few small groups by computation in \magma~\cite{BCP1997}, which confirms Theorem~\ref{ThmLinear} for these groups.

\begin{lemma}\label{ProofLinear-0}
For $L=\SL_4(2)$, $\SL_4(3)$ or $\SL_6(3)$, $G=HK$ if and only if $(G,H,K)$ tightly contains some $(G_0,H_0,K_0)$ in rows~\emph{1--4} and~\emph{10--12} Table~$\ref{TabLinear}$.
\end{lemma}

Next we embark on the case where $L$ is not in Lemma~\ref{ProofLinear-0}. If $n$ is prime, then the triple $(G,H,K)$ has been classified in \cite[Theorem~3.3]{LX}, which implies that Theorem~\ref{ThmLinear} holds for this case. We thus assume that $n$ is not prime, so that $n\geqslant4$. Take $\max^-$ subgroups $A$ and $B$ of $G$ containing $H$ and $K$ respectively.
%From~\cite[Theorem~A]{LPS1990} and~\cite{LPS1996} we obtain the ensuing statement.
%
%\begin{proposition}\label{Prop:Linear-max}
%With the above notation and interchanging $A$ and $B$ if necessary, one of the following holds:
%\begin{enumerate}[{\rm(a)}]
%\item $B\trianglerighteq q^{n-1}{:}\SL_{n-1}(q)$ stabilizes a $1$-space or a hyperplane, and either $A^{(\infty)}=\SL_a(q^b)$ with $ab=n$ and $b$ prime or $A^{(\infty)}=\Sp_n(q)$ with $n\geqslant4$ even;
%\item $B\trianglerighteq q^{n-1}{:}\SL_{n-1}(q)$ stabilizes an antiflag, $n\geqslant4$ is even, and $A^{(\infty)}=\SL_{n/2}(q^2)$ or $\Sp_n(q)$.
%\end{enumerate}
%\end{proposition}
The next lemma determines $K^{(\infty)}$.

\begin{lemma}\label{ProofLinear-2}
$K^{(\infty)}=q^{n-1}{:}\SL_{n-1}(q)$ or $\SL_{n-1}(q)$.
\end{lemma}

\begin{proof}
Let $N=\Rad(B)$. Then $B/N$ is an almost simple group with socle $\PSL_{n-1}(q)$, and $KN/N$ is nonsolvable as $K$ is nonsolvable.
Viewing the irreducible action of $\SL_{n-1}(q)$ on $q^{n-1}$ in $B^{(\infty)}=q^{n-1}{:}\SL_{n-1}(q)$, we only need to prove
\[
KN/N\geqslant\Soc(B/N)=\PSL_{n-1}(q).
\]
Suppose for a contradiction that $KN/N\ngeqslant\Soc(B/N)$. From~\cite[Theorem~A]{LPS1990} and~\cite{LPS1996} we see that $|L|/|A\cap L|$ is divisible by some primitive prime divisor $r$ of $q^{n-1}-1$.
By Lemmas~\ref{LemXia6} and~\ref{LemXia22}, it follows that $|KN/N|$ is divisible by $r$. Suppose further that $n-1$ is prime. Then an inspection of the maximal factorizations of almost simple linear groups in~\cite[Theorem~A]{LPS1990} and~\cite{LPS1996} shows that a group with socle $\PSL_{n-1}(q)$ does not have a factor of order divisible by $r$ in any core-free factorization. This is a contradiction since $KN/N$ should be such a supplement of $(A\cap B)N/N$ in $B/N$. Therefore, $n-1$ is not prime, and as $n$ is not prime either, $n\geqslant9$.
% If $A\cap L$ is in Row~1 of Table~\ref{TabMaxLinear} with $b\geqslant3$, then take $s\in\ppd(q^{n-2}-1)$. If $A\cap L$ is in Row~1 of Table~\ref{TabMaxLinear} with $b=2$ or in Rows~2--4 of Table~\ref{TabMaxLinear}, then $n$ is even and we take $s\in\ppd(q^{n-3}-1)$.
Take $s\in\ppd(q^{n-2}-1)$ if $A^{(\infty)}=\SL_a(q^b)$ with $b\geqslant3$, and take $s\in\ppd(q^{n-3}-1)$ otherwise.
Then $|L|/|A\cap L|$ is divisible by $s$, and so is $|KN/N|$.
However, by~\cite[Theorem~A]{LPS1990} and~\cite{LPS1996}, there is no core-free factor of $B/N$ with order divisible by both $r$ and $s$, which is a contradiction.
\end{proof}

We are now ready to prove Theorem~\ref{ThmLinear}.

\begin{proof}[Proof of Theorem~$\ref{ThmLinear}$]
For convenience, we adopt the notation defined in Lemma~\ref{LemLinear16}. By our assumption before Lemma~\ref{ProofLinear-2}, $n$ is not a prime. If $L=\SL_4(2)$, $\SL_4(3)$ or $\SL_6(3)$, then Theorem~$\ref{ThmLinear}$ is proved by Lemma~\ref{ProofLinear-0}. We thus only need to treat the other cases.
By Lemma~\ref{ProofLinear-2}, either the factor $K\trianglerighteq q^{n-1}{:}\SL_{n-1}(q)$ stabilizes a $1$-space or a hyperplane, or $K\trianglerighteq\SL_{n-1}(q)$ stabilizes an antiflag.

For the former case, conjugating by the transpose-inverse automorphism if necessary, we may assume $K\trianglerighteq L_v$ for some nonzero vector $v\in V$. By~\cite{Hering1985} the triple $(G_0,H_0,K_0)$ is as listed in row~1 of Table~\ref{TabLinear}.

Next, assume that $K\trianglerighteq\SL_{n-1}(q)$ stabilizes an antiflag. Then $H$ is antiflag-transitive. By Theorem~\ref{refined-antiflags}, $H$ is refined-antiflag-transitive, and $H_0$ is explicitly listed in parts~(a)--(c) therein. If $H_0$ is one of the candidates given in part~(a) or (c) of Theorem~\ref{refined-antiflags}, then $G_0=\SL_{2m}(q)$ and $K_0=\SL_{2m-1}(q)$, which is listed in row~2 or~6 of Table~\ref{TabLinear}. Now assume that $H_0$ is as in part~(b) of Theorem~\ref{refined-antiflags}. Then $q\in\{2,4\}$, and by Lemma~\ref{LemLinear16} we have $H_0=H^{(\infty)}{:}\langle\psi\rangle$ or $H^{(\infty)}{:}\langle\psi\gamma\rangle$.
It follows that $G_0=G^{(\infty)}{:}\langle\rho\rangle$ with $\rho\in\{\phi,\phi\gamma\}$ and $K_0=K^{(\infty)}{:}\langle\rho\rangle$.
Thus $(G_0,H_0,K_0)$ lies in rows~3--5 and~7--9 of Table~\ref{TabLinear}.
\end{proof}

%%%%%%%%%%%%%%%%%%%%%%%%%%%%%%%%%%%%%%%%%%%%%%%%%%%%%%%%%%%%%%%%%%%
%%%%%%%%%%%%%%%%%%%%%%%%%%%%%%%%%%%%%%%%%%%%%%%%%%%%%%%%%%%%%%%%%%%

\section{Unitary groups}

%%%%%%%%%%%%%%%%%%%%%%%%%%%%%%%%%%%%%%%%%%%%%%%%%%%%%%%%%%%%%%%%%%%
%%%%%%%%%%%%%%%%%%%%%%%%%%%%%%%%%%%%%%%%%%%%%%%%%%%%%%%%%%%%%%%%%%%

In this section we classify the factorizations of unitary groups. Throughout this section, let $G$ be a classical group with $L=G^{(\infty)}=\SU_n(q)$, where $n\geqslant3$ and $(n,q)\neq(3,2)$, and let $H$ and $K$ be nonsolvable subgroups of $G$ not containing $G^{(\infty)}$.

\begin{theorem}\label{ThmUnitary}
With the above notation, $G=HK$ if and only if either $(G,H,K)$ tightly contains some triple $(G_0,H_0,K_0)$ in Table~$\ref{TabUnitary}$ or $(\overline{G},\overline{H},\overline{K})$ tightly contains some $(G_0,H_0,K_0)$ in Table~$\ref{TabUnitary-2}$, where $\,\overline{\phantom{\varphi}}\,$ is the quotient modulo scalars.
\end{theorem}

Some notation for the whole section is set up in the following Subsection~\ref{SecUnitaryPre}, and the proof of Theorem~\ref{ThmUnitary} will be given in Subsection~\ref{SecUnitaryProof}.

\subsection{Setup}\label{SecUnitaryPre}
\ \vspace{1mm}

In order to determine the pairs $(H,K)$ such that $G=HK$, we let
\[
H\leqslant A\quad\text{and}\quad K\leqslant B
\]
for some $\max^-$ subgroups $A$ and $B$ of $G$. Since $H$ and $K$ are nonsolvable, so are $A$ and $B$.
The $\max^-$ factorizations $G=AB$ are classified by Liebeck-Praeger-Saxl~\cite{LPS1990,LPS1996}, from which we deduce the consequence as follows.

\begin{proposition}\label{Prop:Unitary-max}
If $G=HK$, then with the above notation and interchanging $A$ and $B$ if necessary, one of the following holds:
\begin{enumerate}[{\rm(a)}]
\item $n=2m$ is even with $(m,q)\neq(2,2)$, $B=\N_1[G]$, and $A$ satisfies one of the following:
\begin{itemize}
\item $A=\Pa_m[G]$;
\item $A\cap L=\SL_m(q^2).(q-1).2$ with $q\in\{2,4\}$;
\item $A^{(\infty)}=\Sp_{2m}(q)$;
\item $A^{(\infty)}\in\{3.\PSU_4(3),3.\M_{22}\}$ with $L=\SU_6(2)$;
\item $A^{(\infty)}=3.\Suz$ with $L=\SU_{12}(2)$.
\end{itemize}
\item $L=\SU_4(3)$ or $\SU_9(2)$.
\end{enumerate}
\end{proposition}

For $L\in\{\SU_4(3),\SU_9(2)\}$, computation in \magma~\cite{BCP1997} produces an explicit list of the factorizations, described in the following lemma, where we also include the small group $L=\SU_4(5)$ for later convenience and $\,\overline{\phantom{\varphi}}\,$ denotes the quotient modulo scalars.

\begin{lemma}\label{Lem:UnitaryProof}
For $L=\SU_4(3)$, $\SU_4(5)$ or $\SU_9(2)$, Theorem~$\ref{ThmUnitary}$ holds subject to the following refinements.
\begin{enumerate}[{\rm(a)}]
\item In row~\emph{2} of Table~$\ref{TabUnitary-2}$, the group $G_0$ has a maximal subgroup of the form $3^4{:}(\A_6\times2)$, and $H_0$ is contained in this maximal subgroup.
\item In row~\emph{3} of Table~$\ref{TabUnitary-2}$, the group $G_0$ has a maximal subgroup of the form $\PSL_3(4){:}2$, and $H_0=\Pa_2[G_0]=3^4{:}(\A_6\times2)$ or $3^4{:}\M_{10}$ such that if $H_0=3^4{:}\M_{10}$ then $K_0$ is from one out of two conjugacy classes of maximal subgroups of $G_0$ of the form $\PSL_3(4){:}2$.
\item In row~\emph{5} of Table~$\ref{TabUnitary-2}$, the group $G_0$ has no maximal subgroup of the form $\A_7.2$, the group $H_0=\Pa_2[G_0]$, and $K_0$ is contained in $\N_1[G_0]$.
\end{enumerate}
\end{lemma}

By virtue of Lemma~\ref{Lem:UnitaryProof} we only need to consider case~(a) of Proposition~\ref{Prop:Unitary-max}. Thus we write $n=2m$ with $m\geqslant2$ in the following. Let $q=p^f$ with $p$ prime, and let $V=\bbF_{q^2}^{2m}$ equipped with a nondegenerate Hermitian form $\beta$. Take a standard basis
\[
e_1,f_1,\dots,e_m,f_m
\]
of $V$ with respect to $\beta$, namely, $\beta(e_i,f_j)=\delta_{ij}$ and $\beta(e_i,e_j)=\beta(f_i,f_j)=0$. Let
\[
U=\langle e_1,e_2,\dots,e_m\rangle_{\bbF_{q^2}}\ \text{ and }\ W=\langle f_1,f_2,\dots,f_m\rangle_{\bbF_{q^2}}.
\]
From~\cite[3.6.2]{Wilson2009} we know that $\SU(V)_U=\Pa_m[\SU(V)]$ has solvable residual $R{:}T$ with
\[
R=q^{m^2}\ \text{ and }\ T=\SL_m(q^2),
\]
where $R$ is the kernel of $\SU(V)_U$ acting on $U$,
% By \cite[Theorem~3.2]{Bell1978},
% ${\rm H}^1(\SL_m(q^2),\bbF_q^{m^2})=0$,
% and hence all subgroups of $\SU(V)_U$ isomorphic to $\SL_m(q^2)$ are conjugate.
and $T$ stabilizes both $U$ and $W$ (the action of $T$ on $U$ determines that on $W$, see~\cite[Lemma~2.2.17]{BG2016}).

Let $\calN_1^{(\infty)}=\{u\in V\mid \beta(u,u)=1\}$, take $\lambda\in\bbF_{q^2}$ with $\lambda+\lambda^q=1$ (note that such $\lambda$ exists as the trace of the field extension $\bbF_{q^2}/\bbF_q$ is surjective), and let
\[
v=e_1+\lambda f_1.
\]
Then $\beta(v,v)=\lambda+\lambda^q=1$, and hence $v\in\calN_1^{(\infty)}$. Since $L$ is transitive on $\calN_1^{(\infty)}$ with stabilizer $L_v=\SU(v^\perp)=\SU_{2m-1}(q)$, we can identify $\calN_1^{(\infty)}$ with $L_v\backslash L$. In particular,
\[
|\calN_1^{(\infty)}|=|\SU_{2m}(q)|/|\SU_{2m-1}(q)|=q^{2m-1}(q^{2m}-1).
\]

\subsection{Factorizations with $A=\Pa_m[G]$ and $K^{(\infty)}=\SU_{2m-1}(q)$}\label{SecUnitaryPm}
\ \vspace{1mm}

In this subsection, we classify the factorizations $G=HK$ with $H\leqslant A=\Pa_m[G]$ and $K^{(\infty)}=B^{(\infty)}=\SU_{2m-1}(q)$.
Without loss of generality, let $A=G_U$, where $U=\langle e_1,\dots,e_m\rangle_{\bbF_{q^2}}=\bbF_{q^2}^m$.
For an $\bbF_{p^i}\GL_\ell(p^i)$-module $M$, recall the $\bbF_{p^i}\GL_\ell(p^i)$-module $M^{(p^j)}$ obtained from $M$ twisted by the Frobenius map taking $p^j$th power, as at the end of Subsection~\ref{SecXia2}.
If $m=ab$, then let $U_{(b)}=\bbF_{q^{2b}}^a$ with the same underlying set as $U$, and for $i\in\{1,\dots,\lceil b/2\rceil\}$ denote
\begin{equation}\label{EqnUnitary1}
U_{(b)}^\sharp(i)=
\left\{\begin{array}{ll}
\bigoplus\limits_{j=0}^{2b-1}\left(U_{(b)}\otimes U_{(b)}^{(q^{2i-1})}\right)^{(q^j)}&\text{if }1\leqslant i\leqslant\lfloor b/2\rfloor\vspace{2mm}\\
\bigoplus\limits_{j=0}^{b-1}\left(U_{(b)}\otimes U_{(b)}^{(q^{2i-1})}\right)^{(q^j)}&\text{if }\lfloor b/2\rfloor<i\leqslant\lceil b/2\rceil.
\end{array}\right.
\end{equation}
(Note that $\lfloor b/2\rfloor<i\leqslant\lceil b/2\rceil$ if and only if $i=(b+1)/2$ with $b$ odd.) The $\bbF_{q^{2b}}\GL(U_{(b)})$-module $U_{(b)}^\sharp(i)$ can be realized over $\bbF_q$, and let $U_{(b)}(i)$ denote the realized $\bbF_q\GL(U_{(b)})$-module. Since $U_{(b)}^\sharp(i)=U_{(b)}(i)\otimes\bbF_{q^{2b}}$, the definition of $U_{(b)}^\sharp(i)$ in~\eqref{EqnUnitary1} implies (see \cite[\S26)]{ASCH1993} for instance) that the $\bbF_qS$-module $U_{(b)}(i)$ is isomorphic to $U_{(b)}\otimes U_{(b)}^{(q^{2i-1})}$, where the latter is viewed as an $\bbF_qS$-module.

\begin{lemma}\label{LemUnitaryPm1}
Let $m=ab$, and let $S$ be a subgroup of $\SL(U)$ contained in the group $\SL_a(q^{2b})$ defined over $\bbF_{q^{2b}}$ such that $S=\SL_a(q^{2b})$, $\Sp_a(q^{2b})$ or $\G_2(q^{2b})$ (with $a=6$ and $q$ even). Then the following statements hold.
\begin{enumerate}[{\rm (a)}]
\item For $i\in\{1,\dots,\lceil b/2\rceil\}$, the group $U_{(b)}(i)$ is $q^{2a^2b}$ if $i\neq(b+1)/2$, and is $q^{a^2b}$ if $i=(b+1)/2$.
\item The $\bbF_qS$-module $R$ is the direct sum of pairwise non-isomorphic irreducible submodules
\begin{equation}\label{EqnUnitary2}
U_{(b)}(1),\,\dots,\,U_{(b)}(\lceil b/2\rceil).
\end{equation}
\item The modules in~\eqref{EqnUnitary2} are the irreducible $\bbF_qS$-submodules of $R$.
\item If $H\leqslant G_U=\Pa_m[G]$ with $RH=RS$, then $H=(H\cap R){:}X$ for some $X\leqslant H$ with $X\cong S$.
\end{enumerate}
\end{lemma}

\begin{proof}
Part~(a) follows immediately from the definition of $U_{(b)}^\sharp(i)$ and $U_{(b)}(i)$. Viewing $R$ as an $\bbF_q\GL(U)$-module, we may write $R=U\otimes U^{(q)}$.
%(see~\cite[Page~16]{BL2021})
Since $U_{(b)}$, as an $\bbF_{q^2}S$-module, can be viewed as $U$, we have
\[
\mbox{$U\otimes\bbF_{q^{2b}}=\bigoplus\limits_{i=0}^{b-1}U_{(b)}^{(q^{2i})}$.}
\]
Then straightforward calculation gives
\[\begin{array}{rl}
R\otimes\bbF_{q^{2b}}\,
=&(U\otimes\bbF_{q^{2b}})\otimes(U\otimes\bbF_{q^{2b}})^{(q)}\vspace{2mm}\\
=&\left(\bigoplus\limits_{i=0}^{b-1}U_{(b)}^{(q^{2i})}\right)
\otimes\left(\bigoplus\limits_{i=0}^{b-1}U_{(b)}^{(q^{2i})}\right)^{(q)}\vspace{2mm}\\
=&\bigoplus\limits_{i=0}^{b-1}\bigoplus\limits_{j=0}^{b-1}\left(U_{(b)}\otimes U_{(b)}^{(q^{2i+1})}\right)^{(q^{2j})}\vspace{2mm}\\
%=&\left(\bigoplus\limits_{i=1}^{\lfloor b/2\rfloor}\bigoplus\limits_{j=0}^{2b-1}\left(U_{(b)}\otimes U_{(b)}^{(q^{2i-1})}\right)^{(q^j)}\right)\oplus\left(\bigoplus\limits_{i=\lfloor b/2\rfloor+1}^{\lceil b/2\rceil}\,\bigoplus\limits_{j=0}^{b-1}\left(U_{(b)}\otimes U_{(b)}^{(q^{2i-1})}\right)^{(q^j)}\right).
=&\bigoplus\limits_{i=1}^{\lceil b/2\rceil}U_{(b)}^\sharp(i).
\end{array}\]
The last line realized over $\bbF_q$ shows that the $\bbF_qS$-module $R$ is the direct sum of $U_{(b)}(i)$ with $i$ running over $\{1,\dots,\lceil b/2\rceil\}$.
%\begin{equation}\label{EqnUnitary3}
%R\downarrow S=\bigoplus\limits_{i=1}^{\lceil b/2\rceil}U_{(b)}(i),
%\end{equation}

Since $U_{(b)}\otimes U_{(b)}^{(q^{2i-1})}$ is an irreducible $\bbF_{q^{2b}}S$-module with highest weight $(1+q^{2i-1})\lambda_1$, where $\lambda_1,\dots,\lambda_{a-1}$ are the fundamental dominant weights of $S$, we derive that the $\bbF_{q^{2b}}S$-modules
$U_{(b)}\otimes U_{(b)}^{(q^{2i-1})}$ with $i$ running over $\{1,\dots,\lceil b/2\rceil\}$ are pairwise non-isomorphic (see \cite[\S5.4)]{KL1990} for instance). Note that they are also irreducible as $\bbF_qS$-modules. It follows from part~(a) that $U_{(b)}(1),\dots,U_{(b)}(\lceil b/2\rceil)$ are pairwise non-isomorphic irreducible $\bbF_qS$-modules. This completes the proof of part~(b).

Part~(c) is a consequence of part~(b). To prove part~(d), suppose that $H$ is a subgroup of $G_U$ with $RH=RS=R{:}S$. Write $Q=H\cap R$. Then $H/Q\cong(RH)/R\cong S$ and $Q=\mathbf{O}_p(H)$. By parts~(b) and~(c), we have an $\bbF_q(RS)$-module decomposition $R=P\oplus Q$ for some $P\leqslant R$. In particular, $P$ is normal in $RS$. Since $PH=PQH=RH=RS$ and
\[
|P||H|=|P||Q||S|=|R||S|=|RS|,
\]
it follows that $RS=P{:}H$. Moreover, $RS=PQS$ with $|P||QS|=|P||Q||S|=|RS|$, and so $RS=P{:}(QS)$. This in conjunction with $RS=P{:}H$ leads to $H\cong(RS)/P\cong QS=Q{:}S$, which implies that $H=\mathbf{O}_p(H){:}X$ for some $X\leqslant H$. Thus part~(d) holds.
\end{proof}

For a subset $I=\{i_1,\ldots,i_k\}$ of $\{1,\ldots,\lceil b/2\rceil\}$, denote
\begin{equation}\label{EqnUnitary5}
\gcd(2I-1,b)=\gcd(2i_1-1,\dots,2i_k-1,b)
\end{equation}
and let $U_{(b)}(I)=U_{(b)}(i_1)\cdots U_{(b)}(i_k)\leqslant R$. Then we derive from Lemma~\ref{LemUnitaryPm1} that
\begin{equation}\label{EqnUnitary3}
U_{(b)}(I)=U_{(b)}(i_1)\times\cdots\times U_{(b)}(i_k)=
\begin{cases}
q^{(2k-1)a^2b}&\textup{if }(b+1)/2\in I\\
q^{2ka^2b}&\textup{if }(b+1)/2\notin I.
\end{cases}
\end{equation}
For $H\leqslant G_U=\Pa_m[G]$, the notation $H^U$, as usual, denotes the induced group of $H$ on $U$.

\begin{lemma}\label{LemUnitaryPm2}
Let $K=\N_1[G]$, let $m=ab$, and let $H=U_{(b)}(I).S\leqslant\Pa_m[G]$ with $I\subseteq\{1,\ldots,\lceil b/2\rceil\}$ and $S\leqslant\GaL_a(q^{2b})$ defined over $\bbF_{q^{2b}}$ such that $S^{(\infty)}=\SL_a(q^{2b})$, $\Sp_a(q^{2b})$ or $\G_2(q^{2b})$ (with $a=6$ and $q$ even). If $G=HK$, then $\gcd(2I-1,b)=1$. In particular, if $G=HK$ with $U_{(b)}(I)=1$, then $b=1$.
\end{lemma}

\begin{proof}
Let $I=\{i_1,\ldots,i_k\}$ and $\ell=\gcd(2I-1,b)=\gcd(2i_1-1,\dots,2i_k-1,b)$. Suppose for a contradiction that $G=HK$ with $\ell>1$. Write $b=\ell e$ and $2i_s-1=\ell(2j_s-1)$ for $s\in\{1,\dots,k\}$. Let $M=\Nor_G(\SU_{2ae}(q^\ell))$ be a field-extension subgroup of $G$ over $\bbF_{q^{2\ell}}$. We shall prove that $H$ is contained in $M$ up to conjugation in $G$, which will lead to a contradiction $G=MK$ to Proposition~\ref{Prop:Unitary-max}.

Let $P=q^{\ell(ae)^2}$ be the unipotent radical of $\Pa_{ae}[M]$. Applying Lemma~\ref{LemUnitaryPm1} to the $\bbF_{q^\ell}S^{(\infty)}$ module $P$ we see that $P$ is a direct sum of $P(1),\dots,P(\lceil e/2\rceil)$, where $P(t)$ is the realization of the $\bbF_{q^{2b}}S^{(\infty)}$-module
\[
P^\sharp(t):=
\left\{\begin{array}{ll}
\bigoplus\limits_{s=0}^{2e-1}\left(U_{(b)}\otimes U_{(b)}^{(q^{\ell(2t-1)})}\right)^{(q^{\ell s})}&\text{if }1\leqslant t\leqslant\lfloor e/2\rfloor\vspace{2mm}\\
\bigoplus\limits_{s=0}^{e-1}\left(U_{(b)}\otimes U_{(b)}^{(q^{\ell(2t-1)})}\right)^{(q^{\ell s})}&\text{if }\lfloor e/2\rfloor<t\leqslant\lceil e/2\rceil
\end{array}\right.
\]
over $\bbF_{q^\ell}$. Recall the definition of $U_{(b)}^\sharp(i)$ in~\eqref{EqnUnitary1}. We obtain for $t\in\{1,\dots,\lceil e/2\rceil\}$ that
\[
\mbox{$\bigoplus\limits_{r=0}^{\ell-1}\bigoplus\limits_{s=0}^{2e-1}P(t)^{(q^{\ell s+r})}
=\bigoplus\limits_{r=0}^{\ell-1}\left(\bigoplus\limits_{s=0}^{2e-1}P(t)^{(q^{\ell s})}\right)^{(q^r)}
=\bigoplus\limits_{r=0}^{\ell-1}P^\sharp(t)^{(q^r)}=U_{(b)}^\sharp\Big((2\ell t-\ell+1)/2\Big)$,}
\]
and so $P(t)$ equals $U_{(b)}\big((2\ell t-\ell+1)/2\big)$, the realization of $U_{(b)}^\sharp\big((2\ell t-\ell+1)/2\big)$ over $\bbF_q$.
Since $2i_s-1=\ell(2j_s-1)$ for $s\in\{1,\dots,k\}$, it follows that, up to conjugation in $G$,
\[
U_{(b)}(I)=U_{(b)}(i_1)\cdots U_{(b)}(i_k)=P(j_1)\cdots P(j_k)\leqslant P.
\]
This implies that $H\leqslant\Pa_{ae}[M]$, as desired.
\end{proof}

Recall~\eqref{EqnUnitary5} and~\eqref{EqnUnitary3} for the definition of $\gcd(2I-1,b)$ and $U_{(b)}(I)$. Note that, for a subset $I$ of $\{1,\ldots,\lceil b/2\rceil\}$, we have $U_{(b)}(I)=q^c$ with
\begin{equation}\label{EqnUnitary7}
c=
\begin{cases}
(2|I|-1)a^2b&\textup{if }(b+1)/2\in I\\
2|I|a^2b&\textup{if }(b+1)/2\notin I.
\end{cases}
\end{equation}
Clearly, $0\leqslant c\leqslant m^2$, and $c=m^2=a^2b^2$ if and only if $I=\{1,\ldots,\lceil b/2\rceil\}$ (this can be either seen from Lemma~\ref{LemUnitaryPm1}(b) and the definition of $U_{(b)}(I)$, or confirmed by~\eqref{EqnUnitary7} as $|I|=\lceil b/2\rceil$ if and only if $I=\{1,\ldots,\lceil b/2\rceil\}$). In the following lemma we construct minimal factorizations with $H=U_{(b)}(I){:}S\leqslant\Pa_m[G]$ for nonempty $I$ (so that $c>0$). Minimal factorizations for $I=\emptyset$ (or equivalently, $c=0$) are included in Example~\ref{ex:Unitary-03}.
The special case $b=1$ in Lemma~\ref{LemUnitaryPm3} gives the factorizations for $c=m^2$, which are specifically displayed in the remark after the lemma.

\begin{lemma}\label{LemUnitaryPm3}
Let $G=\SU_{2m}(q)$, let $K=\N_1[G]^{(\infty)}=\SU_{2m-1}(q)$, let $m=ab$, and let $H=U_{(b)}(I){:}S=q^c{:}S\leqslant G_U=\Pa_m[G]$ with $\emptyset\neq I\subseteq\{1,\ldots,\lceil b/2\rceil\}$ such that $\gcd(2I-1,b)=1$ and $S\cong S^U\leqslant\SL_a(q^{2b})$ is defined over $\bbF_{q^{2b}}$.
\begin{enumerate}[{\rm (a)}]
\item If $S=\SL_a(q^{2b})$, then $G=HK$ with $H\cap K=[q^{c-2b+1}].\SL_{a-1}(q^{2b})$.
\item If $S=\Sp_a(q^{2b})$, then $G=HK$ with $H\cap K=[q^{c-2b+1}].\Sp_{a-2}(q^{2b})$.
\item If $S=\G_2(q^{2b})$ with $a=6$ and $q$ even, then $G=HK$ with $H\cap K=[q^{c-2b+1}].\SL_2(q^{2b})$.
\end{enumerate}
\end{lemma}

\begin{proof}
Let $I\cap\{1,\dots,\lfloor b/2\rfloor\}=\{i_1,\ldots,i_k\}$, let $d=\gcd(2i_1-1,\dots,2i_k-1,a)$, and let $e$ be the largest divisor of $a$ coprime to $d$. As $\gcd(2i_1-1,\dots,2i_k-1,b)=\gcd(2I-1,b)=1$, we have
\[
\gcd(2be+2i_1-1,\dots,2be+2i_k-1,m)=1.
\]
Let $C=\langle s\rangle=q^{2m}-1$ be a Singer cycle in $\GL(U)$, and let $P=U_{(b)}(I)$ and $X=PC=P{:}C$.
For each $t\in\{1,\dots,k\}$, by~\cite[Corollary~4.5]{FLWXZ} (with $q$ and $r$ replaced by $q^b$ and $q^{2b}$ respectively), the restriction of $U_{(b)}(i_t)$ to $C$ can be decomposed as
\[
U_{(b)}(i_t,1)\oplus\cdots\oplus U_{(b)}(i_t,a)
\]
such that the value of $\chi_{(b)}(i_t,j)$ at $s$ is
\[
\sum^{2m-1}_{\ell=0}(\omega^{q^{2bj+2i_t-1}+1})^{q^\ell},
\]
where $\chi_{(b)}(i_t,j)$ is the character of the $\bbF_qC$-module $U_{(b)}(i_t,j)$ and $\omega$ is a generator of $\bbF_{q^{2m}}^\times$.
In particular, this character value is completely determined by $2bj+2i_t-1$ up to the action of ${\rm Gal}(\bbF_{q^{2m}}/\bbF_q)$.
By~\cite[Theorem~1,8]{FLWXZ} (our module $U_{(b)}(i_t,j)$ is denoted as $U(bj+i_t)$ in~\cite{FLWXZ}), we have
$\GU_{2m}(q)=(U_{(b)}(i_1,e)\cdots U_{(b)}(i_k,e)C)\GU_{2m}(q)_v$, where
\[
\GU_{2m}(q)_v=\GU_{2m-1}(q)
\]
is the stabilizer of $v\in\calN_1^{(\infty)}$ in $\GU_{2m}(q)$.
%for each $s\in\{1,\dots,a\}$,
%\[
%\mbox{$U_{(b)}(i_t)_s\otimes\bbF_{q^{2m}}=\bigoplus\limits_{j=0}^{a-1}\left(U_{(m)}\otimes U_{(m)}^{(q^{2bs+2i_t-1})}\right)^{(q^{2bj})}$.}
%\]
Since
\[
U_{(b)}(i_1,e)\cdots U_{(b)}(i_k,e)\leqslant U_{(b)}(i_1)\cdots U_{(b)}(i_k)\leqslant P,
\]
it follows that $\GU_{2m}(q)=(PC)\GU_{2m}(q)_v=X\,\GU_{2m}(q)_v$, and so the action of $X$ on $\calN_1^{(\infty)}$ is transitive. Since $X=P{:}C$ and $C=q^{2m}-1$ is a Singer cycle in $\GL(U)$, this implies that the action of $\GL(U)$ on the set $\big(\calN_1^{(\infty)}\big)_P$ of $P$-orbits on $\calN_1^{(\infty)}$ is permutationally equivalent to its natural action on $\bbF_{q^2}^m\Setminus\{0\}$. Then we conclude from Theorem~\ref{ThmLinear} that the action of $S$ on $\big(\calN_1^{(\infty)}\big)_P$ is transitive with point stabilizer $q^{2b(a-1)}{:}\SL_{a-1}(q^{2b})$, $[q^{2b(a-1)}]{:}\Sp_{a-2}(q^{2b})$ or $q^{4b+6b}{:}\SL_2(q^{2b})$ according to $S=\SL_a(q^{2b})$, $\Sp_a(q^{2b})$ or $\G_2(q^{2b})$, respectively. Hence $H=P{:}S$ is transitive on $\calN_1^{(\infty)}$, and so the stabilizer $H_v$ has order $|H|/|\calN_1^{(\infty)}|$. Since $\big(\calN_1^{(\infty)}\big)_P$ is a block system for $\calN_1^{(\infty)}$, the point stabilizer $H_v$ necessarily stabilizes the block containing $v$. Therefore,
\begin{equation}\label{EqnUnitary4}
(H_v)^U\leqslant
\begin{cases}
q^{2b(a-1)}{:}\SL_{a-1}(q^{2b})&\textup{if }S=\SL_a(q^{2b})\\
[q^{2b(a-1)}]{:}\Sp_{a-2}(q^{2b})&\textup{if }S=\Sp_a(q^{2b})\\
q^{4b+6b}{:}\SL_2(q^{2b})&\textup{if }S=\G_2(q^{2b}).
\end{cases}
\end{equation}
As $|H_v|/|(H_v)^U|$ divides $|P|$, we have
\[
|(H_v)^U|_{p'}=|H_v|_{p'}=|H|_{p'}/|\calN_1^{(\infty)}|_{p'}=|H|_{p'}/(q^{2m}-1)=|S|_{p'}/(q^{2m}-1),
\]
which together with~\eqref{EqnUnitary4} implies that $H_v=P.\SL_{a-1}(q^{2b})$, $P.\Sp_{a-2}(q^{2b})$ or $P.\SL_2(q^{2b})$, for some $p$-group $P$, according to $S=\SL_a(q^{2b})$, $\Sp_a(q^{2b})$ or $\G_2(q^{2b})$, respectively. Combined with $|H_v|=|H|/|\calN_1^{(\infty)}|=q^c|S|/|\calN_1^{(\infty)}|$ and $|\calN_1^{(\infty)}|=q^{2m-1}(q^{2m}-1)$, this leads to $|P|=q^{c-2b+1}$. Thus $G=HK$ with $H\cap K=H_v$ as described in parts~(a)--(c) of the lemma.
\end{proof}

%\begin{remark}
%If $a=m$ and $b=1$ in Lemma~\ref{LemUnitaryPm3}, then the conditions $\emptyset\neq I\subseteq\{1,\ldots,\lceil b/2\rceil\}$ and $\gcd(2I-1,b)=1$ turn out to be $I=\{1,\ldots,\lceil b/2\rceil\}=\{1\}$. In this case,~\eqref{EqnUnitary7} gives $c=m^2$, and the factorizations $G=HK$ asserted by Lemma~\ref{LemUnitaryPm3} are with $H=q^{m^2}{:}S$ such that one of the following holds:
%\begin{enumerate}[(a)]
%\item $S=\SL_m(q^2)$, and $H\cap K=[q^{m^2-1}].\SL_{m-1}(q^2)$;
%\item $S=\Sp_m(q^2)$, and $H\cap K=[q^{m^2-1}].\Sp_{m-2}(q^2)$;
%\item $S=\G_2(q^2)$ with $m=6$ and $q$ even, and $H\cap K=[q^{35}].\SL_2(q^2)$.
%\end{enumerate}
%\end{remark}

The proposition below concludes this subsection.

\begin{proposition}\label{PropUnitaryPm}
Let $(m,q)\neq(2,3)$, let $H\leqslant\Pa_m[G]$ with $H^{(\infty)}\neq\SL_m(q^2)$, $\Sp_m(q^2)$ or $\G_2(q^2)$ (with $m=6$ and $q$ even), and let $K^{(\infty)}=\N_1[G]^{(\infty)}=\SU_{2m-1}(q)$. Then $G=HK$ if and only if $(G,H,K)$ tightly contains some triple $(\SU_{2m}(q),q^c{:}S,\SU_{2m-1}(q))$ as in Lemma~$\ref{LemUnitaryPm3}$.
\end{proposition}

\begin{proof}
By Lemma~\ref{LemUnitaryPm3} it suffices to prove the ``only if'' part. Suppose $G=HK$, and let $A=G_U=\Pa_m[G]$ and $B=\N_1[G]$ be maximal subgroups of $G$ containing $H$ and $K$ respectively. Taking $a=m$ in Lemma~\ref{LemUnitaryPm3}(a) we deduce that $(A\cap B)^U$ stabilizes a $1$-space in $\bbF_{q^2}^m$. Then it follows from $A=H(A\cap B)$ that $H^U$ is transitive on the set of $1$-spaces in $\bbF_{q^2}^m$. Hence~\cite{Hering1985} asserts that $H^U\leqslant\GaL_a(q^{2b})$ is defined over $\bbF_{q^{2b}}$ with $m=ab$ and $(H^U)^{(\infty)}=\SL_a(q^{2b})$, $\Sp_a(q^{2b})$ or $\G_2(q^{2b})$ (with $a=6$ and $q$ even). Therefore, by Lemma~\ref{LemUnitaryPm1},
\[
H^{(\infty)}=U_{(b)}(I){:}S
\]
with $I\subseteq\{1,\ldots,\lceil b/2\rceil\}$ and $S=(H^U)^{(\infty)}$. Moreover, by Lemma~\ref{LemUnitaryPm2} we have $\gcd(2I-1,b)=1$.
If $I=\emptyset$, then this implies $b=1$ and hence $H^{(\infty)}=\SL_m(q^2)$, $\Sp_m(q^2)$ or $\G_2(q^2)$ (with $m=6$ and $q$ even), contradicting the assumption of the proposition. Thus $I\neq\emptyset$, and $(G,H,K)$ tightly contains the triple $(G^{(\infty)},H^{(\infty)},K^{(\infty)})=(\SU_{2m}(q),q^c{:}S,\SU_{2m-1}(q))$ as in Lemma~\ref{LemUnitaryPm3}.
\end{proof}

\subsection{Action on vectors of norm $1$}\label{SecUnitaryExample}
\ \vspace{1mm}

Recall that $\calN_1^{(\infty)}=\{u\in V \mid \beta(u,u)=1\}$ with $|\calN_1^{(\infty)}|=q^{2m-1}(q^{2m}-1)$. Finding factorizations $G=HG_v$ is equivalent to finding subgroups $H$ of $G$ that are transitive on $\calN_1^{(\infty)}$. In this section we construct examples of such $H$, other than those given in Proposition~\ref{PropUnitaryPm}.
Let $\gamma$ be the involution in $\GU(V)$ swapping $e_i$ and $f_i$ for all $i\in\{1,\dots,m\}$. If either $m$ or $q$ is even, then $\gamma\in\SU(V)$ as $\det(\gamma)=(-1)^m=1$. Fix a field automorphism $\phi\in\GaU(V)$ such that
\[
\phi\colon a_1e_1+b_1f_1+\dots+a_me_m+b_mf_m\mapsto a_1^pe_1+b_1^pf_1+\dots+a_m^pe_m+b_m^pf_m
\]
for $a_1,b_1\dots,a_m,b_m\in\bbF_{q^2}$ (notice that this definition of $\phi$ is different from that in~\cite[1.7.1]{BHR2013}, see~\cite{BHR2009}).
%By abuse of notation, we also use $\phi$ to denote the corresponding elements in $\Aut(L)$ and $\Out(L)$.

\subsubsection{Subgroups preserving a decomposition into two totally isotropic subspaces}

For the maximal totally isotropic subspaces $U=\langle e_1,\dots,e_m\rangle$ and $W=\langle f_1,\dots,f_m\rangle$, we have
\[
L_{\{U,W\}}=(L_U\cap L_W).2=T.(q-1).2=\SL_m(q^2).(q-1).2
\]
Both $\gamma$ and $\phi$ normalize $T$, and $\gamma$ induces a graph automorphism of the linear group $T$.
Recall $v=e_1+\lambda f_1\in\calN_1^{(\infty)}$ with $\lambda\in\bbF_{q^2}$ such that $\lambda+\lambda^q=1$.

%It follows that
%\begin{equation}\label{Lem:Unitary-Li1}
%|L_{U,W}|/|(L_{U,W})_v|=q^{2m-2}(q^{2m}-1)<q^{2m-1}(q^{2m}-1)=|\calN_1^{(\infty)}|,
%\end{equation}
%and so $L_{U,W}$ is not transitive on $\calN_1^{(\infty)}$.

\begin{lemma}\label{LemUnitary03}
For each group $G$ in the following table, the corresponding subgroup $H$ is transitive on $\calN_1^{(\infty)}$ with stabilizer $H_v=\SL_{m-1}(q^2)$.
\[
\begin{array}{c|cccc}
\hline
G & \SU_{2m}(2) & \SU_{2m}(2){:}\langle\phi\rangle & \SU_{2m}(4){:}\langle\phi\rangle & \SU_{2m}(4){:}\langle\phi\rangle \\
\hline
H & T{:}\langle\gamma\rangle & T{:}\langle\phi\rangle & T{:}\langle\phi\rangle & T{:}\langle\phi\gamma\rangle \\
\hline
\end{array}
\]
\end{lemma}

\begin{proof}
For each $t\in T_v$, we have $e_1^t\in U^t=U$ and $f_1^t\in W^t=W$. Then it follows from
\[
e_1^t+\lambda f_1^t=(e_1+\lambda f_1)^t=v^t=v=e_1+\lambda f_1
\]
that $e_1^t=e_1$ and $f_1^t=f_1$. Hence $T_v=T_{e_1,f_1}$. Moreover, as $U_1=f_1^\perp\cap U$ and $\langle f_1\rangle=U_1^\perp\cap W$, we have $T_{e_1,f_1}=T_{e_1,U_1}$. Therefore,
\[
T_v=T_{e_1,f_1}=T_{e_1,U_1}=\SL_{m-1}(q^2).
\]
Then for the pairs $(G,H)$ in the table we have
\[
|H|/|T_v|=2fq^{2m-2}(q^{2m}-1)=q^{2m-1}(q^{2m}-1)=|\calN_1^{(\infty)}|,
\]
and so it suffices to prove $H_v=T_v$, or equivalently, $H_v\leqslant T$.

First consider the candidate in column~1, namely, $(G,H)=(\SU_{2m}(2),T{:}\langle\gamma\rangle)$. Suppose that $H_v\nleqslant T$. Then $\gamma t$ fixes $v$ for some $t\in T$, and so $e_1+\lambda f_1=v=v^{\gamma t}=\lambda e_1^t+f_1^t$. Since $t$ stabilizes both $U$ and $W$, one has $e_1=\lambda e_1^t$ and $\lambda f_1=f_1^t$. Hence
\[
1=\beta(e_1,f_1)=\beta(e_1^t,f_1^t)=\beta(\lambda^{-1}e_1,\lambda f_1)=\lambda^{q-1},
\]
which implies $\lambda\in\bbF_q$, contradicting the condition that $\lambda+\lambda^q=1$. Thus $H_v\leqslant T$, as desired.

Next consider columns~2--4. Suppose for a contradiction that $H_v\nleqslant T$. Then there exists $t\in T$ such that $\phi^ft$ fixes $v$.
It follows that $e_1+\lambda f_1=v=v^{\phi^ft}=e_1^t+\lambda^qf_1^t$, which implies $e_1=e_1^t$ and $\lambda f_1=\lambda^qf_1^t$. This leads to
\[
1=\beta(e_1,f_1)=\beta(e_1^t,f_1^t)=\beta(e_1,\lambda^{1-q}f_1)=(\lambda^{1-q})^q=\lambda^{q-1},
\]
again contradicting $\lambda+\lambda^q=1$.
%
%Finally, we verify the minimality of these factorizations.
%Let $(G_0,H_0,K_0)$ be a minimal triple with respect to tight containment such that $((G_0)^{(\infty)},(H_0)^{(\infty)},(K_0)^{(\infty)})=(G^{(\infty)},H^{(\infty)},K^{(\infty)})$.
%For the candidates in rows~1--2, since $H=H^{(\infty)}.2$, it follows from~\eqref{Lem:Unitary-Li1} that the triple is indeed minimal.
%Now consider the candidates in rows~3--4, so that $q=|\phi|=4$.
%Suppose that $G_0=L{:}\l\phi^2\r=\SU_{2m}(4).2$.
%Then $H_0\leqslant \N_{G_0}(H^{(\infty)})=\SL_m(16).3.2$, and $K_0\leqslant \N_{G_0}(K^{(\infty)})=\GU_{2m-1}(4).2$.
%Since $H_0\cap K_0\geqslant\SL_{m-1}(16)$, we have
%\[
%\frac{|H_0|_2}{|H_0\cap K_0|_2}\leqslant\frac{|\SL_m(16).3.2|_2}{|\SL_{m-1}(16)|_2}=2^{4m-3}
%<2^{4m-2}=\frac{|\SU_{2m}(4).2|_2}{|\GU_{2m-1}(4).2|_2}\leqslant\frac{|G_0|_2}{|K_0|_2},
%\]
%a contradiction.
%Thus $G_0=L{:}\l\phi\r=\SU_{2m}(4).4$.
%Since $|H_0|_2/|H_0\cap K_0|_2=|G_0|_2/|K_0|_2$, one must have $H_0=\SL_m(16).4$ and $K_0=\SU_{2m-1}(4).4$.
%Hence $(G_0,H_0,K_0)$ lies in rows~3--4.
\end{proof}

For the groups $G$ and $H$ as in Lemma~\ref{LemUnitary03}, each subgroup $X$ of $H$ satisfies $X_v=X\cap H_v$, and $G=XG_v$ if and only if  $H=XH_v$. Then appealing to Theorem~\ref{ThmLinear} we obtain the subsequent example.

\begin{example}\label{ex:Unitary-03}
The triples $(G,H,K)$ in the following table give the minimal factorizations $G=HK$ with respect to tight containment such that $H\leqslant T{:}\langle\phi,\gamma\rangle$ and $K=G_v$, where $\phi$ is the field automorphism of order $2f$ fixing $e_1,f_1,\dots,e_m,f_m$, and $\gamma\in L$ swaps $e_i$ and $f_i$ for all $i$.
\[
\begin{array}{llll}
\hline
G & H & K & H\cap K \\
\hline
\SU_{2m}(2) & \SL_m(4){:}\langle\gamma\rangle & \SU_{2m-1}(2) & \SL_{m-1}(4) \\
 & \Sp_m(4){:}\langle\gamma\rangle & \SU_{2m-1}(2) & \Sp_{m-2}(4) \\
\SU_{12}(2) & \G_2(4){:}\langle\gamma\rangle & \SU_{11}(2) & \SL_2(4) \\
%\hline
\SU_{2m}(2){:}\langle\phi\rangle & \SL_m(4){:}\langle\phi\rangle & \SU_{2m-1}(2){:}2 & \SL_{m-1}(4) \\
 &\Sp_m(4){:}\langle\phi\rangle & \SU_{2m-1}(2).2 & \Sp_{m-2}(4) \\
\SU_{12}(2){:}\langle\phi\rangle & \G_2(4){:}\langle\phi\rangle & \SU_{11}(2){:}2 & \SL_2(4) \\
%\hline
\SU_{2m}(4){:}\langle\phi\rangle & \SL_m(16){:}\langle\phi\rangle,\ \SL_m(16){:}\langle\phi\gamma\rangle & \SU_{2m-1}(4){:}4 & \SL_{m-1}(16) \\
 & \Sp_m(16){:}\langle\phi\rangle,\ \Sp_m(16){:}\langle\phi\gamma\rangle & \SU_{2m-1}(4){:}4 & \Sp_{m-2}(16) \\
\SU_{12}(4){:}\langle\phi\rangle & \G_2(16){:}\langle\phi\rangle,\ \G_2(16){:}\langle\phi\gamma\rangle & \SU_{11}(4){:}4 & \SL_2(16) \\
\hline
\end{array}
\]
Moreover, $H\cap K$ is described in the table for each of the triple $(G,H,K)$.
\end{example}

\subsubsection{Subgroups of $\Sp_{2m}(q)$}

\begin{lemma}\label{LemUnitary09}
Let $G=\SU_{2m}(q)$, let $H=\Sp_{2m}(q)<G$, and let $K=\N_1[G]^{(\infty)}=\SU_{2m-1}(q)$. Then $G=HK$, and $H\cap K=\Sp_{2m-2}(q)$.
\end{lemma}

\begin{proof}
Let $\mu\in\bbF_{q^2}$ such that $\mu^{q-1}=-1$, and let $V_0=\langle\mu e_1,f_1,\dots,\mu e_m,f_m\rangle_{\bbF_q}$.
By~\cite[3.10.6]{Wilson2009}, there is a nondegenerate alternating form $\beta_0$ on $V_0$ such that $\mu e_1,f_1,\dots,\mu e_m,f_m$ is a standard basis of $V_0$ with respect to $\beta_0$ and $H=\Sp(V_0)$. Let $\zeta\in\bbF_{q^2}\Setminus\bbF_q$, and let $u=\mu e_1+\zeta f_1$.
Then $\beta(u,u)=\mu\zeta^q+\zeta\mu^q=\mu(\zeta^q-\zeta)\neq0$, and so we may assume without loss of generality that $K=G_u$.
Let $w=\mu e_1+\zeta^qf_1$. Then $\langle w\rangle_{\bbF_{q^2}}$ is the orthogonal complement of $\langle u\rangle_{\bbF_{q^2}}$ in the unitary space $\langle e_1,f_1\rangle_{\bbF_{q^2}}$, and hence
\[
u^\perp=\langle w,e_2,f_2,\dots,e_m,f_m\rangle_{\bbF_{q^2}}=\langle w,\mu e_2,f_2,\dots,\mu e_m,f_m\rangle_{\bbF_{q^2}}.
\]

Clearly, $V_0\cap u^\perp$ contains $\langle\mu e_2,f_2,\dots,\mu e_m,f_m\rangle_{\bbF_q}$. Suppose that they are not equal.
Then there exists nonzero $\delta\in\bbF_{q^2}$ such that $\delta w\in\langle\mu e_1,f_1\rangle_{\bbF_q}$; that is, $\delta(\mu e_1+\zeta^qf_1)=a\mu e_1+bf_1$ for some $a,b\in\bbF_q$. However, this implies $\zeta^q=b/a\in\bbF_q$, contradicting the condition $\zeta\in\bbF_{q^2}\Setminus\bbF_q$.
Thus, $V_0\cap u^\perp=\langle\mu e_2,f_2,\dots,\mu e_m,f_m\rangle_{\bbF_q}$. Then since $H\cap K$ stabilizes $V_0\cap u^\perp$, it stabilizes
\[
V_1:=(\langle\mu e_2,f_2,\dots,\mu e_m,f_m\rangle_{\bbF_q})^\perp=\langle e_1,f_1\rangle_{\bbF_{q^2}}=\langle\mu e_1,f_1\rangle_{\bbF_{q^2}}.
\]
Let $V_2=V_0\cap V_1=\langle\mu e_1,f_1\rangle_{\bbF_q}$. Then $H\cap K$ also stabilizes $V_2$.
For each $g\in H\cap K$, since $g|_{V_2}\in\Sp(V_2)$, we have $\det(g|_{V_2})=1$ and hence $\det(g|_{V_1})=1$.
This implies that $g|_{V_1}\in\SU(V_1)$. As $g$ fixes the nonsingular vector $u$ in the unitary space $V_1$, it follows that $g|_{V_1}=1$.
In particular, $g$ fixes $e_1$ and $f_1$, and so $H\cap K\leqslant H_{e_1,f_1}$. Conversely, $H_{e_1,f_1}$ is obviously contained in $H_u=H\cap K$. Hence
\[
H\cap K=H_{e_1,f_1}\cong\Sp(\langle\mu e_2,f_2,\dots,\mu e_m,f_m\rangle_{\bbF_q})=\Sp_{2m-2}(q).
\]
It follows that $|H|/|H\cap K|=q^{2m-1}(q^{2m}-1)=|G|/|K|$, and so $G=HK$.
\end{proof}

%\begin{lemma}\label{CorUnitary02}
%The action of $H=\Sp_{2m}(q)$ on $\calN_1^{(\infty)}$ is  permutationally equivalent to the action of $H$ on the set of refined antiflags of $\bbF_q^{2m}$.
%\end{lemma}

A consequence of Lemma~\ref{LemUnitary09} is that the action of $H=\Sp_{2m}(q)$ on $\calN_1^{(\infty)}$ is permutationally equivalent to the transitive action of $H$ on the set of refined antiflags of $\bbF_q^{2m}$. Moreover, it is known that the subgroup $\G_2(q)$ of $\Sp_6(q)$ with even $q$ is transitive on the set of refined antiflags of $\bbF_q^6$ (see Theorem~\ref{refined-antiflags}). Thus we have the next example.

\begin{example}\label{LemUnitary11}
If $q$ is even, then $\SU_6(q)=\G_2(q)\SU_5(q)$ with $\G_2(q)\cap\SU_5(q)=\SL_2(q)$.
\end{example}

\begin{remark}
Similar to Example~\ref{LemUnitary11}, one may construct refined-antiflag-transitive subgroups with solvable residual $\Sp_m(q^2)$ or $\G_2(q^2)$ (arising from part~(b) of Theorem~\ref{refined-antiflags}) that are transitive on $\calN_1^{(\infty)}$. However, these subgroups preserve a decomposition of $V$ into a direct sum of two maximal totally isotropic subspaces, and so the corresponding factorizations are already constructed in Example~\ref{ex:Unitary-03}.
\end{remark}

\subsubsection{Other examples}

The factorizations in row~7 of Table~\ref{TabUnitary} are described in the following example, which is verified by computation in \magma~\cite{BCP1997}.

\begin{example}\label{LemUnitary17}
Let $G=\SU_6(2)$, let $H$ be a subgroup of $G$ of the form $3^{\boldsymbol{\cdot}}\PSU_4(3)$ or $3^{\boldsymbol{\cdot}}\M_{22}$ (see \cite[Table~8.27]{BHR2013}), and let $K=\N_1[G]^{(\infty)}=\SU_5(2)$.
Then $G=HK$ with the intersection
\[
H\cap K=
\begin{cases}
3^5{:}\A_5&\textup{if }H=3^{\boldsymbol{\cdot}}\PSU_4(3)\\
\PSL_2(11)&\textup{if }H=3^{\boldsymbol{\cdot}}\M_{22}.
\end{cases}\qedhere
\]
\end{example}

The next example is known in~\cite[4.4.2]{LPS1990}.

\begin{example}\label{LemUnitary19}
Let $G=\SU_{12}(2)$, let $H=3^{\boldsymbol{\cdot}}\Suz<G$ (see~\cite[Tables~8.79]{BHR2013}), and let $K=\N_1[G]^{(\infty)}=\SU_{11}(2)$. Then $G=HK$ with $H\cap K=3^5.\PSL_2(11)$.
\end{example}

\subsection{Proof of Theorem~\ref{ThmUnitary}}\label{SecUnitaryProof}
\ \vspace{1mm}

By Lemma~\ref{Lem:UnitaryProof} we may assume $L\notin\{\SU_4(3),\SU_4(5),\SU_9(2)\}$. Thus it suffices to deal with case~(a) of Proposition~\ref{Prop:Unitary-max}, namely, when $L=\SU_{2m}(q)$ with $(m,q)\neq(2,2)$ and $K\leqslant B=\N_1[G]$. Under this assumption we have the ensuing lemma.

\begin{lemma}\label{Lem:Unitary-K}
If $G=HK$, then $K^{(\infty)}=B^{(\infty)}=\SU_{2m-1}(q)$, and one of the following holds:
\begin{enumerate}[{\rm (a)}]
\item $H\leqslant\Pa_m[G]$;
\item $H^{(\infty)}=\Sp_{2m}(q)$, or $\G_2(q)'$ with $m=3$ and $q$ even;
\item $q\in\{2,4\}$, and $H^{(\infty)}=\SL_m(q^2)$, $\Sp_m(q^2)$, or $\G_2(q^2)$ with $m=6$;
\item $(L,H^{(\infty)})\in\{(\SU_6(2),3.\PSU_4(3)),(\SU_6(2),3.\M_{22}),(\SU_{12}(2),3.\Suz)\}$.
\end{enumerate}
\end{lemma}

\begin{proof}
Let $G=HK$, and let $N=\Rad(B)$ be the solvable radical of $B$. Then $B/N$ is an almost simple group with socle $\PSU_{2m-1}(q)$.

Suppose that $KN/N\ngeqslant\Soc(B/N)$. Then by Lemma~\ref{LemXia7}, the group $B/N$ has a core-free factorization with a factor $KN/N$, and hence we see from Proposition~\ref{Prop:Unitary-max} that the only possibility is $\Soc(B/N)=\PSU_9(2)$, namely, $(n,q)=(10,2)$. However, for $B/N$ to have such a factorization with $G=AK$, computation in \magma~\cite{BCP1997} shows that this is impossible.

Thus we conclude that $KN/N\trianglerighteq\Soc(B/N)=\PSU_{2m-1}(q)$, and so $K^{(\infty)}=B^{(\infty)}=\SU_{2m-1}(q)$.
Notice that $G$ is transitive on $\mathcal{N}_1[G]$, the set of nonsingular $1$-spaces.
Since $K\leqslant B=\N_1[G]$, we have $G=H\N_1[G]$, which means that $H$ is transitive on $\mathcal{N}_1[G]$.
Then by~\cite[Lemma~4.3]{LPS2010}\footnote{In Lemma~4.3(iv) of~\cite{LPS2010}, $U_4(3)$ should be $3.U_4(3)$, and $M_{22}$ should be $3.M_{22}$.}, either of the lemma holds, or $H^{(\infty)}$ is one of the groups
\begin{equation}\label{EqnUnitary6}
\Sp_{m/2}(q^4)\text{ with $q=2$, }\ \ \SL_{m/2}(q^4)\text{ with $q=2$, }\ \ \G_2(q^4)\text{ with $m=12$ and $q=2$}.
\end{equation}
Suppose on the contrary that $H^{(\infty)}$ is one of these groups. In particular, $q=2$.

First, assume $A^{(\infty)}=\SL_m(q^2)=\SL_m(4)$.
%Let $\gamma$ and $\phi$ be as introduced at the start of Subsection~\ref{SecUnitaryExample}. Then $A\cap L=T{:}\langle\gamma\rangle$, and $A=(A\cap L).\calO$ with $\calO=G/L\leqslant\langle\phi\rangle=2$.
Then $A=\SL_m(4).2.\calO$ with $\calO=G/L\leqslant2$.
From Example~\ref{ex:Unitary-03} we deduce $(A\cap B)^{(\infty)}=\SL_{m-1}(4)$ and $|(A\cap B)/(A\cap B)^{(\infty)}|_2=|\calO|$. However, applying Theorem~\ref{ThmLinear} to the factorization $A=H(A\cap B)$ with $H^{(\infty)}$ described in~\eqref{EqnUnitary6} yields $|(A\cap B)/(A\cap B)^{(\infty)}|_2=|A/A^{(\infty)}|_2$, a contradiction.

Next, assume $A^{(\infty)}=\Sp_{2m}(q)=\Sp_{2m}(2)$. Here subgroups $H$ of $A$ described in~\eqref{EqnUnitary6} preserve a decomposition of $V$ into a direct sum of two maximal totally isotropic subspaces, and hence can be reduced to the previous paragraph, again yielding a contradiction.
\end{proof}

We are now ready to \emph{finish the proof of Theorem~$\ref{ThmUnitary}$}: We discuss the cases for $A$ in Proposition~\ref{Prop:Unitary-max}(a). For the case $A=\Pa_m[G]$, the conclusion of Theorem~\ref{ThmUnitary} is given by Proposition~\ref{PropUnitaryPm} if $H^{(\infty)}\neq\SL_m(q^2)$, $\Sp_m(q^2)$ or $\G_2(q^2)$ (with $m=6$ and $q$ even), and is given by Example~\ref{ex:Unitary-03} otherwise. If $A\cap L=\SL_m(q^2).(q-1).2$ with $q\in\{2,4\}$, then Lemma~\ref{Lem:Unitary-K} and Example~\ref{ex:Unitary-03} together imply the conclusion of the theorem. If $A^{(\infty)}=\Sp_{2m}(q)$, then the conclusion is a consequence of Lemma~\ref{Lem:Unitary-K}, Lemma~\ref{LemUnitary09}, Example~\ref{LemUnitary11} and its remark. Finally, if
\[
(L,A^{(\infty)})\in\{(\SU_6(2),3.\PSU_4(3)),(\SU_6(2),3.\M_{22}),(\SU_{12}(2),3.\Suz)\},
\]
then the conclusion follows from Lemma~\ref{Lem:Unitary-K} combined with Examples~\ref{LemUnitary17} and~\ref{LemUnitary19}. This completes the proof of Theorem~\ref{ThmUnitary}.

%%%%%%%%%%%%%%%%%%%%%%%%%%%%%%%%%%%%%%%%%%%%%%%%%%%%%%%%%%%%%%%%%%%
%%%%%%%%%%%%%%%%%%%%%%%%%%%%%%%%%%%%%%%%%%%%%%%%%%%%%%%%%%%%%%%%%%%

\section{Orthogonal groups in odd dimension}

%%%%%%%%%%%%%%%%%%%%%%%%%%%%%%%%%%%%%%%%%%%%%%%%%%%%%%%%%%%%%%%%%%%
%%%%%%%%%%%%%%%%%%%%%%%%%%%%%%%%%%%%%%%%%%%%%%%%%%%%%%%%%%%%%%%%%%%

In this section we classify the factorizations of classical groups $G$ with $L=G^{(\infty)}=\Omega_{2m+1}(q)$, where $m\geqslant3$ and $q$ is odd.
Such groups $G$ are almost simple with socle $L$. Therefore, throughout this section, we let $L=\Soc(G)=\Omega_{2m+1}(q)$ with $m\geqslant3$ and $q$ odd, and let $H$ and $K$ be nonsolvable core-free subgroups of $G$. Our main result of this section is as follows.

\begin{theorem}\label{ThmOmega}
With the above notation,  $G=HK$ if and only if $(G,H,K)$ tightly contains some triple $(G_0,H_0,K_0)$ in Table~$\ref{TabOmega}$.
\end{theorem}

After the setup in Subsection~\ref{SecOmegaPre} we will embark on the proof of Theorem~\ref{ThmOmega} and finalize it in Subsection~\ref{SecOmegaProof}.

\subsection{Setup}\label{SecOmegaPre}
\ \vspace{1mm}

Let $A$ and $B$ be $\max^-$ subgroups of $G$ containing $H$ and $K$ respectively. Then $A$ and $B$ are nonsolvable, and the $\max^-$ factorizations $G=AB$ are classified by Liebeck-Praeger-Saxl~\cite{LPS1990,LPS1996}, which are summarized in the proposition below.

\begin{proposition}\label{Prop:MaxOmega}
If $G=HK$, then with the above notation and interchanging $A$ and $B$ if necessary, one of the following holds:
\begin{enumerate}[{\rm(a)}]
\item $B=\N_1^-[G]$, and $A$ satisfies one of the following:
\begin{itemize}
\item $A=\Pa_m[G]$;
\item $A^{(\infty)}=\G_2(q)$ with $m=3$;
\item $A^{(\infty)}=\PSp_6(q)$ with $m=6$ and $q=3^f$;
\item $A^{(\infty)}=\F_4(q)$ with $m=12$ and $q=3^f$.
\end{itemize}
\item $m=3$, $B^{(\infty)}=\G_2(q)$, and $A=\Pa_1[G]$, $\N_1^+[G]$, $\N_2^+[G]$ or $\N_2^-[G]$;
\item $L=\Omega_7(3)$.
\end{enumerate}
\end{proposition}

The smallest group $L=\Omega_7(3)$ has more factorizations than the general case, as shown in the next computer-verified lemma.

\begin{lemma}\label{LemXia28}
Let $L=\Omega_7(3)$. Then $G=HK$ if and only if $(G,H,K)$ tightly contains some triple $(G_0,H_0,K_0)$ such that $G_0=\Omega_7(3)$ and $(H_0,K_0)$ lies in either rows~\emph{1--4} of Table~$\ref{TabOmega}$ or the table below, where the remarks restrict the conjugacy classes of valid $K_0$ in $L$ for each $H_0$.
\[
\begin{array}{llll}
\hline
H_0 & K_0 & H_0\cap K_0 & \textup{Remarks}\\
\hline
\G_2(3) & \Omega_5(3) & \SL_2(3) & \textup{$1$ in $2$ conjugacy classes of $K_0<\Sp_6(2)$}\\
 & 3^4{:}\Sy_5  & 3^2 & \textup{$2$ in $4$ conjugacy classes of $K_0$}\\
 & 3^5{:}2^4{:}\A_5 & \ASL_2(3) & \\
%\hline
\A_9,\ \Sp_6(2) & 3^3{:}\SL_3(3) & \Sy_3,\ \GL_2(3) & \\
 & \Omega_6^+(3) & \Sy_5\times2,\ 2^4{:}\Sy_5 & \\
 & \G_2(3) & \PSL_2(7),\ 2^3{}^{\boldsymbol{\cdot}}\PSL_2(7) & \textup{$1$ in $2$ conjugacy classes of $K_0$}\\
%\hline
3^{3+3}{:}\SL_3(3) & 2^6{:}\A_7 & 6^{\boldsymbol{\cdot}}\Sy_4 & \\
 & \Sy_8 & \Sy_3\times\Sy_3 & \\
 & \A_9 & 3^3{:}\Sy_3 & \\
 & 2^{\boldsymbol{\cdot}}\PSL_3(4) & 3^2{:}4 & \\
 & \Sp_6(2) & \SU_3(2){:}\Sy_3 & \\
\hline
%2^6{:}\A_7,\A_8,\A_9,2.\PSL_3(4),\Sp_6(2)&3^{3+3}{:}\SL_3(3) & \\ \hline
\end{array}
\]
\end{lemma}

From now on assume that $L\neq\Omega_7(3)$. Then part~(a) or (b) of Proposition~\ref{Prop:MaxOmega} is satisfied. In particular, we have
\[
B\cap L=\N_1^-[L]=\Omega_{2m}^-(q).2\ \text{ or }\ (L,B\cap L)=(\Omega_7(q),\G_2(q)).
\]
Let $q=p^f$ with an odd prime $p$, let $V=\bbF_q^{2m+1}$ be equipped with a nondegenerate symmetric bilinear form $\beta$, and
let $e_1,f_1,\dots,e_m,f_m,d$ be a standard basis for $V$ as in~\cite[2.2.3]{LPS1990}.
Denote by $\bbF_q^\square$ the subgroup of index $2$ in $\bbF_q^\times$, and $\bbF_q^\boxtimes=\bbF_q^\times\Setminus\bbF_q^\square$.
Fix some
\[
\lambda\in\bbF_q^\boxtimes.
\]
Let $U=\langle e_1,\dots,e_m\rangle$ and $W=\langle f_1,\dots,f_m\rangle$ be a pair of maximal totally isotropic subspaces of $V$, and let $U_1=\langle e_2,\dots,e_m\rangle$. By~\cite[3.7.4]{Wilson2009}, the stabilizer $L_U=\Omega(V)_U$ has the form
\[
L_U=\Pa_m[L]=(q^{m(m-1)/2}.q^m){:}\SL_m(q).\tfrac{q-1}{2}.
\]
Let $R=q^{m(m-1)/2}.q^m$ be the unipotent radical of $L_U$,
and $T=\SL_m(q)$ be a subgroup of $L_U$ stabilizing $U$ and $W$ respectively.
Note that $R$ is the kernel of $L_U$ acting on $U$, and the action of $T$ on $U$ determines that on $W$ (see~\cite[Lemma~2.2.17]{BG2016}).

% For a vector $u\in V$, denote by $r_u$ the reflection with respect to $u$???.

\subsection{Actions on nonzero vectors}
\ \vspace{1mm}

In the orthogonal space $V=\bbF_q^{2m+1}$, there are three types of nonzero vectors, namely, nonzero singular vectors, and nonsingular vectors of plus type and minus type, respectively. Let
\begin{align*}
\mathcal{P}_1^{(\infty)}&=\{u\in V\Setminus\{0\}\mid \beta(u,u)=0\},\\
\mathcal{N}_1^{+(\infty)}&=\{u\in V\mid \beta(u,u)=2\},\\
\mathcal{N}_1^{-(\infty)}&=\{u\in V\mid \beta(u,u)=2\lambda\},
\end{align*}
where, recall that $\lambda\in\bbF_q^\boxtimes$. The group $L=\Omega_{2m+1}(q)$ is transitive on $\mathcal{P}_1^{(\infty)}$, $\mathcal{N}_1^{+(\infty)}$ and $\mathcal{N}_1^{-(\infty)}$, respectively. Each transitive subgroup of $L$ on any of these three sets gives rise to a corresponding factorization of $L$, and we will identify such transitive subgroups in this subsection.

\subsubsection{Nonsingular vectors of minus type}

Take $v=e_1+\lambda f_1$. Then $v\in\mathcal{N}_1^{-(\infty)}$ and $L_v=(\N_1^-[L])^{(\infty)}=\Omega_{2m}^-(q)$. It follows that
\[
|\calN_1^{-(\infty)}|=|\Omega_{2m+1}(q)|/|\Omega_{2m}^-(q)|=q^m(q^m-1).
\]
By Proposition~\ref{Prop:MaxOmega}, if $B=\N_1^-[G]$, then either $A=\Pa_m[G]$, or $A$ is a $\calC_9$-subgroup of $G$ with $m\in\{3,6,12\}$.

We first construct examples of $H\leqslant L_U=\Pa_m[L]$ that are transitive on $\mathcal{N}_1^{-(\infty)}$, which will give rise to factorizations $L=HL_v$. This is based on the following key lemma. Recall that $(L_U)^{(\infty)}=R{:}T$, where $R=q^{m(m-1)/2}.q^m$ is the kernel of $L_U$ acting on $U=\langle e_1,\dots,e_m\rangle$ and $T=\SL_m(q)$ stabilizes $W=\langle f_1,\dots,f_m\rangle$.

\begin{lemma}\label{LemOmega01}
Let $M=(L_U)^{(\infty)}=R{:}T$, let $K=L_v$, where $v=e_1+\lambda f_1\in\mathcal{N}_1^{-(\infty)}$, and let $\big(\mathcal{N}_1^{-(\infty)}\big){}_R$ be the set of orbits of $R$ on $\mathcal{N}_1^{-(\infty)}$. Then the following statements hold:
\begin{enumerate}[{\rm (a)}]
\item The kernel of $M\cap K$ acting on $U$ is the special $p$-group $R\cap K=q^{(m-1)(m-2)/2}.q^{m-1}$.
\item The induced group by the action of $M\cap K$ on $U$ is $\SL(U)_{U_1,e_1+U_1}=q^{m-1}{:}\SL_{m-1}(q)$.
\item The action of $T$ on $\big(\mathcal{N}_1^{-(\infty)}\big){}_R$ is permutationally isomorphic to its action on $\bbF_q^m\Setminus\{0\}$.
\item If $H=R{:}S$ with $S\leqslant T$, then $H\cap K=(R\cap K).S_{U_1,e_1+U_1}$.
\end{enumerate}
\end{lemma}

\begin{proof}
We first calculate $R\cap K$, the kernel of $M\cap K$ acting on $U$. Since each element $r\in R\cap K$ fixes $e_1$ and $v$, it follows that $r$ fixes $\langle e_1,v\rangle=\langle e_1,f_1\rangle$ pointwise. Hence $R\cap K$ is isomorphic to the pointwise stabilizer of $U_1$ in $\Omega(\langle e_2,f_2,\dots,e_m,f_m,d\rangle)$. Then it is shown in~\cite[3.7.4]{Wilson2009} that $R\cap K=q^{(m-1)(m-2)/2}.q^{m-1}$ is a special $p$-group. This proves part~(a).

Since $K=L_v$ stabilizes $v^\perp$, the intersection $M\cap K$ stabilizes $U\cap v^\perp=\langle e_2,\dots,e_m\rangle=U_1$.
For each $h\in M\cap K$, we have $e_1^h=\mu e_1+e$ with $\mu\in\bbF_q$ and $e\in U_1$. Then
\[
\mu\lambda=\beta(\mu e_1+e,e_1+\lambda f_1)=\beta(e_1^h,(e_1+\lambda f_1)^h)=\beta(e_1,e_1+\lambda f_1)=\lambda.
\]
Thus $\mu=1$; that is, $e_1^h=e_1+e$. This means that $M\cap K$ stabilizes $e_1+U_1$, and so the induced group of $M\cap K$ on $U$ is contained in $\SL(U)_{U_1,e_1+U_1}$, with kernel $R\cap K$. Now
\begin{align*}
&(M\cap K)/(R\cap K)\cong(M\cap K)^U\leqslant\SL(U)_{U_1,e_1+U_1}=q^{m-1}{:}\SL_{m-1}(q),\\
&|M\cap K|\geqslant\frac{|M||K|}{|L|}=\frac{|(q^{m(m-1)/2}.q^m){:}\SL_m(q)||\Omega_{2m}^-(q)|}{|\Omega_{2m+1}(q)|}
=|R\cap K||\SL(U)_{U_1,e_1+U_1}|.
\end{align*}
Hence $(M\cap K)^U=\SL(U)_{U_1,e_1+U_1}=q^{m-1}{:}\SL_{m-1}(q)$, proving part~(b).

From parts~(a) and~(b) we deduce $|M|/|M_v|=|M|/|M\cap K|=q^m(q^m-1)=|\calN_1^{-(\infty)}|$, which implies that $M$ is transitive on $\calN_1^{-(\infty)}$. Then $T$ is transitive on $\big(\mathcal{N}_1^{-(\infty)}\big){}_R$ with stabilizer $(M_v)^U$, which is $\SL(U)_{U_1,e_1+U_1}=q^{m-1}{:}\SL_{m-1}(q)$ by part~(b). Thus part~(c) holds.

Finally, let $H=R{:}S$ with $S\leqslant T$. Since $(M\cap K)^U=\SL(U)_{U_1,e_1+U_1}=(M_{U_1,e_1+U_1})^U$, we obtain $M_{U_1,e_1+U_1}=(M\cap K)R$ . This implies that $H_{U_1,e_1+U_1}=(H\cap K)R$, and so
\[
(H\cap K)/(R\cap K)\cong(H\cap K)R/R=H_{U_1,e_1+U_1}/R=S_{U_1,e_1+U_1}R/R\cong S_{U_1,e_1+U_1}.
\]
Therefore, $H\cap K=(R\cap K).S_{U_1,e_1+U_1}$, as part~(d) asserts.
\end{proof}

For $H=R{:}S$ with $S\leqslant T=\SL_m(q)$, Lemma~\ref{LemOmega01} implies that $H$ is transitive on $\mathcal{N}_1^{-(\infty)}$ if and only if $S$ is transitive on $\bbF_q^m\Setminus\{0\}$, in which case, $H_v=(q^{(m-1)(m-2)/2}.q^{m-1}).S_{U_1,e_1+U_1}$. Then appealing to Theorem~\ref{ThmLinear} we obtain the subsequent example.

\begin{example}\label{Ex:OmegaOdd-1}
Let $G=\Omega_{2m+1}(q)$, $K=G_v=\Omega_{2m}^-(q)$ and $H=R{:}S=(q^{m(m-1)/2}.q^m){:}S$ with one of the following subgroup $S$ of $T=\SL_m(q)$:
\begin{enumerate}[{\rm(a)}]
\item $S=\SL_a(q^b)$ or $\Sp_a(q^b)$, defined over $\bbF_{q^b}$, where $m=ab$;
\item $S=2^{\boldsymbol{\cdot}}\Sy_5$, $8^{\boldsymbol{\cdot}}\A_5$ or $2^{1+4}{}^{\boldsymbol{\cdot}}\A_5$, where $G=\Omega_9(3)$;
\item $S=\SL_2(13)$, where $G=\Omega_{13}(3)$.
\end{enumerate}
Then $G=HK$ with $H\cap K$ listed as follows.
\[
\begin{array}{lll}
\hline
G & H & H\cap K \\
\hline
\Omega_{2m+1}(q) & (q^{m(m-1)/2}.q^m){:}\SL_a(q^b) & (q^{(m-1)(m-2)/2}.q^{m-1}.q^{m-b}).\SL_{a-1}(q^b) \\
\textup{($m=ab$)} & (q^{m(m-1)/2}.q^m){:}\Sp_a(q^b) & (q^{(m-1)(m-2)/2}.q^{m-1}.q^{m-b}).\Sp_{a-2}(q^b) \\
%\hline
\Omega_9(3) & 3^{6+4}{:}2.\Sy_5 & 3^{3+3}.3 \\
 & 3^{6+4}{:}8.\A_5 & 3^{3+3}.\Sy_3 \\
 & 3^{6+4}{:}2^{1+4}.\A_5 & 3^{3+3}.\SL_2(3) \\
%\hline
\Omega_{13}(3) & 3^{15+6}{:}\SL_2(13) & 3^{10+5}.3 \\
\hline
\end{array}
\]
\par\vspace{-1\baselineskip}
\qedhere
\end{example}

Besides $\Pa_m[G]$, there are $\calC_9$-subgroups of $G$ as candidates for $A$, namely, $A^{(\infty)}=\G_2(q)$ with $L=\Omega_7(q)$, $A^{(\infty)}=\PSp_6(3^f)$ with $L=\Omega_{13}(3^f)$, and $A^{(\infty)}=\F_4(3^f)$ with $L=\Omega_{25}(3^f)$. In the ensuing three lemmas, we give some subgroups of these candidates that are transitive on $\mathcal{N}_1^{-(\infty)}$, leading to factorizations $L=HL_v$.

\begin{lemma}\label{Lem:Omega7-1}
Let $G=\Omega_7(q)$, and let $K=G_v=\Omega_6^-(q)$. Then the following statements hold:
\begin{enumerate}[{\rm(a)}]
\item If $H=\G_2(q)$, then $G=HK$ with $H\cap K=\SU_3(q)$.
\item If $q=3^f$, then $A=\G_2(q)$ has precisely one conjugacy class of subgroups $H$ isomorphic to $\SL_3(q)$ such that $G=HK$ and $H\cap K=q^2-1$.
\end{enumerate}
\end{lemma}

\begin{proof}
First let $H=\G_2(q)$. Then we see from~\cite[4.3.6]{Wilson2009} that $H\cap K=H_v=\SU_3(q)$. It follows that $|H|/|H_v|=q^3(q^3-1)=|\mathcal{N}_1^{-(\infty)}|$, and so $H$ is transitive on $\mathcal{N}_1^{-(\infty)}$, which implies $G=HG_v=HK$. Thus part~(a) holds.

Next let $q=3^f$ and let $H=\SL_3(q)$ be a subgroup of $A=\G_2(q)$. Part~(a) gives $G=AG_v$ with $A_v=A\cap K=\SU_3(q)$. Then $G=HK$ if and only if $A=HA_v$. According to the proof of~\cite[Proposition~A]{HLS1987}, there is precisely one conjugacy class of subgroups $H$ of $A$ such that $A=HA_v$ and $H\cap A_v=q^2-1$. For each such $H$, we have $H\cap K=H_v=H\cap A_v=q^2-1$.
\end{proof}

\begin{lemma}\label{LemOmega09}
For $G=\Omega_{13}(q)$ with $q=3^f$, let $H=\PSp_6(q)<G$, and let $K=G_v=\Omega_{12}^-(q)$. Then $G=HK$ with $H\cap K=H_v=(\SL_2(q)\times\SL_2(q^2))/\{\pm1\}=(\SL_2(q)\times\SL_2(q^2))/2$.
\end{lemma}

\begin{proof}
We follow the setup in~\cite[4.6.3(a)]{LPS1990}. Let $\beta_0$ be a nondegenerate alternating form on the natural $6$-dimensional module $V_0$ preserved by $X$, let $E_1,F_1,E_2,F_2,E_3,F_3$ be a standard basis for $V_0$ with respect to $\beta_0$, and let $V_1=\langle u\wedge w\mid u,w\in V_0,\,\beta_0(u,w)=0\rangle$.
Then $V_1$ is a $14$-dimensional submodule of the alternating square $\bigwedge^2V_0$, and $V_1$ has a trivial $1$-dimensional $H$-submodule $V_2$ with $V_1/V_2$ irreducible. We may assume $V=V_1/V_2$ and
\[
\beta(u_1\wedge w_1,u_2\wedge w_2)=\beta_0(u_1,u_2)\beta_0(w_1,w_2)-\beta_0(u_1,w_2)\beta_0(u_2,w_1)
\]
for all $u_i,w_i$ among the basis vectors $E_1,F_1,E_2,F_2,E_3,F_3$. Let $S=\langle E_1,F_1,E_2,F_2\rangle$.
As $K=G_v=\Omega_{12}^-(q)$, we may assume $v=E_1\wedge E_2+\mu F_1\wedge F_2$ for some $\mu\in\bbF_q$. Then $B=G_{\langle v\rangle}$ is a maximal subgroup of $G$ containing $K$, and by the proof of Lemma~A in~\cite[4.6.3]{LPS1990},
\begin{equation}\label{EqnOmega01}
H\cap B=H_{\langle v\rangle}=(\SL_2(q)\times\SL_2(q^2).2)/\{\pm1\}
\end{equation}
with $(H\cap B)_S^S=\Sp(S)_{\langle v\rangle}=(2.\Omega_5(q))_{\langle v\rangle}=2.\mathrm{O}_4^-(q)=\SL_2(q^2).2$.
Now $H\cap K=H_v$, and so $(H\cap K)_S^S=\Sp(S)_v=(2.\Omega_5(q))_v=2.\Omega_4^-(q)=\SL_2(q^2)$ has index $2$ in $(H\cap B)_S^S$.
This implies that $H\cap K$ has index $2$ in $H\cap B$. Then we deduce from~\eqref{EqnOmega01} that
\[
H\cap K=(\SL_2(q)\times\SL_2(q^2))/\{\pm1\}=(\SL_2(q)\times\SL_2(q^2))/2.
\]
Hence $|H|/|H_v|=|H|/|H\cap K|=q^6(q^6-1)=|\mathcal{N}_1^{-(\infty)}|$, and so $G=HG_v=HK$.
\end{proof}

\begin{lemma}\label{LemOmega10}
Let $G=\Omega_{25}(q)$ with $q=3^f$, let $H=\F_4(q)<G$, and let $K=G_v=\Omega_{24}^-(q)$. Then $G=HK$ with $H\cap K=\mathrm{Spin}_8^-(q)=2.\Omega_8^-(q)$.
\end{lemma}

\begin{proof}
It suffices to prove $H_v=\mathrm{Spin}_8^-(q)$. By~\cite[Table~2]{CC1988} we have
\begin{equation}\label{EqnOmega02}
H_{\langle v\rangle}=\mathrm{Spin}^-_8(q).2=(2.\Omega_8^-(q)).2.
\end{equation}
Let $M$ be a maximal subgroup of $H$ containing $H_{\langle v\rangle}$. Then $M_{\langle v\rangle}=H_{\langle v\rangle}$, and as $|M|\geqslant|H_{\langle v\rangle}|>q^{24}$,~\cite{LS1987} asserts that $M$ is either parabolic or one of $\mathrm{Spin}_9(q)$, $\mathrm{Spin}^+_8(q).\Sy_3$, $^3{}\D_4(q).3$ or $\F_4(q^{1/2})$. Since $|H_{\langle v\rangle}|$ is divisible by $q^8-1$ and so is $|M|$, it follows that $M=\mathrm{Spin}_9(q)=2.\Omega_9(q)$.
Since $\bbF_q$ is a splitting field for $M$ (see~\cite[Page~241]{Steinberg1968}), we deduce from~\cite[Theorem~1.1]{Liebeck1985} that each irreducible submodule of $M$ acting on $V=\bbF_q^{25}$ is either the natural module of dimension $9$ or the spin module of dimension $16$. Let $S$ be such a submodule. Then $S$ is a nondegenerate subspace of $V$, and so $V=S\oplus S^\perp$ with $S$ and $S^\perp$ both $M$-invariant. Write $v=x+y$ for some $x\in S$ and $y\in S^\perp$. Then $M_{\langle v\rangle}\leqslant M_{\langle x\rangle}\cap M_{\langle y\rangle}$ and $M_v=M_x\cap M_y$.

First assume that $S$ is the natural module of $M$. Then the action of $M$ on $S$ induces $\Omega(S)=\Omega_9(q)$ with kernel $2$.
Hence $(H_v)^S\leqslant(M_v)^S\leqslant(M_x)^S=\Omega(S)_x$, which together with~\eqref{EqnOmega02} implies that $H_v\neq H_{\langle v\rangle}$.
Since $G_v$ has index $2$ in $G_{\langle v\rangle}$, it follows that $H_v$ has index $2$ in $H_{\langle v\rangle}$.
Thus we conclude from~\eqref{EqnOmega02} that $H_v=\mathrm{Spin}^-_8(q)$, as desired.

Next assume that $S$ is the spin module of $M$. Then $S^\perp$ has dimension equal to $25-\dim(S)=9$.
Since $M=\mathrm{Spin}_9(q)=2.\Omega_9(q)$, the induced group of $M$ on $S^\perp$ is either trivial or equal to $\Omega_9(q)$.
For the latter case, replacing $S$ with $S^\perp$ in the previous paragraph gives $H_v=\mathrm{Spin}^-_8(q)$, as desired.
Suppose for the rest of the proof that the induced group of $M$ on $S^\perp$ is trivial. Then $M_{\langle y\rangle}=M$.
From~\eqref{EqnOmega02} we deduce that $M_{\langle x\rangle}\geqslant M_{\langle v\rangle}=H_{\langle v\rangle}=\mathrm{Spin}^-_8(q).2$.
If $x\neq0$, then the Proposition of~\cite[Appendix~3]{LPS1990} shows that $M_{\langle x\rangle}/2$ is one of the groups
\[
\Omega_7(q).2,\quad q^{6+4}{:}\SL_4(q).\tfrac{q-1}{2},\quad q^7.\G_2(q).\tfrac{q-1}{2},
\]
a contradiction.
Therefore, $x=0$. Hence $v=y$, and so
\[
H_{\langle v\rangle}=M_{\langle v\rangle}=M_{\langle y\rangle}=M=\mathrm{Spin}_9(q),
\]
contradicting~\eqref{EqnOmega02}.
\end{proof}

To conclude, we classify the factorizations $G=HK$ with $(\N_1^-[G])^{(\infty)}\leqslant K\leqslant\N_1^-[G]$ in the next proposition. As a corollary, $H$ is transitive on $\calN_1^-$ if and only if it is transitive on $\calN_1^{-(\infty)}$.

\begin{proposition}\label{Thm:gps-on-V^-}
Let $K\leqslant\N_1^-[G]$ with $K^{(\infty)}=\Omega_{2m}^-(q)$. Then $G=HK$ if and only if $(G,H,K)$ tightly contains $(L,H^{(\infty)},K^{(\infty)})$ such that $H^{(\infty)}$ lies in the following table.
\[
\begin{array}{lll}
\hline
H^{(\infty)} & H^{(\infty)}\cap K^{(\infty)} & \textup{Remarks} \\
\hline
(q^{m(m-1)/2}.q^m){:}\SL_a(q^b) & (q^{(m-1)(m-2)/2}.q^{m-1}.q^{m-b}){:}\SL_{a-1}(q^b) & m=ab,\ \textup{as in \ref{Ex:OmegaOdd-1}} \\
(q^{m(m-1)/2}.q^m){:}\Sp_a(q^b) & (q^{(m-1)(m-2)/2}.q^{m-1}.q^{m-b}){:}\Sp_{a-2}(q^b)& m=ab,\ \textup{as in \ref{Ex:OmegaOdd-1}} \\
%\hline
\G_2(q) & \SU_3(q) & L=\Omega_7(q) \\
\SL_3(3^f) & 3^{2f}-1 & L=\Omega_7(3^f),\ \textup{as in \ref{Lem:Omega7-1}(b)} \\
%\hline
3^{6+4}{:}2.\Sy_5 & 3^{3+3}.3 & L=\Omega_9(3),\ \textup{as in \ref{Ex:OmegaOdd-1}} \\
3^{6+4}{:}8.\A_5 & 3^{3+3}.\Sy_3 & L=\Omega_9(3),\ \textup{as in \ref{Ex:OmegaOdd-1}} \\
3^{6+4}{:}2^{1+4}.\A_5 & 3^{3+3}.\SL_2(3) & L=\Omega_9(3),\ \textup{as in \ref{Ex:OmegaOdd-1}} \\
%\hline
\PSp_6(3^f) & (\SL_2(3^f)\times\SL_2(3^{2f}))/2 & L=\Omega_{13}(3^f) \\
3^{15+6}{:}\SL_2(13) & 3^{10+5}.3 & L=\Omega_{13}(3),\ \textup{as in \ref{Ex:OmegaOdd-1}} \\
%\hline
\F_4(3^f) & 2.\Omega_8^-(3^f) & L=\Omega_{25}(3^f) \\
\hline
\end{array}
\]
In this case, $H^{(\infty)}\cap K^{(\infty)}$ is described in the table.
\end{proposition}

\begin{proof}
By Example~\ref{Ex:OmegaOdd-1} and Lemmas~\ref{Lem:Omega7-1}--\ref{LemOmega10}, it suffices to prove the ``only if'' part. Thus suppose that $G=HK$. Then, according to Proposition~\ref{Prop:MaxOmega}, either $H\leqslant\Pa_m[G]$, or
\[
(L,A^{(\infty)})\in\{(\Omega_7(q),\G_2(q)),(\Omega_{13}(3^f),\PSp_6(3^f)),(\Omega_{25}(3^f),\F_4(3^f))\}.
\]
For the latter, since Lemmas~\ref{Lem:Omega7-1}--\ref{LemOmega10} give $A^{(\infty)}\cap K^{(\infty)}=\SU_3(q)$, $(\SL_2(3^f)\times\SL_2(3^{2f}))/2$ or $2.\Omega_8^-(3^f)$ respectively, we apply~\cite{HLS1987},~\cite[Theorem~A]{LPS1990} and~\cite{LPS1996} to the factorization $A=H(A\cap K)$ and conclude that $H^{(\infty)}$ lies in the table of the proposition.

Now let $H\leqslant\Pa_m[G]$. Thus we may assume without loss of generality that $A=G_U$. Let $P=H\cap R$, where, recall that $R$ is the kernel of $L_U$ acting on $U$ and that $R=q^{m(m-1)/2}.q^m$ is a special group with $\Z(R)=\Phi(R)=R'=q^{m(m-1)/2}$. From Lemma~\ref{LemOmega01}(b) we derive that the induced group $(H\cap K)^U\leqslant(A\cap B)^U$ stabilizes a hyperplane of $U$. Since the induced group of $H$ on $U$ is $H/P$, it follows that $H/P$ acts transitively on the set of hyperplanes in $U$. Hence Theorem~\ref{ThmLinear} implies that $H/P$ has a normal subgroup in one of the sets below, where $m=ab$.
\begin{equation}\label{EqnOmega15}
\{\SL_a(q^b),\Sp_a(q^b)\},\ \{2.\Sy_5,8.\A_5,2^{1+4}.\A_5\}\ (L=\Omega_9(3)),\ \{\SL_2(13)\}\ (L=\Omega_{13}(3)).
\end{equation}
In particular, $H/P$ acts irreducibly on $R/R'=q^m$.

Suppose that $P\leqslant R'$. Then $H$ is contained in a subgroup $M$ of $G$ such that
\[
M\cap L=R'{:}T.\tfrac{q-1}{2}=R'{:}\SL_m(q).\tfrac{q-1}{2}.
\]
It follows from Lemma~\ref{LemOmega01} that
\begin{equation}\label{EqnOmega04}
(R{:}T)\cap L_v=(R\cap L_v)'.q^{m-1}.(q^{m-1}{:}\SL_{m-1}(q))
\end{equation}
with $(R\cap L_v)'=q^{(m-1)(m-2)/2}$. Since $T\cap L_v\geqslant\SL(U)_{e_1,U_1}=\SL_{m-1}(q)$, we have
\begin{equation}\label{EqnOmega05}
(R'{:}T)\cap L_v\geqslant(R\cap L_v)'{:}(T\cap L_v)\geqslant(R\cap L_v)'{:}\SL_{m-1}(q).
\end{equation}
Since $G=HK=MK$, by comparing the $p$-parts, we obtain
\[
q^{m-3}|(R\cap L_v)'{:}\SL_{m-1}(q)|_p<|(R'{:}T)\cap L_v|_p<q^{m-1}|(R\cap L_v)'{:}\SL_{m-1}(q)|_p.
\]
However, it is clear from~\eqref{EqnOmega04} and~\eqref{EqnOmega05} that $(R{:}T)\cap L_v$ does not have such a subgroup $(R'{:}T)\cap L_v$.

We thus conclude that $P\nleqslant R'$, and so $PR'/R'=R/R'$ as $H/P$ acts irreducibly on $R/R'$. Consequently, $R=PR'=P\Phi(R)$.
Since $\Phi(R)$ consists of the non-generators of $R$, this implies that $P=R$. Hence $H^{(\infty)}=R{:}S$ such that $S$ belongs to one of the sets in~\eqref{EqnOmega15}, and thus $H^{(\infty)}$ is described in the table of the proposition.
\end{proof}

\subsubsection{Nonzero singular vectors}

Let $v\in\mathcal{P}_1^{(\infty)}$ be a nonzero singular vector of $V$. Then $L=\Omega_{2m+1}(q)$ is a transitive permutation group on $\mathcal{P}_1^{(\infty)}$ with stabilizer $L_v=q^{2m-1}{:}\Omega_{2m-1}(q)$ and degree $|\mathcal{P}_1^{(\infty)}|$, where
%\[
%L_v=q^{2m-1}{:}\Omega_{2m-1}(q)\ \text{ and }\ |\mathcal{P}_1^{(\infty)}|=|L|/|L_v|=q^{2m}-1.
%\]
\[
|\mathcal{P}_1^{(\infty)}|=|L|/|L_v|=q^{2m}-1.
\]
According to Proposition~\ref{Prop:MaxOmega}, to classify factorizations of $G$ with a factor contained in $\Pa_1[G]$, we only need to consider the case when $m=3$.

Later on we will show a factorization of $\Omega_7(q)$ with factors $q^4{:}\Omega_4^-(q)$ and $\G_2(q)$. For this purpose we need to embed $L=\Omega_7(q)$ into $\Omega_8^+(q)$ and construct the two desired factors there.
To start with, we extend the bilinear form $\beta$ on $V=\bbF_q^7$ to a plus type $8$-dimensional orthogonal space over $\bbF_q$ with a basis $x_1,\dots, x_8$ such that $\beta(x_i,x_j)=\delta_{i+j,9}$ for $i,j\in\{1,\dots,8\}$. Let $x=x_4+x_5$ and $y=x_4-x_5$. Then $x$ and $y$ are nonsingular, and
\[
x^\perp=\langle x_1,x_2,x_3,x_6,x_7,x_8,y\rangle.
\]
Thus we may take $L=\Omega(x^\perp)$. Let $W_1=\langle x_2,x_3,x_6,x_7,y\rangle$ and $V_1=\langle x_1\rangle\oplus W_1$.
Then $V_1=x_1^\perp\cap x^\perp$ and hence is stabilized by $L_{\langle x_1\rangle}$. Since $x_1$ is singular, we have $L_{\langle x_1\rangle}=\Pa_1[L]$ and
\[
L_{x_1}=E{:}X,
\]
where $E=q^5$ is the kernel of $L_{x_1}$ acting on $V_1/\langle x_1\rangle$, and $X=L_{x_1,x_8,W_1}=\Omega(W_1)=\Omega_5(q)$.
For $a\in\bbF_q$, let $\sigma(a)$ be the linear transformation on $x^\perp$ satisfying
\begin{equation}\label{EqnOmega06}
\sigma(a)\colon x_3\mapsto x_3-ax_1,\ x_8\mapsto x_8+ax_6
\end{equation}
and fixing the vectors $x_1,x_2,x_6,x_7,y$. Let $F$ be the group generated by $\sigma(a)$ with $a$ running over $\bbF_q$, and let $w=x_3+\lambda x_6$ and $z=x_3-\lambda x_6$ (recall that $\lambda\in\bbF_q^\boxtimes$).

\begin{lemma}\label{Lem:Omega-001}
With the above notation, $E=E_w\times F$ with $E_w=q^4$ and $F=q$, and $X_w=\Omega_4^-(q)$.
\end{lemma}

\begin{proof}
Clearly, $F=q$, and it is straightforward to verify that $\sigma(a)\in\mathrm{O}(x^\perp)$ for all $a\in\bbF_q$. Then since $|\sigma(a)|$ is odd, we deduce that $\sigma(a)\in\Omega(x^\perp)=L$. Moreover, since $\sigma(a)$ fixes $x_1$ and acts trivially on $V_1/\langle x_1\rangle$, we obtain $\sigma(a)\in E$. Note that $E$ stabilizes $w+\langle x_1\rangle$ as $E$ is the kernel of $L_{x_1}$ on $V_1/\langle x_1\rangle$. Then since $F$ is regular on $w+\langle x_1\rangle$, we conclude that $E=E_w\times F$. This together with $E=q^5$ and $F=q$ implies that $E_w=q^4$. Finally, since $w^\perp\cap W_1=\langle x_2,x_7,y,z\rangle$ is a minus type $4$-dimensional orthogonal space, we have $X_w=\Omega_4^-(q)$.
\end{proof}

To consider factorizations $G=HK$ with $L=\Omega_7(q)$, if $H$ is a subgroup of $\Pa_1[L]$ of the form $\Omega_4^-(q)$ or $q^4.\Omega_4^-(q)$, the subsequent lemma enables us to assume $H=X_w$ or $E_w{:}X_w$, respectively.

\begin{lemma}\label{LemOmega20}
Let $A=\Pa_1[L]$ with $m=3$. Then the following statements hold:
\begin{enumerate}[{\em(a)}]
\item $A^{(\infty)}$ has a unique conjugacy class of subgroups isomorphic to $\Omega_4^-(q)$.
\item $A^{(\infty)}$ has a unique conjugacy class of subgroups of the form $q^4.\Omega_4^-(q)$.
\end{enumerate}
\end{lemma}

\begin{proof}
With the above notation we may take $A^{(\infty)}=E{:}X$. From Lemma~\ref{Lem:Omega-001} we deduce that $E$, viewed as an $X_w$-module, is the direct sum of the irreducible $X_w$-modules $E_w$ and $F$. The $X_w$-module $F$ is $1$-dimensional over $\bbF_q$ and thus is trivial. Hence $(E{:}X_w)'=E_w{:}X_w$.

First let $H$ be a subgroup of $A^{(\infty)}$ isomorphic to $\Omega_4^-(q)\cong\PSL_2(q^2)$. Then $H\cap E=1$, and so $HE/E=X_wE/E\cong\PSL_2(q^2)$. It follows that
\begin{equation}\label{EqnOmega11}
E_w{:}H=(E{:}H)'=(E{:}X_w)'=E_w{:}X_w.
\end{equation}
Consider the irreducible $\PSL_2(q^2)$-module $E_w$, which is $4$-dimensional over $\bbF_q$. It can be made a $\SL_2(q^2)$-module via the pullback of the natural homomorphism $\SL_2(q^2)\to\PSL_2(q^2)$. Accordingly, $E_w{:}X_w=\overline{Y}$, where $Y=E_w{:}\SL_2(q^2)$ and $\,\overline{\phantom{\varphi}}\,$ is the homomorphism modulo $\mathbf{Z}(Y)=\mathbf{Z}(\SL_2(q^2))=2$.
By Steinberg's twisted tensor product theorem, the irreducible $\SL_2(q^2)$-module $E_w$ is isomorphic to $M\otimes M^{(q)}$ realized over $\bbF_q$, where $M$ is the natural $\SL_2(q^2)$-module $\bbF_{q^2}^2$. Then by~\cite[Corollary~4.5]{AJL1983} we have $\mathrm{H}^1(\SL_2(q^2),E_w)=0$, which means that subgroups of $Y$ isomorphic to $\SL_2(q^2)$ are conjugate. Consequently, subgroups of $E_w{:}X_w$ isomorphic to $\PSL_2(q^2)$ are conjugate. This together with~\eqref{EqnOmega11} implies part~(a) of the lemma.

Next let $H$ be a subgroup of $A^{(\infty)}$ of the form $q^4.\Omega^-_4(q)$. To prove part~(b) of the lemma, since $\Omega_5(q)$ has a unique conjugacy class of subgroups isomorphic to $\Omega_4^-(q)$, we may assume without loss of generality that $HE/E=X_wE/E$. Then both $H$ and $E_w{:}X_w$ are complements to $F$ in $E{:}X_w$. As a consequence, $H\cong E_w{:}X_w$ is a split extension of $q^4$ by $\Omega_4^-(q)$. By the conclusion of part~(a) we then obtain
\begin{equation}\label{EqnOmega12}
H=\bfO_p(H){:}(X_w)^g=(\bfO_p(H){:}X_w)^g
\end{equation}
for some $g\in A^{(\infty)}$. Now there are two decompositions of the $X_w$-module $E$ into irreducible submodules, namely, $E=E_w\oplus F$ and $E=\bfO_p(H)\oplus F$. Since $E_w$ is not isomorphic to $F$, it follows that $\bfO_p(H)=E_w$, which substituted into~\eqref{EqnOmega12} completes the proof of part~(b).
\end{proof}

We are now in a position to classify the factorizations $G=HK$ with $H\leqslant\Pa_1[G]$. Note from Proposition~\ref{Prop:MaxOmega} that such a factorization only appears when $m=3$.

\begin{proposition}\label{Lem:Omega7-01}
Let $L=\Omega_7(q)$ and $H\leqslant\Pa_1[G]$. Then $G=HK$ if and only if $(G,H,K)$ tightly contains $(L,H^{(\infty)},K^{(\infty)})$ such that $K^{(\infty)}=\G_2(q)$ and $H^{(\infty)}=q^5{:}\Omega_5(q)$, $\Omega_5(q)$ or $q^4{:}\Omega_4^-(q)$. In this case, $H^{(\infty)}\cap K^{(\infty)}$ is given below.
\[
\begin{array}{c|ccc}
\hline
H^{(\infty)} & q^5{:}\Omega_5(q) & \Omega_5(q) & q^4{:}\Omega_4^-(q) \\
H^{(\infty)}\cap K^{(\infty)} & [q^5]{:}\SL_2(q) & \SL_2(q) & [q^3] \\
\hline
\end{array}
\]
\end{proposition}

\begin{proof}
Since $H\leqslant\Pa_1[G]$, we see from Proposition~\ref{Prop:MaxOmega} that $B^{(\infty)}=\G_2(q)$. Hence $|H|$ is divisible by $|G|/|B|=|\Omega_7(q)|/|\G_2(q)|=q^3(q^4-1)$, and $|K|$ is divisible by $|G|/|\Pa_1[G]|=(q^6-1)/(q-1)$.
It follows that $K^{(\infty)}=\G_2(q)$, and $H^{(\infty)}$ is one of the groups
\[
q^5{:}\Omega_5(q),\ \Omega_5(q),\ q^4.\Omega_4^-(q),\ \Omega_4^-(q).
\]
If $H^{(\infty)}=\Omega_4^-(q)$, then by Lemma~\ref{LemOmega20}(a) we have $H\leqslant\N_3[G]$, which leads to $G=\N_3[G]K$, contradicting Proposition~\ref{Prop:MaxOmega}. To finish the proof, we show that the first three candidates for $H^{(\infty)}$ satisfy $L=H^{(\infty)}K^{(\infty)}$ with $H^{(\infty)}\cap K^{(\infty)}$ described in the table. Thus we assume $H=H^{(\infty)}$ and $K=K^{(\infty)}=\G_2(q)$ in the remainder of the proof.

First let $H=q^5{:}\Omega_5(q)=L_v$. Then by~\cite[4.3.5]{Wilson2009} we have $H\cap K=K_v=[q^5]{:}\SL_2(q)$. It follows that $|K|/|K_v|=q^6-1=|\calP_1^{(\infty)}|$, and so $L=HK$.

Next let $H=\Omega_5(q)$. Since $\mathrm{H}^1(\Omega_5(q),\mathbb{F}_q^5)=0$ (see~\cite[Table~(4.5)]{CPS1975}), the subgroups of $\Pa_1[G]$ isomorphic to $\Omega_5(q)$ are all conjugate. Thus $H$ fixes pointwise a nondegenerate $2$-dimensional subspace of $V$.
Hence by~\cite[4.3.6]{Wilson2009} we have $H\cap K=\SL_2(q)$, which implies $L=HK$.

Finally let $H=q^4.\Omega^-_4(q)$. Adopt the notation before Lemma~\ref{Lem:Omega-001}. By the virtue of Lemma~\ref{LemOmega20}(b) we may assume $H=E_w{:}X_w$. In particular, $H=q^4{:}\Omega^-_4(q)$. For $a\in\bbF_q$, let $h(a)$ and $k(a)$ be the linear transformations on $x^\perp$ satisfying
\begin{align*}
h(a)&\colon x_7\mapsto x_7+ay+a^2x_2,\ y\mapsto y+2ax_2\\
k(a)&\colon x_3\mapsto x_3-ax_1,\ x_7\mapsto x_7-ay+a^2x_2,\ x_8\mapsto x_8+ax_6,\ y\mapsto y-2ax_2
\end{align*}
and fixing the remaining vectors in the basis $x_1,x_2,x_3,x_6,x_7,x_8,y$. It is straightforward to verify that $h(a)k(a)=\sigma(a)$ (recall the definition of $\sigma(a)$ in~\eqref{EqnOmega06}) and $h(a)\in\mathrm{O}(x^\perp)_{w,x_1,x_8,W_1}$. Then since $|h(a)|$ is odd, we deduce that $h(a)\in\Omega(x^\perp)_{w,x_1,x_8,W_1}=X_w<H$.
Note that there are precisely two conjugacy classes of subgroups $\G_2(q)$ in $L$, fused in $\mathrm{O}(x^\perp)$ (see~\cite[Table~8.40]{BHR2013}). We may assume that $K=\G_2(q)$ is the subgroup of $L$ described in~\cite[4.3.4]{Wilson2009}, so that by~\cite[Equation~(4.34)]{Wilson2009} we have $k(a)\in K$. Then as $k(a)$ fixes $x_1$, we obtain $k(a)\in K_{x_1}$. Now it follows from $\sigma(a)=h(a)k(a)$ that $F\subseteq HK_{x_1}$, and so
\begin{equation}\label{EqnOmega03}
E=E_wF\subseteq H(HK_{x_1})=HK_{x_1}.
\end{equation}
Let $\overline{\phantom{x}}\colon L_{x_1}\to L_{x_1}/E$ be the quotient modulo $E$. Then $\overline{L_{x_1}}=\Omega_5(q)$ and $\overline{H}=\Omega_4^-(q)$. From~\cite[4.3.5]{Wilson2009} we see that $\overline{K_{x_1}}=q^{1+2}{:}\SL_2(q)$. Consequently, $\overline{L_{x_1}}=\overline{H}\,\overline{K_{x_1}}$ by Example~\ref{Ex:OmegaOdd-1} (it works the same for $m=2$). This together with~\eqref{EqnOmega03} yields $L_{x_1}=HK_{x_1}$. Since $K$ is transitive on the set of singular vectors in $x^\perp$, it follows that $L=L_{x_1}K=(HK_{x_1})K=HK$. This also implies
\[
|H\cap K|=|H||K|/|L|=|q^4.\Omega_4^-(q)||\G_2(q)|/|\Omega_7(q)|=q^3,
\]
that is, $H\cap K=[q^3]$.
\end{proof}

\subsubsection{Nonsingular vectors of plus type}

To close this subsection, we classify the factorizations $G=HK$ with $K\leqslant\N_1^+[G]$. By Proposition~\ref{Prop:MaxOmega}, such a factorization only occurs when $m=3$. The classification is given in the coming proposition, which also implies that $H$ is transitive on $\calN_1^+$ if and only if it is transitive on $\calN_1^{+(\infty)}$. Recall that $\Omega_{2m+1}(q)$ is transitive on $\mathcal{N}_1^{+(\infty)}$ with stabilizer $\Omega_{2m}^+(q)$. Hence
\[
|\calN_1^{+(\infty)}|=|\Omega_{2m+1}(q)|/|\Omega_{2m}^+(q)|=q^m(q^m+1).
\]

\begin{proposition}\label{LemOmega07}
Let $L=\Omega_7(q)$ and $K\leqslant\N_1^+[G]$. Then $G=HK$ if and only if $(G,H,K)$ tightly contains $(L,H^{(\infty)},K^{(\infty)})$ such that $(H^{(\infty)},K^{(\infty)})$ is listed in the following table.
\[
\begin{array}{llll}
\hline
H^{(\infty)} & K^{(\infty)} & H^{(\infty)}\cap K^{(\infty)} & \textup{Remarks} \\
\hline
\G_2(q) & \Omega_5(q) & \SL_2(q) & K^{(\infty)}\textup{ not necessarily in }\Pa_1[L] \\
\G_2(q) & \Omega_6^+(q) & \SL_3(q) &  \\
\SU_3(q) & \Omega_6^+(q) & q^2-1 & q=3^f,\ H^{(\infty)}\textup{ in one conjugacy class of }\SU_3(q) \\
{^2}\G_2(q) & \Omega_6^+(q) & \tfrac{q-1}{2}.2 & q=3^f\textup{ with $f$ odd} \\
\hline
\end{array}
\]
In this case, $H^{(\infty)}\cap K^{(\infty)}$ is described in the table.
\end{proposition}

\begin{proof}
By Proposition~\ref{Prop:MaxOmega} we have $A^{(\infty)}=\G_2(q)$. Let $B=\N_1^+[G]$ and take $v\in\calN_1^{+(\infty)}$. We may assume $B=G_{\langle v\rangle}$, so that $B^{(\infty)}=L_v=\Omega_6^+(q)$. From~\cite[4.3.6]{Wilson2009} we see that
\begin{equation}\label{EqnOmega13}
A^{(\infty)}\cap B^{(\infty)}=(A^{(\infty)})_v=\SL_3(q).
\end{equation}
Consequently, $|A^{(\infty)}|/|A^{(\infty)}\cap B^{(\infty)}|=q^3(q^3+1)=|L|/|B^{(\infty)}|$, and so $L=A^{(\infty)}B^{(\infty)}$.
If $K^{(\infty)}=B^{(\infty)}$, then $G=AK$ with $A\cap K\geqslant\SL_3(q)$. In this case, $G=HK$ if and only if $A=H(A\cap K)$, and~\cite{HLS1987} shows that this occurs if and only if $(G,H,K)$ tightly contains $(L,H^{(\infty)},K^{(\infty)})$ with $H^{(\infty)}$ in the last three rows of the table. In what follows we assume $K^{(\infty)}<B^{(\infty)}$.

Since $B^{(\infty)}=\Omega_6^+(q)\cong\SL_4(q)/2$, it follows from Theorem~\ref{ThmLinear} and~\eqref{EqnOmega13} that $B=(A\cap B)K$ if and only if $K^{(\infty)}=\Omega_5(q)\cong\PSp_4(q)$, and in this case we have
\begin{equation}\label{EqnOmega14}
A^{(\infty)}\cap K^{(\infty)}=(A^{(\infty)}\cap B^{(\infty)})\cap K^{(\infty)}=\SL_2(q).
\end{equation}
This implies that the first row of the table gives rise to $L=H^{(\infty)}K^{(\infty)}$ and that a necessary condition for $G=HK$ is $K^{(\infty)}=\Omega_5(q)$. Moreover, if $G=HK$, then $A=H(A\cap K)$, which together with~\cite{HLS1987} and~\eqref{EqnOmega14} implies $H^{(\infty)}=A^{(\infty)}=\G_2(q)$, as in the first row of the table.
\end{proof}

\subsection{Proof of Theorem~\ref{ThmOmega}}\label{SecOmegaProof}
\ \vspace{1mm}

By Lemma~\ref{LemXia28} we may assume $L\neq\Omega_7(3)$. Then by Propositions~\ref{Prop:MaxOmega},~\ref{Lem:Omega7-01} and~\ref{LemOmega07}, there are two cases to deal with:
\begin{itemize}
\item $m=3$, $B^{(\infty)}=\G_2(q)$, and $A=\N_1^-[G]$, $\N_2^+[G]$ or $\N_2^-[G]$;
\item $A^{(\infty)}\neq\G_2(q)$, and $B=\N_1^-[G]$.
\end{itemize}

The first case is handled in the lemma below.

\begin{lemma}\label{LemOmegaRow1,3,4}
Let $L=\Omega_7(q)$, let $A=\N_1^-[G]$, $\N_2^+[G]$ or $\N_2^-[G]$, and let $B^{(\infty)}=\G_2(q)$. Then $G=HK$ if and only if $(G,H,K)$ tightly contains $(L,H^{(\infty)},K^{(\infty)})$ such that $(H^{(\infty)},K^{(\infty)})$ is listed in the following table.
\[
\begin{array}{llll}
\hline
H^{(\infty)} & K^{(\infty)} & H^{(\infty)}\cap K^{(\infty)} & \textup{Remarks} \\
\hline
q^4{:}\Omega_4^-(q) & \G_2(q) & [q^3] &  \\
\Omega_5(q) & \G_2(q) & \SL_2(q) & H^{(\infty)}\textup{ not necessarily in }\Pa_1[L] \\
\Omega_6^-(q) & \G_2(q) & \SU_3(q) &  \\
\Omega_6^-(q) & \SL_3(q) & q^2-1 & q=3^f,\ K^{(\infty)}\textup{ in one conjugacy class of }\SL_3(q) \\
\hline
\end{array}
\]
In this case, $H^{(\infty)}\cap K^{(\infty)}$ is described in the table.
\end{lemma}

\begin{proof}
Since the index $|L|/|B\cap L|=|\Omega_7(q)|/|\G_2(q)|$ is divisible by $\ppd(q^4-1)$, it suffices to consider $H\leqslant A$ such that $|H|$ is divisible by $\ppd(q^4-1)$. Since $A\in\{\N_1^-[G],\N_2^+[G],\N_2^-[G]\}$, such a subgroup $H$ satisfies either $H\leqslant\Pa_1[\Omega_6^-(q)]\leqslant\Pa_1[G]$ or $H^{(\infty)}\in\{\Omega_4^-(q),\Omega_5(q),\Omega_6^-(q)\}$.
If $H\leqslant\Pa_1[\Omega_6^-(q)]$ or $H^{(\infty)}\in\{\Omega_4^-(q),\Omega_5(q)\}$, then the lemma follows from Propositions~\ref{Lem:Omega7-01},~\ref{Prop:MaxOmega} and~\ref{LemOmega07}. Now let $H^{(\infty)}=\Omega_6^-(q)$. Then we see from~\cite[4.3.6]{Wilson2009} that $H^{(\infty)}\cap B^{(\infty)}=\SU_3(q)$, whence $L=H^{(\infty)}B^{(\infty)}$.
It follows that $G=HB$ with $H\cap B\geqslant\SU_3(q)$. Consequently, $G=HK$ if and only if $B=(H\cap B)K$. Then~\cite{HLS1987} shows that this occurs if and only if $(G,H,K)$ tightly contains $(L,H^{(\infty)},K^{(\infty)})$ with $K^{(\infty)}$ in the last two rows of the table.
\end{proof}

The second case is done by the subsequent lemma in conjunction with Proposition~\ref{Thm:gps-on-V^-}.

\begin{lemma}
Let $G=HK$ with $A^{(\infty)}\neq\G_2(q)$ and $B=\N_1^-[G]$. Then $K^{(\infty)}=\Omega_{2m}^-(q)$.
\end{lemma}

\begin{proof}
In this case, $A^{(\infty)}$ is one of the groups
\[
(q^{m(m-1)/2}.q^m){:}\SL_m(q),\ \ \PSp_6(q)\text{ with }(m,q)=(6,3^f),\ \ \F_4(q)\text{ with }(m,q)=(12,3^f).
\]
For each of these candidates of $A^{(\infty)}$, direct calculation shows that $|L|/|A^{(\infty)}|$ is divisible by
\[
r\in\ppd(q^{2m}-1)\ \text{ and }\ s\in\ppd(q^{2m-2}-1).
\]
Hence $|K|$ is divisible by $rs$. Let $N=\Rad(B)$, and let $\,\overline{\phantom{\varphi}}\,$ be the quotient from $B$ to $B/N$. Then $\overline{B}$ is an almost simple group with socle $\POm_{2m}^-(q)$, and $\overline{K}$ is a factor of order divisible by $rs$ in the factorization $\overline{B}=(\overline{A\cap B})\overline{K}$. An inspection of~\cite[Theorem~A]{LPS1990} and~\cite{LPS1996} shows that $\overline{B}$ has no $\max^-$ factorizations with the order of a factor divisible by $rs$. Thus $\overline{K}\geqslant\Soc(\overline{B})=\POm_{2m}^-(q)$, and so $K^{(\infty)}=B^{(\infty)}=\Omega_{2m}^-(q)$, as required.
\end{proof}

%%%%%%%%%%%%%%%%%%%%%%%%%%%%%%%%%%%%%%%%%%%%%%%%%%%%%%%%%%%%%%%%%%%
%%%%%%%%%%%%%%%%%%%%%%%%%%%%%%%%%%%%%%%%%%%%%%%%%%%%%%%%%%%%%%%%%%%

\section{Orthogonal groups of minus type}

%%%%%%%%%%%%%%%%%%%%%%%%%%%%%%%%%%%%%%%%%%%%%%%%%%%%%%%%%%%%%%%%%%%
%%%%%%%%%%%%%%%%%%%%%%%%%%%%%%%%%%%%%%%%%%%%%%%%%%%%%%%%%%%%%%%%%%%

In this section we embark on orthogonal groups of minus type. Throughout this section, let $G$ be a classical group with $L=G^{(\infty)}=\Omega_{2m}^-(q)$, where $m\geqslant4$, and let $H$ and $K$ be nonsolvable subgroups of $G$ not containing $G^{(\infty)}$.

\begin{theorem}\label{ThmOmegaMinus}
With the above notation, $G=HK$ if and only if $(G,H,K)$ tightly contains some triple $(G_0,H_0,K_0)$ in Table~$\ref{TabOmegaMinus}$.
\end{theorem}

Some notation throughout this section is set up in the following Subsection~\ref{SecOmegaMinus01}, and the proof of Theorem~\ref{ThmOmegaMinus} will be given in Subsection~\ref{SecOmegaMinusProof}.

\subsection{Setup}\label{SecOmegaMinus01}
\ \vspace{1mm}

Let $H$ and $K$ be nonsolvable subgroups of $G$ such that neither $H$ nor $K$ contains $L$, and let $A$ and $B$ be $\max^-$ subgroups of $G$ containing $H$ and $K$ respectively. Then $A$ and $B$ are nonsolvable.
The $\max^-$ factorizations $G=AB$ are classified by Liebeck-Praeger-Saxl~\cite{LPS1990,LPS1996} and presented in the following proposition.

\begin{proposition}\label{Prop:MaxOmegaminu}
If $G=HK$, then with the notation defined above and interchanging $A$ and $B$ if necessary, one of the following holds:
\begin{enumerate}[{\rm(a)}]
\item $B=\Pa_1[G]$ with $m$ odd, and $A^{(\infty)}=\SU_m(q)$, or $\A_{12}$ with $L=\Omega_{10}^-(2)$;
\item $B=\N_1[G]$, and $A^{(\infty)}=\SU_m(q)$ with $m$ odd, or $\Omega^-_m(q^2)$ with $q\in\{2,4\}$ and $G=\GaO^-_{2m}(q)$;
\item $B=\N^+_2[G]$ with $m$ odd and $G=\GaO^-_{2m}(4)$, and $A=\GU_m(4).4$.
\end{enumerate}
\end{proposition}

%{\small\begin{table}[htbp]
%\caption{Maximal factorizations with socle $L=\POm_{2m}^-(q)$, $m\geqslant4$}\label{TabMaxOmegaMinus}
%\centering
%\begin{tabular}{|l|l|l|l|}
%\hline
%Row & $A\cap L$ & $B\cap L$ & Remark\\
%\hline
%1 & $\lefthat\GU_m(q)$ & $\Pa_1$ & $m$ odd \\
%2 & $\lefthat\GU_m(q)$ & $\N_1$ & $m$ odd \\
%3 & $\Omega_m^-(q^2).2$ & $\N_1$ & $m$ even, $q\in\{2,4\}$, $G=\Aut(L)$ \\
%4 & $\GU_m(4)$ & $\N_2^+$ & $m$ odd, $q=4$, $G=\Aut(L)$ \\
%5 & $\A_{12}$ & $\Pa_1$ & $m=5$, $q=2$ \\
%\hline
%\end{tabular}
%\end{table}}

We observe from Proposition~\ref{Prop:MaxOmegaminu} that there are three possibilities for $B\cap L$, namely, $\Pa_1[L]$, $\N_1[L]$ and $\N_2^+[L]$. Our analysis will proceed in the coming three subsections accordingly.

% For the case $L=\Omega_{10}^-(2)$, computation in \magma~\cite{BCP1997} verifies the following:
%
% \begin{lemma}\label{Omega^-(10,2)}
% For $L=\Omega_{10}^-(2)$, $G=HK$ if and only if $(G,H,K)$ tightly contains some $(G_0,H_0,K_0)$ in the following table:
% \[\begin{array}{|l|l|l|}\hline
% H_0 & K_0 & H_0\cap K_0 \\
% \hline
% \SU_5(2),\ \M_{12},\ \A_{12} & 2^8{:}\Omega_8^-(2) & 2^{1+6}{:}\SU_3(2),\ 2_+^{1+4}.\Sy_3,\ (\A_4\times\A_8).2 \\
% \SU_5(2) & \Omega_8^-(2) & \SU_3(2) \\
% ???
% \hline
% %3.\J_3&2^{2+4}.(3\times\Sy_3) &(m,q)=(9,2) \\
% %\hline
% \end{array}\]
% \end{lemma}

Throughout this section, let $q=p^f$ be a power of a prime $p$, and let $V=\bbF_q^{2m}$ be equipped with a nondegenerate quadratic form $Q$ of minus type, associated with bilinear form $\beta$.
For $w\in V$ with $Q(w)\neq0$, let $r_w$ be the reflection with respect to $w$ defined by
\[
r_w\colon V\to V,\quad v\mapsto v-\frac{\beta(v,w)}{Q(w)}w.
\]
We fix the notation for some field-extension subgroups of $\mathrm{O}(V)$ as follows, according to the parity of $m$. Let $V_\sharp$ be a vector space of dimension $m$ over $\bbF_{q^2}$ with the same underlying set as $V$, and let $\Tr$ be the trace of the field extension $\bbF_{q^2}/\bbF_q$.
Take $\mu\in\bbF_{q^2}$ such that the polynomial $x^2+x+\mu$ is irreducible over $\bbF_{q^2}$.

First, assume that $m=2\ell+1$ is odd. Equip $V_\sharp$ with a nondegenerate unitary form $\beta_\sharp$ such that $Q(v)=\beta_\sharp(v,v)$ for all $v\in V$, and then $\GU(V_\sharp)<\mathrm{O}(V)$. Take a standard basis $E_1,F_1,\dots,E_\ell,F_\ell,D$ for $V_\sharp$, so that
\[
\beta_\sharp(E_i,E_j)=\beta_\sharp(F_i,F_j)=\beta_\sharp(E_i,D)=\beta_\sharp(F_i,D)=0,\quad\beta_\sharp(E_i,F_j)=\delta_{i,j},\quad\beta_\sharp(D,D)=1
\]
for all $i,j\in\{1,\dots,\ell\}$. Let $\psi\in\GaU(V_\sharp)$ such that
\[
\psi\colon a_1E_1+b_1F_1+\dots+a_\ell E_\ell+b_\ell F_\ell+cD\mapsto a_1^pE_1+b_1^pF_1+\dots+a_\ell^pE_\ell+b_\ell^pF_\ell+c^pD
\]
for $a_1,b_1\dots,a_\ell,b_\ell,c\in\bbF_{q^2}$, and let $\lambda\in\bbF_{q^2}$ such that $\lambda+\lambda^q=1$. Then $(\lambda E_1,F_1)$ is a hyperbolic pair with respect to $Q$. Thus there exists a standard basis $e_1,f_1,\dots,e_{m-1},f_{m-1},d,d'$ for $V$ as in~\cite[2.2.3]{LPS1990} such that
\[
e_1=\lambda E_1\ \text{ and }\ f_1=F_1.
\]
It follows that $Q(e_1+f_1)=\beta_\sharp(e_1+f_1,e_1+f_1)=1$, and so $e_1+f_1$ is a nonsingular vector in both $V$ (with respect to $Q$) and $V_\sharp$ (with respect to $\beta_\sharp$). Let $\phi\in\GaO(V)$ be as defined in~\cite[\S2.8]{KL1990} with respect to the basis $e_1,f_1,\dots,e_{m-1},f_{m-1},d,d'$.
Then $\phi$ has order $2f$, fixes $e_1,f_1,\dots,e_{m-1},f_{m-1},d$ and commutes with $r_{e_1+f_1}$.

Next, assume that $m=2\ell$ is even. Equip $V_\sharp$ with a nondegenerate quadratic form $Q_\sharp$ of minus type such that $Q(v)=\Tr(Q_\sharp(v))$ for all $v\in V$, and then $\mathrm{O}(V_\sharp)<\mathrm{O}(V)$. Let $\beta_\sharp$ be the associated symmetric bilinear form of $Q_\sharp$, and take a standard basis $E_1,F_1,\dots,E_{\ell-1},F_{\ell-1},D,D'$ for $V_\sharp$, so that for all $i,j\in\{1,\dots,\ell-1\}$,
\begin{align*}
&Q_\sharp(E_i)=Q_\sharp(F_i)=0,\quad Q_\sharp(D)=1=\beta_\sharp(D,D'),\quad Q_\sharp(D')=\mu,\quad\beta_\sharp(E_i,F_j)=\delta_{i,j},\\
&\beta_\sharp(E_i,E_j)=\beta_\sharp(F_i,F_j)=\beta_\sharp(E_i,D)=\beta_\sharp(F_i,D)=\beta_\sharp(E_i,D')=\beta_\sharp(F_i,D')=0.
\end{align*}

\subsection{Actions on $\calP_1$}
\ \vspace{1mm}

Let $\mathcal{P}_1^{(\infty)}$ be the set of nonzero singular vectors of $V$. By Witt's Lemma, $L=\Omega_{2m}^-(q)$ is transitive on $\mathcal{P}_1^{(\infty)}$. For any $v\in\mathcal{P}_1^{(\infty)}$ we have $L_v=q^{2m-2}{:}\Omega_{2m-2}^-(q)<\Pa_1[L]$, and hence
\[
|\mathcal{P}_1^{(\infty)}|=|L|/|L_v|=(q^m+1)(q^{m-1}-1).
\]
Each transitive subgroup $H$ of $G$ on $\mathcal{P}_1^{(\infty)}$ gives rise to a factorization $G=HG_v$.

\begin{lemma}\label{LemOmegaMinusRow1,5}
Let $G=HK$ with $B=\Pa_1[G]$. Then $m$ is odd, and one of the following holds:
\begin{enumerate}[{\rm (a)}]
\item $K^{(\infty)}=q^{2m-2}{:}\Omega_{2m-2}^-(q)$, and $H^{(\infty)}$ satisfies one of:
\begin{itemize}
\item $H^{(\infty)}=\SU_m(q)$;
\item $H^{(\infty)}\in\{\A_{12},\M_{12}\}$ with $(m,q)=(5,2)$;
\item $H^{(\infty)}=3^{\boldsymbol{\cdot}}\J_3<\SU_9(2)$ with $(m,q)=(9,2)$.
\end{itemize}
\item $K^{(\infty)}=\Omega_{2m-2}^-(q)$, and $H^{(\infty)}=\SU_m(q)$;
\end{enumerate}
\end{lemma}

\begin{proof}
For $(m,q)=(5,2)$, computation in \magma~\cite{BCP1997} verifies the lemma. Thus we assume $(m,q)\neq(5,2)$. By Proposition~\ref{Prop:MaxOmegaminu} we have
$m\geqslant5$ is odd, $A^{(\infty)}=\SU_m(q)$ and
\begin{equation}\label{EqnXia7}
B^{(\infty)}=q^{2m-2}{:}\Omega^-_{2m-2}(q).
\end{equation}
It follows that $|L|/|A\cap L|$ is divisible by $r\in\ppd(q^{2m-2}-1)$ and $s\in\ppd(q^{m-2}-1)$, whence $|K|$ is divisible by $rs$. Moreover, $|L|/|B\cap L|$ is divisible by $t\in\ppd(q^{2m}-1)$, and so is $|H|$.

We first describe $K$. Let $N=\Rad(B)$. Then $B/N$ is an almost simple group with socle $\POm_{2m-2}^-(q)$.
%By Lemmas~\ref{LemXia3} and~\ref{LemXia7}, the group $KN/N$ is a factor of $B/N$.
Since $|K|$ is divisible by $rs$, so is $|KN/N|$ by Lemma~\ref{LemXia22}. Inspecting the $\max^-$ factorizations of groups with socle $\POm_{2m-2}^-(q)$, we conclude that there is no factor with order divisible by $rs$ in any core-free factorization of $B/N$. Consequently, $KN/N\trianglerighteq\Soc(B/N)=\POm_{2m-2}^-(q)$. Note in~\eqref{EqnXia7} that the conjugation action of $\Omega_{2m-2}^-(q)$ on $\mathbf{O}_p(B^{(\infty)})=q^{2m-2}$ is irreducible. We then obtain $K^{(\infty)}=q^{2m-2}{:}\Omega_{2m-2}^-(q)$ or $\Omega^-_{2m-2}(q)$.
%\begin{itemize}
%\item $K\cap\mathbf{O}_p(B^{(\infty)})\neq1$, and $K^{(\infty)}=q^{2m-2}{:}\Omega_{2m-2}^-(q)$, or
%\item $K\cap\mathbf{O}_p(B^{(\infty)})=1$, and $K^{(\infty)}=\Omega^-_{2m-2}(q)$.
%\end{itemize}

We now consider $H$ in its overgroup $A$. Let $M=\Rad(A)$. Then $A/M$ is an almost simple group with socle $\PSU_m(q)$. Since $|H|$ is divisible by $t\in\ppd(q^{2m}-1)$, so is $|HM/M|$. Then since $m\geqslant5$ is odd and $A/M$ is the product of $HM/M$ and the core-free subgroup $(A\cap B)M/M$, we derive from Theorem~\ref{ThmUnitary} and~\cite[Theorem~1.1]{LX} that one of the following occurs:
\begin{enumerate}[{\rm (i)}]
\item $HM/M\trianglerighteq\Soc(A/M)=\PSU_m(q)$;
\item $(m,q)=(9,2)$ and $(HM/M)\cap\Soc(A/M)=\J_3$.
\end{enumerate}
For~(i), it follows that $H^{(\infty)}=\SU_m(q)$, and the conclusion of the lemma holds. Next, assume that~(ii) occurs. Then $H^{(\infty)}=3^{\boldsymbol{\cdot}}\J_3$ is a maximal subgroup of $A^{(\infty)}=\SU_9(2)$ (see~\cite[Table~8.57]{BHR2013}). Since $H$ is an irreducible subgroup of $G$, we deduce that $|H|_2\leqslant|H^{(\infty)}|_2|\Out(L)|_2=2|\J_3|_2=2^8$. As $G=HK$ implies $|H|_2|K|_2\geqslant|G|_2\geqslant2^{72}$,
it follows that $|K|_2\geqslant2^{64}$. Hence $K^{(\infty)}\neq\Omega^-_{16}(2)$, and so $K^{(\infty)}=2^{16}{:}\Omega_{16}^-(2)$, as in part~(a) of the lemma.
\end{proof}

To give examples of factorizations $G=HK$ satisfying Lemma~\ref{LemOmegaMinusRow1,5}(a), we adopt the notation defined after Proposition~\ref{Prop:MaxOmegaminu}. In particular, for odd $m$, we have singular vectors $e_1\in V=\bbF_q^{2m}$ and $E_1\in V_\sharp=\bbF_{q^2}^m$ such that $e_1=\lambda E_1$, where $\lambda\in\bbF_{q^2}$ with $\lambda+\lambda^q=1$.

\begin{example}\label{ex:OmegaMinus-1}
Let $H=\SU(V_\sharp)=\SU_m(q)<L$ with $m$ odd. Since $e_1$ and $E_1=\lambda^{-1}e_1$ are singular vectors in $V$ and $V_\sharp$, respectively, we have
\[
L_{e_1}=q^{2m-2}{:}\Omega_{2m-2}^-(q)\ \text{ and }\ H_{e_1}=H_{E_1}=\SU(V_\sharp)_{E_1}=q^{1+(2m-4)}{:}\SU_{m-2}(q).
\]
Hence $|H|/|H_{e_1}|=(q^m+1)(q^{m-1}-1)=|\mathcal{P}_1^{(\infty)}|$, and so $H$ is transitive on $\mathcal{P}_1^{(\infty)}$. This gives
\[
L=HL_{e_1}=\SU_m(q)(q^{2m-2}{:}\Omega_{2m-2}^-(q))
\]
with the intersection $H\cap L_{e_1}=H_{e_1}=q^{1+(2m-4)}{:}\SU_{m-2}(q)$.
\end{example}

Since $\Omega_{2m-2}^-(q)$ has trivial first cohomology on the natural module (see~\cite[Theorem~2.14]{KantorLiebler1982}),
subgroups of $\Pa_1[L]$ isomorphic to $\Omega_{2m-2}^-(q)$ are all conjugate in $L$. Thus in part~(b) of Lemma~\ref{LemOmegaMinusRow1,5}, $K$ stabilizes a nonsingular $2$-subspace of plus type, which will be treated in Subsection~\ref{nonsingular 2-spaces}. In the ensuing proposition we classify the factorizations $G=HK$ described in part~(a) of Lemma~\ref{LemOmegaMinusRow1,5}, namely, those with $\Pa_1[G]^{(\infty)}\leqslant K\leqslant\Pa_1[G]$. It turns out that $H$ is transitive on $\calP_1$ if and only if it is transitive on $\calP_1^{(\infty)}$.

\begin{proposition}\label{prop:OmegaMinus-1}
Let $K\leqslant\Pa_1[G]$ with $K^{(\infty)}=q^{2m-2}{:}\Omega_{2m-2}^-(q)$. Then $G=HK$ if and only if $(G,H,K)$ tightly contains $(L,H^{(\infty)},K^{(\infty)})$ such that $H^{(\infty)}$ lies in the following table. In this case, $H^{(\infty)}\cap K^{(\infty)}$ is described in the table.
\[
\begin{array}{lll}
\hline
H^{(\infty)} & H^{(\infty)}\cap K^{(\infty)} & \textup{Remarks} \\
\hline
\SU_m(q) & q^{1+(2m-4)}{:}\SU_{m-2}(q) & m~\textup{odd} \\
\A_{12} & (\A_8\times\A_4){:}2 & L=\Omega_{10}^-(2) \\
\M_{12} & 2_+^{1+4}{:}\Sy_3 & L=\Omega_{10}^-(2) \\
3^{\boldsymbol{\cdot}}\J_3 & 2^{2+4}{:}(\Sy_3\times3) & L=\Omega_{18}^-(2),\ H^{(\infty)}<\SU_9(2) \\
\hline
\end{array}
\]
\end{proposition}

\begin{proof}
Computation in \magma~\cite{BCP1997} confirms the proposition for $(m,q)=(5,2)$. Thus we assume that $(m,q)\neq(5,2)$.
By Lemma~\ref{LemOmegaMinusRow1,5} and Example~\ref{ex:OmegaMinus-1}, we only need to show that, when $(m,q)=(9,2)$, the subgroup $H=3^{\boldsymbol{\cdot}}\J_3$ of $\SU_9(2)$ is transitive on $\mathcal{P}_1^{(\infty)}$ with stabilizer $2^{2+4}{:}(\Sy_3\times3)$.

Let $L=\Omega_{18}^-(2)$, $M=\SU_9(2)$ and $H=3^{\boldsymbol{\cdot}}\J_3<M$. We see from Example~\ref{ex:OmegaMinus-1} that $M_{e_1}=2^{1+14}{:}\SU_7(2)$. Then since $M_{e_1}\ngeqslant\Z(M)=3$, it follows from Lemma~\ref{Lem:UnitaryProof} that
\[
H_{e_1}=H\cap M_{e_1}\cong(H/\Z(M))\cap(M_{e_1}\Z(M)/\Z(M))=2^{2+4}{:}(\Sy_3\times3).
\]
Hence $|H|/|H_{e_1}|=(2^9+1)(2^8-1)=|\mathcal{P}_1^{(\infty)}|$, and so $H$ is transitive on $\mathcal{P}_1^{(\infty)}$.
\end{proof}

\subsection{Actions on $\calN_1$}
\ \vspace{1mm}

Let $\mathcal{N}_1^{(\infty)}=\{u\in V\mid Q(u)=1\}$. Then $L=\Omega_{2m}^-(q)$ is transitive on $\mathcal{N}_1^{(\infty)}$ by Witt's Lemma, and the stabilizer is $\Omega_{2m-1}(q)$. Hence
\[
|\mathcal{N}_1^{(\infty)}|=|\Omega_{2m}^-(q)|/|\Omega_{2m-1}(q)|=q^{m-1}(q^m+1).
\]
The next lemma is proved by the same quotient and primitive prime divisor argument as Lemma~\ref{LemOmegaMinusRow1,5}, except that~\cite[Lemma~4.4]{LPS2010} also comes into play when determining $H^{(\infty)}$.

\begin{lemma}\label{LemOmegaMinusRow2,3}
Let $G=HK$ with $B=\N_1[G]$. Then one of the following holds:
\begin{enumerate}[{\rm(a)}]
\item $H^{(\infty)}=\SU_m(q)$ with $m$ odd, and $K^{(\infty)}=\Omega_{2m-1}(q)$ or $\Omega_{2m-2}^-(q)$;
\item $G=\GaO^-_{2m}(q)$ with $q\in\{2,4\}$, $H^{(\infty)}=\SU_{m/2}(q^2)$ ($m/2$ odd) or $\Omega_m^-(q^2)$, and $K^{(\infty)}=\Omega_{2m-1}(q)$;
\item $G=\mathrm{O}_{2m}^-(2)$, $H^{(\infty)}=\SU_{m/4}(q^4)$ ($m/4$ odd) or $\Omega_{m/2}^-(q^4)$, and $K^{(\infty)}=\Omega_{2m-1}(q)$.
\end{enumerate}
\end{lemma}

Next we construct examples of factorizations $G=HK$ corresponding to the candidates in Lemma~\ref{LemOmegaMinusRow2,3}. Note that $e_1+f_1\in\mathcal{N}_1^{(\infty)}$.

\begin{example}\label{ex:OmegaMinus01}
Let $H=\SU(V_\sharp)=\SU_m(q)$ with $m$ odd, and let $v=e_1+f_1$. Since $v$ is a nonsingular vector in $V_\sharp$ (with respect to $\beta_\sharp$), we have $H_v=\SU(V_\sharp)_v=\SU_{m-1}(q)$. It follows that $|H|/|H_v|=q^{m-1}(q^m+1)=|\mathcal{N}_1^{(\infty)}|$, and thus $H$ is transitive on $\mathcal{N}_1^{(\infty)}$. This leads to a factorization $L=HL_v$ with $L_v=\Omega_{2m-1}(q)$ and $H\cap L_v=\SU_{m-1}(q)$.
\end{example}

When $m$ is even, recall that we equip $V_\sharp$ with a nondegenerate quadratic form $Q_\sharp$ of minus type such that $Q(v)=\Tr(Q_\sharp(v))$ for all $v\in V$, and thus $\mathrm{O}(V_\sharp)<\mathrm{O}(V)$. Also recall the associated symmetric bilinear form $\beta_\sharp$ and the standard $\bbF_{q^2}$-basis $E_1,F_1,\dots,E_{\ell-1},F_{\ell-1},D,D'$.
The next lemma gives a maximal factorization $G=HK$ described in~\cite[Theorem~A]{LPS1990}.

\begin{lemma}\label{ex:OmegaMinus02}
For $q\in\{2,4\}$ and $m$ even, let $G=\GaO(V)=\GaO_{2m}^-(q)$, $H=\GaO(V_\sharp)=\GaO_m^-(q^2)$ and $K=G_{D'}=\Omega_{2m-1}(q){:}(2f)$. Then $G=HK$ with $H\cap K=\mathrm{O}(V_\sharp)\cap K=\Omega_{m-1}(q^2){:}2$.
\end{lemma}

For $q\in\{2,4\}$, combining the factorizations $\GaO_{4\ell}^-(q)=\GaO_{2\ell}^-(q^2)(\Omega_{4\ell-1}(q).(2f))$ in Lemma~\ref{ex:OmegaMinus02} and $\GaO_{2\ell}^-(q^2)=(\SU_\ell(q^2){:}(4f))(\Omega_{2\ell-1}(q^2){:}2)$ implied by Example~\ref{ex:OmegaMinus01}, we obtain the factorization
\[\GaO_{4\ell}^-(q)=(\SU_\ell(q^2){:}(4f))\,(\Omega_{4\ell-1}(q).(2f))
\]
in the example as shown below.

\begin{example}\label{ex:OmegaMinus05}
For $q\in\{2,4\}$ and $m=2\ell$ with odd $\ell$, let $G=\GaO(V)=\GaO_{2m}^-(q)$, let $H=\SU_\ell(q^2){:}(4f)<\GaO_m^-(q^2)<G$, and let $K=G_{D'}=\Omega_{2m-1}(q){:}(2f)$. Then $G=HK$ with $H\cap K=\SU_{\ell-1}(q^2){:}2$.
\end{example}

We now classify the factorizations $G=HK$ with $K^{(\infty)}=\Omega_{2m-1}(q)$.

\begin{proposition}\label{prop:OmegaMinus-02}
Let $K\leqslant\N_1[G]$ with $K^{(\infty)}=\Omega_{2m-1}(q)$. Then $G=HK$ if and only if $(G,H,K)$ tightly contains some $(G_0,H_0,K_0)$ in the following table. In this case, $H_0\cap K_0$ is described in the table.
\[
\begin{array}{lllll}
\hline
G_0 & H_0 & K_0 & H_0\cap K_0 & \textup{Remarks} \\
\hline
\Omega_{2m}^-(q) &\SU_m(q) & \Omega_{2m-1}(q) & \SU_{m-1}(q) & m\textup{ odd}  \\
\GaO_{2m}^-(q) & \GaO_m^-(q^2) & \Sp_{2m-2}(q){:}(2f) & \Sp_{m-2}(q^2){:}2 & q\in\{2,4\} \\
\GaO_{2m}^-(q) & \SU_{m/2}(q^2){:}(4f) & \Sp_{2m-2}(q){:}(2f) & \SU_{m/2-1}(q^2){:}2 & m/2\textup{ odd},\ q\in\{2,4\} \\
\hline
\end{array}
\]
\end{proposition}

\begin{proof}
By Example~\ref{ex:OmegaMinus01}, Lemma~\ref{ex:OmegaMinus02} and Example~\ref{ex:OmegaMinus05}, it suffices to prove the ``only if'' part. Let $G=HK$ with $\Omega_{2m-1}(q)\leqslant K\leqslant\N_1[G]$. Then the candidates for $H^{(\infty)}$ are given in Lemma~\ref{LemOmegaMinusRow2,3}. For $H^{(\infty)}=\SU_m(q)$ ($m$ odd) as in Lemma~\ref{LemOmegaMinusRow2,3}(a), the triple $(G,H,K)$ tightly contains $(G^{(\infty)},H^{(\infty)},K^{(\infty)})=(\Omega_{2m}^-(q),\SU_m(q),\Omega_{2m-1}(q))$, which is the triple $(G_0,H_0,K_0)$ in the first row of the table.

Next, assume that Lemma~\ref{LemOmegaMinusRow2,3}(b) holds. Then $G=\GaO^-_{2m}(q)$ with $q\in\{2,4\}$, and $H^{(\infty)}=\SU_{m/2}(q^2)$ ($m/2$ odd) or $\Omega_m^-(q^2)$. Since $G=\Aut(L)=L.(2f)$, it follows that the triples $(G_0,H_0,K_0)$ in the second and third rows of the table give rise to minimal factorizations $G_0=H_0K_0$ with respect to tight containment.

Finally, we prove that Lemma~\ref{LemOmegaMinusRow2,3}(c) cannot occur. Suppose for a contradiction that it occurs. Then $G=\GaO_{2m}^-(2)$, and $H^{(\infty)}=\SU_{m/4}(q^4)$ ($m/4$ odd) or $\Omega_{m/2}^-(q^4)$. In particular, $H\leqslant\Nor_G(H^{(\infty)})\leqslant\GaO^-_{m/2}(16)$. To derive a contradiction from $G=HK$, we may assume that $H=\GaO^-_{m/2}(16)$ and $K=\N_1[G]$. It follows that $H\cap K=\N_1[H]\geqslant\GO_{m/2-1}(16)$ (see \cite[Proposition~4.1.7]{KL1990}), and so
\[
|H|_2/|H\cap K|_2\leqslant|\GaO^-_{m/2}(16)|_2/|\GO_{m/2-1}(16)|_2=2^{m-2}<2^{m-1}=|G|_2/|K|_2,
\]
contradicting the requirement $|H|/|H\cap K|=|G|/|K|$ from the factorization $G=HK$. %This completes the proof.
\end{proof}

\subsection{Actions on $\calN_2^+$}\label{nonsingular 2-spaces}
\ \vspace{1mm}

In this subsection we handle the cases identified in Proposition~\ref{Prop:MaxOmegaminu}(c) and Lemma~\ref{LemOmegaMinusRow1,5}(b). The following lemma gives the candidates of $H^{(\infty)}$ and $K^{(\infty)}$ for Proposition~\ref{Prop:MaxOmegaminu}(c), namely, the case when $G=\GaO^-_{2m}(4)$ with $m$ odd, $A=\GU_m(4).4$ and $B=\N_2^+[G]$.

\begin{lemma}\label{LemOmegaMinusRow4}
Let $G=\GaO^-_{2m}(4)$ with $m$ odd, and let $G=HK$ with $A=\GU_m(4).4$ and $B=\N_2^+[G]$. Then $H^{(\infty)}=A^{(\infty)}=\SU_m(4)$, and $K^{(\infty)}=B^{(\infty)}=\Omega^-_{2m-2}(4)$.
\end{lemma}

\begin{proof}
Since $|L|/|A\cap L|$ is divisible by $r\in\ppd(4^{2m-2}-1)$ and $s\in\ppd(4^{m-2}-1)$, so is $|K|$. Let $N=\Rad(B)$. Then $B/N$ is almost simple with socle $\Omega_{2m-2}^-(4)$ and is the product of $(A\cap B)N/N$ and $KN/N$, where $|KN/N|$ is divisible by $rs$.
By the classification of $\max^+$ factorizations of almost simple groups with socle $\Omega_{2m-2}^-(4)$, which only appears in part~(b) of Proposition~\ref{Prop:MaxOmegaminu} with $m$ replaced by $m-1$ (as $m$ is odd), $B/N$ has no core-free factor of order divisible by $rs$. Hence $KN/N\trianglerighteq\Soc(B/N)=\Omega_{2m-2}^-(4)$, and so $K^{(\infty)}=B^{(\infty)}=\Omega^-_{2m-2}(4)$.

Let $M=\Rad(A)$. Then $A/M$ is an almost simple group with socle $\PSU_m(4)$ and is the product of $HM/M$ and $(A\cap B)M/M$. By~\cite[Theorem~A]{LPS1990}, there is no factorization of $A/M$ with core-free factors. Therefore, $HM/M\trianglerighteq\Soc(A/M)=\PSU_m(4)$, and hence $H^{(\infty)}=A^{(\infty)}=\SU_m(4)$.
\end{proof}

The next two lemmas show the existence of factorizations $G=HK$ with $H^{(\infty)}=\SU_m(q)$ and $K^{(\infty)}=\Omega^-_{2m-2}(q)$ for $q\in\{2,4\}$. These are desired by Lemma~\ref{LemOmegaMinusRow4} and also give examples in Lemma~\ref{LemOmegaMinusRow1,5}(b) and Lemma~\ref{LemOmegaMinusRow2,3}(a). Recall the notation of the field-extension subgroups of $\mathrm{O}(V)$ for odd $m$ given in Subsection~\ref{SecOmegaMinus01}.

\begin{lemma}\label{LemOmegaMinus03}
Let $G=\Omega_{2m}^-(2)$ with $m$ odd, let $H=\SU(V_\sharp)=\SU_m(2)$, and let $K=G_{\{e_1,f_1\}}=\Omega_{2m-2}^-(2){:}2$. Then $G=HK$ with $H\cap K=\SU_{m-2}(2)$.
\end{lemma}

\begin{proof}
Since $(e_1,f_1)$ is a hyperbolic pair with respect to $Q$, we have $G_{e_1,f_1}=\Omega_{2m-2}^-(2)$ and $K=G_{\{e_1,f_1\}}=G_{e_1,f_1}{:}\langle s\rangle$ for some involution $s$ swapping $e_1$ and $f_1$.
Suppose that $H_{\{e_1,f_1\}}\neq H_{e_1,f_1}$; that is, there exists $t\in H$ swapping $e_1=\lambda E_1$ and $f_1=F_1$. Then as $t\in H=\SU(V_\sharp)$,
\[
\lambda=\beta_\sharp(\lambda E_1,F_1)=\beta_\sharp((\lambda E_1)^t,F_1^t)=\beta_\sharp(F_1,\lambda E_1)=\lambda^q,
\]
contradicting the condition $\lambda+\lambda^q=1$. Thus we conclude that $H_{\{e_1,f_1\}}=H_{e_1,f_1}$, and so
\[
H\cap K=G_{\{e_1,f_1\}}\cap H=H_{\{e_1,f_1\}}=H_{e_1,f_1}=\SU(V_\sharp)_{e_1,f_1}=\SU(V_\sharp)_{\lambda E_1,F_1}=\SU_{m-2}(2).
\]
which yields $G=HK$.
\end{proof}

\begin{lemma}\label{LemOmegaMinus04}
Let $G=\GaO_{2m}^-(q)$ with $q\in\{2,4\}$ and $m$ odd, let $H=\SU(V_\sharp){:}\langle\psi\rangle=\SU_m(q){:}(2f)$, and let $K=\Omega(V)_{e_1,f_1}{:}\langle\rho\rangle=\Omega_{2m-2}^-(q){:}(2f)$ with $\rho\in\{\phi,(r_{e_1+f_1})^{2/f}\phi\}$. Then $G=HK$ with $H\cap K=\SU_{m-2}(q)$.
\end{lemma}

\begin{proof}
Suppose that $H\cap K\nleqslant\Omega(V)$. Then there exist $s\in\Omega(V)_{e_1,f_1}$ and $t\in\SU(V_\sharp)$ such that $\phi^fs=\psi^ft$. Since $\phi$ and $s$ both fix $e_1$ and $f_1$, it follows that $\phi^fs$ fixes $e_1=\lambda E_1$ and $f_1=F_1$. This together with $t\in\SU(V_\sharp)$ implies that
\begin{align*}
\beta_\sharp(\lambda E_1,F_1)=\beta_\sharp((\lambda E_1)^{\phi^fs},F_1^{\phi^fs})=\beta_\sharp((\lambda E_1)^{\psi^ft},F_1^{\psi^ft})=\beta_\sharp((\lambda^qE_1)^t,F_1^t)=\beta_\sharp(\lambda^qE_1,F_1),
\end{align*}
contradicting the condition $\lambda+\lambda^q=1$. Thus $H\cap K\leqslant\Omega(V)$. Accordingly,
\[
H\cap K=H\cap(K\cap\Omega(V))=H\cap\Omega(V)_{e_1,f_1}=\SU(V_\sharp)_{e_1,f_1}=\SU(V_\sharp)_{\lambda E_1,F_1}=\SU_{m-2}(q),
\]
and so $G=HK$.
\end{proof}

To finish this subsection, we classify the factorizations $G=HK$ with $H^{(\infty)}=\SU_m(q)$ ($m$ odd) and $K^{(\infty)}=\Omega_{2m-2}^-(q)$.

\begin{proposition}\label{LemOmegaMinus09}
Let $H^{(\infty)}=\SU_m(q)$ with $m$ odd, and let $K^{(\infty)}=\Omega_{2m-2}^-(q)$. Then $G=HK$ if and only if $(G,H,K)$ tightly contains some $(G_0,H_0,K_0)$ in the following table. In this case, $H_0\cap K_0$ is described in the table.
\[
\begin{array}{lllll}
\hline
G_0 & H_0 & K_0 & H_0\cap K_0 & \textup{Remarks}\\
\hline
\Omega_{2m}^-(2) & \SU_m(2) & \Omega_{2m-2}^-(2){:}2 & \SU_{m-2}(2) & \textup{as in~\ref{LemOmegaMinus03}}\\
\GO_{2m}^-(2) & \SU_m(2){:}2 & \Omega_{2m-2}^-(2){:}2 & \SU_{m-2}(2) & \textup{as in~\ref{LemOmegaMinus04}}\\
\GaO_{2m}^-(4) & \SU_m(4){:}4 & \Omega_{2m-2}^-(4){:}4 & \SU_{m-2}(4) & \textup{as in~\ref{LemOmegaMinus04}}\\
\hline
\end{array}
\]
\end{proposition}

\begin{proof}
By Lemmas~\ref{LemOmegaMinus03} and~\ref{LemOmegaMinus04}, it suffices to prove the ``only if'' part. Since $K\leqslant \Nor_G(K^{(\infty)})=\N_2^+[G]$, we have $G=H\N_2^+[G]$. Then by Proposition~\ref{Prop:MaxOmegaminu} we obtain $q\leqslant4$ (note that $\N_2^+[G]$ is not maximal in $G$ for $q<4$, see~\cite[Table~3.5.F]{KL1990}). Without loss of generality, we may assume $H^{(\infty)}=\SU(V_\sharp)$ and $K^{(\infty)}=\Omega(V)_{e_1,f_1}$. Then
\[
H\cap K\geqslant H^{(\infty)}\cap K^{(\infty)}=\SU(V_\sharp)_{e_1,f_1} =\SU(V_\sharp)_{E_1,F_1}=\SU_{m-2}(q).
\]
Let $\calO=G/L$, $\calO_1=H/H^{(\infty)}$, and $\calO_2=K/K^{(\infty)}$.
Comparing the $p$-part, we obtain $|\calO_2|_p\geqslant q$.

First, assume $q=2$. Then $|\calO_2|_2\geqslant2$, and so $K\geqslant K_0$ for some $K_0=\Omega_{2m-2}^-(2).2$. Applying this conclusion to the factorization $KL\cap HL=(K\cap HL)(H\cap KL)$ we also obtain $K\cap HL\geqslant\Omega_{2m-2}^-(2).2$. If $K_0\leqslant L$, then $(H,K)$ tightly contains the pair $(H_0,K_0)$ with $H_0=H^{(\infty)}$, as in the first row of the table. If $K_0\nleqslant L$ and $H\leqslant L$, then $K\cap L=K\cap HL\geqslant\Omega_{2m-2}^-(2).2$ and so we are led to the previous case where $K_0\leqslant L$. If $K_0\nleqslant L$ and $H\nleqslant L$, then $(H,K)$ tightly contains the pair $(H_0,K_0)$ with $H_0=\SU_m(2).2$, as in the second row of the table.

Next, assume $q=3$. Then $|K|_3\leqslant|\N_2^+[G]|_3=|\Omega_{2m-2}^-(3)|_3=|K^{(\infty)}|_3$, which implies that $|\calO_2|_3=1$, contradicting the conclusion $|\calO_2|_p\geqslant q$.

Finally assume $q=4$. By~Proposition~\ref{Prop:MaxOmegaminu}, we have $G=\Aut(L)=\Omega_{2m}^-(4).4$. This conclusion applied to the factorization $KL\cap HL=(K\cap HL)(H\cap KL)$ leads to $KL=HL=\Omega_{2m}^-(4).4$. It follows that $K=(K\cap L).4$ and $H=(H\cap L).4$. Hence $(G,H,K)$ tightly contains the triple $(G_0,H_0,K_0)$ in the third row of the table.
\end{proof}

\subsection{Proof of Theorem~\ref{ThmOmegaMinus}}\label{SecOmegaMinusProof}
\ \vspace{1mm}

As factorizations $G_0=H_0K_0$ in Table~\ref{TabOmegaMinus} have been constructed in this section, we are left to prove the ``only if'' part. Let $G=HK$. Then Proposition~\ref{Prop:MaxOmegaminu} shows that $K\leqslant\Pa_1[G]$, $\N_1[G]$ or $\N_2^+[G]$. If $K^{(\infty)}=q^{2m-2}{:}\Omega_{2m-2}^-(q)$ or $K^{(\infty)}=\Omega_{2m-1}(q)$, then by Propositions~\ref{prop:OmegaMinus-1} and~\ref{prop:OmegaMinus-02}, the triple $(G,H,K)$ tightly contains some $(G_0,H_0,K_0)$ in rows~1--8 of Table~\ref{TabOmegaMinus}. Now assume that $K^{(\infty)}\neq q^{2m-2}{:}\Omega_{2m-2}^-(q)$ or $\Omega_{2m-1}(q)$. Combining Proposition~\ref{Prop:MaxOmegaminu} with Lemmas~\ref{LemOmegaMinusRow1,5}, \ref{LemOmegaMinusRow2,3} and \ref{LemOmegaMinusRow4}, we obtain $H^{(\infty)}=\SU_m(q)$ with $m$ odd and $K^{(\infty)}=\Omega_{2m-2}^-(q)$. Then by Proposition~\ref{LemOmegaMinus09}, the triple $(G,H,K)$ tightly contains some $(G_0,H_0,K_0)$ in rows~9--11 of Table~\ref{TabOmegaMinus}, completing the proof.

%%%%%%%%%%%%%%%%%%%%%%%%%%%%%%%%%%%%%%%%%%%%%%%%%%%%%%%%%%%%%%%%%%%
%%%%%%%%%%%%%%%%%%%%%%%%%%%%%%%%%%%%%%%%%%%%%%%%%%%%%%%%%%%%%%%%%%%

\section{Orthogonal groups of plus type}\label{SecOmegaPlus05}

%%%%%%%%%%%%%%%%%%%%%%%%%%%%%%%%%%%%%%%%%%%%%%%%%%%%%%%%%%%%%%%%%%%
%%%%%%%%%%%%%%%%%%%%%%%%%%%%%%%%%%%%%%%%%%%%%%%%%%%%%%%%%%%%%%%%%%%

In this section we deal with factorizations of orthogonal groups of plus type. The main result is the theorem below.

\begin{theorem}\label{ThmOmegaPlus}
Let $G$ be a classical group with $G^{(\infty)}=\Omega_{2m}^+(q)$, where $m\geqslant4$, and let $H$ and $K$ be subgroups of $G$ not containing $G^{(\infty)}$ such that both $H$ and $K$ have a unique nonsolvable composition factor. Then $G=HK$ if and only if either $(G,H,K)$ tightly contains some $(G_0,H_0,K_0)$ in Table~$\ref{TabOmegaPlus}$ or $(\overline{G},\overline{H},\overline{K})$ tightly contains some $(G_0,H_0,K_0)$ in Table~$\ref{TabOmegaPlus2}$, where $\,\overline{\phantom{\varphi}}\,$ is the quotient modulo scalars.
\end{theorem}

For convenience, we collect in Subsection~\ref{SecOmegaPlus01} the notation that will be used throughout this section. The proof of Theorem~\ref{ThmOmegaPlus} will start in Subsection~\ref{SecOmegaPlus02} and conclude in Subsection~\ref{SecOmegaPlus03}.

\subsection{Notation}\label{SecOmegaPlus01}
\ \vspace{1mm}

Throughout this section, let $q=p^f$ be a power of a prime $p$, let $m\geqslant4$ be an integer,
% let $\,\overline{\phantom{\varphi}}\,$ be the homomorphism from $\GaO_{2m}^+(q)$ to $\mathrm{P\Gamma O}_{2m}^+(q)$ modulo scalars,
let $\tau$ be a triality automorphism of $\Omega_8^+(q)$, let $V$ be a vector space of dimension $2m$ over $\bbF_q$ equipped with a nondegenerate quadratic form $Q$ of plus type, whose associated bilinear form is $\beta$, let $\perp$ denote the perpendicularity with respect to $\beta$, let $e_1,f_1,\dots,e_m,f_m$ be a standard basis for $V$ as in~\cite[2.2.3]{LPS1990}, let $\phi\in\GaO(V)$ such that
\[
\phi\colon a_1e_1+b_1f_1+\dots+a_me_m+b_mf_m\mapsto a_1^pe_1+b_1^pf_1+\dots+a_m^pe_m+b_m^pf_m
\]
for $a_1,b_1\dots,a_m,b_m\in\bbF_q$, let $\gamma$ be the involution in $\GO(V)$ swapping $e_i$ and $f_i$ for all $i\in\{1,\dots,m\}$, let $U=\langle e_1,\dots,e_m\rangle_{\bbF_q}$, and let $W=\langle f_1,\dots,f_m\rangle$. From~\cite[3.7.4~and~3.8.2]{Wilson2009} we know that $\Omega(V)_U=\Pa_m[\Omega(V)]$ has a subgroup $R{:}T$, where
\[
R=q^{m(m-1)/2}
\]
is the kernel of $\Omega(V)_U$ acting on $U$, and
\[
T=\SL_m(q)
\]
stabilizes both $U$ and $W$. The action of $T$ on $U$ determines that on $W$ in the way described in~\cite[Lemma~2.2.17]{BG2016}, from which we can see that $\gamma$ normalizes $T$ and induces the transpose inverse automorphism of $T$.
% Also note that $\gamma\in\Omega(V)$ if $m$ is even.
Let $w$ be nonsingular and define the reflection with respect to $w$ as the map $r_w$ given by
\[
r_w\colon V\to V,\quad v\mapsto v-\frac{\beta(v,w)}{Q(w)}w.
\]
Take $\mu\in\bbF_q$ such that the polynomial $x^2+x+\mu$ is irreducible over $\bbF_q$, and let
\begin{align*}
&u=e_1+e_2+\mu f_2,\\
&u'=(1-\mu^2)e_1+(\mu^2-1+\mu^{-1})e_2+\mu^2f_1+\mu^2f_2.
\end{align*}
Then $Q(e_1+f_1)=\beta(e_1+f_1,u)=\beta(e_1+f_1,u')=1$ and $Q(u)=Q(u')=\mu$. Hence $\langle e_1+f_1,u\rangle$ and $\langle e_1+f_1,u'\rangle$ are both nondegenerate $2$-subspaces of minus type in $V$, and there exists $g\in\Omega(V)$ (respectively, $g\in\mathrm{O}(V)\Setminus\Omega(V)$) such that
\[
(e_1+f_1)^g=e_1+f_1\ \text{ and }\ u^g=u'.
\]

If $m=2\ell$ is even, then let $V_\sharp$ be a vector space of dimension $m$ over $\bbF_{q^2}$ with the same underlying set as $V$, let $\Tr$ be the trace of the field extension $\bbF_{q^2}/\bbF_q$, and let $\lambda\in\bbF_{q^2}$ such that
\[
\Tr(\lambda)=\lambda+\lambda^q=1.
\]
In this case, we fix the notation for some field-extension subgroups of $\mathrm{O}(V)$ as follows.

Firstly, we may equip $V_\sharp$ with a nondegenerate unitary form $\beta_\sharp$ such that $Q(v)=\beta_\sharp(v,v)$ for all $v\in V$, and thus $\GU(V_\sharp)<\mathrm{O}(V)$. Take a standard $\bbF_{q^2}$-basis $E_1,F_1,\dots,E_\ell,F_\ell$ for $V_\sharp$, so that
\[
\beta_\sharp(E_i,E_j)=\beta_\sharp(F_i,F_j)=0,\quad\beta_\sharp(E_i,F_j)=\delta_{i,j}
\]
for all $i,j\in\{1,\dots,\ell\}$. Then we have
\begin{align}
&Q(\lambda E_1)=\beta_\sharp(\lambda E_1,\lambda E_1)=0,\label{EqnOmegaPlus01}\\
&Q(F_1)=\beta_\sharp(F_1,F_1)=0,\nonumber\\
&Q(\lambda E_1+F_1)=\beta_\sharp(\lambda E_1+F_1,\lambda E_1+F_1)=\lambda+\lambda^q=1,\nonumber
\end{align}
and hence
\[
\beta(\lambda E_1,F_1)=Q(\lambda E_1+F_1)-Q(\lambda E_1)-Q(F_1)=1.
\]
From these we see that $(\lambda E_1,F_1)$ is a hyperbolic pair with respect to $Q$.
Thus we may assume without loss of generality that
\[
e_1=\lambda E_1\ \text{ and }\ f_1=F_1.
\]
Then $Q(e_1+f_1)=\beta_\sharp(e_1+f_1,e_1+f_1)=1$, which shows that $e_1+f_1$ is a nonsingular vector in both $V$ (with respect to $Q$) and $V_\sharp$ (with respect to $\beta_\sharp$). Moreover,~\eqref{EqnOmegaPlus01} shows that $e_1$ is a singular vector in both $V$ and $V_\sharp$.
Let $\xi\in\GaU(V_\sharp)$ such that
\begin{equation}\label{EqnOmegaPlus02}
\xi\colon a_1E_1+b_1F_1+\dots+a_\ell E_\ell+b_\ell F_\ell\mapsto a_1^pE_1+b_1^pF_1+\dots+a_\ell^pE_\ell+b_\ell^pF_\ell
\end{equation}
for $a_1,b_1\dots,a_\ell,b_\ell\in\bbF_{q^2}$.

Secondly, we may equip $V_\sharp$ with a nondegenerate quadratic form $Q_\sharp$ of plus type such that $Q(v)=\Tr(Q_\sharp(v))$ for all $v\in V$, and thus $\mathrm{O}(V_\sharp)<\mathrm{O}(V)$. For $w\in V_\sharp$ with $Q_\sharp(w)\neq0$, let $r'_w$ be the reflection with respect to $w$ defined by
\begin{equation}\label{EqnOmegaPlus09}
r'_w\colon V_\sharp\to V_\sharp,\quad v\mapsto v-\frac{Q_\sharp(v+w)-Q_\sharp(v)-Q_\sharp(w)}{Q_\sharp(w)}w.
\end{equation}
Take a standard $\bbF_{q^2}$-basis $E_1,F_1,\dots,E_\ell,F_\ell$ for $V_\sharp$. Then we have
\begin{align*}
&Q(\lambda E_1)=\Tr(Q_\sharp(\lambda E_1))=\Tr(0)=0,\\
&Q(F_1)=\Tr(Q_\sharp(F_1))=\Tr(0)=0,\\
&Q(\lambda E_1+F_1)=\Tr(Q_\sharp(\lambda E_1+F_1))=\Tr(\lambda)=\lambda+\lambda^q=1,
\end{align*}
and hence
\[
\beta(\lambda E_1,F_1)=Q(\lambda E_1+F_1)-Q(\lambda E_1)-Q(F_1)=1.
\]
From these we see that $(\lambda E_1,F_1)$ is a hyperbolic pair with respect to $Q$.
Thus we may assume without loss of generality that
\[
e_1=\lambda E_1\ \text{ and }\ f_1=F_1.
\]
It follows that $Q(e_1+f_1)=1$ and $Q_\sharp(e_1+f_1)=\lambda$, whence $e_1+f_1$ is a nonsingular vector in both $V$ (with respect to $Q$) and $V_\sharp$ (with respect to $Q_\sharp$). Let $\psi\in\GaO(V_\sharp)$ such that
\begin{equation}\label{EqnOmegaPlus10}
\psi\colon a_1E_1+b_1F_1+\dots+a_\ell E_\ell+b_\ell F_\ell\mapsto a_1^pE_1+b_1^pF_1+\dots+a_\ell^pE_\ell+b_\ell^pF_\ell
\end{equation}
for $a_1,b_1\dots,a_\ell,b_\ell\in\bbF_{q^2}$.

\subsection{Factorizations with socle $\Sp_6(q)$ and a factor normalizing $\mathrm{G}_2(q)$}
\ \vspace{1mm}

When dealing with factorizations of almost simple groups with socle $\POm_8^+(q)$ in Subsection~\ref{SecOmegaPlus02}, where $q$ is even, we will need to apply results on certain factorizations of almost simple groups with socle $\Sp_6(q)$. These results will be needed in Section~\ref{SecSymplectic02} again. So we single out this part in the current subsection.

% We follow the convention as in~\cite{LPS1990} to denote $\SL_n^+(q)=\SL_n(q)$ and $\SL_n^-(q)=\SU_n(q)$.

\begin{lemma}\label{LemSymplectic07}
Let $G=\Sp_6(q)$ with $q$ even, let $H=\Omega_6^\varepsilon(q)<G$ with $\varepsilon\in\{+,-\}$, and let $K=\G_2(q)<G$. Then $G=HK$ with $H\cap K=\SL_3^\varepsilon(q)$.
\end{lemma}

\begin{proof}
By~\cite[5.2.3(b)]{LPS1990}, $\mathrm{O}_6^\varepsilon(q)\cap K=\SL_3^\varepsilon(q).2$. Since $H$ has no subgroup of the form $\SL_3^\varepsilon(q).2$, it follows that $H\cap K=\SL_3^\varepsilon(q)$, and then the order identity gives $G=HK$.
\end{proof}

\begin{remark}
If we let $K=\G_2(q)'$ in Lemma~\ref{LemSymplectic07}, then the conclusion $G=HK$ would not hold for $(q,\varepsilon)=(2,+)$.
\end{remark}

\begin{lemma}\label{LemSymplectic008}
Let $G=\Sp_6(q)$ with $q\geqslant4$ even, let $H=\Sp_4(q)<\N_2[G]$, and let $K=\G_2(q)<G$. Then $G=HK$ with $H\cap K=\SL_2(q)$.
\end{lemma}

\begin{proof}
This follows from~\cite[4.3.6]{Wilson2009} and the order identity.
\end{proof}

The next lemma is verified by computation in \magma~\cite{BCP1997}.

\begin{lemma}\label{LemSymplectic10}
Let $G=\GaSp_6(4)$ and $K=\GaG_2(4)<G$. Then $G$ has the following two conjugacy classes of subgroups $H=\SiL_2(16)$ such that $G=HK$ with $H\cap K=1$.
\begin{enumerate}[{\rm(a)}]
\item $H<(\SL_2(16)\times5){:}4<\GaO_6^+(4)<G$;
\item $H<\SiL_2(16)\times3<\GaO_6^-(4)<G$.
\end{enumerate}
\end{lemma}

We now characterize factorizations of $6$-dimensional symplectic groups with a factor contained in $\GaG_2(q)$.

\begin{proposition}\label{LemSymplectic12}
Let $\Sp_6(q)\leqslant G\leqslant\GaSp_6(q)$ with $q\geqslant4$ even, let $H$ be a core-free subgroup of $G$ with a unique nonsolvable composition factor, and let $K$ be a subgroup of $G$ such that $K\leqslant\GaG_2(q)$. Then $G=HK$ if and only if $(G,H,K)$ tightly contains some $(G_0,H_0,K_0)$ in the following table. In this case, $H_0\cap K_0$ is described in the table.
\[
\begin{array}{lllll}
\hline
G_0 & H_0 & K_0 & H_0\cap K_0 & \textup{Remarks} \\
\hline
\Sp_6(q) & \Omega^+_6(q),\ \Omega^-_6(q) & \G_2(q) & \SL_3(q),\ \SU_3(q) &  \\
\Sp_6(q) & \Sp_4(q) & \G_2(q) & \SL_2(q) & H_0\textup{ not necessarily in }\N_2[G] \\
\Sp_6(q) & q^5{:}\Sp_4(q) & \G_2(q) & q^{2+3}{:}\SL_2(q) & H_0=\Pa_1[G]^{(\infty)} \\
\Sp_6(q) & q^4{:}\Omega^-_4(q) & \G_2(q) & [q^3] & H_0=\Pa_1[\Omega_6^-(q)]^{(\infty)} \\
\Sp_6(4) & \Omega_6^-(4) & \J_2 & 5^2{:}\Sy_3 &  \\
\GaSp_6(4) & \GaO_6^-(4) & \SU_3(3)\times2 & \Sy_3 &  \\
\GaSp_6(4) & \GaO_6^-(4) & \G_2(2) & \Sy_3 & K_0\nleqslant\Sp_6(4) \\
\GaSp_6(4) & \SiL_2(16) & \GaG_2(4) & 1 & \textup{as in~\ref{LemSymplectic10}} \\
\hline
\end{array}
\]
\end{proposition}

\begin{proof}
The case $q=4$ is verified in \magma~\cite{BCP1997}. For $q\geqslant8$, by~\cite{HLS1987} and~\cite[Theorem~A]{LPS1990}, we only need to consider the case where $K^{(\infty)}=\G_2(q)$ and $H$ is contained in a maximal subgroup $M$ of $G$ such that either $M^{(\infty)}\in\{\Omega^+_6(q),\Omega_6^-(q)\}$ or $M\in\{\Pa_1[G],\N_2[G]\}$.

If $M^{(\infty)}=\Omega^\varepsilon_6(q)$ with $\varepsilon\in\{+,-\}$, then Lemma~\ref{LemSymplectic07} together with Theorems~\ref{ThmLinear} and~\ref{ThmUnitary} verifies the proposition.

If $M=\N_2[G]$, then $M$ acts on a nondegenerate $4$-dimensional subspace with kernel, say, $N$. By~\cite[4.3.6]{Wilson2009}, $(M^{(\infty)}\cap K)N/N=\SL_2(q)$. Then considering the quotient of $M$ modulo $N$, we obtain the conclusion of the proposition by Lemma~\ref{LemSymplectic008}.

Similarly, for $M=\Pa_1[G]$, we see from~\cite[4.3.7]{Wilson2009} that $\bfO_2(M)\cap K^{(\infty)}=q^2$. Then considering the quotient of $M$ modulo $\bfO_2(M)$, the proposition follows from Example~\ref{ex:Linear} and Lemma~\ref{LemSymplectic07}.
\end{proof}

\subsection{Reduction}\label{SecOmegaPlus02}
\ \vspace{1mm}

Let $L=G^{(\infty)}=\Omega(V)=\Omega_{2m}^+(q)$ with $m\geqslant4$, and let $H$ and $K$ be subgroups of $G$ such that both $H$ and $K$ have a unique nonsolvable composition factor and neither $H$ nor $K$ contains $L$. Take $\max^-$ subgroups $A$ and $B$ of $G$ containing $H$ and $K$ respectively. It is known from~\cite[Theorem~A]{LPS1990} and~\cite{LPS1996} that, in most cases, one of the two factors in a $\max^-$ factorization of $G$ stabilizes a subspace of dimension $1$ or $2$. More precisely, we have the proposition below.

\begin{proposition}\label{Prop:MaxO+subgps}
Let $(m,q)\not=(4,2)$, $(4,3)$ or $(4,4)$. If $G=HK$, then with the above notation and interchanging $A$ and $B$ if necessary, one of the following holds:
\begin{enumerate}[{\rm(a)}]
\item $B=\Pa_1[G]$, and $A^{(\infty)}=\SU_m(q)$ with $m$ even;
\item $B=\N_2^+[G]$, and $A^{(\infty)}=\SU_m(q)$ with $m$ even and $q=4$;
\item $B=\N_2^-[G]$, and either $A=\Pa_{m-1}[G]$ or $\Pa_m[G]$, or $A^{(\infty)}=\SL_m(q)$ with $q\in\{2,4\}$;
\item $B=\N_1[G]$, and $A^{(\infty)}$ is one of:
\begin{align*}
&q^{m(m-1)/2}{:}\SL_m(q),\,\ \SL_m(q),\,\ \SU_m(q)\text{ with }m\text{ even},\,\ \Omega_m^+(q^2)\text{ with }q\in\{2,4\},\\
&(\Sp_2(q)\circ\Sp_m(q))^{(\infty)}\text{ with }q>2,\,\ \mathrm{Spin}_9(q)\text{ with }m=8,\,\ \Co_1\text{ with }(m,q)=(12,2);
\end{align*}
\item $m=4$, $\overline{B}^{(\infty)}=\Omega_7(q)$, and $\overline{A}^{(\infty)}=\Omega_7(q)$, $\Omega^{\pm}_6(q)$ or $\Omega_8^-(q^{1/2})$, where $\,\overline{\phantom{\varphi}}\,$ is the quotient modulo scalars.
\item $m=4$, $q=16$, $B\cap L=(17\times\Omega_6^-(16)).2$, and $A\cap L=\Omega_8^-(4)$.
\end{enumerate}
\end{proposition}

\begin{remark}
Part~(f) of Proposition~\ref{Prop:MaxO+subgps} is missing in~\cite[Table~4]{LPS1990}, see~\cite[Remark~1.5]{GGP}.
\end{remark}

In view of the triality automorphism of $\POm_8^+(q)$, the case $(\overline{A}^{(\infty)},\overline{B}^{(\infty)})=(\Omega^{\pm}_6(q),\Omega_7(q))$ in part~(e) of Proposition~\ref{Prop:MaxO+subgps} can be reduced to part~(d) with $A^{(\infty)}=\SL_4^\pm(q)$, and in part~(f) we may assume $B=\N_2^-[G]$.

\begin{lemma}\label{PropOmegaPlusO+8}
Let $m=4$, $q\notin\{2,3,4\}$, $\overline{K}^{(\infty)}=\Omega_7(q)$, and $\overline{A}^{(\infty)}=\Omega_7(q)$ or $\Omega_8^-(q^{1/2})$, where $\,\overline{\phantom{\varphi}}\,$ is the quotient modulo scalars. Then $G=HK$ if and only if $(\overline{G},\overline{H},\overline{K})$ tightly contains the triple $(\POm_8^+(q),H_0,\Omega_7(q))$ for some $H_0$ in the following table. In this case, $H_0\cap\overline{K}^{(\infty)}$ is described in the table.
\[
\begin{array}{c|cccccc}
\hline
H_0 & \Omega_7(q) & \Omega_6^\pm(q) & q^5{:}\Omega_5(q) & \Omega_5(q) & q^4{:}\Omega^-_4(q) & \Omega_8^-(q^{1/2}) \\
H_0\cap\overline{K}^{(\infty)} & \G_2(q) & \SL_3^\pm(q) & [q^5]{:}\SL_2(q) & \SL_2(q) & [q^3] & \G_2(q^{1/2}) \\
\hline
\end{array}
\]
\end{lemma}

\begin{proof}
By triality we may take $\overline K^{(\infty)}=\Omega_7(q)$ in class $\calC_9$. Its intersections with $\Omega_7(q)$ and $\Omega_8^-(q^{1/2})$ are respectively $\G_2(q)$ and $\G_2(q^{1/2})$ by~\cite[Proposition~3.1.1(iv)]{Kleidman1987} and~\cite[Page~105]{LPS1990}. This gives $\POm_8^+(q)=\overline{A}^{(\infty)}\overline{K}^{(\infty)}$. The result now follows inside $\overline A$ from Theorem~\ref{ThmOmega}, Proposition~\ref{LemSymplectic12} and Theorem~\ref{ThmOmegaMinus}.
\end{proof}

We treat the special cases with $L\in\{\Omega_8^+(2),\POm_8^+(3),\Omega_8^+(4)\}$ by computation in \magma~\cite{BCP1997}, which leads to the following result on the factorization $G=HK$.

\begin{lemma}\label{LemXia29}
Let $L=\Omega_8^+(2)$, $\Omega_8^+(3)$ or $\Omega_8^+(4)$. Then $G=HK$ if and only if either $(G,H,K)$ tightly contains some $(G_0,H_0,K_0)$ in rows~\emph{1--3},~\emph{5--7} or~\emph{15--23} of Table~$\ref{TabOmegaPlus}$ or $(\overline{G},\overline{H},\overline{K})$ tightly contains some $(G_0,H_0,K_0)$ in Table~$\ref{TabOmegaPlus2}$, where they satisfy the following extra conditions and $\,\overline{\phantom{\varphi}}\,$ is the quotient modulo scalars.
\begin{enumerate}[{\rm(a)}]
\item In rows~\emph{2--5} of Table~$\ref{TabOmegaPlus2}$, for each $K_0$ and each isomorphism type in the corresponding $H_0$ columns, $H_0$ lies in one of the two conjugacy classes of subgroups of $G_0$ such that $G_0=H_0K_0$ and that the intersection $H_0\cap K_0$ is as described.
\item In rows~\emph{6--8} of Table~$\ref{TabOmegaPlus2}$, for each $K_0$ and each isomorphism type in the corresponding $H_0$ columns, $H_0$ lies in one of the $k$ conjugacy classes of subgroups of $G_0$ such that $G_0=H_0K_0$, where $k$ is given in the following three tables for $K_0=\Omega_7(3)$, $3^6{:}\PSL_4(3)$ and $\Omega_8^+(2)$ respectively.
\[
\begin{array}{c|ccccccc}
\hline
H_0 & 3^4{:}\Sy_5 & 3^4{:}(\A_5\times4) & (3^5{:}2^4){:}\A_5 & \A_9 & \SU_4(2) & \Sp_6(2) & \Omega_8^+(2) \\
k & 4 & 4 & 4 & 4 & 8 & 4 & 2 \\
\hline
\end{array}
\]
\[
\begin{array}{c|ccccccc}
\hline
H_0 & 2^6{:}\A_7 & \A_8 & 2^6{:}\A_8 & \A_9 & 2^{\boldsymbol{\cdot}}\PSL_3(4) & \Sp_6(2) & \Omega_6^-(3)\\
k & 4 & 8 & 4 & 8 & 2 & 8 & 2\\
\hline
\end{array}
\]
\[
\begin{array}{c|cccc}
\hline
H_0  & 3^6{:}\SL_3(3) & 3^{3+3}{:}\SL_3(3) & 3^{6+3}{:}\SL_3(3) & 3^6{:}\PSL_4(3) \\
k  & 3 & 3 & 3 & 3 \\
\hline
\end{array}
\]
\item In row~\emph{9} of Table~$\ref{TabOmegaPlus2}$, $G_0=\Omega_8^+(4){:}\langle\phi\rangle$, $K_0=\N_1[G_0]$, and
%\[
%H_0<(\SL_2(16)\times5){:}4<\GaO_6^+(4)<(K_0)^\tau\ \text{ or }\ H_0<\GaL_2(16)\times3<\GaO_6^-(4)<(K_0)^\tau
%\]
either $H_0<(\SL_2(16)\times5){:}4<\GaO_6^+(4)<(K_0)^\tau$ or $H_0<\GaL_2(16)\times3<\GaO_6^-(4)<(K_0)^\tau$
for some triality automorphism $\tau$ of $\Omega_8^+(4)$.
\end{enumerate}
\end{lemma}

Our analysis will proceed by different cases for $K$ categorized in Proposition~\ref{Prop:MaxO+subgps}. Notice
\[
\Pa_1[G]^{(\infty)}=q^{2m-2}{:}\Omega_{2m-2}^+(q),\,\ \N_1[G]^{(\infty)}=\Omega_{2m-1}(q),\,\ \N_2^{\pm}[G]^{(\infty)}=\Omega_{2m-2}^{\pm}(q).
\]
Also note that $\Omega_{2m-2}^+(q)$ has trivial first cohomology on the natural $(2m-2)$-dimensional module (see~\cite[Table~4.5]{CPS1975}) and that the Levi factor of $\Pa_1[L]$ is a subgroup of $\N_2^+[L]$. In the ensuing subsections we first handle the cases $B=\Pa_1[G]$ and $B=\N_2^+[G]$ together, then the case $B=\N_2^-[G]$, and finally $B=\N_1[G]$.

We close this subsection with two technical lemmas that will be applied multiple times in this section. Recall the notation $r'$ defined in~\eqref{EqnOmegaPlus09} and $\psi$ in~\eqref{EqnOmegaPlus10}.

\begin{lemma}\label{LemOmegaPlus12}
Let $G=L{:}\langle\phi\rangle=\Omega_{2m}^+(q){:}f$ with $q\in\{2,4\}$ and $m$ even, let $H=\Omega(V_\sharp){:}\langle\rho\rangle=\Omega_m^+(q^2){:}(2f)$ with $\rho\in\{\psi,\psi r'_{E_1+F_1}\}$, and let $K=G_{e_1+f_1}=\GaSp_{2m-2}(q)$. Then $G=HK$ with $H\cap K=\Sp_{m-2}(q^2)$.
\end{lemma}

\begin{proof}
Let $S=\Omega(V_\sharp)$. Since $e_1+f_1$ is nonsingular in both $V$ and $V_\sharp$, we have
\[
S\cap K=\Omega(V_\sharp)\cap G_{e_1+f_1}=\Omega(V_\sharp)_{e_1+f_1}=\Omega_{m-1}(q^2)\cong\Sp_{m-2}(q^2).
\]
Now as $q\in\{2,4\}$ we have $q=2f$, and so $|H|/|S\cap K|=|G|/|K|$. Thus it suffices to prove $H\cap K=S\cap K$.
Suppose for a contradiction that there exists $y\in(H\Setminus S)\cap K$. Then there exists $s\in S$ such that $\rho^fs\in K$.
Since $K=G_{e_1+f_1}=G_{\lambda E_1+F_1}$, it follows that
\[
\lambda E_1+F_1=(\lambda E_1+F_1)^{\rho^fs}=(\lambda E_1+F_1)^{\psi^f\psi^f\rho^fs}=(\lambda^qE_1+F_1)^{\psi^f\rho^fs}.
\]
Then as $\psi^f\rho^fs\in\mathrm{O}(V_\sharp)$, we obtain $Q_\sharp(\lambda E_1+F_1)=Q_\sharp(\lambda^qE_1+F_1)$, contradicting $\lambda^q+\lambda=1$.
\end{proof}

\begin{lemma}\label{OmegaPlusClaim}
Let $G=HK$ with $K\leqslant B=\N_1[G]$. Then the following statements hold.
\begin{enumerate}[{\rm (a)}]
\item If $q$ is even, then $H\nleqslant\mathrm{O}(V_\sharp){:}\langle\psi^2\rangle$;
\item If $H$ is in a field-extension subgroup $M$ of $G$ defined over $\bbF_{q^b}$, then $b=2$ and $q\in\{2,4\}$.
\end{enumerate}
\end{lemma}

\begin{proof}
For~(a), if $H\leqslant M:=\mathrm{O}_m^+(q^2).f$, then $G=MB$ and
\[
q^{m-1}\leqslant |M|_2/|M\cap B|_2\leqslant
|\mathrm{O}_m^+(q^2).f|_2/|\Sp_{m-2}(q^2)\times2|_2=q^{m-2}f,
\]
a contradiction. To prove part~(b), by Proposition~\ref{Prop:MaxO+subgps}, it suffices to exclude the possibility for $A^{(\infty)}=\Omega_m^+(4)$ with $q=2$ such that $H\cap A^{(\infty)}$ is contained in $\Omega_{m/2}^+(16).2^2$. Suppose for a contradiction that this is the case. From Lemma~\ref{LemOmegaPlus12} we see that $(\Omega(V_\sharp){:}\langle\psi\rangle)\cap B$ and $(\Omega(V_\sharp){:}\langle\psi r'_{E_1+F_1}\rangle)\cap B$ are both contained in $\Omega(V_\sharp)$. Then since $A\leqslant\Omega(V_\sharp){:}\langle\psi,r'_{E_1+F_1}\rangle$,
it follows that $A\cap B\leqslant \Omega(V_\sharp){:}\langle r'_{E_1+F_1}\rangle=\mathrm{O}(V_\sharp)$. Hence the factorization $A=H(A\cap B)$ yields $A\cap\mathrm{O}(V_\sharp)=(H\cap\mathrm{O}(V_\sharp))(A\cap B)$. However, the conclusion of part~(a) applied to the factorization $A\cap\mathrm{O}(V_\sharp)=(H\cap\mathrm{O}(V_\sharp))(A\cap B)$ shows that it is impossible.
\end{proof}

\subsection{Actions on $\calP_1$ and $\calN_2^+$}
\ \vspace{1mm}

The stabilizer of the singular vector $e_1$ in $L$ has the form
\[
L_{e_1}=q^{2m-2}{:}\Omega_{2m-2}^+(q)\leqslant\Pa_1[L].
\]
By Witt's Lemma, $L=\Omega_{2m}^+(q)$ is transitive on the set of nonzero singular vectors, whence the number of singular vectors is
\[
|\calP_1|=|L|/|L_{e_1}|=(q^m-1)(q^{m-1}+1).
\]
For the hyperbolic pair $(e_1,f_1)$, we have the inclusion of stabilizers
\[
\Omega_{2m-2}^+(q)=L_{e_1,f_1}\leqslant L_{\{e_1,f_1\}}\leqslant L_{\langle e_1,f_1\rangle}=\N_2^+[L].
\]
%and moreover,
%\[\begin{array}{l}
%L_{e_1,f_1}=\Omega(\l e_2,f_2,\dots,e_m,f_m\r)=\Omega_{2m-2}^+(q),\\
%L_{\{e_1,f_1\}}=L_{e_1,f_1}{:}\l s\r=\Omega_{2m-2}^+(q){:}2,\\
%\end{array}
%\]
%where $s$ interchanges $e_i$ and $f_i$ for $1\leqslant i\leqslant m-1$ if $m$ is odd, and $1\leqslant i\leqslant m$ if $m$ is even.

\subsubsection{Candidates for $H^{(\infty)}$ and $K^{(\infty)}$}

Below is a lemma determining the pair $(H^{(\infty)},K^{(\infty)})$ for the case where one of the factors stabilizes a singular 1-space, or a nonsingular 2-subspace of plus type.

\begin{lemma}\label{P_1N_2^+<O+}
Let $G=HK$ with $B=\Pa_1[G]$ or $\N_2^+[G]$. Then $m$ is even, $H^{(\infty)}=\SU_m(q)$, and either $K^{(\infty)}=q^{2m-2}{:}\Omega_{2m-2}^+(q)$, or $K^{(\infty)}=\Omega_{2m-2}^+(q)$ with $q\in\{2,4\}$.
\end{lemma}

\begin{proof}
Proposition~\ref{Prop:MaxO+subgps} first gives $A^{(\infty)}=\SU_m(q)$ with $m$ even. Primitive prime divisors of $q^{2m-2}-1$ and $q^m-1$, together with Proposition~\ref{Prop:Unitary-max} applied to $A/\Rad(A)$, force $H^{(\infty)}=A^{(\infty)}$. Similarly, considering $\ppd(q^{m-1}-1)$ and $\ppd(q^{2m-4}-1)$ and using Proposition~\ref{Prop:MaxO+subgps} and~\cite[Tables~8.8~and~8.9]{BHR2013} on factorizations of $B/\Rad(B)$, gives $K^{(\infty)}$ as stated in the lemma.
\end{proof}

\subsubsection{Examples}

We now give examples of factorizations $G=HK$ with $K\leqslant\Pa_1[G]$ or $\N_2^+[G]$. Their proofs are standard: intersect $\SU(V_\sharp)$ with the indicated stabilizer and use $|G||H\cap K|=|H||K|$. The argument as in Lemma~\ref{LemOmegaPlus12} shows that no element outside $\SU(V_\sharp)$ lies in $H\cap K$.

\begin{lemma}\label{Ex:OmegaPlus17}
Let $G=L=\Omega_{2m}^+(q)$ with $m$ even, let $H=\SU(V_\sharp)=\SU_m(q)$, and let $K=G_{e_1}=q^{2m-2}{:}\Omega_{2m-2}^+(q)$. Then $G=HK$ with $H\cap K=q^{1+(2m-4)}{:}\SU_{m-2}(q)$.
\end{lemma}

The three subsequent lemmas describe examples of factorizations satisfying $K^{(\infty)}=\Omega_{2m-2}^+(q)$ with $q\in\{2,4\}$ in Lemma~\ref{P_1N_2^+<O+}. The examples for $q=2$ are given in Lemmas~\ref{LemOmegaPlus18} and~\ref{LemOmegaPlus19}, while those for $q=4$ lie in Lemmas~\ref{LemOmegaPlus18} and~\ref{LemOmegaPlus22}.

\begin{lemma}\label{LemOmegaPlus18}
Let $G=L{:}\langle\phi\rangle=\Omega_{2m}^+(q){:}f$ with $q\in\{2,4\}$ and $m$ even, let $H=\SU(V_\sharp){:}\langle\xi\rangle=\SU_m(q){:}(2f)$, and let $K=G_{e_1,f_1}=\Omega_{2m-2}^+(q){:}f$. Then $G=HK$ with $H\cap K=H^{(\infty)}\cap K^{(\infty)}=\SU_{m-2}(q)$.
\end{lemma}

For $q=2$, the stabilizer in $L$ of the unordered pair $\{e_1,f_1\}$ is
\[
L_{\{e_1,f_1\}}=L_{e_1,f_1}{:}\langle r_{e_1+f_1}r_{e_2+f_2}\rangle=L_{e_1,f_1}{:}2=\Omega_{2m-2}^+(2){:}2.
\]

\begin{lemma}\label{LemOmegaPlus19}
Let $G=L=\Omega_{2m}^+(2)$ with $m$ even, let $H=\SU(V_\sharp)=\SU_m(2)$, and let $K=L_{\{e_1,f_1\}}=\Omega_{2m-2}^+(2){:}2$. Then $G=HK$ with $H\cap K=\SU_{m-2}(2)$.
\end{lemma}

For $q=4$, since $r_{e_1+f_1}r_{e_2+f_2}$ and $\phi$ are two commutative involutions that normalize $L_{e_1,f_1}$, the product $r_{e_1+f_1}r_{e_2+f_2}\phi$ is an involution normalizing $L_{e_1,f_1}$, and so
\[
L_{e_1,f_1}{:}\langle r_{e_1+f_1}r_{e_2+f_2}\phi\rangle=L_{e_1,f_1}{:}2=\Omega_{2m-2}^+(4){:}2.
\]

\begin{lemma}\label{LemOmegaPlus22}
Let $G=L{:}\langle\phi\rangle=\Omega_{2m}^+(4){:}2$ with $m$ even, let $H=\SU(V_\sharp){:}\langle\xi\rangle=\SU_m(4){:}4$, and let
$K=L_{e_1,f_1}{:}\langle r_{e_1+f_1}r_{e_2+f_2}\phi\rangle=\Omega_{2m-2}^+(4){:}2$. Then $G=HK$ with $H\cap K=\SU_{m-2}(4)$.
\end{lemma}

\subsubsection{Classifying factorizations with $K\leqslant\Pa_1[G]$ or $\N_2^+[G]$}

We present a classification of factorizations $G=HK$ with $K\leqslant\Pa_1[G]$ or $\N_2^+[G]$ in the following proposition.

\begin{proposition}\label{prop:O+N_2^+}
Let $K\leqslant\Pa_1[G]$ or $\N_2^+[G]$. Then $G=HK$ if and only if $(G,H,K)$ tightly contains some $(G_0,H_0,K_0)$ in the following table with $m$ even. In this case, $H_0\cap K_0$ is described in the table.
\[
\begin{array}{lllll}
\hline
G_0 & H_0 & K_0 & H_0\cap K_0 & \textup{Remarks} \\
\hline
\Omega_{2m}^+(q) & \SU_m(q) & q^{2m-2}{:}\Omega_{2m-2}^+(q) & q^{1+(2m-4)}{:}\SU_{m-2}(q) & \textup{as in \ref{Ex:OmegaPlus17}} \\
\Omega_{2m}^+(2) & \SU_m(2){:}\calO & \Omega_{2m-2}^+(2){:}(2/\calO) & \SU_{m-2}(2) & \calO\leqslant2,\textup{ as in \ref{LemOmegaPlus18} or \ref{LemOmegaPlus19}}\\
\Omega_{2m}^+(4){:}\langle\phi\rangle & \SU_m(4){:}4 & \Omega_{2m-2}^+(4){:}2 & \SU_{m-2}(4) & \textup{as in \ref{LemOmegaPlus18} or \ref{LemOmegaPlus22}}\\
\hline
\end{array}
\]
\end{proposition}

\begin{proof}
By the lemmas shown in the last column of the table, it suffices to prove the ``only if'' part. Suppose that $G=HK$. By Lemma~\ref{P_1N_2^+<O+}, $H^{(\infty)}=\SU_m(q)$ with $m$ even, and either $K^{(\infty)}=q^{2m-2}{:}\Omega_{2m-2}^+(q)$, or $K^{(\infty)}=\Omega_{2m-2}^+(q)$ with $q\in\{2,4\}$. For the former case, $(G,H,K)$ tightly contains $(G_0,H_0,K_0)$ in the first row of the table in the proposition. Now let $K^{(\infty)}=\Omega_{2m-2}^+(q)$ with $q\in\{2,4\}$.
By~\cite[Table~4.5]{CPS1975}, $\Omega_{2m-2}^+(q)$ has trivial first cohomology on the natural $(2m-2)$-dimensional module, and hence $K\leqslant\N_2^+[G]$ in this case. Since $H^{(\infty)}=\SU_m(q)$, we have $H\leqslant A<L{:}\langle\phi\rangle$ (see~\cite[Table~3.5.E]{KL1990}). Note that $\GU_m(q).(2f)$ has a unique conjugacy class of subgroups isomorphic to $\SU_m(q).(2f)$ as $q\in\{2,4\}$.

Assume first $q=2$. Then $H<L$, and so $L=H(K\cap L)$. From Lemma~\ref{LemOmegaPlus18} we derive $2|L||H^{(\infty)}\cap K^{(\infty)}|=|H^{(\infty)}||K^{(\infty)}|$. It follows that either $H\leqslant H^{(\infty)}.2$ or $K\leqslant K^{(\infty)}.2$. Thus $(G,H,K)$ tightly contains $(G_0,H_0,K_0)$ in the second row of the table.

Now assume $q=4$. By~\cite[Theorem~A]{LPS1990}, such a factorization does not exist for $G=L$. This together with $H<L{:}\langle\phi\rangle$ implies $G\geqslant L{:}\langle\phi\rangle$ and hence $L\langle\phi\rangle=(H\cap L\langle\phi\rangle)K$ with $H\cap L\langle\phi\rangle\nleqslant L$ and $K\nleqslant L$. Thus $(G,H,K)$ tightly contains $(G_0,H_0,K_0)$ in the third row of the table.
\end{proof}

\subsection{Actions on ${\mathcal N}_2^-$}
\ \vspace{1mm}

In this subsection, we handle the candidate with $K\leqslant B=\N_2^-[G]$. Recall that, for the vector
\[
u=e_1+e_2+\mu f_2
\]
introduced in Subsection~\ref{SecOmegaPlus01}, the 2-subspace $\langle e_1+f_1,u\rangle_{\bbF_q}$ is an orthogonal space of minus type. It follows that $\langle e_1+f_1,u\rangle^\perp=\langle e_1-f_1+f_2,-e_2+\mu f_2,e_3,f_3,\dots,e_m,f_m\rangle$ is an orthogonal space of minus type, and the pointwise stabilizer of $\{e_1+f_1,u\}$ in $L$ is
\[
L_{e_1+f_1,u}=\Omega(\langle e_1+f_1,u\rangle^\perp)=\Omega_{2m-2}^-(q).
\]
Thus, $\Omega_{2m-2}^-(q)=L_{e_1+f_1,u}\leqslant L_{\{e_1+f_1,u\}}\leqslant L_{\langle e_1+f_1,u\rangle}\leqslant\N_2^-[L]$.
Since $\SL_m(q)$ has trivial first cohomology on the alternating square of its natural module (see~\cite[Page~324]{JonPar1976}), there is a unique conjugacy class of subgroups isomorphic to $\SL_m(q)$ in $\Pa_m[L]^{(\infty)}=q^{m(m-1)/2}{:}\SL_m(q)$.

\subsubsection{Candidates for $H^{(\infty)}$ and $K^{(\infty)}$}

We first identify the possibilities for $H^{(\infty)}$ and $K^{(\infty)}$.

\begin{lemma}\label{LemOmegaPlus1Row4,9,10}
Let $G=HK$ with $B=\N_2^-[G]$. Then $K^{(\infty)}=\Omega^-_{2m-2}(q)$, and one of the following holds:
\begin{enumerate}[{\rm (a)}]
\item $H^{(\infty)}=A^{(\infty)}=q^{m(m-1)/2}{:}\SL_m(q)$;
\item $H^{(\infty)}=A^{(\infty)}=\SL_m(q)$ with $q\in\{2,4\}$;
\item $H^{(\infty)}=A^{(\infty)}=\Omega_8^-(4)$ with $(m,q)=(4,16)$.
\end{enumerate}
\end{lemma}

\begin{proof}
Proposition~\ref{Prop:MaxO+subgps} gives the three displayed candidates for $A^{(\infty)}$. Primitive prime divisors of $q^{2m-2}-1$ and $q^{2m-4}-1$, followed by Proposition~\ref{Prop:MaxOmegaminu} (or Proposition~\ref{Prop:Unitary-max} when $m=4$), force $K^{(\infty)}=B^{(\infty)}$. The remaining primitive divisors of $q^m-1$ and $q^{m-1}-1$ force $H^{(\infty)}=A^{(\infty)}$, except that $H^{(\infty)}=A^{(\infty)}=\SL_m(q)$ with $q\in\{2,4\}$, where the restriction $q\in\{2,4\}$ is obtained by comparing the $p$-parts of $|H|/|H\cap K|$ and $|G|/|K|$.
\end{proof}

\subsubsection{Examples}

Next, we describe examples of factorizations $G=HK$ satisfying Lemma~\ref{LemOmegaPlus1Row4,9,10}, which are verified by the standard argument to determine the intersection $H\cap K$. We record the index
\[
|L|/|L_{e_1+f_1,u}|=|\Omega_{2m}^+(q)|/|\Omega_{2m-2}^-(q)|=q^{2m-2}(q^m-1)(q^{m-1}-1).
\]
In part~(a) of Lemma~\ref{LemOmegaPlus1Row4,9,10}, the group $A$ is the stabilizer in $G$ of an totally singular subspace of dimension $m$, and with the notation in Subsection~\ref{SecOmegaPlus01}, $A^{(\infty)}=R{:}T$, where $R=q^{m(m-1)/2}$ is the kernel of $L_U$ acting on $U$, and $T=\SL_m(q)$ stabilizes $U=\langle e_1,\dots,e_m\rangle$ and $W=\langle f_1,\dots,f_m\rangle$.

\begin{lemma}\label{LemOmegaPlus13}
Let $G=L=\Omega_{2m}^+(q)$, let $H=R{:}T=q^{m(m-1)/2}{:}\SL_m(q)$, and let $K=L_{e_1+f_1,u}=\Omega_{2m-2}^-(q)$. Then $G=HK$ with $H\cap K=(q^{(m-2)(m-3)/2}.q^{2m-4}){:}\SL_{m-2}(q)$ and $T\cap K=\SL_{m-2}(q)$.
\end{lemma}

\begin{proof}
Let $U_1=\langle e_2,\dots,e_m\rangle$ and $U_2=\langle e_3,\dots,e_m\rangle$. Standard calculation shows that $T\cap K=T_{e_1,e_2,f_1,f_2}=\SL_{m-2}(q)$. As the kernel of $H\cap K$ acting on $U$, the group $R\cap K$ is isomorphic to the pointwise stabilizer of $U_2$ in $\Omega(\langle e_3,f_3,\dots,e_m,f_m\rangle)$, and so $R\cap K=q^{(m-2)(m-3)/2}$. Moreover,
\[
(H\cap K)^U\leqslant\left\{
\begin{pmatrix}
1&&A_{13}\\
&1&A_{23}\\
&&A_{33}
\end{pmatrix}
\,\middle|\,A_{13},A_{23}\in\bbF_q^{m-2},\A_{33}\in\SL_{m-2}(q)\right\}=q^{2m-4}{:}\SL_{m-2}(q).
\]
This forces $H\cap K=(R\cap K).(q^{2m-4}{:}\SL_{m-2}(q))$, and the factorization $G=HK$ follows.
\end{proof}

In the next two lemmas, we give examples of factorizations $G=HK$ satisfying part~(b) of Lemma~\ref{LemOmegaPlus1Row4,9,10} with $q=2$. In this case, since $Q(e_1+f_1+u)=\mu=Q(u)$, there are elements in $L$ swapping $e_1+f_1+u$ and $u$, whence $L_{\{e_1+f_1+u,u\}}=L_{e_1+f_1+u,u}.2$.

\begin{lemma}\label{LemOmegaPlus15}
Let $G=L=\Omega_{2m}^+(2)$, let $H=T=\SL_m(2)$, and let $K=L_{\{e_1+f_1+u,u\}}=\Omega_{2m-2}^-(2).2$. Then $G=HK$ with $H\cap K=\SL_{m-2}(2)$.
\end{lemma}

For $q\in\{2,4\}$, recall from Subsection~\ref{SecOmegaPlus01} that $\gamma=r_{e_1+f_1}\cdots r_{e_m+f_m}$ is a product of $m$ reflections. Hence $\gamma\in\Omega(V)$ if and only if $m$ is even. Consequently, $\langle L,\gamma\rangle=\Omega_{2m}^+(q){:}(2,m-1)$. Moreover, as $r_{e_1+f_1}r_{e_2+f_2}\gamma$ fixes both $e_1+f_1$ and $u$, we have $\langle L,\gamma\rangle_{e_1+f_1,u}=\Omega_{2m-2}^-(q).(2,m-1)$.

\begin{lemma}\label{LemOmegaPlus14}
Let $G=\langle L,\gamma\rangle=\Omega_{2m}^+(2){:}(2,m-1)$, let $H=T{:}\langle\gamma\rangle=\SL_m(2){:}2$, and let $K=G_{e_1+f_1,u}=\Omega_{2m-2}^-(2).(2,m-1)$. Then $G=HK$ with $H\cap K=\SL_{m-2}(2)$.
\end{lemma}

In the next lemma, we construct a factorization $G=HK$ satisfying part~(b) of Lemma~\ref{LemOmegaPlus1Row4,9,10} with $q=4$. Recall from Subsection~\ref{SecOmegaPlus01} that, with
\[
u'=(1-\mu^2)e_1+(\mu^2-1+\mu^{-1})e_2+\mu^2f_1+\mu^2f_2,
\]
we have some $g\in\mathrm{O}(V)$ with $(e_1+f_1)^g=e_1+f_1$ and $u^g=u'$, and such an element $g$ exists both in $\Omega(V)$ and $\mathrm{O}(V)\Setminus\Omega(V)$.

\begin{lemma}\label{LemOmegaPlus16}
Let $G=L{:}\langle\rho\rangle=\Omega_{2m}^+(4){:}2$ with $\rho\in\{\phi,\gamma\phi\}$, let $H=T{:}\langle\rho\rangle=\SL_m(4){:}2$, and let $K=\langle\Omega(V)_{e_1+f_1,u},g\phi\rangle$ such that $g\phi\in G$, $(e_1+f_1)^g=e_1+f_1$ and $u^g=u'$. Then $G=HK$ with $K=\Omega_{2m-2}^-(4).4$ and $H\cap K=\SL_{m-2}(4)$.
\end{lemma}

Finally, we give an example below that satisfies part~(c) of Lemma~\ref{LemOmegaPlus1Row4,9,10}. Although we have assumed $(m,q)\neq(4,4)$, it is natural to still include $q=4$ to give a uniform description of the examples for $q=4$ and $q=16$. Note here that the triality automorphism of $\Omega_8^+(q)$ maps the $\mathcal{C}_3$-subgroup $\SU_4(q)$ to $\Omega_6^-(q)=\N_2^-[L]^{(\infty)}$ and maps the $\mathcal{C}_5$-subgroup $\Omega_8^-(q^{1/2})$ to a $\mathcal{C}_9$-subgroup.

\begin{lemma}\label{LemOmegaPlus23}
Let $G=\Omega_8^+(q){:}\langle\phi\rangle=\Omega_8^+(q){:}f$ with $q\in\{4,16\}$, let $H=\Omega_8^-(q^{1/2}).f$ be a $\mathcal{C}_5$-subgroup of $G$, and let $K=\SU_4(q).(2f)<\GU_4(q).(2f)$ be a $\mathcal{C}_3$-subgroup of $G$. Then $G=HK$ with $H\cap K=\SL_2(q^{1/2})\times2$.
\end{lemma}

\subsubsection{Classifying factorizations with $K\leqslant\N_2^-[G]$}

As a conclusion of this subsection, we classify the factorizations $G=HK$ with $K\leqslant\N_2^-[G]$.

\begin{proposition}\label{prop:O+N_2^-}
Let $K\leqslant\N_2^-[G]$. Then $G=HK$ if and only if $(G,H,K)$ tightly contains some $(G_0,H_0,K_0)$ in the following table, where $P=q^{(m-2)(m-3)/2}.q^{2m-4}$ and $d=(2,m-1)$. In this case, $H_0\cap K_0$ is described in the table.
\[
\begin{array}{lllll}
\hline
G_0 & H_0 & K_0 & H_0\cap K_0 & \textup{Remarks} \\
\hline
\Omega_{2m}^+(q) & q^{m(m-1)/2}{:}\SL_m(q)&\Omega_{2m-2}^-(q) & P{:}\SL_{m-2}(q) & \textup{as in \ref{LemOmegaPlus13}} \\
\Omega_{2m}^+(2) & \SL_m(2) & \Omega_{2m-2}^-(2).2 & \SL_{m-2}(2) & \textup{as in \ref{LemOmegaPlus15}} \\
\Omega_{2m}^+(2){:}d & \SL_m(2){:}2 & \Omega_{2m-2}^-(2).d & \SL_{m-2}(2) & \textup{as in \ref{LemOmegaPlus14}} \\
\Omega_{2m}^+(4){:}2 & \SL_m(4){:}2 & \Omega_{2m-2}^-(4).4 & \SL_{m-2}(4) & \textup{as in \ref{LemOmegaPlus16}} \\
\Omega_8^+(16){:}4 & \Omega_8^-(4).4 & \SU_4(16).8 & \SL_2(4)\times2 & \textup{as in \ref{LemOmegaPlus23}} \\
\hline
\end{array}
\]
\end{proposition}

\begin{proof}
By the lemmas shown in the last column of the table, it remains to prove the ``only if'' part.
Lemma~\ref{LemOmegaPlus1Row4,9,10} first reduces the possibilities for $(H^{(\infty)},K^{(\infty)})$ to the parabolic case,
the case $H^{(\infty)}=\SL_m(q)$ with $q\in\{2,4\}$, and the exceptional $\Omega_8^-(4)$ case. The parabolic and exceptional cases yield the first and last rows, respectively, the latter by~\cite[Remark~1.5]{GGP}.
In the $\SL_m(q)$ case, Lemmas~\ref{LemOmegaPlus15} and~\ref{LemOmegaPlus14} together with the subgroup structure of
$\N_2^-[G]$ determine the extensions $G_0/L$, $H_0/H^{(\infty)}$ and $K_0/K^{(\infty)}$ for $q=2$, while~\cite[Page~69,~(i)]{LPS1990} forces the required extension for $q=4$, yielding rows~2--4.
\end{proof}

\subsection{Factorizations with $A=\Pa_m[G]$ and $K^{(\infty)}=\Omega_{2m-1}(q)$}\label{SecOmegaPlus04}
\ \vspace{1mm}

Without loss of generality, let $A=G_U$, where $U=\langle e_1,\dots,e_m\rangle=\bbF_q^m$.
If $m=ab$, then let $U_{(b)}=\bbF_{q^b}^a$ with the same underlying set as $U$, and for $i\in\{1,\dots,\lfloor b/2\rfloor\}$ denote
\begin{equation}\label{EqnOrthonal1}
U_{(b)}^\sharp(i)=
\left\{\begin{array}{ll}
\bigoplus\limits_{j=0}^{b-1}\left(U_{(b)}\otimes U_{(b)}^{(q^i)}\right)^{(q^j)}&\text{if }1\leqslant i\leqslant\lfloor(b-1)/2\rfloor\vspace{2mm}\\
\bigoplus\limits_{j=0}^{b/2-1}\left(U_{(b)}\otimes U_{(b)}^{(q^i)}\right)^{(q^j)}&\text{if }\lfloor(b-1)/2\rfloor<i\leqslant\lfloor b/2\rfloor.
\end{array}\right.
\end{equation}
(Note that $\lfloor(b-1)/2\rfloor<i\leqslant\lfloor b/2\rfloor$ if and only if $i=b/2$ with $b$ even.) The $\bbF_{q^b}\GL(U_{(b)})$-module $U_{(b)}^\sharp(i)$ can be realized over $\bbF_q$, and let $U_{(b)}(i)$ denote the realized $\bbF_q\GL(U_{(b)})$-module. For each $i\in\{1,\dots,\lfloor b/2\rfloor\}$, since $U_{(b)}^\sharp(i)=U_{(b)}(i)\otimes\bbF_{q^b}$, the definition of $U_{(b)}^\sharp(i)$ in~\eqref{EqnOrthonal1} implies (see \cite[\S26)]{ASCH1993} for instance) that the $\bbF_qS$-module $U_{(b)}(i)$ is isomorphic to $U_{(b)}\otimes U_{(b)}^{(q^i)}$, where the latter is viewed as an $\bbF_qS$-module. Denote
\[
\mbox{$U_{(b)}(0)=\bfO_p\Big(\Pa_a[\Omega_{2a}^+(q^b)]\Big)=\bigwedge^2U_{(b)}=q^{ba(a-1)/2}$},
\]
where $\Omega_{2a}^+(q^b)$ is a field-extension subgroup of $L$ defined over $\bbF_{q^b}$.
%The action of $\SL(U)=\SL_a(q^b)$ is irreducible on $U_{(b)}(0)$, while the action of the $\calC_8$-subgroup $\Sp_a(q^b)$ of $\SL(U)$ on $U_{(b)}(0)$ is not.

\begin{lemma}\label{LemOmegaPlusPm1}
Let $m=ab$, and let $S$ be a subgroup of $\SL(U)$ contained in the group $\SL_a(q^b)$ defined over $\bbF_{q^b}$ such that $S=\SL_a(q^b)$, $\Sp_a(q^b)$ or $\G_2(q^b)$ (with $a=6$ and $q$ even). Then the following statements hold.
\begin{enumerate}[{\rm (a)}]
\item For $i\in\{1,\dots,\lfloor b/2\rfloor\}$, the group $U_{(b)}(i)$ is $q^{a^2b}$ if $i\neq b/2$, and is $q^{a^2b/2}$ if $i=b/2$.
\item The $\bbF_qS$-module $R$ is the direct sum of pairwise non-isomorphic submodules
\[
U_{(b)}(0),\,U_{(b)}(1),\,\dots,\,U_{(b)}(\lfloor b/2\rfloor).
\]
\item The irreducible $\bbF_qS$-submodules of $R$ are $U_{(b)}(1),\,\dots,\,U_{(b)}(\lfloor b/2\rfloor)$ and the irreducible $\bbF_qS$-submodules of $U_{(b)}(0)$.
\item If $H\leqslant G_U=\Pa_m[G]$ with $RH=RS$, then $H=(H\cap R){:}X$ for some $X\leqslant H$ with $X\cong S$.
\end{enumerate}
\end{lemma}
\begin{proof}
Part~(a) follows immediately from the definition of $U_{(b)}^\sharp(i)$.
Viewing $R$ as an $\bbF_q\GL(U)$-module, we may write $R=\bigwedge^2(U)$.
%(see~\cite[Page~22]{BL2021})
Since $U_{(b)}$, as an $\bbF_qS$-module, can be viewed as $U$, we have
\[
\mbox{$U\otimes\bbF_{q^b}=\bigoplus\limits_{i=0}^{b-1}U_{(b)}^{(q^i)}$.}
\]
It then follows that
\[\begin{array}{rl}
R\otimes\bbF_{q^b}=\bigwedge^2\left(U\otimes\bbF_{q^b}\right)
\,=&\left(\bigoplus\limits_{i=0}^{b-1}U_{(b)}^{(q^i)}\right)\wedge\left(\bigoplus\limits_{i=0}^{b-1}U_{(b)}^{(q^i)}\right)\vspace{2mm}\\
=&\bigoplus\limits_{i=0}^{b-1}\left(\bigwedge^2(U_{(b)})\right)^{(q^i)}
\oplus\bigoplus\limits_{0\leqslant i<j\leqslant b-1}U_{(b)}^{(q^i)}\wedge U_{(b)}^{(q^j)}.
\end{array}\]
If $b$ is odd, then
\[
\mbox{$\bigoplus\limits_{0\leqslant i<j\leqslant b-1}U_{(b)}^{(q^i)}\wedge U_{(b)}^{(q^j)}
=\bigoplus\limits_{i=1}^{(b-1)/2}\,\bigoplus\limits_{j=0}^{b-1}\left(U_{(b)}\otimes U_{(b)}^{(q^i)}\right)^{(q^j)}
=\bigoplus\limits_{i=1}^{(b-1)/2}U_{(b)}^\sharp(i)$}
\]
If $b$ is even, then
\[
\mbox{$\bigoplus\limits_{0\leqslant i<j\leqslant b-1}U_{(b)}^{(q^i)}\wedge U_{(b)}^{(q^j)}
=\bigoplus\limits_{i=1}^{b/2-1}\,\bigoplus\limits_{j=0}^{b-1}\left(U_{(b)}\otimes U_{(b)}^{(q^i)}\right)^{(q^j)}\oplus\bigoplus\limits_{j=0}^{b/2-1}\left(U_{(b)}\otimes U_{(b)}^{(q^{b/2})}\right)^{(q^j)}=\bigoplus\limits_{i=1}^{b/2}U_{(b)}^\sharp(i)$.}
\]
Thus, for all $b$,
\[
\mbox{$R\otimes\bbF_{q^b}=\bigoplus\limits_{i=0}^{b-1}\left(\bigwedge^2(U_{(b)})\right)^{(q^i)}
\oplus\bigoplus\limits_{0\leqslant i<j\leqslant b-1}U_{(b)}^{(q^i)}\wedge U_{(b)}^{(q^j)}
=\bigoplus\limits_{i=0}^{b-1}\left(\bigwedge^2(U_{(b)})\right)^{(q^i)}
\oplus\bigoplus\limits_{i=1}^{\lfloor b/2\rfloor}U_{(b)}^\sharp(i)$}.
\]
This realized over $\bbF_q$ shows that the $\bbF_qS$-module $R$ is the direct sum of $U_{(b)}(i)$ with $i$ running over $\{0,1,\dots,\lfloor b/2\rfloor\}$, which proves part~(b).

Since $U_{(b)}\otimes U_{(b)}^{(q^i)}$ is an irreducible $\bbF_{q^b}S$-module with highest weight $(1+q^i)\lambda_1$, where $\lambda_1,\dots,\lambda_{a-1}$ are the fundamental dominant weights of $S$, we derive that the $\bbF_{q^b}S$-modules
$U_{(b)}\otimes U_{(b)}^{(q^i)}$ with $i$ running over $\{1,\dots,\lfloor b/2\rfloor\}$ are pairwise non-isomorphic (see \cite[\S5.4)]{KL1990} for instance). Since they are also irreducible as $\bbF_qS$-modules, it follows from part~(a) that $U_{(b)}(1),\dots,U_{(b)}(\lfloor b/2\rfloor)$ are pairwise non-isomorphic irreducible $\bbF_qS$-modules. As the dimension of $U_{(b)}(0)$ is smaller than $U_{(b)}(i)$ for $i\in\{1,\dots,\lfloor b/2\rfloor\}$, the irreducible $\bbF_qS$-submodules of $U_{(b)}(0)$ are not isomorphic to any of $U_{(b)}(1),\dots,U_{(b)}(\lfloor b/2\rfloor)$. This in conjunction with part~(b) implies part~(c).

Let $R_1=(H\cap U_{(b)}(0))U_{(b)}(1)\cdots U_{(b)}(\lfloor b/2\rfloor)$. Then along the same lines as the proof of Lemma~\ref{LemUnitaryPm1}(d) with $R$ replaced by $R_1$, one obtains $H\cong\bfO_p(H){:}S$, proving part~(d) here.
\end{proof}

\begin{remark}
For the group $S$ in Lemma~\ref{LemOmegaPlusPm1}, the action of $S$ on $U_{(b)}(0)=\bigwedge^2U_{(b)}$ is well known (for example, see~\cite[Section~1]{Liebeck1985} and~\cite[Section~2]{Liebeck1987}). The group $S=\SL_a(q^b)$ is irreducible on $U_{(b)}(0)$. For $S=\Sp_a(q^b)$, there is a unique maximal submodule $\bbF_{q^b}^{a(a-1)/2-1}$ of $U_{(b)}(0)$, which has an irreducible quotient $\bbF_{q^b}^{a(a-1)/2-2}$ or $\bbF_{q^b}^{a(a-1)/2-1}$ according to whether $p$ divides $a/2$ or not. In particular, both $\Sp_a(q^b)$ and $\G_2(q^b)$ are reducible on $U_{(b)}(0)$.
\end{remark}

For a subset $I=\{i_1,\ldots,i_k\}$ of $\{1,\ldots,\lfloor b/2\rfloor\}$, denote
\begin{equation}\label{EqnOmegaPlus5}
\gcd(I,b)=\gcd(i_1,\dots,i_k,b)
\end{equation}
and let $U_{(b)}(I)=U_{(b)}(i_1)\cdots U_{(b)}(i_k)\leqslant R$. Then we derive from Lemma~\ref{LemOmegaPlusPm1} that
\[
U_{(b)}(I)=U_{(b)}(i_1)\times\cdots\times U_{(b)}(i_k)=
\begin{cases}
q^{(2k-1)a^2b/2}&\textup{if }b/2\in I\\
q^{ka^2b}&\textup{if }b/2\notin I.
\end{cases}
\]

\begin{lemma}\label{LemOmegaPlusPm2}
Let $m=ab=ae\ell$, let $M=\Nor_G(\Omega_{2ae}^+(q^\ell))$ be a field-extension subgroup of $G$ over $\bbF_{q^\ell}$, let $P=q^{\ell(ae)(ae-1)/2}$ be the unipotent radical of $\Pa_{ae}[M]$, and let $S\leqslant\GaL_a(q^b)$ defined over $\bbF_{q^b}$ such that $S^{(\infty)}=\SL_a(q^b)$, $\Sp_a(q^b)$ or $\G_2(q^b)$ (with $a=6$ and $q$ even). Then the following statements hold.
\begin{enumerate}[{\rm (a)}]
\item The $\bbF_{q^\ell}S^{(\infty)}$ module $P$ is a direct sum of $U_{(\ell)}(0),U_{(b)}(\ell),U_{(b)}(2\ell),\dots,U_{(b)}(\ell\lfloor e/2\rfloor)$.
\item If $H=(U_{(b)}(0)\cap H)U_{(b)}(I).S\leqslant\Pa_m[G]$ with $I\subseteq\{1,\ldots,\lfloor b/2\rfloor\}$ such that $\ell$ divides $\gcd(I,b)$, then $H\leqslant\Pa_{ae}[M]$ up to conjugation in $G$.
\end{enumerate}
\end{lemma}

\begin{proof}
Applying Lemma~\ref{LemOmegaPlusPm1} to the $\bbF_{q^\ell}S^{(\infty)}$ module $P$ we see that $P$ is a direct sum of $\bigwedge^2U_{(\ell)},P(1),\dots,P(\lfloor e/2\rfloor)$, where $P(t)$ is the realization of the $\bbF_{q^b}S^{(\infty)}$-module
\[
P^\sharp(t):=
\left\{\begin{array}{ll}
\bigoplus\limits_{s=0}^{e-1}\left(U_{(b)}\otimes U_{(b)}^{(q^{\ell t})}\right)^{(q^{\ell s})}&\text{if }1\leqslant t\leqslant\lfloor(e-1)/2\rfloor\vspace{2mm}\\
\bigoplus\limits_{s=0}^{e/2-1}\left(U_{(b)}\otimes U_{(b)}^{(q^{\ell t})}\right)^{(q^{\ell s})}&\text{if }\lfloor(e-1)/2\rfloor<t\leqslant\lfloor e/2\rfloor
\end{array}\right.
\]
over $\bbF_{q^\ell}$. Recall the definition of $U_{(b)}^\sharp(i)$ in~\eqref{EqnOrthonal1}. We obtain for $t\in\{1,\dots,\lfloor e/2\rfloor\}$ that
\[
\mbox{$\bigoplus\limits_{r=0}^{\ell-1}\bigoplus\limits_{s=0}^{e-1}P(t)^{(q^{\ell s+r})}
=\bigoplus\limits_{r=0}^{\ell-1}\left(\bigoplus\limits_{s=0}^{e-1}P(t)^{(q^{\ell s})}\right)^{(q^r)}
=\bigoplus\limits_{r=0}^{\ell-1}P^\sharp(t)^{(q^r)}=U_{(b)}^\sharp(\ell t)$,}
\]
and so $P(t)$ equals $U_{(b)}(\ell t)$, the realization of $U_{(b)}^\sharp(\ell t)$ over $\bbF_q$. This proves part~(a).

Now let $H=(U_{(b)}(0)\cap H)U_{(b)}(I).S\leqslant\Pa_m[G]$ with $I=\{i_1,\ldots,i_k\}\subseteq\{1,\ldots,\lfloor b/2\rfloor\}$ such that $\ell$ divides $\gcd(I,b)$. Write $i_s=\ell j_s$ for $s\in\{1,\dots,k\}$. It follows that, up to conjugation in $G$,
\begin{align*}
(U_{(b)}(0)\cap H)U_{(b)}(I)&=(U_{(b)}(0)\cap H)U_{(b)}(\ell i_1)\cdots U_{(b)}(\ell i_k)\\
&\leqslant U_{(b)}(0)U_{(b)}(\ell i_1)\cdots U_{(b)}(\ell i_k)\leqslant U_{(\ell)}(0)U_{(b)}(\ell i_1)\cdots U_{(b)}(\ell i_k)\leqslant P
\end{align*}
as $U_{(b)}(0)=\bfO_p\big(\Pa_a[\Omega_{2a}^+(q^b)]\big)\leqslant\bfO_p\big(\Pa_{ae}[\Omega_{2ae}^+(q^\ell)]\big)=U_{(\ell)}(0)$. Hence $H\leqslant\Pa_{ae}[M]$, as part~(b) asserts.
\end{proof}

For $H\leqslant G_U=\Pa_m[G]$, the notation $H^U$, as usual, denotes the induced group of $H$ on $U$. By Lemma~\ref{LemOmegaPlusPm1}, if $H^U=S\leqslant\GaL_a(q^b)$ defined over $\bbF_{q^b}$ such that $S^{(\infty)}=\SL_a(q^b)$, $\Sp_a(q^b)$ or $\G_2(q^b)$ (with $a=6$ and $q$ even), then we may write
\begin{equation}\label{EqnOmegaPlus3}
H=\big((U_{(b)}(0)\cap H)\times U_{(b)}(I)\big){:}S=(E\times U_{(b)}(I)){:}S=q^c{:}S
\end{equation}
for some $\bbF_{q^b}S$-submodule $E$ of $U_{(b)}(0)=\bigwedge^2U_{(b)}=\bigwedge^2\bbF_{q^b}^a$ and subset $I$ of $\{1,\ldots,\lfloor b/2\rfloor\}$, where
\begin{equation}\label{EqnOmegaPlus7}
c=\dim_{\bbF_q}(E)+\dim_{\bbF_q}(U_{(b)}(I))=
\begin{cases}
\log_q|E|+(2|I|-1)a^2b/2&\textup{if }b/2\in I\\
\log_q|E|+|I|a^2b&\textup{if }b/2\notin I.
\end{cases}
\end{equation}
The next two lemmas give some subgroups of $\Pa_m[G]$ that are transitive on $\calN_1^{(\infty)}$.
Note that, if $c=0$, then $I=\emptyset$ and so $\gcd(I,b)=b$.

\begin{lemma}\label{LemOmegaPlusPm3}
Let $G=\Omega_{2m}^+(q)$ with $m=ab$, let $K=\N_1[G]^{(\infty)}=\Omega_{2m-1}(q)$, and let $H\leqslant G_U=\Pa_m[G]$ satisfying~\eqref{EqnOmegaPlus3} such that $\gcd(I,b)=1$ and $S\leqslant\SL_a(q^b)$ is defined over $\bbF_{q^b}$.
\begin{enumerate}[{\rm (a)}]
\item If $S=\SL_a(q^b)$, then $G=HK$ with $H\cap K=[q^{c-b+1}].\SL_{a-1}(q^b)$.
\item If $S=\Sp_a(q^b)$, then $G=HK$ with $H\cap K=[q^{c-b+1}].\Sp_{a-2}(q^b)$.
\item If $S=\G_2(q^b)'$ with $a=6$ and $q$ even, then $G=HK$ with $H\cap K=[(q^{c-b+1},q^{c+1}/4)].\SL_2(q^b)$.
\end{enumerate}
\end{lemma}

\begin{proof}
Let $I=\{i_1,\ldots,i_k\}$, let $d=\gcd(i_1,\dots,i_k,a)$, and let $e$ be the largest divisor of $a$ coprime to $d$. It follows from $\gcd(i_1,\dots,i_k,b)=\gcd(I,b)=1$ that
\[
\gcd(be+i_1,\dots,be+i_k,m)=1.
\]
Let $C=\langle s\rangle=q^m-1$ be a Singer cycle in $\GL(U)$, let $P=\bfO_p(H)$, and let $X=PC=P{:}C$.
For each $t\in\{1,\dots,k\}$, by~\cite[Corollary~4.5]{FLWXZ} (with both $r$ and $q$ replaced by $q^b$), the restriction of $U_{(b)}(i_t)$ to $C$ can be decomposed as
\[
U_{(b)}(i_t,1)\oplus\cdots\oplus U_{(b)}(i_t,a)
\]
such that the value of $\chi_{(b)}(i_t,j)$ at $s$ is
\[
\sum^{m-1}_{\ell=0}(\omega^{q^{bj+i_t}+1})^{q^\ell},
\]
where $\chi_{(b)}(i_t,j)$ is the character of the $\bbF_qC$-module $U_{(b)}(i_t,j)$ and $\omega$ is a generator of $\bbF_{q^m}^\times$.
In particular, this character value is completely determined by $bj+i_t$ up to the action of ${\rm Gal}(\bbF_{q^m}/\bbF_q)$.
By~\cite[Theorem~1.9]{FLWXZ} (our module $U_{(b)}(i_t,j)$ is denoted as $U(bj+i_t)$ in \cite{FLWXZ}), we have $\mathrm{O}_{2m}^+(q)=(U_{(b)}(i_1,e)\cdots U_{(b)}(i_k,e)C)\mathrm{O}_{2m}^+(q)_v$, where $\mathrm{O}_{2m}^+(q)_v$ is the stabilizer of some $v\in\calN_1^{(\infty)}$ in $\mathrm{O}_{2m}^+(q)$. Since $U_{(b)}(i_1,e)\cdots U_{(b)}(i_k,e)\leqslant P$, it follows that
\[
\mathrm{O}_{2m}^+(q)=(PC)\mathrm{O}_{2m}^+(q)_v=X\,\mathrm{O}_{2m}^+(q)_v,
\]
and so the action of $X$ on $\calN_1^{(\infty)}$ is transitive. Since $X=P{:}C$ and $C=q^m-1$ is a Singer cycle in $\GL(U)$, this implies that the action of $\GL(U)$ on the set $\big(\calN_1^{(\infty)}\big)_P$ of $P$-orbits on $\calN_1^{(\infty)}$ is permutationally equivalent to its natural action on $\bbF_q^m\Setminus\{0\}$. Then we conclude from Theorem~\ref{ThmLinear} that the action of $S$ on $\big(\calN_1^{(\infty)}\big)_P$ is transitive with point stabilizer $q^{b(a-1)}{:}\SL_{a-1}(q^b)$, $[q^{b(a-1)}]{:}\Sp_{a-2}(q^b)$ or $[(q^{5b},q^{6b}/4)]{:}\SL_2(q^b)$ according to $S=\SL_a(q^b)$, $\Sp_a(q^b)$ or $\G_2(q^b)'$, respectively. Hence $H=P{:}S$ is transitive on $\calN_1^{(\infty)}$, and so $H_v$ has order $|H|/|\calN_1^{(\infty)}|$. Since $\big(\calN_1^{(\infty)}\big)_P$ is a block system for $\calN_1^{(\infty)}$, the point stabilizer $H_v$ necessarily stabilizes the block containing $v$. It follows that
\begin{equation}\label{EqnOmegaPlus4}
(H_v)^U\leqslant
\begin{cases}
q^{b(a-1)}{:}\SL_{a-1}(q^b)&\textup{if }S=\SL_a(q^b)\\
[q^{b(a-1)}]{:}\Sp_{a-2}(q^b)&\textup{if }S=\Sp_a(q^b)\\
[(q^{5b},q^{6b}/4)]{:}\SL_2(q^b)&\textup{if }S=\G_2(q^b)'.
\end{cases}
\end{equation}
As $|H_v|/|(H_v)^U|$ divides $|P|$, we have
\[
|(H_v)^U|_{p'}=|H_v|_{p'}=|H|_{p'}/|\calN_1^{(\infty)}|_{p'}=|H|_{p'}/(q^m-1)=|S|_{p'}/(q^m-1),
\]
which in conjunction with~\eqref{EqnOmegaPlus4} implies that $H_v=\bfO_p(H_v).\SL_{a-1}(q^b)$, $\bfO_p(H_v).\Sp_{a-2}(q^b)$ or $\bfO_p(H_v).\SL_2(q^b)$, according to $S=\SL_a(q^b)$, $\Sp_a(q^b)$ or $\G_2(q^b)'$, respectively. Combined with $|H_v|=|H|/|\calN_1^{(\infty)}|=q^c|S|/|\calN_1^{(\infty)}|$ and $|\calN_1^{(\infty)}|=q^{m-1}(q^m-1)$, this leads to
\[
|\bfO_p(H_v)|=
\begin{cases}
q^{c-b+1}&\textup{if }S=\SL_a(q^b)\textup{ or }\Sp_a(q^b)\\
(q^{c-b+1},q^{c+1}/4)&\textup{if }S=\G_2(q^b)'.
\end{cases}
\]
Thus $G=HK$ with $H\cap K=H_v$ as described in parts~(a)--(c) of the lemma.
\end{proof}

Recall that $b_2$ denotes the $2$-part of $b$. The next result follows from Lemmas~\ref{LemOmegaPlus12},~\ref{LemOmegaPlusPm2} and~\ref{LemOmegaPlusPm3}.

\begin{lemma}\label{LemOmegaPlusPm4}
Let $G=\Omega(V){:}\langle\phi\rangle=\Omega_{2m}^+(q){:}f$ with $q\in\{2,4\}$, let $K=\N_1[G]=\GaSp_{2m-2}(q)$, let $m=ab$ with $b$ even, and let $H\leqslant G_U=\Pa_m[G]$ satisfying~\eqref{EqnOmegaPlus3} such that $\gcd(I,b)=2$ and $S=S^{(\infty)}{:}\langle\theta\rangle\leqslant\SiL_a(q^b)$ is defined over $\bbF_{q^b}$, where $\theta$ is the field automorphism of $\SL_a(q^b)$ of order $fb_2$.
\begin{enumerate}[{\rm (a)}]
\item If $S^{(\infty)}=\SL_a(q^b)$, then $G=HK$ with $H\cap K=[q^{c-b+2}].\SL_{a-1}(q^b).(b_2/2)$.
\item If $S^{(\infty)}=\Sp_a(q^b)$, then $G=HK$ with $H\cap K=[q^{c-b+2}].\Sp_{a-2}(q^b).(b_2/2)$.
\item If $S^{(\infty)}=\G_2(q^b)$ with $a=6$, then $G=HK$ with $H\cap K=[q^{c-b+2}].\SL_2(q^b).(b_2/2)$.
\end{enumerate}
\end{lemma}

We conclude this subsection by the following proposition.

\begin{proposition}\label{PropOmegaPlusPm}
Let $H\leqslant\Pa_m[G]$ with $(m,q)\neq(4,2)$ or $(4,3)$, and let $K^{(\infty)}=\N_1[G]^{(\infty)}=\Omega_{2m-1}(q)$. Then $G=HK$ if and only if $(G,H,K)$ tightly contains some $(G_0,H_0,K_0)$ in the following table. In this case, $H_0\cap K_0$ is described in the table.
\[
\begin{array}{lllll}
\hline
G_0 & H_0 & K_0 & H_0\cap K_0 & \textup{Remarks} \\
\hline
\Omega_{2ab}^+(q) & q^c{:}\SL_a(q^b) & \Omega_{2ab-1}(q) & [q^{c-b+1}].\SL_{a-1}(q^b) & \textup{as in \ref{LemOmegaPlusPm3}} \\
 & q^c{:}\Sp_a(q^b) &  & [q^{c-b+1}].\Sp_{a-2}(q^b) & \textup{as in \ref{LemOmegaPlusPm3}} \\
 & q^c{:}\G_2(q^b)'  &  & [(q^{c-b+1},q^{c+1}/4)].\SL_2(q^b) & a=6,\ q\textup{ even}, \\
 &  &  &  & \textup{as in \ref{LemOmegaPlusPm3}} \\
\Omega_{2ab}^+(2) & 2^c{:}\SL_a(2^b){:}b_2 & \Sp_{2ab-2}(2) & [2^{c-b+2}].\SL_{a-1}(2^b).(b_2/2) & \textup{as in \ref{LemOmegaPlusPm4}} \\
\textup{($b$ even)} & 2^c{:}\Sp_a(2^b){:}b_2 &  & [2^{c-b+2}].\Sp_{a-2}(2^b).(b_2/2) & \textup{as in \ref{LemOmegaPlusPm4}} \\
 & 2^c{:}\G_2(2^b){:}b_2 &  & [2^{c-b+2}].\SL_2(2^b).(b_2/2) & a=6,\ \textup{as in \ref{LemOmegaPlusPm4}} \\
\Omega_{2ab}^+(4){:}\langle\phi\rangle & 2^{2c}{:}\SL_a(4^b){:}2b_2 & \GaSp_{2ab-2}(4) & [4^{c-b+2}].\SL_{a-1}(4^b).(b_2/2) & \textup{as in \ref{LemOmegaPlusPm4}} \\
\textup{($b$ even)} & 2^{2c}{:}\Sp_a(4^b){:}2b_2 &  & [4^{c-b+2}].\Sp_{a-2}(4^b).(b_2/2) & \textup{as in \ref{LemOmegaPlusPm4}} \\
 & 2^{2c}{:}\G_2(4^b){:}2b_2 &  & [4^{c-b+2}].\SL_2(4^b).(b_2/2) & a=6,\ \textup{as in \ref{LemOmegaPlusPm4}} \\
\Omega_{12}^+(3) & 3^{14}{:}\SL_2(13) & \Omega_{11}(3) & 3^{9+1} & \SL_2(13)<T \\
\hline
\end{array}
\]
\end{proposition}

\begin{proof}
For $L=\Omega_{12}^+(3)$, computation in \magma~\cite{BCP1997} verifies the conclusion of the proposition. By Lemmas~\ref{LemOmegaPlusPm3} and~\ref{LemOmegaPlusPm4}, we are left to prove the ``only if'' part. Suppose $G=HK$, and let $A=G_U=\Pa_m[G]$ and $B=\N_1[G]$ be maximal subgroups of $G$ containing $H$ and $K$ respectively.

Taking $a=m$ in Lemma~\ref{LemOmegaPlusPm3}(a) we deduce that $(A\cap B)^U$ stabilizes a $1$-space in $\bbF_q^m$. Thus it follows from $A=H(A\cap B)$ that $H^U$ is transitive on the set of $1$-spaces in $\bbF_q^m$, and so~\cite{Hering1985} asserts that $H^U\leqslant\GaL_a(q^b)$ is defined over $\bbF_{q^b}$ with $m=ab$ and $(H^U)^{(\infty)}=\SL_a(q^b)$, $\Sp_a(q^b)$ or $\G_2(q^b)'$ (with $a=6$ and $q$ even). Let $S=(H^U)^{(\infty)}$. Then by Lemma~\ref{LemOmegaPlusPm1} we may write
\[
H=(E\times U_{(b)}(I)){:}H^U=q^c{:}H^U
\]
for some $\bbF_{q^b}S$-submodule $E$ of $U_{(b)}(0)=\bigwedge^2U_{(b)}$ and some subset $I$ of $\{1,\ldots,\lfloor b/2\rfloor\}$. Let $\ell=\gcd(I,b)$. If $\ell=1$, then $(G,H,K)$ tightly contains the triple $(G^{(\infty)},H^{(\infty)},K^{(\infty)})$ as in Lemma~\ref{LemOmegaPlusPm3}. Now assume $\ell>1$. Then Lemma~\ref{LemOmegaPlusPm2} implies that $H$ is contained in some field-extension subgroup $M$ of $G$ defined over $\bbF_{q^\ell}$, and so Lemma~\ref{OmegaPlusClaim}(b) shows $\ell=2$ and $q\in\{2,4\}$. Moreover, Lemma~\ref{OmegaPlusClaim}(a) implies $\Omega(V)H\geqslant\Omega(V){:}\langle\phi\rangle$. As a consequence, $G\geqslant\Omega(V){:}\langle\phi\rangle=\Omega_{2m}^+(q){:}f$ and $H^U\geqslant S{:}\langle\theta\rangle=S{:}(fb)_2$, where $\theta$ is the field automorphism of $\SL_a(q^b)$ of order $(fb)_2=fb_2$. Hence $(G,H,K)$ tightly contains the triple $(G_0,H_0,K_0)$ as in Lemma~\ref{LemOmegaPlusPm4}.
\end{proof}

\subsection{Actions on ${\mathcal N}_1$}
\ \vspace{1mm}

In this subsection, we assume that $K\leqslant B=\N_1[G]$, which is the case when $K$ fixes a nonsingular vector. In view of the triality automorphism of $\POm_8^+(q)$, we note that this covers both cases~(d) and~(e) of Proposition~\ref{Prop:MaxO+subgps}.

\subsubsection{Candidates for $K^{(\infty)}$}

We first determine $K^{(\infty)}$ for the factorizations $G=HK$, as given in the following lemma.

\begin{lemma}\label{LemOmegaPlus1Row1--3,6--8,12--13}
Let $G=HK$ with $B=\N_1[G]$. Then either $K^{(\infty)}=\Omega_{2m-1}(q)$ or $\Omega^\pm_{2m-2}(q)$, or $H^{(\infty)}=\mathrm{Spin}_7(q)$ with $m=4$.
\end{lemma}

\begin{proof}
Apply Proposition~\ref{Prop:MaxO+subgps} to obtain candidates of $A^{(\infty)}$. If $K^{(\infty)}<B^{(\infty)}$, consider the factorization $B=(A\cap B)K$ (and its image under the quotient modulo $\Rad(B)$) and use primitive prime divisors of $q^{2m-2}-1$, $q^m-1$, $q^{m-1}-1$ and $q^{2m-4}-1$.
\end{proof}

Let $\,\overline{\phantom{\varphi}}\,$ be the quotient modulo scalars. Among the candidates given in the conclusion of Lemma~\ref{LemOmegaPlus1Row1--3,6--8,12--13}, the case $K^{(\infty)}=\Omega^\pm_{2m-2}(q)$ is treated in Propositions~\ref{prop:O+N_2^+} and \ref{prop:O+N_2^-}, and the case $H^{(\infty)}=\mathrm{Spin}_7(q)$ with $m=4$ can be reduced to that $K^{(\infty)}=\Omega_{2m-1}(q)$ by switching $\overline{H}$ and $\overline{K}$ and applying the triality automorphism of $\POm_8^+(q)$. Thus we only need to consider the case $K^{(\infty)}=\Omega_{2m-1}(q)$. In this case, since $e_1+f_1$ is a nonsingular vector, we may assume (applying an automorphism of $L$ if necessary)
\[
K^{(\infty)}=L_{e_1+f_1}.
\]
By Proposition~\ref{Prop:MaxO+subgps} and the triality automorphism of $\POm_8^+(q)$, either $A^{(\infty)}$ is one of the groups
\begin{align*}
&q^{m(m-1)/2}{:}\SL_m(q),\,\ \SL_m(q),\,\ \SU_m(q)\text{ with }m\text{ even},\,\ \Omega_m^+(q^2)\text{ with }q\in\{2,4\},\\
&(\Sp_2(q)\circ\Sp_m(q))^{(\infty)}\text{ with }q>2,\,\ \mathrm{Spin}_9(q)\text{ with }m=8,\,\ \Co_1\text{ with }(m,q)=(12,2),
\end{align*}
or $m=4$ and $\overline{A}^{(\infty)}\in\{\Omega_7(q),\Omega_8^-(q^{1/2})\}$. The latter is determined by Lemma~\ref{PropOmegaPlusO+8}, while the case $A^{(\infty)}=q^{m(m-1)/2}{:}\SL_m(q)$ is done in Proposition~\ref{PropOmegaPlusPm}. In the following, we deal with the remaining cases. Since their proofs follow the same general strategy as above---reducing to a factorization inside $A$, applying the relevant earlier classification, and then determining the outer extensions---we omit them and state only the resulting propositions.

\subsubsection{The case $A^{(\infty)}=\SL_m(q)$}

In this case, we may assume $A^{(\infty)}=T$ (applying an involutory graph automorphism of $L$ if necessary, see~\cite[Table~3.5.E]{KL1990}). Recall for even $q$ that $\gamma=r_{e_1+f_1}\cdots r_{e_m+f_m}$ normalizes $T$ and induces the graph automorphism on $T$ by swapping $e_i$ and $f_i$ for all $i\in\{1,\dots,m\}$, and that $\gamma$ is in $L$ if and only if $m$ is even.

\begin{proposition}\label{Prop:O^+=(SL,N_1)}
Let $A^{(\infty)}=T=\SL_m(q)$ with $(m,q)\neq(4,2)$, and let $K^{(\infty)}=L_{e_1+f_1}=\Omega_{2m-1}(q)$. Then $G=HK$ if and only if $(G,H,K)$ tightly contains some $(G_0,H_0,K_0)$ in the following table. In this case, $H_0\cap K_0$ is described in the table.
\[
\begin{array}{lllll}
\hline
G_0&H_0 &K_0& H_0\cap K_0 & \textup{Remarks}\\
\hline
\Omega^+_{2m}(q) & \SL_m(q) & \Omega_{2m-1}(q) & \SL_{m-1}(q) & \\
&\Sp_m(q) &&\Sp_{m-2}(q) & \\
&\G_2(q) &&\SL_2(q)& m=6,\ q\textup{ even} \\
\Omega^+_{2m}(2)&\SL_{m/2}(4){:}2& \Sp_{2m-2}(2) &\SL_{m/2-1}(4)& H_0<T\textup{ or }T{:}\langle\gamma\rangle\\
&\GaSp_{m/2}(4)& & \Sp_{m/2-2}(4) & H_0<T\textup{ or }T{:}\langle\gamma\rangle \\
&\GaG_2(4) &&\SL_2(4)&m=12,\ H_0<T\textup{ or }T{:}\langle\gamma\rangle \\
\Omega^+_{2m}(4){:}\langle\phi\rangle&\SL_{m/2}(16){:}4& \GaSp_{2m-2}(4)&\SL_{m/2-1}(16) & H_0<T{:}\langle\phi\rangle\textup{ or }T{:}\langle\phi\gamma\rangle\\
&\GaSp_{m/2}(16)& & \Sp_{m/2-2}(16)& H_0<T{:}\langle\phi\rangle\textup{ or }T{:}\langle\phi\gamma\rangle\\
&\GaG_2(16) & &\SL_2(16) &m=12,\ H_0<T{:}\langle\phi\rangle\textup{ or }T{:}\langle\phi\gamma\rangle\\
\hline
\end{array}
\]
\end{proposition}

\subsubsection{The case $A^{(\infty)}=\SU_m(q)$}

Recall $V_\sharp$ and $\beta_\sharp$ introduced in Subsection~\ref{SecOmegaPlus01}. Here we may assume $A^{(\infty)}=\SU(V_\sharp)$ (applying an involutory graph automorphism of $L$ if necessary, see~\cite[Table~3.5.E]{KL1990}). Also, by the results in Propositions~\ref{PropOmegaPlusPm} and~\ref{Prop:O^+=(SL,N_1)} we may assume that $H$ is not contained in any maximal subgroup of $G$ with solvable residual $q^{m(m-1)/2}{:}\SL_m(q)$ or $\SL_m(q)$.

\begin{proposition}\label{Prop:O^+=(SU,N_1)}
Let $A^{(\infty)}=\SU(V_\sharp)=\SU_m(q)$ with $m$ even and $(m,q)\neq(4,2)$, and let $K^{(\infty)}=L_{e_1+f_1}=\Omega_{2m-1}(q)$.
Suppose that $H$ is not contained in any maximal subgroup of $G$ with solvable residual $q^{m(m-1)/2}{:}\SL_m(q)$ or $\SL_m(q)$.
Then $G=HK$ if and only if $(G,H,K)$ tightly contains some $(G_0,H_0,K_0)$ in the following table. In this case, $H_0\cap K_0$ is described in the table.
\[
\begin{array}{llll}
\hline
G_0&H_0 &K_0& H_0\cap K_0 \\
\hline
\Omega^+_{2m}(q) & \SU_m(q) & \Omega_{2m-1}(q) & \SU_{m-1}(q) \\
\Omega^+_{12}(2) & 3^{\boldsymbol{\cdot}}\PSU_4(3) & \Sp_{10}(2) & 3^5{:}\A_5 \\
 & 3^{\boldsymbol{\cdot}}\M_{22} &  &  \PSL_2(11) \\
\Omega^+_{24}(2) &3^{\boldsymbol{\cdot}}\Suz & \Sp_{22}(2) & 3^5{:}\PSL_2(11) \\
\hline
\end{array}
\]
\end{proposition}

\subsubsection{The case $A^{(\infty)}=(\Sp_2(q)\circ\Sp_m(q))^{(\infty)}$}

Recall from Proposition~\ref{Prop:MaxO+subgps} that $q>2$ in this case. Note also that $(\Sp_2(q)\circ\Sp_m(q))^{(\infty)}=\Sp_2(q)\circ\Sp_m(q)$ for $q\geqslant4$.

\begin{proposition}\label{Prop:O^+=(tensor,N_1)}
Let $A^{(\infty)}=(\Sp_2(q)\circ\Sp_m(q))^{(\infty)}$ with $q>2$, and let $K^{(\infty)}=L_{e_1+f_1}=\Omega_{2m-1}(q)$.
Then $G=HK$ if and only if $(G,H,K)$ tightly contains some $(G_0,H_0,K_0)$ in the following table. In this case, $H_0\cap K_0$ is described in the table.
\[
\begin{array}{lllll}
\hline
G_0 & H_0 & K_0 & H_0\cap K_0 & \textup{Remarks}\\ \hline
\Omega^+_{2m}(q) & \Sp_m(q) & \Omega_{2m-1}(q) & \Sp_{m-2}(q) & \\
 &\G_2(q) & &\SL_2(q) & m=6,\ q\textup{ even} \\
%\Omega^+_{2m}(2) & \GaSp_{m/2}(4) & \Sp_{2m-2}(2) & \Sp_{m/2-2}(4) & H_0<T\textup{ or }T{:}\langle\gamma\rangle \\
% &\GaG_2(4) &&\SL_2(4)&m=12,\ H_0<T\textup{ or }T{:}\langle\gamma\rangle \\
\Omega^+_{2m}(4){:}\langle\phi\rangle & \GaSp_{m/2}(16) & \GaSp_{2m-2}(4) & \Sp_{m/2-2}(16) & H_0<T{:}\langle\phi\rangle\textup{ or }T{:}\langle\phi\gamma\rangle\\
 &\GaG_2(16) & &\SL_2(16) &m=12,\ H_0<T{:}\langle\phi\rangle\textup{ or }T{:}\langle\phi\gamma\rangle\\
\hline
\end{array}
\]
\end{proposition}

\subsubsection{$\calC_9$-subgroups}

As illustrated after Lemma~\ref{LemOmegaPlus1Row1--3,6--8,12--13}, the candidates to consider for $A$ being a $\calC_9$-subgroup of $G$ are $\mathrm{Spin}_9(q)$\ with $m=8$ and $\Co_1$ with $(m,q)=(12,2)$.

\begin{proposition}\label{C_9-subgroups}
Let $A$ be a $\calC_9$-subgroup of $G$ with $m>4$, and let $K^{(\infty)}=L_{e_1+f_1}=\Omega_{2m-1}(q)$. Then $G=HK$ if and only if $(G,H,K)$ tightly contains some $(G_0,H_0,K_0)$ in the following table. In this case, $H_0\cap K_0$ is described in the table.
\[
\begin{array}{lllll}
\hline
G_0 & H_0 & K_0 & H_0\cap K_0 & \textup{Remarks}\\
\hline
\Omega^+_{16}(q) & \mathrm{Spin}_9(q) & \Omega_{15}(q) & \mathrm{Spin}_7(q) & \\
\Omega^+_{16}(2) & \GaSp_4(4) & \Sp_{14}(2) & \Sp_2(4) & H_0<\mathrm{Spin}_9(2)\cong\Sp_8(2)\\
\Omega^+_{16}(4){:}\langle\phi\rangle & \GaSp_4(16) & \GaSp_{14}(4) & \Sp_2(16) & H_0<\mathrm{Spin}_9(4)\cong\Sp_8(4)\\
\Omega^+_{24}(2) & \Co_1 & \Sp_{22}(2) & \Co_3 & \\
 & 3^{\boldsymbol{\cdot}}\Suz,\ \G_2(4){:}2 & & 3^5{:}\PSL_2(11),\ \A_5 & H_0<\Co_1\\
\hline
\end{array}
\]
\end{proposition}

\subsubsection{The case $A^{(\infty)}=\Omega_m^+(q^2)$}

For $G=L{:}\langle\phi\rangle=\Omega_{2m}^+(q){:}f$, $M=\Omega(V_\sharp){:}\langle\psi\rangle=\Omega_m^+(q^2){:}(2f)$ and $K=G_{e_1+f_1}=\GaSp_{2m-2}(q)$, where $q\in\{2,4\}$ and $m$ is even, Lemma~\ref{LemOmegaPlus12} shows that $G=MK$ with $M\cap K=\Sp_{m-2}(q^2)$. Hence for each subgroup $H$ of $M$, we have $G=HK$ if and only if $M=H(M\cap K)$. This together with Propositions~\ref{Prop:O^+=(SU,N_1)} and~\ref{C_9-subgroups} leads to the example below.

\begin{example}\label{LemOmegaPlus35}
Let $G=L{:}\langle\phi\rangle=\Omega_{2m}^+(q){:}f$ with $q\in\{2,4\}$ and $m/2$ even, let $M=\Omega(V_\sharp){:}\langle\psi\rangle=\Omega_m^+(q^2){:}(2f)$, and let $K=G_{e_1+f_1}=\GaSp_{2m-2}(q)$.
\begin{enumerate}[{\rm (a)}]
\item If $H=\SU_{m/2}(q^2){:}(4f)<M$, then $G=HK$ with $H\cap K=\SU_{m/2-1}(q^2).2$.
\item If $H=\Sp_{m/2}(q^2).(4f)<\SU_{m/2}(q^2){:}(4f)<M$, then $G=HK$ with $H\cap K=\Sp_{m/2-2}(q^2).2$.
\item For $m=8$, if $H$ is a $\calC_9$-subgroup of $M$ isomorphic to $\GaSp_6(q^2)$ or $\Omega_8^-(q).(2f)$, then $G=HK$ with
\[
H\cap K=
\begin{cases}
\G_2(q^2)&\textup{if }H=\GaSp_6(q^2)\\
\G_2(q)&\textup{if }H=\Omega_8^-(q).(2f).
\end{cases}
\]
\item For $m=8$, if $M_1=\Omega_7(q^2){:}(2f)$ is a $\calC_9$-subgroup of $M$ and $H$ is a subgroup of $M_1$ isomorphic to $\Omega_6^+(q^2){:}(2f)$, $\Omega_6^-(q^2){:}(2f)$ or $\GaSp_4(q^2)$, then $G=HK$ with
\[
H\cap K=
\begin{cases}
\SL_3(q^2)&\textup{if }H=\Omega_6^+(q^2){:}(2f)\\
\SU_3(q^2)&\textup{if }H=\Omega_6^-(q^2){:}(2f)\\
\SL_2(q^2)&\textup{if }H=\GaSp_4(q^2).
\end{cases}
\]
\item For $m=12$, if $H=\G_2(q^2).(4f)<\Sp_6(q^2).(4f)<\SU_6(q^2){:}(4f)<M$, then $G=HK$ with $H\cap K=\SL_2(q^2).2$.
\item For $m=16$, if $H=\GaSp_8(q^2)=\mathrm{Spin}_9(q^2){:}(2f)$ is a $\calC_9$-subgroup of $M$, then $G=HK$ with $H\cap K=\Sp_6(q^2)$.\qedhere
\end{enumerate}
\end{example}

Now we classify the factorizations with $A^{(\infty)}=\Omega_m^+(q^2)$ and $K^{(\infty)}=\Omega_{2m-1}(q)$. By~\cite[Table~3.5.E]{KL1990}, applying an involutory graph automorphism of $L$ if necessary, we may assume that $A^{(\infty)}=\Omega(V_\sharp)$.

\begin{proposition}\label{prop:O+N1}
Let $A^{(\infty)}=\Omega(V_\sharp)=\Omega_m^+(q^2)$ with $q\in\{2,4\}$ and $m\geqslant6$, and let $K^{(\infty)}=L_{e_1+f_1}=\Omega_{2m-1}(q)$.
Suppose that $H\nleqslant\Pa_{m-1}[G]$ or $\Pa_m[G]$ and that $(G,H,K)$ is not as described in Proposition~$\ref{Prop:O^+=(SL,N_1)}$, $\ref{Prop:O^+=(SU,N_1)}$ or~$\ref{Prop:O^+=(tensor,N_1)}$. Then $G=HK$ if and only if $(G,H,K)$ tightly contains some $(G_0,H_0,K_0)$ in the following table. In this case, $H_0\cap K_0$ is described in the table.
\[
\begin{array}{lllll}
\hline
G_0 & H_0 & K_0 & H_0\cap K_0 & \textup{Remarks}\\
\hline
\Omega_{2m}^+(q){:}\langle\phi\rangle & \Omega_m^+(q^2){:}(2f) & \GaSp_{2m-2}(q) & \Sp_{m-2}(q^2) & \textup{as in \ref{LemOmegaPlus12}}\\
 & \SU_{m/2}(q^2){:}(4f) & & \SU_{m/2-1}(q^2).2 & m/2\textup{ even, as in \ref{LemOmegaPlus35}(a)}\\
 & \Sp_{m/2}(q^2).(4f) & & \Sp_{m/2-2}(q^2).2 & \textup{as in \ref{LemOmegaPlus35}(b)}\\
\Omega^+_{16}(q){:}\langle\phi\rangle & \Omega_8^-(q).(2f),\ \GaSp_6(q^2) & \GaSp_{14}(q) & \G_2(q),\ \G_2(q^2) & \textup{as in \ref{LemOmegaPlus35}(c)}\\
 & \Omega_6^\pm(q^2){:}(2f) & & \SL_3^\pm(q^2) & \textup{as in \ref{LemOmegaPlus35}(d)}\\
 & \GaSp_4(q^2) & & \SL_2(q^2) & \textup{as in \ref{LemOmegaPlus35}(d)}\\
\Omega^+_{24}(q){:}\langle\phi\rangle & \G_2(q^2).(4f) & \GaSp_{22}(q) & \SL_2(q^2).2 & \textup{as in \ref{LemOmegaPlus35}(e)}\\
\Omega^+_{32}(q){:}\langle\phi\rangle & \GaSp_8(q^2) & \GaSp_{30}(q) & \Sp_6(q^2) & \textup{as in \ref{LemOmegaPlus35}(f)}\\
\hline
\end{array}
\]
\end{proposition}

\subsection{Proof of Theorem~\ref{ThmOmegaPlus}}\label{SecOmegaPlus03}
\ \vspace{1mm}

Since the triples $(G_0,H_0,K_0)$ in Tables~\ref{TabOmegaPlus} and~\ref{TabOmegaPlus2} are already shown to give rise to factorizations $G_0=H_0K_0$, it suffices to prove the ``only if'' part. Let $G=HK$. By Proposition~\ref{Prop:MaxO+subgps} and Lemmas~\ref{PropOmegaPlusO+8} and~\ref{LemXia29}, and in view of the triality automorphism of $\POm_8^+(q)$, we only need to consider parts~(a)--(d) of Proposition~\ref{Prop:MaxO+subgps}. If part~(a) or~(b) appears, then $K\leqslant\Pa_1[G]$ or $\N_2^+[G]$, and Proposition~\ref{prop:O+N_2^+} asserts that $(G,H,K)$ tightly contains some $(G_0,H_0,K_0)$ in rows~15--18 of Table~\ref{TabOmegaPlus}. For part~(c), we see from Proposition~\ref{prop:O+N_2^-} that $(G,H,K)$ tightly contains some $(G_0,H_0,K_0)$ in rows~19--22 or~24 of Table~\ref{TabOmegaPlus}. Now assume that part~(d) appears. Then by Lemma~\ref{LemOmegaPlus1Row1--3,6--8,12--13} and the paragraph thereafter, we may assume that $K^{(\infty)}=L_{e_1+f_1}=\Omega_{2m-1}(q)$ and $A^{(\infty)}$ is one of
\begin{align*}
&q^{m(m-1)/2}{:}\SL_m(q),\,\ \SL_m(q),\,\ \SU_m(q)\text{ with }m\text{ even},\,\ \Omega_m^+(q^2)\text{ with }q\in\{2,4\},\\
&(\Sp_2(q)\circ\Sp_m(q))^{(\infty)}\text{ with }q>2,\,\ \mathrm{Spin}_9(q)\text{ with }m=8,\,\ \Co_1\text{ with }(m,q)=(12,2).
\end{align*}
For these candidates, Propositions~\ref{PropOmegaPlusPm}, \ref{Prop:O^+=(SL,N_1)}, \ref{Prop:O^+=(SU,N_1)}, \ref{Prop:O^+=(tensor,N_1)}, \ref{C_9-subgroups} and \ref{prop:O+N1} together show that $(G,H,K)$ tightly contains some $(G_0,H_0,K_0)$ in rows~1--14 of Table~\ref{TabOmegaPlus}, where the reference labels are listed in the last column of the table. This completes the proof.

%%%%%%%%%%%%%%%%%%%%%%%%%%%%%%%%%%%%%%%%%%%%%%%%%%%%%%%%%%%%%%%%%%%
%%%%%%%%%%%%%%%%%%%%%%%%%%%%%%%%%%%%%%%%%%%%%%%%%%%%%%%%%%%%%%%%%%%

\section{Symplectic groups}\label{SecSymplectic02}

%%%%%%%%%%%%%%%%%%%%%%%%%%%%%%%%%%%%%%%%%%%%%%%%%%%%%%%%%%%%%%%%%%%
%%%%%%%%%%%%%%%%%%%%%%%%%%%%%%%%%%%%%%%%%%%%%%%%%%%%%%%%%%%%%%%%%%%

In this section we tackle factorizations symplectic groups. The main result is the following theorem.

\begin{theorem}\label{ThmSymplectic}
Let $G$ be a classical group with $L=G^{(\infty)}=\Sp_{2m}(q)$, where $m\geqslant2$, and let $H$ and $K$ be subgroups of $G$ not containing $L$ such that both $H$ and $K$ have a unique nonsolvable composition factor. Then $G=HK$ if and only if $(G,H,K)$ tightly contains some $(G_0,H_0,K_0)$ in Table~$\ref{TabSymplectic1}$ or~$\ref{TabSymplectic2}$.
\end{theorem}

For convenience, we collect in Subsection~\ref{SecSymplectic01} the notation that will be used throughout this section. The remaining subsections are then devoted to the proof of Theorem~\ref{ThmSymplectic}.

\subsection{Notation}\label{SecSymplectic01}
\ \vspace{1mm}

Throughout this section, let $q=p^f$ be a power of a prime $p$, let $m\geqslant2$ be an integer such that $(m,q)\neq(2,2)$, let $V$ be a vector space of dimension $2m$ over $\bbF_q$ equipped with a nondegenerate alternating form $\beta$, let $e_1,f_1,\dots,e_m,f_m$ be a standard basis for $V$ as in~\cite[2.2.3]{LPS1990}, let $\gamma$ be the involution in $\Sp(V)$ swapping $e_i$ and $f_i$ for all $i\in\{1,\dots,m\}$, and let $\phi\in\GaSp(V)$ such that
\[
\phi\colon a_1e_1+b_1f_1+\dots+a_me_m+b_mf_m\mapsto a_1^pe_1+b_1^pf_1+\dots+a_m^pe_m+b_m^pf_m
\]
for $a_1,b_1\dots,a_m,b_m\in\bbF_q$, let $U=\langle e_1,\dots,e_m\rangle$, let $U_1=\langle e_1,\dots,e_{m-1}\rangle$, and let $W=\langle f_1,\dots,f_m\rangle$. From~\cite[3.5.4]{Wilson2009} we know that $\Sp(V)_U=\Pa_m[\Sp(V)]$ has a subgroup $R{:}T$, where
\[
R=q^{m(m+1)/2}
\]
is the kernel of $\Sp(V)_U$ acting on $U$, and
\[
T=\SL_m(q)
\]
stabilizes both $U$ and $W$ (the action of $T$ on $U$ determines that on $W$ in the way described in~\cite[Lemma~2.2.17]{BG2016}).
Take $\mu\in\bbF_q$ such that the polynomial $x^2+x+\mu$ is irreducible over $\bbF_q$.

If $m=ab$, then let $V_{(b)}$ be a vector space of dimension $2a$ over $\bbF_{q^b}$ with the same underlying set as $V$, let $\Tr_{(b)}$ be the trace of the field extension $\bbF_{q^b}/\bbF_q$. In this case, we may equip $V_{(b)}$ with a nondegenerate alternating form $\beta_{(b)}$ such that $\beta(v,w)=\Tr_{(b)}(\beta_{(b)}(v,w))$ for all $v,w\in V$, and thus $\GaSp(V_{(b)})<\GaSp(V)$ (see~\cite[\S4.3]{KL1990}).

If $m=2\ell$, then let $V_1=\langle e_1,f_1,\dots,e_\ell,f_\ell\rangle$, let $V_2=\langle e_{\ell+1},f_{\ell+1},\dots,e_m,f_m\rangle$, and let $\sigma$ be the linear transformation swapping $e_i$ with $e_{\ell+i}$ and $f_i$ with $f_{\ell+i}$ for $i\in\{1,\dots,\ell\}$.
In this case, there is a subgroup $(\Sp(V_1)\times\Sp(V_2)){:}\langle\sigma\rangle=\Sp_m(q)\wr\Sy_2$ of $\Sp(V)$.

Now we give some notation assuming that $q$ is even.
For $\varepsilon\in\{+,-\}$, let $Q^\varepsilon$ be a nondegenerate quadratic form of type $\varepsilon$ on $V$ with associated bilinear form $\beta$ such that $e_1,f_1,\dots,e_m,f_m$ is a standard basis for the orthogonal space $(V,Q^\varepsilon)$. Thus
\[
Q^\varepsilon(e_i)=Q^\varepsilon(f_i)=Q^+(e_m)=Q^+(f_m)=0,\quad Q^-(e_m)=1,\quad Q^-(f_m)=\mu
\]
for all $\varepsilon\in\{+,-\}$ and $i\in\{1,\dots,m-1\}$. If $m=ab$, then let $Q_{(b)}^\varepsilon$ be a nondegenerate quadratic form of type $\varepsilon$ on $V_{(b)}$ with associated bilinear form $\beta_{(b)}$ such that $Q^\varepsilon(v)=\Tr_{(b)}(Q_{(b)}^\varepsilon(v))$ for all $v\in V$.

% Sometimes it is convenient to regard $\Sp_{2m}(q)$ as the subgroup $\Omega_{2m+1}(q)$ fixing a nonsingular vector in $\Omega_{2m+2}^+(q)$.
% Thus let $V'$ be a vector space of dimension $2m+2$ over $\bbF_q$, extend the alternating form $\beta$ with standard basis $e_1,f_1,\dots,e_m,f_m$ to an alternating form on $V'$ with standard basis $e_1,f_1,\dots,e_{m+1},f_{m+1}$, equip $V'$ with a nondegenerate quadratic form $Q$ of plus type with associated bilinear form $\beta$ such that $e_1,f_1,\dots,e_{m+1},f_{m+1}$ is a standard basis for the orthogonal space $(V',Q)$, and let
% \[
% d=e_{m+1}+f_{m+1}.
% \]
% Then $Q|_V=Q^+$, and $\Omega(V')_d\cong\Sp_{2m}(q)$.
% % Extend $\phi$ to $\phi\in\GaO(V')$ such that
% % \[
% % \phi\colon a_1e_1+b_1f_1+\dots+a_{m+1}e_{m+1}+b_{m+1}f_{m+1}\mapsto a_1^2e_1+b_1^2f_1+\dots+a_{m+1}^2e_{m+1}+b_{m+1}^2f_{m+1}
% % \]
% % for $a_1,b_1\dots,a_{m+1},b_{m+1}\in\bbF_q$.
% Let
% \[
% u=e_m+\mu f_m+e_{m+1},
% \]
% let $U=\langle e_1,\dots,e_m\rangle_{\bbF_q}$, let $U_1=\langle e_1,\dots,e_{m-1}\rangle_{\bbF_q}$, and let $W=\langle f_1,\dots,f_m\rangle_{\bbF_q}$.
% Then $\Omega(V')_{d,U}=\Pa_m[\Omega(V')_d]$ has a subgroup $R{:}T$, where
% \[
% R=q^{m(m+1)/2}
% \]
% is the kernel of $\Omega(V')_d$ acting on $U$, and
% \[
% T=\SL_m(q)
% \]
% stabilizes both $U$ and $W$.

\subsection{Reduction}
\ \vspace{1mm}

Let $L=G^{(\infty)}=\Sp_{2m}(q)$ with $m\geqslant2$, and let $H$ and $K$ be subgroups of $G$ such that both $H$ and $K$ have a unique nonsolvable composition factor and neither $H$ nor $K$ contains $L$. If $G=HK$, then $H$ and $K$ are contained in some subgroups $A$ and $B$ of $G$ described in~\cite[Theorem~A]{LPS1990} and~\cite{LPS1996}. (Note that the exceptional maximal factorization in~\cite[Table~3]{LPS1990} for $(m,q)=(2,3)$ does not have $B$ nonsolvable.) We list the candidates in the following proposition according to our analysis procedure.

\begin{proposition}\label{lem:K<B-Sp}
Let $(m,q)\neq(3,3)$ or $(4,2)$. If $G=HK$, then with the above notation and interchanging $A$ and $B$ if necessary, one of the following holds:
\begin{enumerate}[{\rm(a)}]
\item $B=\Pa_1[G]$, and $H$ is transitive on the set of $1$-subspaces of $\bbF_q^{2m}$;
\item $B=\N_2[G]$ with $q$ even, and $A^{(\infty)}=\Sp_{m}(q^2)$ with $q\in\{2,4\}$, or $\G_2(q)'$ with $m=3$;
%\item $\{A^{(\infty)},B^{(\infty)}\}=\{\Omega_{2m}^+(q),\Omega_{2m}^-(q)\}$ with $q\in\{2,4\}$;
\item $B^{(\infty)}=\Omega_{2m}^+(q)$ with $q$ even, and $A^{(\infty)}=\Omega_{2m}^-(q)$ with $q\in\{2,4\}$, $\Sp_{2a}(q^b)$ with $m=ab$ and $b$ prime, $\Sz(q)$ with $m=2$, or $\G_2(q)'$ with $m=3$;
\item $B^{(\infty)}=\Omega_{2m}^-(q)$ with $q$ even, and $A^{(\infty)}=\Sp_m(q)\times\Sp_m(q)$, $\Sp_{2m}(q^{1/2})'$ with $q\in\{4,16\}$, $\Sp_{2a}(q^b)$ with $m=ab$ and $b$ prime, $q^{m(m+1)/2}{:}\SL_m(q)$, or $\G_2(q)'$ with $m=3$.
\end{enumerate}
\end{proposition}

Before discussing cases~(a)--(d) in Proposition~\ref{lem:K<B-Sp}, we first treat the small groups
\[
L\in\{\Sp_4(3),\Sp_4(4),\Sp_6(2),\Sp_6(3),\Sp_6(4),\Sp_8(2)\}
\]
by computation in \magma~\cite{BCP1997}.

\begin{lemma}\label{LemSymplectic52}
Let $L=\Sp_4(3)$, $\Sp_4(4)$, $\Sp_6(2)$, $\Sp_6(3)$, $\Sp_6(4)$ or $\Sp_8(2)$. Then $G=HK$ if and only if $(G,H,K)$ tightly contains some $(G_0,H_0,K_0)$ in Table~$\ref{TabSymplectic1}$ or rows~\emph{1},~\emph{8--15} and~\emph{17--19} in Table~$\ref{TabSymplectic2}$, where they satisfy the following extra conditions.
\begin{enumerate}[{\rm(a)}]
\item In rows~\emph{1},~\emph{9},~\emph{10} and~\emph{19} of Table~$\ref{TabSymplectic2}$, for each isomorphism type of $H_0$ in the following table, the number $h$ of conjugacy classes of $H_0$ in $G_0$ satisfying $G_0=H_0K_0$ is given in the table.
\[
\begin{array}{c|ccccccccccccc}
\hline
\textup{row} & 1 & 9 & 9 & 10 & 19 & 19 & 19 & 19 & 19 & 19 & 19 & 19 & 19 \\
H_0 & \Sy_6 & \Sy_5\times2 & \Sy_6 & \Sy_5 & \Sy_5 & \A_5{:}4 & \A_6 & 2^5{:}\A_6 & \A_7 & 2^6{:}\A_7 & \A_8 & \A_9 & \Sp_6(2) \\
h & 1 & 1 & 2 & 2 & 1 & 2 & 1 & 2 & 1 & 1 & 2 & 1 & 1 \\
\hline
\end{array}
\]
\item In row~\emph{19} of Table~$\ref{TabSymplectic2}$, for each of the isomorphism types $\Sy_5\times2$, $\A_6\times2$ and $\Sy_6$ of $H_0$, the group $H_0$ lies in the unique conjugacy class of subgroups of $G_0$ such that $(H_0)'\cap K_0=H_0\cap K_0$ is isomorphic to $2$, $\Sy_3$ and $\Sy_3$, respectively.
\end{enumerate}
\end{lemma}

Lemma~\ref{LemSymplectic52} enables us to assume in the following that $L=G^{(\infty)}=\Sp_{2m}(q)$ with $m\geqslant2$ and
\[
(m,q)\neq(2,3),\ (2,4),\ (3,2),\ (3,3),\ (3,4)\text{ or }(4,2).
\]
We remark that $\Sp_{2m}(q)$ has a unique conjugacy class of subgroups isomorphic to $\Sp_{2m-2}(q)$ that are contained in $\N_2[\Sp_{2m}(q)]$. If $q$ is even, then $\Sp_{2m}(q)$ has another conjugacy class of subgroups isomorphic to $\Sp_{2m-2}(q)$, which are contained in $\N_1[\mathrm{O}_{2m}^+(q)]$ as well as in $\N_1[\mathrm{O}_{2m}^-(q)]$. To sum up, $\Sp_{2m}(q)$ has exactly $(2,q)$ conjugacy classes of subgroups isomorphic to $\Sp_{2m-2}(q)$.

\subsection{Preparation}
\ \vspace{1mm}

We will need the following lemma several times. As the convention before, throughout this subsection, $L=G^{(\infty)}=\Sp_{2m}(q)$ with $m\geqslant2$, while $H$ and $K$ are subgroups of $G$ such that both $H$ and $K$ have a unique nonsolvable composition factor and neither $H$ nor $K$ contains $L$.

\begin{lemma}\label{LemXia13}
Let $H$ be a field-extension subgroup of $G$ such that $H^{(\infty)}=\Sp_{2d}(q^e)$ with $m=de$, and let $K^{(\infty)}=\Sp_{2m-2}(q)$ such that $G=HK$. Then $H^{(\infty)}=\Sp_m(q^2)$ with $q\in\{2,4\}$.
\end{lemma}

\begin{proof}
Let $X=\Nor_G(H^{(\infty)})$. Then $G=XK$ and $X\cap L=\Sp_{2d}(q^e).e$. By the preceding description of the conjugacy classes of subgroups of $L$ isomorphic to $\Sp_{2m-2}(q)$, either $K\leqslant\N_2[G]$, or $q$ is even and $K\leqslant\N_1[B]$ for some $B$ with
$B\cap L=\GO_{2m}^{\pm}(q)$.
In the first case, $G=X\N_2[G]$, and the argument in~\cite[3.2.1(a)]{LPS1990} yields $e=2$ and $q\in\{2,4\}$. In the second case,
\[
X\cap\N_1[B]\geqslant
\GO_{2d-1}(q^e)=\Sp_{2d-2}(q^e)\times2.
\]
Comparing $2$-parts in $G=X\N_1[B]$ again forces $e=2$ and $q\in\{2,4\}$.
\end{proof}

\begin{remark}
In the condition of Lemma~\ref{LemXia13} we do not require $e$ to be prime.
\end{remark}

The next lemma gives candidates for $H^{(\infty)}$ and $K^{(\infty)}$ if $G=HK$ with $K\leqslant\Pa_1[G]$. Recall that the factorizations $G=HK$ with $G^{(\infty)}=\Sp_6(q)$ and one of $H$ or $K$ contained in $\GaG_2(q)$ are classified in Proposition~\ref{LemSymplectic12}.

\begin{lemma}\label{P_1-infty<Sp}
Let $G=HK$ with $K\leqslant B=\Pa_1[G]$ and $(m,q)\neq(3,2)$. Then one of the following holds:
\begin{enumerate}[{\rm(a)}]
\item $H^{(\infty)}=\Sp_{2a}(q^b)$ with $m=ab$ or $\G_2(q^b)$ with $m=3b$ and $q$ even, defined over $\bbF_{q^b}$, and $K^{(\infty)}=\Pa_1[G]^{(\infty)}=[q^{2m-1}]{:}\Sp_{2m-2}(q)$;
\item $H^{(\infty)}=\Sp_m(q^2)$ with $q\in\{2,4\}$, and $K^{(\infty)}=\Sp_{2m-2}(q)$;
\item $L=\Sp_4(q)$, $H^{(\infty)}=\Sp_2(q^2)$, and $K^{(\infty)}=q^{1+2}{:}\SL_2(5)$ with $q\in\{9,11,19,29,59\}$;
\item $L=\Sp_6(q)$ with $q\geqslant4$ even, and $H^{(\infty)}=\G_2(q)$;
\item $L=\Sp_{12}(q)$ with $q\in\{2,4\}$, $H^{(\infty)}=\G_2(q^2)$, and $K^{(\infty)}=\Sp_{10}(q)$.
\end{enumerate}
\end{lemma}

\begin{proof}
Since $G=H\Pa_1[G]$, the group $H$ is transitive on $\calP_1[G]$. Thus~\cite{Hering1985} gives $H^{(\infty)}=\Sp_{2a}(q^b)$ with $m=ab$, or $H^{(\infty)}=\G_2(q^b)$ with $m=3b$ and $q$ even.
Put $N=\Rad(B)$. If the image of $K$ in $B/N$, whose socle is $\PSp_{2m-2}(q)$, contains the socle, then $K^{(\infty)}$ is either $[q^{2m-1}]{:}\Sp_{2m-2}(q)$ or $\Sp_{2m-2}(q)$. This gives~(a), or, by Lemma~\ref{LemXia13}, the cases in~(b) and~(e).

If the image of $K$ is core-free, Proposition~\ref{lem:K<B-Sp} applied to $B/N$ leaves, for $m=2$, precisely the groups
$q^{1+2}{:}\SL_2(5)$ in~(c). For $m=3$, the remaining possibilities are excluded by the order divisibility in the
factorization of $B/N$ (the cases $q\in\{5,8,11\}$ are checked in \magma~\cite{BCP1997}); the $\G_2(q)$ alternative is~(d).
For $m\geqslant4$, primitive prime divisors of $q^{2m-2}-1$ and $q^{2m-4}-1$ exclude every core-free factor listed in Proposition~\ref{lem:K<B-Sp}. This proves the lemma.
\end{proof}

The candidates for $H^{(\infty)}$ and $K^{(\infty)}$ in a factorization $G=HK$ with $K\leqslant\N_2[G]$ are given in the subsequent lemma.

\begin{lemma}\label{LemSymplecticRow10,11}
Let $G=HK$ with $K\leqslant B=\N_2[G]$ and $(m,q)\neq(3,2)$ or $(3,4)$. Then $q$ is even, and one of the following holds:
\begin{enumerate}[{\rm (a)}]
\item $H^{(\infty)}=\Sp_m(q^2)$ with $q\in\{2,4\}$, and $K^{(\infty)}=\Sp_{2m-2}(q)$;
\item $L=\Sp_6(q)$, and $H^{(\infty)}=\G_2(q)$;
\item $L=\Sp_{12}(q)$ with $q\in\{2,4\}$, $H^{(\infty)}=\G_2(q^2)$ and $K^{(\infty)}=\Sp_{10}(q)$.
\end{enumerate}
\end{lemma}

\begin{proof}
The proof follows the same quotient argument as that of Lemma~\ref{P_1-infty<Sp}: Proposition~\ref{lem:K<B-Sp} gives the initial possibilities, and after passing to the almost simple quotient of $B$ with socle $\Sp_{2m-2}(q)$, the remaining cases follow from
Proposition~\ref{LemSymplectic12}, \cite[Lemma~4.2]{LPS2010} and Lemma~\ref{LemXia13}.
\end{proof}

We take the chance to give a corrected version of~\cite[Lemma~5.1]{LX2019}, which will be used in the proof of Lemma~\ref{LemSymplectic1Row2,3} below.

\begin{lemma}\label{LemXia9}
Let $L=\Sp_4(2^f)$ with $f\geqslant3$, and let $X$ be a subgroup of $G$ such that $G=X\N_2[G]$. Then one of the following holds:
\begin{enumerate}[{\rm (a)}]
\item $X\geqslant L$;
\item $X=\GaSp_2(16)$.
\end{enumerate}
\end{lemma}

\begin{remark}
Part~(b) of Lemma~\ref{LemXia9} is missing in~\cite[Lemma~5.1]{LX2019}. This is caused by the missing case $X\cap L\leqslant\Sp_2(4^2).2$ with $q=4$ in the proof of~\cite[Lemma~5.1]{LX2019} (here $\Sp_2(4^2).2$ is a $\calC_3$-subgroup of $L$), and leads to the triple $(L,H\cap L,K\cap L)=(\Sp_{6}(4),(\Sp_2(4)\times\Sp_{2}(16)).2,\G_2(4))$ missing in~\cite[Table~1]{LX2019}. The proof of~\cite[Lemma~5.1]{LX2019} is still valid under the assumption $f\geqslant3$, which together with computation in \magma~\cite{BCP1997} for $f=2$ yields the above Lemma~\ref{LemXia9}.
\end{remark}

The proofs of the next three lemmas are omitted, as they all use the standard induced factorizations $A=H(A\cap B)$ and $B=(A\cap B)K$,
together with primitive prime divisor arguments.

\begin{lemma}\label{LemSymplectic1Row2,3}
Let $G=HK$ such that $A^{(\infty)}=\Sp_{2a}(q^b)$ with $q$ even, $m=ab$ and $b$ prime, and $B^{(\infty)}=\Omega_{2m}^\varepsilon(q)$ with $\varepsilon\in\{+,-\}$. Then one of the following holds:
\begin{enumerate}[{\rm (a)}]
\item $K^{(\infty)}=B^{(\infty)}=\Omega_{2m}^\varepsilon(q)$;
\item $H^{(\infty)}=A^{(\infty)}=\Sp_m(q^2)$ with $b=2$ and $q\in\{2,4\}$, and $K^{(\infty)}=\Omega_{2m-1}(q)$;
\item $L=\Sp_{12}(q)$ with $q\in\{2,4\}$, $H^{(\infty)}=\G_2(q^2)$, and $K^{(\infty)}=\Omega_{11}(q)$.
\end{enumerate}
\end{lemma}

The next lemma shows that, in a factorization $G=HK$ with $B^{(\infty)}=\Omega_{2m}^+(q)$, either $H^{(\infty)}=A^{(\infty)}=\Omega_{2m}^-(q)$ or $K^{(\infty)}\in\{\Omega_{2m}^+(q),\Omega_{2m-1}(q)\}$.

\begin{lemma}\label{LemSymplectic01}
Let $G=HK$ with $B^{(\infty)}=\Omega_{2m}^+(q)$ and $(m,q)\neq(2,4)$, $(3,2)$ or $(3,4)$. Then one of the following holds:
\begin{enumerate}[{\rm (a)}]
\item $K^{(\infty)}=B^{(\infty)}=\Omega_{2m}^+(q)$;
\item $K^{(\infty)}=\N_1[B]^{(\infty)}=\Omega_{2m-1}(q)$, and $A^{(\infty)}\neq\Omega_{2m}^-(q)$;
\item $q\in\{2,4\}$, and $H^{(\infty)}=A^{(\infty)}=\Omega_{2m}^-(q)$.
\end{enumerate}
\end{lemma}

The next result is parallel to Lemma~\ref{LemSymplectic01}.

\begin{lemma}\label{LemSymplectic1Rows4--9}
Let $G=HK$ with $B^{(\infty)}=\Omega_{2m}^-(q)$ and $(m,q)\neq(2,4)$, $(3,2)$ or $(3,4)$. Then one of the following holds:
\begin{enumerate}[{\rm (a)}]
\item $K^{(\infty)}=B^{(\infty)}=\Omega_{2m}^-(q)$;
\item $K^{(\infty)}=\N_1[B]^{(\infty)}=\Omega_{2m-1}(q)$, and $A^{(\infty)}\neq\Omega_{2m}^+(q)$;
\item $m=3$, $K^{(\infty)}=\Pa_1[B]^{(\infty)}=q^4{:}\Omega_4^-(q)$, and $A^{(\infty)}=\G_2(q)$;
\item $q\in\{2,4\}$, and $H^{(\infty)}=A^{(\infty)}=\Omega_{2m}^+(q)$.
\end{enumerate}
\end{lemma}

Recall the notation $V_{(b)}$ and $Q_{(b)}^\varepsilon$ defined in Subsection~\ref{SecSymplectic01}. Based on the Lemma on Page~23 of~\cite{LPS1990}, we have the following observation.

\begin{lemma}\label{LemSymplectic03}
Let $q$ be even, $m=ab$ and $\varepsilon\in\{+,-\}$. Then the following statements hold:
\begin{enumerate}[{\rm (a)}]
\item $\Sp(V_{(b)})\cap\mathrm{O}(V,Q^\varepsilon)=\mathrm{O}(V_{(b)},Q_{(b)}^\varepsilon)$;
\item if $b$ is odd then $\Sp(V_{(b)})\cap\Omega(V,Q^\varepsilon)=\Omega(V_{(b)},Q_{(b)}^\varepsilon)$;
\item if $b$ is even then $\Sp(V_{(b)})\cap\Omega(V,Q^\varepsilon)=\mathrm{O}(V_{(b)},Q_{(b)}^\varepsilon)$.
\end{enumerate}
\end{lemma}

From Lemma~\ref{LemSymplectic03} we immediately obtain the example below.

\begin{example}\label{LemSymplectic04}
Let $G=\Sp(V)=\Sp_{2m}(q)$ with $q$ even, let $H=\Sp(V_{(b)})=\Sp_{2a}(q^b)$ with $m=ab$, and let $K=\Omega_{2m}^\varepsilon(q){:}(2,b)<G$ with $\varepsilon\in\{+,-\}$. Then $G=HK$ with
\[
H\cap K=\Omega_{2a}^\varepsilon(q^b){:}(2,b)=
\begin{cases}
\Omega_{2a}^\varepsilon(q^b)&\textup{if }b\textup{ is odd}\\
\mathrm{O}_{2a}^\varepsilon(q^b)&\textup{if }b\textup{ is even}.
\end{cases}\qedhere
\]
\end{example}

We also need the following observation.

\begin{lemma}\label{LemSymplectic02}
Let $G=HK$ with $K^{(\infty)}=\Omega_{2m}^+(q)$, and let $H$ be contained in some field-extension subgroup $M$ of $G$ with $M^{(\infty)}=\Sp_{2a}(q^b)$, where $m=ab$. Then $K\cap L\geqslant\Omega_{2m}^+(q){:}(2,b)$.
\end{lemma}

\subsection{Actions on $\calP_1$ and $\calN_2$}
\ \vspace{1mm}

In this subsection, we determine the factorizations $G=HK$ with $K\leqslant B=\Pa_1[G]$ or $\N_2[G]$.
Recall that Lemmas~\ref{P_1-infty<Sp} and~\ref{LemSymplecticRow10,11} give the candidates for $H^{(\infty)}$ and $K^{(\infty)}$ for such factorizations. Now we construct the minimal (with respect to tight containment) examples of factorizations in Lemmas~\ref{P_1-infty<Sp} and~\ref{LemSymplecticRow10,11}. Also recall that the factorizations in Lemma~\ref{P_1-infty<Sp}(d) and Lemma~\ref{LemSymplecticRow10,11}(b) are classified in Proposition~\ref{LemSymplectic12}. Before introducing the minimal ones for Lemma~\ref{P_1-infty<Sp}(a), note that
\[
\Pa_1[\Sp_{2m}(q)]^{(\infty)}=[q^{2m-1}]{:}\Sp_{2m-2}(q)'=
\begin{cases}
q^{1+(2m-2)}{:}\Sp_{2m-2}(q)&\textup{if }q\textup{ is odd}\\
q^{2m-1}{:}\Sp_{2m-2}(q)'&\textup{if }q\textup{ is even}
\end{cases}
\]
and that $\Sp_4(2)'=\A_6$ while $\Sp_{2m-2}(q)'=\Sp_{2m-2}(q)$ if $(m,q)\neq(3,2)$.

\begin{example}\label{ex:K<P1<Sp}
Let $G=\Sp_{2m}(q)$ with $m=ab$, and let $K=\Pa_1[G]^{(\infty)}$. If $(m,q)\neq(3,2)$, then Example~\ref{ex:Linear} implies $G=HK$ with $H=\Sp_{2a}(q^b)$ or $\G_2(q^b)$ ($a=3$, $q$ even), where the intersection of the two factors is $[q^{2ab-b}]{:}\Sp_{2a-2}(q^b)$ or $q^{2b+3b}{:}\SL_2(q^b)$, respectively. If $(m,q)=(3,2)$, then $G=\Sp_2(2^3)K=\SU_3(3)K$ with
\[
\Sp_2(2^3)\cap K=2^2\ \text{ and }\ \SU_3(3)\cap K=4^2{:}3.
\]
To summarize, we have $G=\Sp_{2a}(q^b)K$ with $\Sp_{2a}(q^b)\cap K=[q^d]{:}\Sp_{2a-2}(q^b)$, where
\begin{equation}\label{EqnSymplectic13}
d=
\begin{cases}
2ab-b&\textup{if }(a,b,q)\neq(1,3,2)\\
2&\textup{if }(a,b,q)=(1,3,2),
\end{cases}
\end{equation}
and $G=\G_2(q)'K$ with $\G_2(q)'\cap K=[(q^{5b},q^{6b}/4)]{:}\SL_2(q^b)'$.
\end{example}

For $G=\Sp_4(q)$ with odd $q$, Example~\ref{ex:K<P1<Sp} implies the factorization $G=\Pa_1[G]\,\Sp_2(q^2)$. Based on this, we construct some new factorizations, which are examples for Lemma~\ref{P_1-infty<Sp}(c). Let $N=q^{1+2}$ be the unipotent radical of $\Pa_1[G]$, let $K=\Sp(V_{(2)})=\Sp_2(q^2)$, and let $\Delta=K\backslash G$ be the set of right cosets of $K$ in $G$, so that $|\Delta|=|G|/|K|=q^2(q^2-1)$. Then $\Pa_1[G]$ is transitive on $\Delta$, and so $\GL_1(q)\times\SL_2(q)=\Pa_1[G]/N$ is transitive on the set $\Delta_N$ of orbits of $N$ on $\Delta$, which is of size $q^2-1$. Moreover, the action of $\GL_1(q)\times\SL_2(q)$ is permutationally equivalent to its natural action on $\bbF_q^2\Setminus\{0\}$. Then~\cite{Hering1985} implies:

\begin{example}\label{lem:K<P1-Sp(4,q)}
Let $G=\Sp_4(q)$ with $q\in\{9,11,19,29,59\}$, and $K=\Sp_2(q^2)$ be a field-extension subgroup of $G$. For $H=q^{1+2}{:}S$ with $S$ in the following table, where $q^{1+2}$ is the unipotent radical of $\Pa_1[G]$ and $S$ is a subgroup of the Levi complement $\GL_1(q)\times\Sp_2(q)$ of $\Pa_1[G]$, we have $G=HK$ with $H\cap K$ described in the table.
\[
\begin{array}{c|ccccc}
\hline
q & 9&11&19&29&59 \\
S & \SL_2(5)\times8 & \SL_2(5) & \SL_2(5)\times9 & \SL_2(5)\times7 & \SL_2(5)\times29 \\
H\cap K & 3^3{:}4 & 11 & 19{:}3 & 29 & 59\\
\hline
\end{array}
\]
\par\vspace{-1\baselineskip}
\qedhere
\end{example}

Computation in \magma \cite{BCP1997} shows that the following is a minimal example for Lemma~\ref{P_1-infty<Sp}(c) besides those in Example~\ref{lem:K<P1-Sp(4,q)}.

\begin{example}\label{ex:S5<S6<Sp(4,9).2}
Let $G=\Sp_4(9){:}\langle\phi\rangle$, $H=3^{2+4}{:}\SL_2(5){:}\langle\phi\rangle<\Pa_1[G]$, and $K=\SiL_2(81)$ be a maximal subgroup of $G$. Then $G=HK$ with $H\cap K=3^3{:}2$.
\end{example}

Recall the following example from the list of refined-antiflag-transitive groups in Theorem~\ref{refined-antiflags}.

\begin{example}\label{LemSymplectic13}
Let $G=\GaSp_{2m}(q)$ with $m$ even and $q\in\{2,4\}$, let $H=\GaSp_m(V_{(2)})=\GaSp_m(q^2)$, and let $K=\Sp_{2m-2}(q){:}f<\N_2[G]$ such that $K\cap L=\Sp_{2m-2}(q)$. Then $G=HK$ with $H\cap K=\Sp_{m-2}(q^2)<\N_2[H^{(\infty)}]$.
\end{example}

Example~\ref{LemSymplectic13} together with the ensuing lemma gives factorizations for Lemma~\ref{P_1-infty<Sp}(b) and Lemma~\ref{LemSymplecticRow10,11}(a).

\begin{lemma}\label{LemSymplectic19}
Let $G=\GaSp_{2m}(q)$ with $m$ even, and let $H=\GaSp(V_{(2)})=\GaSp_m(q^2)$.
\begin{enumerate}[{\rm (a)}]
\item If $q=2$ and $K=\Omega_{2m-1}(2)\times2<\mathrm{O}(V,Q^+)$, then $G=HK$ with $H\cap K=\Omega_{m-1}(4)\times2$.
\item If $q=4$ and $K=(\Omega_{2m-1}(4)\times2).2<\GaO(V,Q^+)$ such that $K\cap\mathrm{O}(V,Q^+)=\Omega_{2m-1}(4)\times2$, then $G=HK$ with $K=(\Omega_{2m-1}(4).2)\times2=\Omega_{2m-1}(4).2^2$ and $H\cap K=\Omega_{m-1}(16)\times2$.
\item If $q=4$ and $K=\Omega_{2m-1}(4).4<\GaO(V,Q^-)$ such that $K\cap\Omega(V,Q^-)=\Omega_{2m-1}(4)$, then $G=HK$ with $H\cap K=\Omega_{m-1}(16).2$.
\end{enumerate}
\end{lemma}

\begin{proof}
For (a) and (b), put $M=\GaO(V,Q^+)$. Example~\ref{LemSymplectic04} gives $G=HM$ with $H\cap M=\GaO_m^+(q^2)$, and Lemma~\ref{LemOmegaPlus12} gives $M=(H\cap M)K$ with $(H\cap M)\cap K=\Omega_{m-1}(q^2)\times2$. For (c), use instead
$M=\GaO(V,Q^-)$. Then Example~\ref{LemSymplectic04} and Lemma~\ref{ex:OmegaMinus02} give the asserted intersection
$\Omega_{m-1}(16).2$. In each case the two factorizations combine to give $G=HK$.
\end{proof}

The final example in this subsection, as shown below, constructs factorizations for Lemma~\ref{P_1-infty<Sp}(e) and Lemma~\ref{LemSymplecticRow10,11}(c). Part~(a) of the example follows from Example~\ref{LemSymplectic13} and Lemma~\ref{LemSymplectic008}, while parts~(b) and~(c) follow from Lemmas~\ref{LemSymplectic19} and~\ref{LemSymplectic008}.

\begin{example}\label{LemSymplectic36}
Let $G=\GaSp_{12}(q)$ with $m=6$ and $q\in\{2,4\}$, and let $H=\GaG_2(q^2)<\GaSp_6(q^2)=\GaSp(V_{(2)})$.
\begin{enumerate}[{\rm (a)}]
\item If $K$ is as in Example~$\ref{LemSymplectic13}$, then $G=HK$ with $H\cap K=\SL_2(q^2)$.
\item If $K$ is as in part~(a) or~(b) of Lemma~$\ref{LemSymplectic19}$, then $G=HK$ with $H\cap K=\SL_2(q^2)\times2$.
\item If $K$ is as in part~(c) of Lemma~$\ref{LemSymplectic19}$, then $G=HK$ with $H\cap K=\SL_2(16).2$.\qedhere
\end{enumerate}
\end{example}

%\subsubsection{Classifying factorizations with $K\leqslant\Pa_1[G]$ or $\N_2[G]$}\

We are now ready to present a classification of $G=HK$ with $K\leqslant B=\Pa_1[G]$ or $\N_2[G]$.

\begin{proposition}\label{prop:K=P1-Sp}
Let $(m,q)\neq(2,4)$, $(3,2)$ or $(3,4)$, and let $K\leqslant\Pa_1[G]$ or $\N_2[G]$. Then $G=HK$ if and only if $(G,H,K)$ tightly contains some $(G_0,H_0,K_0)$ in the following table. In this case, $H_0\cap K_0$ is described in the table.
{\small
\[
\begin{array}{lllll}
\hline
G_0 & H_0 & K_0 & H_0\cap K_0 & \textup{Remarks} \\
\hline
\Sp_{2ab}(q) & \Sp_{2a}(q^b) & [q^{2ab-1}]{:}\Sp_{2ab-2}(q) & [q^{2ab-b}]{:}\Sp_{2a-2}(q^b) & \textup{as in \ref{ex:K<P1<Sp}}\\
\Sp_{6b}(2^f) & \G_2(2^{fb}) & 2^{6fb-f}{:}\Sp_{6b-2}(2^f) & 2^{2fb+3fb}{:}\SL_2(2^{fb}) & \textup{as in \ref{ex:K<P1<Sp}}\\
\Sp_{2m}(2) & \GaSp_m(4) & \Sp_{2m-2}(2) & \Sp_{m-2}(4) & \textup{as in \ref{LemSymplectic13}}\\
 & & \Omega_{2m-1}(2)\times2 & \Omega_{m-1}(4)\times2 & \textup{as in \ref{LemSymplectic19}}\\
\GaSp_{2m}(4) & \GaSp_m(16) & \Sp_{2m-2}(4){:}2 & \Sp_{m-2}(16) & \textup{as in \ref{LemSymplectic13}}\\
 & & \Omega_{2m-1}(4).2^2 & \Omega_{m-1}(16)\times2 & \textup{as in \ref{LemSymplectic19}}\\
 & & \Omega_{2m-1}(4).4 & \Omega_{m-1}(16){:}2 & \textup{as in \ref{LemSymplectic19}}\\
\Sp_6(2^f) & \Sp_4(2^f) & \G_2(2^f) & \SL_2(2^f) & \textup{as in \ref{LemSymplectic12}}\\
 & 2^{5f}{:}\Sp_4(2^f) & & 2^{2f+3f}{:}\SL_2(2^f) & \textup{as in \ref{LemSymplectic12}}\\
 & 2^{4f}{:}\Omega^-_4(2^f) & & [2^{3f}] & \textup{as in \ref{LemSymplectic12}}\\
\Sp_4(9) & 3^{2+4}{:}(\SL_2(5)\times8) & \Sp_2(81) & 3^3{:}4 & \textup{as in \ref{lem:K<P1-Sp(4,q)}}\\
\Sp_4(9){:}\langle\phi\rangle & (3^{2+4}{:}\SL_2(5)){:}\langle\phi\rangle & \SiL_2(81) & 3^3{:}2 & \textup{as in \ref{ex:S5<S6<Sp(4,9).2}}\\
\Sp_4(11) & 11^{1+2}{:}\SL_2(5) & \Sp_2(11^2) & 11 & \textup{as in \ref{lem:K<P1-Sp(4,q)}}\\
\Sp_4(19) & 19^{1+2}{:}(\SL_2(5)\times9) & \Sp_2(19^2) & 19{:}3 & \textup{as in \ref{lem:K<P1-Sp(4,q)}}\\
\Sp_4(29) & 29^{1+2}{:}(\SL_2(5)\times7) & \Sp_2(29^2) & 29 & \textup{as in \ref{lem:K<P1-Sp(4,q)}}\\
\Sp_4(59) & 59^{1+2}{:}(\SL_2(5)\times29) & \Sp_2(59^2) & 59 & \textup{as in \ref{lem:K<P1-Sp(4,q)}}\\
\Sp_{12}(2) & \GaG_2(4) & \Sp_{10}(2) & \SL_2(4) & \textup{as in \ref{LemSymplectic36}(a)}\\
 & & \Omega_{11}(2)\times2 & \SL_2(4)\times2 & \textup{as in \ref{LemSymplectic36}(b)}\\
\GaSp_{12}(4) & \GaG_2(16) & \Sp_{10}(4){:}2 & \SL_2(16) & \textup{as in \ref{LemSymplectic36}(a)}\\
 & & \Omega_{11}(4).2^2 & \SL_2(16)\times2 & \textup{as in \ref{LemSymplectic36}(b)}\\
 & & \Omega_{11}(4).4 & \SL_2(16){:}2 & \textup{as in \ref{LemSymplectic36}(c)}\\
\hline
\end{array}
\]
}
\end{proposition}

\begin{proof}
The ``if'' part follows from the results recorded in the last column of the table. For the ``only if'' part, suppose that $G=HK$. If $K\leqslant\Pa_1[G]$, then Lemma~\ref{P_1-infty<Sp} applies; if $K\leqslant\N_2[G]$, then Lemma~\ref{LemSymplecticRow10,11} gives the cases corresponding to parts~(b), (d) and~(e) of Lemma~\ref{P_1-infty<Sp}. Thus one of parts~(a)--(e) of that lemma occurs.

In part~(a), $(G,H,K)$ tightly contains $(G^{(\infty)},H^{(\infty)},K^{(\infty)})$, as in the first two rows of
the table. In part~(c), we have $L=\Sp_4(q)$ with $q\in\{9,11,19,29,59\}$, and computation in \magma~\cite{BCP1997} yields precisely the corresponding minimal triples in the table. In part~(d), Proposition~\ref{LemSymplectic12} shows that, up to interchanging the two factors, $(G,H,K)$ tightly contains the listed triple with
\[
K^{(\infty)}\in
\{\Sp_4(q),q^5{:}\Sp_4(q),q^4{:}\Omega_4^-(q)\}.
\]
In parts~(b) and~(e), the required minimal triples are exactly those in Lemmas~\ref{LemSymplectic13} and~\ref{LemSymplectic19}, and Example~\ref{LemSymplectic36}, where Lemma~\ref{LemSymplectic02} is applied to determine the necessary outer extension.
\end{proof}

\subsection{Factorizations with $\{A^{(\infty)},B^{(\infty)}\}=\{\Omega_{2m}^+(q),\Omega_{2m}^-(q)\}$}
\ \vspace{1mm}

In this subsection, we discuss factorizations $G=HK$ with $\{A^{(\infty)},B^{(\infty)}\}=\{\Omega_{2m}^+(q),\Omega_{2m}^-(q)\}$ and $q\in\{2,4\}$, as a special case in part~(c) of Proposition~\ref{lem:K<B-Sp}. In this case, Lemma~\ref{LemSymplectic1Rows4--9} shows that, under the assumption $(m,q)\neq(2,4)$ and interchanging $H$ and $K$ if necessary, either $A^{(\infty)}=\Omega_{2m}^-(q)$ and $K^{(\infty)}=\Omega_{2m}^+(q)$, or $A^{(\infty)}=\Omega_{2m}^+(q)$ and $K^{(\infty)}=\Omega_{2m}^-(q)$.

\begin{example}\label{LemSymplectic14}
According to~\cite[3.2.4(e)]{LPS1990}, the maximal factorization $\Sp_{2m}(2)=\mathrm{O}_{2m}^-(2)\mathrm{O}_{2m}^+(2)$ has the intersection $\mathrm{O}_{2m}^-(2)\cap\mathrm{O}_{2m}^+(2)=\Omega_{2m-1}(2)\times2$ maximal in both the factors $\mathrm{O}_{2m}^-(2)$ and $\mathrm{O}_{2m}^+(2)$. This leads to the factorizations
\[
\Sp_{2m}(2)=\mathrm{O}_{2m}^-(2)\Omega_{2m}^+(2)=\Omega_{2m}^-(2)\mathrm{O}_{2m}^+(2)
\]
with $\mathrm{O}_{2m}^-(2)\cap\Omega_{2m}^+(2)=\Omega_{2m}^-(2)\cap\mathrm{O}_{2m}^+(2)=\Omega_{2m-1}(2)$.
\end{example}

For each $w\in V$ with $Q^+(w)\neq0$, let $r_w\in\mathrm{O}(V,Q^+)$ be the reflection defined in Subsection~\ref{SecOmegaPlus01}.

\begin{lemma}\label{LemSymplectic16}
Let $G=\GaSp_{2m}(4)$, let $H=\GaO_{2m}^-(4)$, and let $K=\Omega(V,Q^+){:}\langle\rho\rangle=\Omega_{2m}^+(4){:}2$ with $\rho\in\{\phi,r_{e_1+f_1}\phi\}$. Then $G=HK$ with $H\cap K=\Omega_{2m-1}(4)$.
\end{lemma}

\begin{proof}
Let $B=\GaO_{2m}^+(4)$. By~\cite[3.2.4(e)]{LPS1990}, $H\cap B=\Omega_{2m-1}(4)\times2$, while $K$ has index $2$ in $B$. Thus $H\cap K$ has index at most $2$ in $H\cap B$.
It cannot equal $H\cap B$, since $K\cap L=\Omega_{2m}^+(4)$ contains no subgroup $\Omega_{2m-1}(4)\times2$. Hence $H\cap K=\Omega_{2m-1}(4)$, and so $G=HK$.
\end{proof}

The first proposition of this subsection classifies factorizations $G=HK$ with $A^{(\infty)}=\Omega_{2m}^-(q)$ and $K^{(\infty)}=\Omega_{2m}^+(q)$.

\begin{proposition}\label{prop:Sp(2)=O^-O^+}
Let $q\in\{2,4\}$ and $(m,q)\neq(2,4)$, let $A^{(\infty)}=\Omega_{2m}^-(q)$, and let $K^{(\infty)}=\Omega_{2m}^+(q)$.
Then $G=HK$ if and only if $(G,H,K)$ tightly contains some $(G_0,H_0,K_0)$ in the following table. In this case, $H_0\cap K_0$ is described in the table.
\[
\begin{array}{lllll}
\hline
G_0 & H_0 & K_0 & H_0\cap K_0 & \textup{Remarks}\\ \hline
\Sp_{2m}(2) & \mathrm{O}_{2m}^-(2) & \Omega_{2m}^+(2) & \Omega_{2m-1}(2) & \\
 & \SU_m(2){:}2 & & \SU_{m-1}(2) & m\textup{ odd}\\
\Sp_{2m}(2) & \Omega_{2m}^-(2) & \mathrm{O}_{2m}^+(2) & \Omega_{2m-1}(2) & \\
 & \SU_m(2) & & \SU_{m-1}(2) & m\textup{ odd}\\
 & \GaO_m^-(4) & & \Omega_{m-1}(4){:}2 & \\
 & \SU_{m/2}(4).4 & & \SU_{m/2-1}(4){:}2 & m/2\textup{ odd}\\
\GaSp_{2m}(4) & \GaO_{2m}^-(4) & \Omega_{2m}^+(4){:}2 & \Omega_{2m-1}(4) & K_0\neq\mathrm{O}_{2m}^+(4)\\
 & \SU_m(4).4 & & \SU_{m-1}(4) & m\textup{ odd},\ K_0\neq\mathrm{O}_{2m}^+(4)\\
\hline
\end{array}
\]
\end{proposition}

\begin{proof}
The ``if'' part follow by combining Example~\ref{LemSymplectic14} or Lemma~\ref{LemSymplectic16} with the appropriate factorization in
Theorem~\ref{ThmOmegaMinus}. Conversely, suppose $G=HK$. For $q=2$, the conclusion follows by considering the factorization $A=H(A\cap K)$ and applying Example~\ref{LemSymplectic14} and Theorem~\ref{ThmOmegaMinus}.
For $q=4$, we consider $A=H(A\cap B)$ and apply~\cite[Theorem~A]{LPS1990} and Lemma~\ref{LemSymplectic16}.
\end{proof}

Recall the definition of $T$ and $\gamma$ in Subsection~\ref{SecSymplectic01}. Note that, for odd $m$ and even $q$, the group $T{:}\langle\gamma\rangle=\SL_m(q){:}2$ is contained in $\mathrm{O}(V,Q^+)$ such that $(T{:}\langle\gamma\rangle)\cap\Omega(V,Q^+)=T=\SL_m(q)$.
%To be as general as possible, we include the previously excluded pairs $(2,4)$, $(3,2)$, $(3,4)$ and $(4,2)$ for $(m,q)$ in the following Example~\ref{LemSymplectic30}.
Combining Example~\ref{LemSymplectic16} with Theorem~\ref{ThmOmegaPlus}, we obtain the next example (the conclusion for $m\in\{2,3\}$ can be directly verified by computation in \magma~\cite{BCP1997}).

\begin{example}\label{LemSymplectic30}
Let $G=\GaSp_{2m}(4)$, let $M=\Omega(V,Q^+){:}\langle\rho\rangle=\Omega_{2m}^+(4){:}2$ with $\rho\in\{\phi,r_{e_1+f_1}\phi\}$, and let $K=\GaO_{2m}^-(4)$.
\begin{enumerate}[{\rm(a)}]
\item If $H=N{:}2<M$ such that $N=H\cap\Omega(V,Q^+)$ is in the following table, then $G=HK$ with $H\cap K$ shown in the same table.
\[
\begin{array}{c|ccccc}
\hline
N & \SL_m(4) & \SU_m(4) & \Sp_m(4) & \G_2(4) & \Omega_9(4)\\
H\cap K & \SL_{m-1}(4) & \SU_{m-1}(4) & \Sp_{m-2}(4) & \SL_2(4) & \Omega_7(4)\\
\textup{Condition} & & m\text{ even} & & m=6 & m=8\\
\hline
\end{array}
\]
\item If $m=4$, then for each isomorphism type in the $H$ row in the following table, there are precisely $k$ conjugacy classes of subgroups $H$ of $M$ such that $M=H(M\cap K)$, and for each $H$, we have $G=HK$ with $H\cap K$ shown in the same table.
\[
\begin{array}{c|ccccc}
\hline
H & \Omega_8^-(2).2 & \Sp_6(4).2 & \Omega_6^+(4).2 & \Omega_6^-(4).4 & \Sp_4(4).2\\
k & 2 & 2 & 4 & 2 & 2 \\
H\cap K & \G_2(2) & \G_2(4) & \SL_3(4) & \SU_3(4).2 & \SL_2(4) \\
\hline
\end{array}
\]
\item If $H=T{:}\langle\gamma\rangle=\SL_m(4){:}2$ with $m$ odd, then $G=HK$ with $H\cap K=\SL_{m-1}(4)$.\qedhere
\end{enumerate}
\end{example}

Before stating the next proposition, we remark that, in the case $K^{(\infty)}=\Omega_{2m}^-(q)$, the candidates with $H^{(\infty)}=\Omega_{2m}^+(q)$ have been dealt with in Proposition~\ref{prop:Sp(2)=O^-O^+} (with $H$ and $K$ interchanged), and the candidates with $H\leqslant\Pa_m[G]$ will be handled in Proposition~\ref{PropSymplecticPm}.

\begin{proposition}\label{prop:Sp(2)=O^+O^-}
Let $q\in\{2,4\}$, let $A^{(\infty)}=\Omega_{2m}^+(q)$, and let $K^{(\infty)}=\Omega_{2m}^-(q)$. Suppose $H^{(\infty)}\neq\Omega_{2m}^+(q)$ and $H\nleqslant\Pa_m[G]$. Then $G=HK$ if and only if $(G,H,K)$ tightly contains some $(G_0,H_0,K_0)$ in the following table. In this case, $H_0\cap K_0$ is described in the table.
{\small
\[
\begin{array}{lllll}
\hline
G_0&H_0 &K_0& H_0\cap K_0 & \textup{Remarks}\\
\hline
\Sp_{2m}(2) & \SL_m(2){:}2 & \Omega_{2m}^-(2) & \SL_{m-1}(2) & m\textup{ odd}\\
\Sp_{2m}(2) & \SL_m(2) & \mathrm{O}_{2m}^-(2) & \SL_{m-1}(2) & \\
 & \SU_m(2) & & \SU_{m-1}(2) & m\textup{ even}\\
 &\Sp_m(2) & & \Sp_{m-2}(2) & H_0<T\\
 & \Omega_m^+(4){:}2 & & \Sp_{m-2}(4) & \textup{as in \ref{LemOmegaPlus12}}\\
 & \SL_{m/2}(4){:}2 & & \SL_{m/2-1}(4) & H_0<T\textup{ or }T{:}\langle\gamma\rangle\\
 & \SU_{m/2}(4){:}4 & & \SU_{m/2-1}(4){:}2 & m/2\textup{ even},\ H_0<\Omega_m^+(4){:}2\\
 & \GaSp_{m/2}(4) & & \Sp_{m/2-2}(4) & H_0<T\textup{ or }T{:}\langle\gamma\rangle \\
 & \Sp_{m/2}(4).4 & & \Sp_{m/2-2}(4){:}2 & \textup{as in \ref{LemOmegaPlus35}(b)}\\
\Sp_{12}(2) &\G_2(2) & \mathrm{O}_{12}^-(2) & \SL_2(2) & H_0<T \\
 & 3^{\boldsymbol{\cdot}}\PSU_4(3),\ 3^{\boldsymbol{\cdot}}\M_{22} & & 3^5{:}\A_5,\ \PSL_2(11) & H_0<\SU_6(2)\\
\Sp_{16}(2) & \Omega_9(2) & \mathrm{O}_{16}^-(2) & \Omega_7(2) & \\
 & \GaSp_4(4) & & \Sp_2(4) & H_0<\Omega_9(2)\cong\Sp_8(2)\\
 & \Omega_8^-(2){:}2,\ \GaSp_6(4) &  & \G_2(2),\ \G_2(4) & \textup{as in \ref{LemOmegaPlus35}(c)}\\
 & \Omega_6^\pm(4){:}2,\ \GaSp_4(4) & & \SL_3^\pm(4),\ \SL_2(4) & \textup{as in \ref{LemOmegaPlus35}(d)}\\
\Sp_{24}(2) & \GaG_2(4) & \mathrm{O}_{24}^-(2) & \SL_2(4) & H_0<T\textup{ or }T{:}\langle\gamma\rangle \\
 & \G_2(4).4 & & \SL_2(4){:}2 & \textup{as in \ref{LemOmegaPlus35}(e)}\\
 & 3^{\boldsymbol{\cdot}}\Suz & & 3^5{:}\PSL_2(11) & H_0<\SU_{12}(2)\\
 & \Co_1 & & \Co_3 & \\
 & 3^{\boldsymbol{\cdot}}\Suz,\ \G_2(4){:}2 & & 3^5{:}\PSL_2(11),\ \A_5 & H_0<\Co_1\\
\Sp_{32}(2) & \GaSp_8(4) & \mathrm{O}_{32}^-(2) & \Sp_6(4) & \textup{as in \ref{LemOmegaPlus35}(f)}\\
\GaSp_{2m}(4) & \SL_m(4){:}2 & \GaO_{2m}^-(4) & \SL_{m-1}(4) & \textup{as in \ref{LemSymplectic30}(a)(c)}\\
 & \SU_m(4){:}2 & & \SU_{m-1}(4) & m\textup{ even, as in \ref{LemSymplectic30}(a)}\\
  & \Sp_m(4){:}2 & & \Sp_{m-2}(4) & \textup{as in \ref{LemSymplectic30}(a)}\\
\GaSp_8(4) & \Omega_8^-(2){:}2 & \GaO_8^-(4) & \G_2(2) & \textup{as in \ref{LemSymplectic30}(b)}\\
 & \Sp_6(4){:}2,\ \Omega_6^+(4){:}2 & & \G_2(4),\ \SL_3(4) & \textup{as in \ref{LemSymplectic30}(b)}\\
 & \Omega_6^-(4).4,\ \Sp_4(4){:}2 & & \SU_3(4){:}2,\ \SL_2(4) & \textup{as in \ref{LemSymplectic30}(b)}\\
\GaSp_{12}(4) & \G_2(4){:}2 & \GaO_{12}^-(4) & \SL_2(4) & \textup{as in \ref{LemSymplectic30}(a)}\\
\GaSp_{16}(4) & \Omega_9(4){:}2 & \GaO_{16}^-(4) & \Omega_7(4) & \textup{as in \ref{LemSymplectic30}(a)}\\
\hline
\end{array}
\]
}
\end{proposition}

\begin{proof}
By Lemma~\ref{LemSymplectic52}, we may assume that $(m,q)\notin\{(2,4),(3,2),(3,4),(4,2)\}$.
For $q=2$, the conclusion follows from Example~\ref{LemSymplectic14} and Theorem~\ref{ThmOmegaPlus}.
For $q=4$ we see from~\cite[Theorem~A]{LPS1990} that $G=\GaSp_{2m}(4)$. Then the conclusion follows from Example~\ref{LemSymplectic16}  and Theorem~\ref{ThmOmegaPlus}.
\end{proof}

\subsection{Actions on $\mathcal{N}_1^+[\Omega_{2m+1}(q)]$}
\ \vspace{1mm}

For even $q$, the action of $\Omega_{2m+1}(q)$ on $\mathcal{N}_1^+[\Omega_{2m+1}(q)]$ is equivalent to the action of $\Sp_{2m}(q)$ on $\mathrm{O}_{2m}^+(q)\backslash\Sp_{2m}(q)$. In this subsection, we determine the factorizations $G=HK$ with
\[
K^{(\infty)}=\Omega_{2m}^+(q).
\]
Here $m\geqslant3$ as $K$ has a unique nonsolvable composition factor. Then Proposition~\ref{lem:K<B-Sp} shows that $A^{(\infty)}$ is one of
\[
\mbox{$\Sp_{2d}(q^e)$ with $m=de$ and $e$ prime,\ \ $\Omega_{2m}^-(q)$ with $q\in\{2,4\}$,\ \ $\G_2(q)'$ with $m=3$.}
\]
By Lemma~\ref{LemSymplectic52} and Propositions~\ref{LemSymplectic12} and~\ref{prop:Sp(2)=O^-O^+}, it suffices to deal with the first case. To start with, we construct two examples.

For $\varepsilon\in\{+,-\}$, Lemma~\ref{LemSymplectic07} gives the factorization $\Sp_6(q^b)=\G_2(q^b)\Omega_6^\varepsilon(q^b)$, where the intersection of factors $\G_2(q^b)\cap\Omega_6^\varepsilon(q^b)=\SL_3^\varepsilon(q^b)$. Hence
\[
\Sp_6(q^b)=\G_2(q^b)(\Omega_6^\varepsilon(q^b){:}(2,b))\ \text{ with }\ \G_2(q^b)\cap(\Omega_6^\varepsilon(q^b){:}(2,b))=\SL_3^\varepsilon(q^b).(2,b).
\]
This together with Example~\ref{LemSymplectic04} leads to the next example.

\begin{example}\label{LemSymplectic49}
Let $G=\Sp_{2m}(q)$ with $q$ even, let $H=\G_2(q^b)<\Sp(V_{(b)})=\Sp_6(q^b)$ with $m=3b$, and let $K=\Omega_{2m}^\varepsilon(q){:}(2,b)$ with $\varepsilon\in\{+,-\}$. Then $G=HK$ with $H\cap K=\SL_3^\varepsilon(q^b).(2,b)$.
%Moreover, when $(m,q)=(3,2)$, computation in \magma~\cite{BCP1997} shows $\Sp_6(2)=\G_2(2)'\,\Omega_6^-(2)$ with $\G_2(2)'\cap\Omega_6^-(2)=3^2{:}\Sy_3.2$.
\end{example}

If $m=2b$ with $b$ odd, then Example~\ref{LemSymplectic04} implies that $\Sp_{2m}(2^f)=\Sp_4(2^{fb})\mathrm{O}_{2m}^+(2^f)$ with $\Sp_4(2^{fb})\cap\mathrm{O}_{2m}^+(2^f)=\mathrm{O}_4^+(2^{fb})$. If in addition $f$ is odd, then it is shown in~\cite[5.1.7(b)]{LPS1990} that $\Sp_4(2^{fb})=\Sz(2^{fb})\mathrm{O}_4^+(2^{fb})$, where $\Sz(q^b)\cap\mathrm{O}_4^+(q^b)=\Sz(q^b)\cap\Omega_4^+(q^b)=\D_{2(q^b-1)}$. Thus we obtain a new example.

\begin{example}\label{LemSymplectic21}
Let $G=\Sp_{2m}(q)$ with $q=2^f$ and $m=2b$ such that $fb>1$ is odd, let $H=\Sz(q^b)<\Sp(V_{(b)})=\Sp_4(q^b)$, and let $K=\mathrm{O}_{2m}^+(q)<G$. Then $G=HK$ with $H\cap K=H\cap K^{(\infty)}=\D_{2(q^b-1)}$.
\end{example}

Now we present a classification of factorizations $G=HK$ such that $K^{(\infty)}=\Omega_{2m}^+(q)$ and $H^{(\infty)}\leqslant A^{(\infty)}=\Sp_{2d}(q^e)$ with $m=de$ and $e$ prime. By Proposition~\ref{prop:Sp(2)=O^-O^+} we may also assume that $H$ is not contained in any maximal subgroup of $G$ with solvable residual $\Omega_{2m}^-(q)$.

\begin{proposition}\label{prop:SpaO+<Sp}
Let $(m,q)\neq(3,2)$, let $A^{(\infty)}=\Sp_{2d}(q^e)$ with $m=de$ and $e$ prime, and let $K^{(\infty)}=\Omega_{2m}^+(q)$. Suppose that $(G,H,K)$ is not as described in Proposition~$\ref{prop:Sp(2)=O^-O^+}$. Then $G=HK$ if and only if $(G,H,K)$ tightly contains some $(G_0,H_0,K_0)$ in the following table. In this case, $H_0\cap K_0$ is described in the table.
\[
\begin{array}{lllll}
\hline
G_0&H_0 &K_0& H_0\cap K_0 & \textup{Remarks}\\ \hline
\Sp_{2ab}(q) & \Sp_{2a}(q^b) & \Omega_{2ab}^+(q){:}(2,b) & \Omega_{2a}^+(q^b){:}(2,b) & \textup{as in \ref{LemSymplectic04}}\\
\Sp_{6b}(q) & \G_2(q^b) & \Omega_{6b}^+(q){:}(2,b) & \SL_3(q^b).(2,b) & \textup{as in \ref{LemSymplectic49}}\\
\Sp_{4b}(q) & \Sz(q^b) & \mathrm{O}_{4b}^+(q) & \D_{2(q^b-1)} & \textup{as in \ref{LemSymplectic21}}\\
\hline
\end{array}
\]
\end{proposition}

\begin{proof}
It is shown in Examples~\ref{LemSymplectic04},~\ref{LemSymplectic49} and~\ref{LemSymplectic21} that each triple $(G_0,H_0,K_0)$ in the table gives rise to the factorization $G_0=H_0K_0$. In the following, suppose that $G=HK$.
Since $(G,H,K)$ is not as described in Proposition~\ref{prop:Sp(2)=O^-O^+}, the group $H$ is not contained in any maximal subgroup of $G$ with solvable residual $\Omega_{2m}^-(q)$. Repeatedly applying Proposition~\ref{lem:K<B-Sp} we derive that $H$ is contained in some field-extension subgroup $M=\Sp_{2a}(q^b).(b\ell)$ of $G$ with $m=ab$ and $\ell=|G/L|$ dividing $f$, and one of the following holds:
\begin{enumerate}[{\rm (i)}]
\item $H^{(\infty)}=\Sp_{2a}(q^b)$;
\item $a=3$, and $H^{(\infty)}=\G_2(q^b)$;
\item $a=2$, $bf$ is odd, and $H^{(\infty)}=\Sz(q^b)$.
\end{enumerate}
If~(i) or~(ii) appears, then Lemma~\ref{LemSymplectic02} asserts that $K\cap L\geqslant\Omega_{2m}^+(q){:}(2,b)$. If~(iii) appears, then analysis of the $2$-parts gives $K\cap L\geqslant\mathrm{O}_{2m}^+(q)$. Hence $(G,H,K)$ tightly contains some triple $(G_0,H_0,K_0)$ in the table of the proposition.
\end{proof}

\subsection{Actions on $\mathcal{N}_1^-[\Omega_{2m+1}(q)]$}
\ \vspace{1mm}

For even $q$, the action of $\Omega_{2m+1}(q)$ on $\mathcal{N}_1^-[\Omega_{2m+1}(q)]$ is equivalent to the action of $\Sp_{2m}(q)$ on $\mathrm{O}_{2m}^-(q)\backslash\Sp_{2m}(q)$. In this subsection, we determine the factorizations $G=HK$ with
\[
K^{(\infty)}=\Omega_{2m}^-(q).
\]
According to Proposition~\ref{lem:K<B-Sp}, we have the following cases for $A^{(\infty)}$:
\begin{align*}
&\mbox{$\Sp_m(q)\times\Sp_m(q)$,\ \ $\Sp_{2m}(q^{1/2})$ with $q\in\{4,16\}$,\ \ $\Sp_{2d}(q^e)$ with $m=de$ and $e$ prime,}\\
&\mbox{$q^{m(m+1)/2}{:}\SL_m(q)$,\ \ $\Omega_{2m}^+(q)$ with $q\in\{2,4\}$,\ \ $\G_2(q)$ with $m=3$.}
\end{align*}
By Lemma~\ref{LemSymplectic52} and Propositions~\ref{LemSymplectic12} and~\ref{prop:Sp(2)=O^-O^+}, we only need to deal with the first four cases.

\subsubsection{The case $A^{(\infty)}=\Sp_m(q)\times\Sp_m(q)$}

In this case, $m$ is even and $A\cap L=(\Sp(V_1)\times\Sp(V_2)){:}\langle\sigma\rangle$, where $V_1$, $V_2$ and $\sigma$ are defined in Subsection~\ref{SecSymplectic01}.
%To be as general as possible, we include the previously excluded pairs $(2,4)$ and $(4,2)$ for $(m,q)$ in this case.

\begin{lemma}\label{LemSymplectic18}
Let $G=\GaSp(V)=\GaSp_{2m}(q)$ with $m$ even and $q\in\{2,4\}$, let $\rho\in\{1,\phi\}$, let
\[
N=\{(x,x^\rho)\mid x\in\Sp_m(q)\}<\Sp_m(q)\times\Sp_m(q)=\Sp(V_1)\times\Sp(V_2),
\]
let $H=N{:}(\langle\sigma\rangle\times\langle\phi\rangle)=\Sp_m(q){:}(2\times f)$, and let $K=\GaO(V,Q^-)=\GaO_{2m}^-(q)$. Then $G=HK$ with $H\cap K=N\cap K=\{(x,x^\rho)\mid x\in\Omega_{m-1}(q)\times2\}=\Omega_{m-1}(q)\times2$.
\end{lemma}

\begin{proof}
The case $(m,q)=(2,4)$ is verified in \magma~\cite{BCP1997}. Now assume $m\geqslant4$.
Since the restriction of $Q^-$ on $V_1$ and $V_2$ is a nondegenerate quadratic form of type $+$ and $-$, respectively, we obtain by using $\rho^2=1$ that
\[
H\cap K=\{(x,x^\rho)\mid x\in\GaO_m^+(q)\cap(\GaO_m^-(q))^\rho\}=\{(x,x^\rho)\mid x\in\GaO_m^+(q)\cap\GaO_m^-(q)\}.
\]
By~\cite[3.2.4(e)]{LPS1990}, $\GaO_m^+(q)\cap\GaO_m^-(q)=\mathrm{O}_m^+(q)\cap\mathrm{O}_m^-(q)=\Omega_{m-1}(q)\times2$. So the lemma follows.
\end{proof}

By Lemmas~\ref{LemSymplectic18} and~\ref{LemSymplectic19} we obtain the subsequent example.

\begin{example}\label{LemSymplectic15}
Let $G=\Sp_{2m}(2)$ with $m/2$ even, let $N=\{(x,x)\mid x\in\GaSp_{m/2}(4)\}<\Sp_m(2)\times\Sp_m(2)=\Sp(V_1)\times\Sp(V_2)$, let $H=N\times\langle\sigma\rangle=\GaSp_{m/2}(4)\times2=\Sp_{m/2}(4){:}2^2$, and let $K=\mathrm{O}_{2m}^-(2)$. Then $G=HK$ with $H\cap K=N\cap K=\Omega_{m/2-1}(4)\times2$.
\end{example}

For even $q$, according to~\cite{LMM2003}, $\Aut(\Sp_4(q))$ is a split extension of $\Sp_4(q)$ by $\Out(\Sp_4(q))$, and note that $\Out(\Sp_4(q))$ is a cyclic group of order $2f$ and $\GaSp_4(q)=\Sp_4(q){:}f$. Similar to Lemma~\ref{LemSymplectic18}, we have the following factorization.

\begin{lemma}\label{LemSymplectic41}
Let $G=\Sp_8(q)$ with $q$ even and $f$ odd, let $\rho$ be an involution in $\Aut(\Sp_4(q))\Setminus\GaSp_4(q)$, let $N=\{(x,x^\rho)\mid x\in\Sp_4(q)\}<\Sp_4(q)\times\Sp_4(q)=\Sp(V_1)\times\Sp(V_2)$, let $H=N{:}\langle\sigma\rangle=\Sp_4(q){:}2$, and let $K=\Omega(V,Q^-)=\Omega_8^-(q)$. Then $G=HK$ with $H\cap K=\mathrm{O}_2^+(q^2)$.
\end{lemma}

The next example follows from Example~\ref{LemSymplectic04} and Lemma~\ref{LemSymplectic41}.

\begin{example}\label{LemSymplectic42}
Let $G=\Sp_{2m}(q)$ with $q$ even, $f$ odd and $m=4b$ for some odd $b$, let $H<\Sp(V_{(b)})=\Sp_8(q^b)$ be the group $H$ defined in Lemma~$\ref{LemSymplectic41}$ with $q$ replaced by $q^b$, and let $K=\Omega_{2m}^-(q)$. Then $G=HK$ with $H=\Sp_4(q^b){:}2$ and $H\cap K=\mathrm{O}_2^+(q^{2b})$.
\end{example}

For $q\in\{2,4\}$ we have $\Sp_6(q)=\G_2(q)'(\Omega_5(q)\times2)$ with $\G_2(q)'\cap(\Omega_5(q)\times2)=\SL_2(q)\times(q/2)$, which can be verified by computation in \magma~\cite{BCP1997} for $q=2$ and be read off from Proposition~\ref{LemSymplectic12} for $q=4$. This combined with Lemma~\ref{LemSymplectic18} yields the ensuing example.

\begin{example}\label{LemSymplectic51}
Let $G=\GaSp_{12}(q)$ with $q\in\{2,4\}$, let $N=\{(x,x^\rho)\mid x\in\G_2(q)'\}<\Sp_6(q)\times\Sp_6(q)=\Sp(V_1)\times\Sp(V_2)$ with $\rho\in\{1,\phi\}$, let $H=N{:}(\langle\sigma\rangle\times\langle\phi\rangle)=\G_2(q)'{:}(2\times f)$, and let $K=\GaO_{12}^-(q)$. Then $G=HK$ with $H\cap K=\SL_2(q)\times(q/2)$.
\end{example}

Below is the last example for the case $A^{(\infty)}=\Sp_m(q)\times\Sp_m(q)$, which follows from Example~\ref{LemSymplectic15} and Lemma~\ref{LemSymplectic07}.

\begin{example}\label{LemSymplectic47}
Let $G=\Sp_{24}(2)$, let $N=\{(x,x)\mid x\in\GaG_2(4)\}<\Sp_{12}(2)\times\Sp_{12}(2)=\Sp(V_1)\times\Sp(V_2)$, let $H=N\times\langle\sigma\rangle=\GaG_2(4)\times2=\G_2(4){:}2^2$, and let $K=\mathrm{O}_{24}^-(2)$. Then $G=HK$ with $H\cap K=N\cap K=\SL_2(4)\times2$.
\end{example}

We now classify the factorizations with $A^{(\infty)}=\Sp_m(q)\times\Sp_m(q)$.

\begin{proposition}\label{C2-subgroup}
Let $(m,q)\neq(4,2)$, let $A\cap L=\Sp_m(q)\wr\Sy_2$, and let $K^{(\infty)}=\Omega_{2m}^-(q)$. Then $G=HK$ if and only if $(G,H,K)$ tightly contains some $(G_0,H_0,K_0)$ in the following table. In this case, $H_0\cap K_0$ is described in the table.
\[
\begin{array}{lllll}
\hline
G_0 & H_0 & K_0 & H_0\cap K_0 & \textup{Remarks}\\
\hline
\Sp_{8b}(2^f) & \Sp_4(2^{fb}){:}2 & \Omega_{8b}^-(2^f) & \D_{2(4^{fb}-1)} & \textup{$fb$ odd, as in \ref{LemSymplectic42}}\\
\Sp_{2m}(2) & \Sp_m(2)\times2 & \mathrm{O}_{2m}^-(2) & \Sp_{m-2}(2)\times2 & \textup{as in \ref{LemSymplectic18}}\\
 & \Sp_{m/2}(4){:}2^2 &  & \Sp_{m/2-2}(4)\times2 & \textup{as in \ref{LemSymplectic15}}\\
\GaSp_{2m}(4) & \Sp_m(4){:}2^2 & \GaO_{2m}^-(4) & \Sp_{m-2}(4)\times2 & \textup{as in \ref{LemSymplectic18}}\\
\Sp_{12}(2) & \SU_3(3)\times2 & \mathrm{O}_{12}^-(2) & \SL_2(2) & \textup{as in \ref{LemSymplectic51}}\\
\GaSp_{12}(4) & \G_2(4){:}2^2 & \GaO_{12}^-(4) & \SL_2(4)\times2 & \textup{as in \ref{LemSymplectic51}}\\
\Sp_{24}(2) & \G_2(4){:}2^2 & \mathrm{O}_{24}^-(2) & \SL_2(4)\times2 & \textup{as in \ref{LemSymplectic47}}\\
\hline
\end{array}
\]
\end{proposition}

\begin{proof}
It suffices to prove the ``only if'' part. Thus let $G=HK$.
Write $\calO=G/L$ and $A\cap L=(\Sp(V_1)\times\Sp(V_2)){:}\langle\sigma\rangle=\Sp_m(q)\wr\Sy_2$. Let $N=(\Sp(V_1)\times\Sp(V_2)).\calO<A$
such that $N\cap L=\Sp(V_1)\times\Sp(V_2)$, and let $\pi_1$ and $\pi_2$ be the projections of $N$ onto $\Sp(V_1).\calO$ and $\Sp(V_2).\calO$, respectively. Then $A\cap B\cap L=\mathrm{O}_m^+(q)\times\mathrm{O}_m^-(q)<N\cap L$, and we have factorizations
\begin{equation}\label{EqnSymplectic10}
N^{\pi_1}=(H\cap N)^{\pi_1}(A\cap B)^{\pi_1}\ \text{ and }\ N^{\pi_2}=(H\cap N)^{\pi_2}(A\cap B)^{\pi_2}
\end{equation}
of $N^{\pi_1}\cong N^{\pi_2}=\Sp_m(q).\calO$.
If one of $(H\cap N)^{\pi_1}$ or $(H\cap N)^{\pi_2}$ contains $\Sp_m(q)$, then the uniqueness of its nonsolvable composition factor makes
\[
(H\cap N)^{(\infty)}=\{(x,x^\rho)\mid x\in\Sp_m(q)\}<\Sp_m(q)\times\Sp_m(q)=\Sp(V_1)\times\Sp(V_2).
\]
for some automorphism $\rho$ of $\Sp_m(q)$, and $(G,H,K)$ tightly contains some $(G_0,H_0,K_0)$ in the first, second or fourth row of the table. If neither $(H\cap N)^{\pi_1}$ nor $(H\cap N)^{\pi_2}$ contains $\Sp_m(q)$, then considering the two factorizations in~\eqref{EqnSymplectic10} simultaneously shows that $(G,H,K)$ tightly contains some $(G_0,H_0,K_0)$ in the first, third or the last three rows of the table.
\end{proof}

\subsubsection{The case $A^{(\infty)}=\Sp_{2m}(q^{1/2})$ with $q\in\{4,16\}$}

We first recall the following example of maximal factorization given in~\cite[3.2.4(d)]{LPS1990}.

\begin{example}\label{LemSymplectic20}
Let $G=\GaSp_{2m}(q)$ with $q\in\{4,16\}$, let $H$ be the maximal subgroup of $G$ such that
\[
H=\Sp_{2m}(q^{1/2}){:}f=
\begin{cases}
\Sp_{2m}(2)\times2&\textup{if }q=4\\
\Sp_{2m}(4){:}4&\textup{if }q=16,
\end{cases}
\]
and let $K=\GaO_{2m}^-(q)$. Then $G=HK$ with $H\cap K=H\cap K\cap L=\Omega_{2m-1}(q^{1/2})\times2$.
\end{example}

The next example follows from Example~\ref{LemSymplectic20} and Lemma~\ref{LemSymplectic07}.

\begin{example}\label{LemSymplectic32}
Let $G=\GaSp_6(16)$, let $H=\G_2(4){:}4<\Sp_6(4){:}4<G$, and let $K=\GaO_6^-(16)<G$. Then $G=HK$ with $H\cap K=\SL_2(4)\times2$.
\end{example}

It turns out that Examples~\ref{LemSymplectic20} and~\ref{LemSymplectic32} are the only new examples in the case $A^{(\infty)}=\Sp_{2m}(q^{1/2})$, as shown below.

\begin{proposition}\label{prop:Sp(q)=Sp(q^{1/2})Q-}
Let $A^{(\infty)}=\Sp_{2m}(q^{1/2})$ with $q\in\{4,16\}$, and let $K^{(\infty)}=\Omega_{2m}^-(q)$.
Suppose that $(G,H,K)$ is not as described in Proposition~$\ref{C2-subgroup}$.
Then $G=HK$ if and only if $(G,H,K)$ is as described in Example~$\ref{LemSymplectic20}$ or~$\ref{LemSymplectic32}$.
\end{proposition}

\begin{proof}
It suffices to prove the ``only if'' part. Suppose that $G=HK$.
By~\cite[Theorem~A]{LPS1990}, $A=\Sp_{2m}(q^{1/2}){:}f$ and $K=B=\GaO_{2m}^-(q)$. If $H^{(\infty)}=A^{(\infty)}$, then the triple $(G,H,K)$ is described in Example~\ref{LemSymplectic20}. Otherwise, applying Proposition~\ref{lem:K<B-Sp} to the induced factorization of $A/\Rad(A)=\GaSp_{2m}(q^{1/2})$ leads to either Example~\ref{LemSymplectic32} or that $H$ is contained in a subgroup $M$ of $A$ such that $M^{(\infty)}=\Sp_m(q)$ is a field-extension subgroup of $A^{(\infty)}$. For the latter, $M$ is contained in a $\calC_2$-subgroup of type $\Sp_m(q)\wr\Sy_2$ in $G$, and so is $H$. This contradicts the assumption that $(G,H,K)$ is not as described in Proposition~\ref{C2-subgroup}.
\end{proof}

\subsubsection{The case $A^{(\infty)}=\Sp_{2d}(q^e)$ with $m=de$ and $e$ prime}

%To be as general as possible, we include the pairs $(2,4)$, $(3,2)$, $(3,4)$ and $(4,2)$ for $(m,q)$ in the following proposition.
Recall that, for a positive integer $b$, we denote the $2$-part of $b$ by $b_2$. Lemma~\ref{LemSymplectic03} gives $\Sp(V_{(b)})\cap\Omega(V,Q^\varepsilon)=\mathrm{O}(V_{(b)},Q_{(b)}^\varepsilon)$ for even $b$. Hence we have the following result.

\begin{lemma}\label{LemSymplectic05}
Let $G=\Sp_{2m}(q)$ with $q$ even and $m=ab$ such that $b$ is even, let $H=\Sp_{2a}(q^b){:}b_2\leqslant\GaSp(V_{(b)})$, and let $K=\Omega_{2m}^-(q)<G$. Then $G=HK$ with $H\cap K=\mathrm{O}_{2a}^-(q^b).(b_2/2)$.
\end{lemma}

Lemma~\ref{LemSymplectic05} together with Lemma~\ref{LemSymplectic07} yields the following example.

\begin{example}\label{LemSymplectic50}
Let $G=\Sp_{6b}(q)$ with $q$ and $b$ even, let $H=\G_2(q^b){:}b_2<\Sp_6(q^b){:}b_2\leqslant\GaSp(V_{(b)})$, and let $K=\Omega_{6b}^-(q)$. Then $G=HK$ with $H\cap K=\SU_3(q^b).b_2$.
\end{example}
%
%The next example follows from Lemma~\ref{LemSymplectic41} and Example~\ref{LemSymplectic04}.
%
%\begin{example}\label{LemSymplectic42}
%Let $G=\Sp_{2m}(q)$ with $m=4b$ for some odd $b$ and $q=2^f$ for some odd $f$, let $M=\Sp(V_{(b)})=\Sp_8(q^b)$, let $H<M$ be the group $H$ defined in Lemma~$\ref{LemSymplectic41}$ with $q$ replaced by $q^b$, and let $K=\Omega_{2m}^-(q)<G$. Then $G=HK$ with $H=\Sp_4(q^b){:}2$ and $H\cap K=\mathrm{O}_2^+(q^{2b})$.
%\end{example}

In what follows we show that all the minimal (with respect to tight containment) factorizations in the case $A^{(\infty)}=\Sp_{2d}(q^e)$ have been given. To achieve this result, we need the ensuing lemma, which also obtained based on Lemma~\ref{LemSymplectic03}.

\begin{lemma}\label{LemSymplectic48}
Let $G=HK$ such that $H^{(\infty)}=\Sp_{2a}(q^b)$ is a field-extension subgroup of $L$ with $q$ even and $m=ab$, and $K^{(\infty)}=\Omega_{2m}^-(q)$. Then either $H\cap L\geqslant\Sp_{2a}(q^b).b_2$ with $b$ even, or $K\cap L\geqslant\Omega_{2m}^-(q){:}(2,b)$.
\end{lemma}

We are now ready to give the proposition determining the factorizations $G=HK$ in the case $A^{(\infty)}=\Sp_{2d}(q^e)$. Note that the factorizations with $H\leqslant\Pa_m[G]$ will be classified in Proposition~\ref{PropSymplecticPm}.

\begin{proposition}\label{prop:SpaO-<Sp}
Let $A^{(\infty)}=\Sp_{2d}(q^e)$ with $m=de$ and $e$ prime, and let $K^{(\infty)}=\Omega_{2m}^-(q)$. Suppose that $H\nleqslant\Pa_m[G]$ and that $(G,H,K)$ is not as described in Proposition~$\ref{prop:Sp(2)=O^+O^-}$,~$\ref{C2-subgroup}$ or~$\ref{prop:Sp(q)=Sp(q^{1/2})Q-}$. Then $G=HK$ if and only if $(G,H,K)$ tightly contains some $(G_0,H_0,K_0)$ in the following table. In this case, $H_0\cap K_0$ is described in the table.
\[
\begin{array}{lllll}
\hline
G_0 & H_0 & K_0 & H_0\cap K_0 & \textup{Remarks}\\ \hline
\Sp_{2ab}(q) & \Sp_{2a}(q^b){:}b_2 & \Omega_{2ab}^-(q) & \mathrm{O}_{2a}^-(q^b).(b_2/2) & \textup{as in \ref{LemSymplectic05}}\\
\textup{($b$ even)} & \G_2(q^b){:}b_2 &  & \SU_3(q^b).b_2 & \textup{$a=3$, as in \ref{LemSymplectic50}}\\
%\Sp_{8b}(2^f) & \Sp_4(2^{fb}){:}2 & \Omega_{8b}^-(2^f) & \mathrm{O}_2^+(q^{2b}) & \textup{$fb$ odd, $H_0$ as in Example~\ref{LemSymplectic42}}\\
\Sp_{2ab}(q) & \Sp_{2a}(q^b) & \Omega_{2ab}^-(q){:}(2,b) & \Omega_{2a}^-(q^b){:}(2,b) & \textup{as in \ref{LemSymplectic04}}\\
\Sp_{6b}(q) & \G_2(q^b) & \Omega_{6b}^-(q){:}(2,b) & \SU_3(q^b).(2,b) & \textup{as in \ref{LemSymplectic49}}\\
\Sp_{12}(2) & \J_2{:}2 & \Omega_{12}^-(2) & 5^2{:}\Sy_3.2 & H_0<\GaG_2(4)<\GaSp_6(4)\\
\Sp_{12}(2) & \J_2 & \mathrm{O}_{12}^-(2) & 5^2{:}\Sy_3.2 & H_0<\G_2(4)<\Sp_6(4)\\
\hline
\end{array}
\]
\end{proposition}

\begin{proof}
This follows by the same field-extension reduction as Proposition~\ref{prop:SpaO+<Sp}: choosing a field-extension subgroup
containing $H$ with maximal defining field and applying the $\max^-$ factorization classification.
\end{proof}

\subsubsection{The case $A^{(\infty)}=q^{m(m+1)/2}{:}\SL_m(q)$}

Here $A=G_U=\Pa_m[G]$. Recall $R=\bfO_2(A)=q^{m(m+1)/2}$ and $T=\SL_m(q)$ in Subsection~\ref{SecSymplectic01}.
The following result is parallel to Lemma~\ref{LemOmega01}.

\begin{lemma}\label{LemSymplectic27}
Let $G=\Sp(V)=\Sp_{2m}(q)$ with $q$ even, let $M=A^{(\infty)}=R{:}T$, and let $K=\Omega(V,Q_-)=\Omega_{2m}^-(q)$. Then the following statements hold:
\begin{enumerate}[{\rm (a)}]
\item The kernel of $M\cap K$ acting on $U$ is $R\cap K=q^{m(m-1)/2}$.
\item The induced group by the action of $M\cap K$ on $U$ is $\SL(U)_{U_1,e_m+U_1}=q^{m-1}{:}\SL_{m-1}(q)$.
% \item $T\cap K=\SL_{m-1}(q)$;
\item If $H=R{:}S$ with $S\leqslant T$, then $H\cap K=(R\cap K).S_{U_1,e_m+U_1}$.
\item If $S$ is a transitive subgroup of $T$ on $U\Setminus\{0\}$, then $G=(R{:}S)K$.
\end{enumerate}
\end{lemma}

If $m=ab$, then let $U_{(b)}=\bbF_{q^b}^a$ be a vector space of dimension $a$ over $\bbF_{q^b}$ with the same underlying set as $U$, and define $U_{(b)}^\sharp(1),\dots,U_{(b)}^\sharp(\lfloor b/2\rfloor)$ as in~\eqref{EqnOrthonal1}; that is,
\[
U_{(b)}^\sharp(i)=
\left\{\begin{array}{ll}
\bigoplus\limits_{j=0}^{b-1}\left(U_{(b)}\otimes U_{(b)}^{(q^i)}\right)^{(q^j)}&\text{if }1\leqslant i\leqslant\lfloor(b-1)/2\rfloor\vspace{2mm}\\
\bigoplus\limits_{j=0}^{b/2-1}\left(U_{(b)}\otimes U_{(b)}^{(q^i)}\right)^{(q^j)}&\text{if }\lfloor(b-1)/2\rfloor<i\leqslant\lfloor b/2\rfloor.
\end{array}\right.
\]
Again as in Subsection~\ref{SecOmegaPlus04}, for $i\in\{1,\dots,\lfloor b/2\rfloor\}$, the $\bbF_q\GL(U_{(b)})$-module $U_{(b)}(i)$ is defined to be the realization of the $\bbF_{q^b}\GL(U_{(b)})$-modules $U_{(b)}^\sharp(i)$ over $\bbF_q$. Denote
\[
U_{(b)}(0)=\bfO_p\Big(\Pa_a[\Sp_{2a}(q^b)]\Big)=q^{ba(a+1)/2},
\]
where $\Sp_{2a}(q^b)$ is a field-extension subgroup of $L$ defined over $\bbF_{q^b}$. Then $U_{(b)}(0)$ is the symmetric square of $U_{(b)}$ with a submodule
\[
\mbox{$\Lambda_{(b)}:=\bfO_p\Big(\Pa_a[\Omega_{2a}^+(q^b)]\Big)=\bigwedge^2(U_{(b)})=q^{ba(a-1)/2}$},
\]
where $\Omega_{2a}^+(q^b)$ is a $\calC_8$-subgroup of $\Sp_{2a}(q^b)$, such that the quotient $U_{(b)}(0)/\Lambda_{(b)}$ is isomorphic to the irreducible $\bbF_{q^b}\GL(U_{(b)})$-module $U_{(b)}$.
%As an $\bbF_{q^b}\GL(U_{(b)})$-module, $U_{(b)}(0)$ is indecomposable (see~\cite[Page~13]{BL2021} for instance).

\begin{lemma}\label{LemSymplecticPm1}
Let $q$ be even and $m=ab$, and let $S$ be a subgroup of $\SL(U)$ contained in the group $\SL_a(q^b)$ defined over $\bbF_{q^b}$ such that $S=\SL_a(q^b)$, $\Sp_a(q^b)$ or $\G_2(q^b)$ (with $a=6$). Then the following statements hold.
\begin{enumerate}[{\rm (a)}]
\item The $\bbF_qS$-module $R$ is the direct sum of pairwise non-isomorphic submodules
\[
U_{(b)}(0),\,U_{(b)}(1),\,\dots,\,U_{(b)}(\lfloor b/2\rfloor).
\]
\item The irreducible $\bbF_qS$-submodules of $R$ are $U_{(b)}(1),\,\dots,\,U_{(b)}(\lfloor b/2\rfloor)$ and the irreducible $\bbF_qS$-submodules of $\Lambda_{(b)}$.
\item Every proper $\bbF_qS$-submodule of $U_{(b)}(0)$ is a submodule of $\Lambda_{(b)}$.
\end{enumerate}
\end{lemma}

\begin{proof}
From Lemma~\ref{LemOmegaPlusPm1} we see that $U_{(b)}(1),\,\dots,\,U_{(b)}(\lfloor b/2\rfloor)$ are pairwise non-isomorphic irreducible $\bbF_qS$-submodules of $\bfO_2(\Pa_m[\Omega_{2m}^+(q)])$ and hence of $R$. Since $\Lambda_{(b)}=\bbF_q^{ba(a-1)/2}$ is an $\bbF_qS$-submodule of $U_{(b)}(0)$ with quotient $\bbF_q^m$, the composition factors of $U_{(b)}(0)$ have dimensions less than those of $U_{(b)}(1),\,\dots,\,U_{(b)}(\lfloor b/2\rfloor)$.
Thus there is an $\bbF_qS$-submodule of $R$ as the direct sum of pairwise non-isomorphic submodules $U_{(b)}(0),U_{(b)}(1),\,\dots,\,U_{(b)}(\lfloor b/2\rfloor)$. Since this direct sum has dimension
\[
\frac{ba(a+1)}{2}+\sum_{i=1}^{\lfloor(b-1)/2\rfloor}ba^2+\sum_{i=\lfloor(b-1)/2\rfloor+1}^{\lfloor b/2\rfloor}\frac{ba^2}{2}
=\frac{ba(a+1)}{2}+ba^2\cdot\frac{b-1}{2}=\frac{m(m+1)}{2}
\]
as an $\bbF_qS$-module, which is the same as the dimension of $R$, it follows that
\[
R=U_{(b)}(0)\oplus U_{(b)}(1)\oplus\dots\oplus U_{(b)}(\lfloor b/2\rfloor).
\]
Hence part~(a) is confirmed.

Part~(a) reduces the verification of part~(b) to proving that the $\bbF_qS$-submodule $\Lambda_{(b)}$ has no complement $\bbF_qS$-submodule in $U_{(b)}(0)$.
This conclusion is well known for $S=\SL_a(q^b)$. In fact, suppose for a contradiction that there is such a complement $\bbF_q\SL_a(q^b)$-submodule $X$. Then $X=\bbF_q^{ab}=\bbF_q^m$. Let $E_1,\dots,E_a$ be a basis of $U_{(b)}$. Identify $U_{(b)}(0)$ as the symmetric square of $U_{(b)}$ and identify
\[
\mbox{$\Lambda_{(b)}=\bigwedge^2(U_{(b)})=\bigoplus\limits_{i=1}^{a-1}\bigoplus\limits_{j=i+1}^a\bbF_{q^b}(E_i\otimes E_j+E_j\otimes E_i)$}.
\]
Since $U_{(b)}(0)=X\oplus\bigwedge^2(U_{(b)})$, there exists $x\in X$ such that $x=E_1\otimes E_1+y$ for some $y\in\bigwedge^2(U_{(b)})$. As an $\bbF_q\SL_a(q^b)$-module, $X$ is spanned by $x^{\SL_a(q^b)}$. By taking appropriate elements in $\SL_a(q^b)$, elementary calculation leads to nonzero vectors in $X\cap\bigwedge^2(U_{(b)})$, a contradiction.
Next let $S=\Sp_a(q^b)$ or $\G_2(q^b)$ (with $a=6$). Then there exists a cyclic subgroup $C$ of order $(q^m-1)/(q^b-1)$ in $\SL_a(q^b)$ such that $\langle S,C\rangle=\SL_a(q^b)$ and $|S\cap C|$ is divisible by $\ppd(q^m-1)$. Suppose for a contradiction that the $\bbF_qS$-submodule $\Lambda_{(b)}$ has a complement $\bbF_qS$-submodule $X=\bbF_q^m$ in $U_{(b)}(0)$.
In particular, $X$ is $(S\cap C)$-invariant, which implies that $X$ is $C$-invariant, as $S\cap C$ is a subgroup of order divisible by $\ppd(q^m-1)$ in the cyclic group $C$ while $C$ is irreducible on $X\cong U_{(b)}(0)/\Lambda_{(b)}$. However, this yields that $X$ is invariant under $\langle S,C\rangle=\SL_a(q^b)$, contradicting the conclusion from the previous paragraph that $\Lambda_{(b)}$ has no complement $\bbF_q\SL_a(q^b)$-submodule in $U_{(b)}(0)$. This confirms part~(b).

Now we embark on part~(c). Since both $\Lambda_{(b)}$ and $U_{(b)}(0)/\Lambda_{(b)}$ are irreducible $\bbF_q\SL_a(q^b)$-submodules, the conclusion for $S=\SL_a(q^b)$ follows from part~(b). For the rest of the proof, assume $S=\Sp_a(q^b)$ or $\G_2(q^b)$ (with $a=6$), and suppose for a contradiction that $U_{(b)}(0)$ has an $\bbF_qS$-submodule $X$ not contained in $\Lambda_{(b)}$. Then the projection of $X$ to $U_{(b)}(0)/\Lambda_{(b)}$ is surjective, as $S$ is irreducible on $U_{(b)}(0)/\Lambda_{(b)}\cong\bbF_q^m$.

\textsf{Case}~1: $S=\Sp_a(q^b)$. Let $S_0=\SL_2(q^{m/2})=\Sp_2(q^{m/2})$ be a field-extension subgroup of $S$. Since the underlying vector space of $X\cap\Lambda_{(b)}$ is $\bbF_{q^b}$ or $\bbF_{q^b}^{a(a-1)/2-1}$ (see the remark after Lemma~\ref{LemOmegaPlusPm1}), we have $X\cong\bbF_q^{m+b}$ or $\bbF_q^{m+ba(a-1)/2-b}$. Applying the conclusions of part~(a) of the lemma and Lemma~\ref{LemOmegaPlusPm1}(b) (both with $(q,m,a,b,S)$ replaced by $(q^b,a,2,a/2,S_0)$) we obtain the decomposition of $\bbF_{q^b}S_0$-modules $U_{(b)}(0)$ and $\Lambda_{(b)}$ into pairwise non-isomorphic submodules
\begin{align*}
U_{(b)}(0)&=U_{(m/2)}(0)\oplus W(1)\oplus\dots\oplus W(\lfloor a/4\rfloor),\\
\Lambda_{(b)}&=\Lambda_{(m/2)}\oplus W(1)\oplus\dots\oplus W(\lfloor a/4\rfloor),
\end{align*}
where $W(t)=\bbF_{q^b}^{2a}$ if $1\leqslant t\leqslant\lfloor(a-2)/4\rfloor$ and $W(t)=\bbF_{q^b}^a$ if $\lfloor(a-2)/4\rfloor<t\leqslant\lfloor a/4\rfloor$. Then the condition $X\nleqslant\Lambda_{(b)}$ implies $X\cap U_{(m/2)}(0)\nleqslant\Lambda_{(m/2)}$, which contradicts the above conclusion that every proper $\bbF_q\SL_2(q^{m/2})$-submodule of $U_{(m/2)}(0)$ is a submodule of $\Lambda_{(m/2)}$.

\textsf{Case}~2: $S=\G_2(q^b)$ with $a=6$. Take a cyclic subgroup $C_0$ of order $q^{3b}+1$ in $\Sp_6(q^b)$ such that $\langle S,C_0\rangle=\Sp_6(q^b)$ and $|S\cap C_0|=(q^{3b}+1)/(q^b+1)$. Since $X$ is $(S\cap C_0)$-invariant and~\cite[Corollary~4.5]{FLWXZ} implies that each irreducible $\bbF_qC_0$-submodule of $U_{(b)}(0)$ is either the irreducible module $\bbF_q^{6b}$ or the trivial module $\bbF_q^{3b}$, we derive that $X$ is $C_0$-invariant. However, it follows that $X$ is invariant under $\langle S,C_0\rangle=\Sp_6(q^b)$, contradicting the conclusion from Case~1. This completes the proof.
\end{proof}

To construct subgroups of $\Pa_m[G]$ that are transitive on $\mathcal{N}_1^-[\Omega_{2m+1}(q)]$, we need the two subsequent lemmas.

\begin{lemma}\label{LemSymplecticPm2}
Let $G=\Sp_{2m}(q)$ with $q$ even, let $K=\Omega_{2m}^-(q)$, let $m=ab$ with $b$ even, and let $H=U_{(b)}(0){:}S=q^c{:}S^{(\infty)}{:}b_2$ with $c=ba(a+1)/2$ such that $S=S^{(\infty)}{:}\langle\theta\rangle\leqslant\SiL_a(q^b)$ is defined over $\bbF_{q^b}$, where $\theta$ is the field automorphism of $\SL_a(q^b)$ of order $b_2$.
\begin{enumerate}[{\rm (a)}]
\item If $S^{(\infty)}=\SL_a(q^b)$, then $G=HK$ with $H\cap K=[q^{c-b}].\SL_{a-1}(q^b).b_2$.
\item If $S^{(\infty)}=\Sp_a(q^b)$, then $G=HK$ with $H\cap K=[q^{c-b}].\Sp_{a-2}(q^b).b_2$.
\item If $S^{(\infty)}=\G_2(q^b)$ with $a=6$, then $G=HK$ with $H\cap K=[q^{c-b}].\SL_2(q^b).b_2$.
\end{enumerate}
\end{lemma}

\begin{proof}
Since $U_{(b)}(0)=\bfO_2(\Pa_a[\Sp_{2a}(q^b)])$, the group $H=U_{(b)}(0){:}S$ is contained in a field-extension subgroup $M=\Sp_{2a}(q^b){:}b_2$ of $G$. Then the lemma follows from Lemmas~\ref{LemSymplectic05} and~\ref{LemSymplectic27}, together with~\cite{Hering1985}.
\end{proof}

Similarly as Lemma~\ref{LemSymplecticPm2}, it follows from Example~\ref{LemSymplectic04}, Lemma~\ref{LemSymplectic27} and~\cite{Hering1985} that the next lemma holds.

\begin{lemma}\label{LemSymplecticPm3}
Let $G=\Sp_{2m}(q)$ with $q$ even and $m=ab$, let $K=\Omega_{2m}^-(q){:}(2,b)$, and let $H=U_{(b)}(0){:}S=q^c{:}S$ with $c=ba(a+1)/2$ such that $S\leqslant\SL_a(q^b)$ is defined over $\bbF_{q^b}$.
\begin{enumerate}[{\rm (a)}]
\item If $S=\SL_a(q^b)$, then $G=HK$ with $H\cap K=[q^{c-b}].\SL_{a-1}(q^b).(2,b)$.
\item If $S=\Sp_a(q^b)$, then $G=HK$ with $H\cap K=[q^{c-b}].\Sp_{a-2}(q^b).(2,b)$.
\item If $S=\G_2(q^b)'$ with $a=6$, then $G=HK$ with $H\cap K=[(q^{c-b},q^c/4)].\SL_2(q^b).(2,b)$.
\end{enumerate}
\end{lemma}

Recall from Subsection~\ref{SecOmegaPlus04} that, for a subset $I=\{i_1,\ldots,i_k\}$ of $\{1,\ldots,\lfloor b/2\rfloor\}$, we denote
\[
U_{(b)}(I)=U_{(b)}(i_1)\cdots U_{(b)}(i_k)\ \text{ and }\,\gcd(I,b)=\gcd(i_1,\dots,i_k,b).
\]
For $H\leqslant G_U=\Pa_m[G]$, the notation $H^U$, as usual, denotes the induced group of $H$ on $U$. By Lemma~\ref{LemSymplecticPm1}, if $q$ is even and $H^U=S\leqslant\GaL_a(q^b)$ defined over $\bbF_{q^b}$ such that $S^{(\infty)}=\SL_a(q^b)$, $\Sp_a(q^b)$ or $\G_2(q^b)$ (with $a=6$), then we may write
\begin{equation}\label{EqnSymplectic1}
H=\Big((U_{(b)}(0)\cap H)\times U_{(b)}(I)\Big){:}S=(E\times U_{(b)}(I)){:}S=q^c{:}S
\end{equation}
for some $\bbF_{q^b}S$-submodule $E$ of $U_{(b)}(0)=q^{ba(a+1)/2}$ and some subset $I$ of $\{1,\ldots,\lfloor b/2\rfloor\}$, where
\begin{equation}\label{EqnSymplectic2}
c=\dim_{\bbF_q}(E)+\dim_{\bbF_q}(U_{(b)}(I))=
\begin{cases}
\log_q|E|+(2|I|-1)a^2b/2&\textup{if }b/2\in I\\
\log_q|E|+|I|a^2b&\textup{if }b/2\notin I.
\end{cases}
\end{equation}
In Lemmas~\ref{LemSymplecticPm2} and~\ref{LemSymplecticPm3}, changing $H$ to its overgroups in the form of~\eqref{EqnSymplectic1} with $E=U_{(b)}(0)$, one immediately establishes the following two examples.

\begin{example}\label{LemSymplecticPm4}
Let $G=\Sp_{2m}(q)$ with $q$ even, let $K=\Omega_{2m}^-(q)$, let $m=ab$ with $b$ even, and let $H\leqslant G_U=\Pa_m[G]$ satisfying~\eqref{EqnSymplectic1} such that $E=U_{(b)}(0)=q^{ba(a+1)/2}$ and $S=S^{(\infty)}{:}\langle\theta\rangle\leqslant\SiL_a(q^b)$ is defined over $\bbF_{q^b}$, where $\theta$ is the field automorphism of $\SL_a(q^b)$ of order $b_2$.
\begin{enumerate}[{\rm (a)}]
\item If $S^{(\infty)}=\SL_a(q^b)$, then $G=HK$ with $H\cap K=[q^{c-b}].\SL_{a-1}(q^b).b_2$.
\item If $S^{(\infty)}=\Sp_a(q^b)$, then $G=HK$ with $H\cap K=[q^{c-b}].\Sp_{a-2}(q^b).b_2$.
\item If $S^{(\infty)}=\G_2(q^b)$ with $a=6$, then $G=HK$ with $H\cap K=[q^{c-b}].\SL_2(q^b).b_2$.
\end{enumerate}
\par\vspace{-1.2\baselineskip}
\qedhere
\end{example}

\begin{example}\label{LemSymplecticPm5}
Let $G=\Sp_{2m}(q)$ with $q$ even and $m=ab$, let $K=\Omega_{2m}^-(q){:}(2,b)$, and let $H\leqslant G_U=\Pa_m[G]$ satisfying~\eqref{EqnSymplectic1} such that $E=U_{(b)}(0)=q^{ba(a+1)/2}$ and $S\leqslant\SL_a(q^b)$ is defined over $\bbF_{q^b}$.
\begin{enumerate}[{\rm (a)}]
\item If $S=\SL_a(q^b)$, then $G=HK$ with $H\cap K=[q^{c-b}].\SL_{a-1}(q^b).(2,b)$.
\item If $S=\Sp_a(q^b)$, then $G=HK$ with $H\cap K=[q^{c-b}].\Sp_{a-2}(q^b).(2,b)$.
\item If $S=\G_2(q^b)'$ with $a=6$, then $G=HK$ with $H\cap K=[(q^{c-b},q^c/4)].\SL_2(q^b).(2,b)$.
\end{enumerate}
\par\vspace{-1.2\baselineskip}
\qedhere
\end{example}

For a group $H$ in the form~\eqref{EqnSymplectic1}, if $E\neq U_{(b)}(0)$, then Lemma~\ref{LemSymplecticPm1} implies that $E\leqslant\bfO_2\Big(\Pa_a[\Omega_{2a}^+(q^b)]\Big)=\bigwedge^2(U_{(b)})$ and hence $H\leqslant\Pa_m[M]$ or $\Pa_{m-1}[M]$ for some subgroup $M$ of $G$ with $M\cap L=\Omega_{2m}^+(q)$. Thus we obtain the next two examples from Example~\ref{LemSymplectic14}, Lemma~\ref{LemSymplectic16} and Proposition~\ref{PropOmegaPlusPm}.
Note from~\eqref{EqnSymplectic2} that, if $c=0$, then $I=\emptyset$ and so $\gcd(I,b)=b$.

\begin{example}\label{LemSymplecticPm7}
Let $G=\GaSp_{2m}(q)$ with $q\in\{2,4\}$ and $m=ab$, let $K=\GaO_{2m}^-(q)$, and let $H\leqslant G_U=\Pa_m[G]$ satisfying~\eqref{EqnSymplectic1} such that $E\leqslant\bigwedge^2U_{(b)}=q^{ba(a-1)/2}$, $\gcd(I,b)=2/f$ and $S=S^{(\infty)}{:}\langle\theta\rangle\leqslant\SiL_a(q^b)$ is defined over $\bbF_{q^b}$, where $\theta$ is the field automorphism of $\SL_a(q^b)$ of order $fb_2$.
\begin{enumerate}[{\rm (a)}]
\item If $S^{(\infty)}=\SL_a(q^b)$, then $G=HK$ with $H\cap K=[4q^{c-b}].\SL_{a-1}(q^b).(fb_2/2)$.
\item If $S^{(\infty)}=\Sp_a(q^b)$, then $G=HK$ with $H\cap K=[4q^{c-b}].\Sp_{a-2}(q^b).(fb_2/2)$.
\item If $S^{(\infty)}=\G_2(q^b)$ with $a=6$, then $G=HK$ with $H\cap K=[4q^{c-b}].\SL_2(q^b).(fb_2/2)$.
\end{enumerate}
\par\vspace{-1.2\baselineskip}
\qedhere
\end{example}

\begin{example}\label{LemSymplecticPm6}
Let $G=\Sp_{2m}(2)$ with $m=ab$, let $K=\mathrm{O}_{2m}^-(2)$, and let $H\leqslant G_U=\Pa_m[G]$ satisfying~\eqref{EqnSymplectic1} such that $E\leqslant\bigwedge^2U_{(b)}=2^{ba(a-1)/2}$, $\gcd(I,b)=1$ and $S\leqslant\SL_a(2^b)$ is defined over $\bbF_{2^b}$.
\begin{enumerate}[{\rm (a)}]
\item If $S=\SL_a(2^b)$, then $G=HK$ with $H\cap K=[2^{c-b+1}].\SL_{a-1}(2^b)$.
\item If $S=\Sp_a(2^b)$, then $G=HK$ with $H\cap K=[2^{c-b+1}].\Sp_{a-2}(2^b)$.
\item If $S=\G_2(2^b)'$ with $a=6$, then $G=HK$ with $H\cap K=[(2^{c-b+1},2^{c-1})].\SL_2(2^b)$.
\end{enumerate}
\par\vspace{-1.2\baselineskip}
\qedhere
\end{example}

Now we prove that the above four examples exhaust the minimal (with respect to tight containment) factorizations with $H\leqslant\Pa_m[G]$ and $K^{(\infty)}=\Omega_{2m}^-(q)$.

\begin{proposition}\label{PropSymplecticPm}
Let $(m,q)\neq(2,4)$, $(3,2)$, $(3,4)$ or $(4,2)$, let $H\leqslant\Pa_m[G]$, and let $K^{(\infty)}=\Omega_{2m}^-(q)$. Then $G=HK$ if and only if $(G,H,K)$ tightly contains some $(G_0,H_0,K_0)$ in the following table. In this case, $H_0\cap K_0$ is described in the table.
\[
\begin{array}{lllll}
\hline
G_0 & H_0 & K_0 & H_0\cap K_0 & \textup{Remarks} \\
\hline
\Sp_{2ab}(q) & q^c{:}\SL_a(q^b){:}b_2 & \Omega_{2ab}^-(q) & [q^{c-b}].\SL_{a-1}(q^b).b_2 & \textup{as in \ref{LemSymplecticPm4}} \\
\textup{($b$ even)} & q^c{:}\Sp_a(q^b){:}b_2 &  & [q^{c-b}].\Sp_{a-2}(q^b).b_2 & \textup{as in \ref{LemSymplecticPm4}} \\
 & q^c{:}\G_2(q^b){:}b_2  &  & [q^{c-b}].\SL_2(q^b).b_2 & a=6,\textup{ as in \ref{LemSymplecticPm4}} \\
\Sp_{2ab}(q) & q^c{:}\SL_a(q^b) & \Omega_{2ab}^-(q){:}(2,b) & [q^{c-b}].\SL_{a-1}(q^b).(2,b) & \textup{as in \ref{LemSymplecticPm5}} \\
 & q^c{:}\Sp_a(q^b) &  & [q^{c-b}].\Sp_{a-2}(q^b).(2,b) & \textup{as in \ref{LemSymplecticPm5}} \\
 & q^c{:}\G_2(q^b)'  &  & [(q^{c-b},q^c/4)].\SL_2(q^b).(2,b) & a=6,\textup{ as in \ref{LemSymplecticPm5}} \\
\Sp_{2ab}(2) & 2^c{:}\SL_a(2^b){:}b_2 & \mathrm{O}_{2ab}^-(2) & [2^{c-b+2}].\SL_{a-1}(2^b).(b_2/2) & \textup{as in \ref{LemSymplecticPm7}} \\
\textup{($b$ even)} & 2^c{:}\Sp_a(2^b){:}b_2 &  & [2^{c-b+2}].\Sp_{a-2}(2^b).(b_2/2) & \textup{as in \ref{LemSymplecticPm7}} \\
 & 2^c{:}\G_2(2^b){:}b_2  &  & [2^{c-b+2}].\SL_2(2^b).(b_2/2) & a=6,\textup{ as in \ref{LemSymplecticPm7}} \\
\Sp_{2ab}(2) & 2^c{:}\SL_a(2^b) & \mathrm{O}_{2ab}^-(2) & [2^{c-b+1}].\SL_{a-1}(2^b) & \textup{as in \ref{LemSymplecticPm6}} \\
 & 2^c{:}\Sp_a(2^b) &  & [2^{c-b+1}].\Sp_{a-2}(2^b) & \textup{as in \ref{LemSymplecticPm6}} \\
 & 2^c{:}\G_2(2^b)'  &  & [(2^{c-b+1},2^{c-1})].\SL_2(2^b) & a=6,\textup{ as in \ref{LemSymplecticPm6}} \\
\GaSp_{2ab}(4) & 2^{2c}{:}\SL_a(4^b){:}2b_2 & \GaO_{2ab}^-(4) & [4^{c-b+1}].\SL_{a-1}(4^b).b_2 & \textup{as in \ref{LemSymplecticPm7}} \\
\textup{($b$ even)} & 2^{2c}{:}\Sp_a(4^b){:}2b_2 &  & [4^{c-b+1}].\Sp_{a-2}(4^b).b_2 & \textup{as in \ref{LemSymplecticPm7}} \\
 & 2^{2c}{:}\G_2(4^b){:}2b_2  &  & [4^{c-b+1}].\SL_2(4^b).b_2 & a=6,\textup{ as in \ref{LemSymplecticPm7}} \\
\hline
\end{array}
\]
\end{proposition}

\begin{proof}
By Examples~\ref{LemSymplecticPm4}--\ref{LemSymplecticPm7}, it remains to prove the ``only if'' part. Suppose $G=HK$, and let $A=G_U=\Pa_m[G]$. The condition $K^{(\infty)}=\Omega_{2m}^-(q)$ implies that $q$ is even, and the group $B$ is a maximal subgroup of $G$ with $B\cap L=\mathrm{O}_{2m}^-(q)$. According to Lemma~\ref{LemSymplectic27}, $(A\cap B)^U$ stabilizes a $1$-space in $\bbF_q^m$. Thus it follows from $A=H(A\cap B)$ that $H^U$ is transitive on the set of $1$-spaces in $\bbF_q^m$, and so~\cite{Hering1985} asserts that $H^U\leqslant\GaL_a(q^b)$ is defined over $\bbF_{q^b}$ with $m=ab$ and $(H^U)^{(\infty)}=\SL_a(q^b)$, $\Sp_a(q^b)$ or $\G_2(q^b)'$ (with $a=6$). Let $S=(H^U)^{(\infty)}$. Then by Lemma~\ref{LemSymplecticPm1},
\[
H=(E\times U_{(b)}(I)){:}H^U=q^c{:}H^U
\]
for some $\bbF_{q^b}S$-submodule $E$ of $U_{(b)}(0)$ and some subset $I$ of $\{1,\ldots,\lfloor b/2\rfloor\}$. In particular, $H^{(\infty)}=(E\times U_{(b)}(I)){:}S$.

First, assume that $E$ is a proper submodule of $U_{(b)}(0)$. It follows from Lemma~\ref{LemSymplecticPm1} that $E$ is a submodule of $\bfO_2(\Pa_a[\Omega_{2a}^+(q^b)])$ and hence $H\leqslant\Pa_m[M]$ or $\Pa_{m-1}[M]$ for some subgroup $M$ of $G$ with $M\cap L=\Omega_{2m}^+(q)$. Then by Proposition~\ref{lem:K<B-Sp} we have $q\in\{2,4\}$, and applying Proposition~\ref{prop:Sp(2)=O^-O^+} to the factorization $G=KM$ we conclude that the triple $(G,K,M)$ tightly contains $(\GaSp_{2m}(q),\GaO_{2m}^-(q),\Omega_{2m}^+(q){:}\langle\phi\rangle)$ with $M\cap K=\Sp_{2m-2}(q)$. Now consider the factorization $M=H(M\cap K)$. Since $H\leqslant\Pa_m[M]$ or $\Pa_{m-1}[M]$, Proposition~\ref{PropOmegaPlusPm} implies that $(G,H,K)$ tightly contains some triple $(G_0,H_0,K_0)$ as in Example~\ref{LemSymplecticPm6} or~\ref{LemSymplecticPm7}.

Next, assume that $E=U_{(b)}(0)$. If $b$ is odd, then $(G,H,K)$ tightly contains the triple $(G_0,H_0,K_0)=(L,H^{(\infty)},K^{(\infty)})$ as in Example~\ref{LemSymplecticPm5}. Thus assume for the rest of the proof that $b$ is even.
Suppose for a contradiction that there exists a factorization $L=H_1K_1$ with $H_1=H^{(\infty)}{:}(b_2/2)=(E\times U_{(b)}(I)){:}S{:}(b_2/2)<(E\times U_{(b)}(I)){:}\SiL_a(q^b)$ and $K_1=\Omega_{2m}^-(q)$. Let
\[
H_2=H^{(\infty)}{:}b_2=(E\times U_{(b)}(I)){:}S{:}b_2\leqslant(E\times U_{(b)}(I)){:}\SiL_a(q^b),
\]
such that $H_1<H_2$, and let $H_3=E{:}S{:}(b_2/2)\leqslant H_1$ and $H_4=E{:}S{:}b_2\leqslant H_2$ such that $H_3<H_4$. From Lemma~\ref{LemSymplecticPm2} we derive that
$L=H_4K_1$ with $(H_4\cap K_1)^U=S_{U_1,e_m+U_1}.b_2$. Moreover, Lemma~\ref{LemSymplectic27} implies that $(H_2\cap K_1)^U\leqslant(R{:}S{:}b_2\cap K)^U=S_{U_1,e_m+U_1}.b_2$. Then we deduce from $H_4\cap K_1\leqslant H_2\cap K_1$ that $(H_2\cap K_1)^U=(H_4\cap K_1)^U$.
On the other hand, since $H_3$ is contained in the field-extension subgroup $\Sp_{2a}(q^b){:}(b/2)$ of $L$, Proposition~\ref{prop:SpaO-<Sp} shows that $L\neq H_3K_1$. As $H_3$ has index $2$ in $H_4$, this implies that $H_3\cap K_1=H_4\cap K_1$. In particular,
\[
(H_3\cap K_1)^U=(H_4\cap K_1)^U=(H_2\cap K_1)^U.
\]
Since $H_3\cap K_1\leqslant H_1\cap K_1\leqslant H_2\cap K_1$, it follows that $(H_1\cap K_1)^U=(H_2\cap K_1)^U$. For $i\in\{1,2\}$, the kernel of $H_i\cap K_1$ acting on $U$ is
\[
H_i\cap K_1\cap R=(H_i\cap R)\cap K_1=(E\times U_{(b)}(I))\cap K_1.
\]
This in conjunction with $(H_1\cap K_1)^U=(H_2\cap K_1)^U$ leads to $|H_1\cap K_1|=|H_2\cap K_1|$. However, the factorizations $L=H_1K_1=H_2K_1$ require $|H_2\cap K_1|/|H_1\cap K_1|=|H_2|/|H_1|=2$, a contradiction.
If $|H|_2\leqslant|H_1|_2$ and $|K|_2\leqslant|K_1|_2$, then this argument leads to a contradiction again. Hence $(G,H,K)$ tightly contains the triple $(L,H^{(\infty)}{:}b_2,\Omega_{2m}^-(q))$ or $(L,H^{(\infty)},\mathrm{O}_{2m}^-(q))$, as in Example~\ref{LemSymplecticPm4} or~\ref{LemSymplecticPm5}.
\end{proof}

\subsection{Proof of Theorem~\ref{ThmSymplectic}}
\ \vspace{1mm}

Since the triples $(G_0,H_0,K_0)$ in Tables~\ref{TabSymplectic1} and~\ref{TabSymplectic2} are already shown to give rise to factorizations $G_0=H_0K_0$, it remains to prove the ``only if'' part. Let $G=HK$. By Lemma~\ref{LemSymplectic52} we may assume $(m,q)\neq(2,3)$, $(2,4)$, $(3,2)$, $(3,3)$, $(3,4)$ or $(4,2)$. If $K\leqslant\Pa_1[G]$ or $\N_2[G]$, then Proposition~\ref{prop:K=P1-Sp} shows that $(G,H,K)$ tightly contains some $(G_0,H_0,K_0)$ in Table~\ref{TabSymplectic1} or Table~\ref{TabSymplectic2}. Now assume that $K$ is not contained in $\Pa_1[G]$ or $\N_2[G]$. Then part~(c) or~(d) of Proposition~\ref{lem:K<B-Sp} holds, and in particular we have $B^{(\infty)}\in\{\Omega_{2m}^+(q),\Omega_{2m}^-(q)\}$. It then follows from Lemmas~\ref{LemSymplectic01} and~\ref{LemSymplectic1Rows4--9} that, interchanging $H$ and $K$ if necessary, $K^{(\infty)}\in\{\Omega_{2m}^+(q),\Omega_{2m}^-(q)\}$.
For these candidates, Propositions~\ref{LemSymplectic12}, \ref{prop:Sp(2)=O^-O^+}, \ref{prop:Sp(2)=O^+O^-}, \ref{prop:SpaO+<Sp}, \ref{C2-subgroup}, \ref{prop:Sp(q)=Sp(q^{1/2})Q-}, \ref{prop:SpaO-<Sp} and \ref{PropSymplecticPm} together show that $(G,H,K)$ tightly contains some $(G_0,H_0,K_0)$ in Table~\ref{TabSymplectic1} or Table~\ref{TabSymplectic2}, where the reference labels are listed in the last column of the table. This completes the proof.

%%%%%%%%%%%%%%%%%%%%%%%%%%%%%%%%%%%%%%%%%%%%%%%%%%%%%%%%%%%%%%%%%%%
%%%%%%%%%%%%%%%%%%%%%%%%%%%%%%%%%%%%%%%%%%%%%%%%%%%%%%%%%%%%%%%%%%%

\section{Biperfect bicrossproduct Hopf algebras}\label{sec:Application}

%%%%%%%%%%%%%%%%%%%%%%%%%%%%%%%%%%%%%%%%%%%%%%%%%%%%%%%%%%%%%%%%%%%
%%%%%%%%%%%%%%%%%%%%%%%%%%%%%%%%%%%%%%%%%%%%%%%%%%%%%%%%%%%%%%%%%%%

We call a triple $(G,H,K)$ \emph{biperfect} if $H$ and $K$ are perfect and self-normalizing subgroups of the group $G$ such that $G=HK$ is an exact factorization. The following example of a biperfect triple, which can be confirmed in \magma~\cite{BCP1997}, has $|G|<|\M_{24}|$.

\begin{example}\label{Ex:Hopf}
There are two non-isomorphic groups $G$ of the form
\[
2^{14}{:}\A_7=2^4.(2^6.(2^4{:}\A_7))=2^4.(2^4.(2^6{:}\A_7))
\]
that have exact factorizations $G=HK$ with both $H$ and $K$ perfect and self-normalizing. One has precisely three conjugacy classes of normal subgroups of the form $2^4$, and the other has only two. For the first $G$, interchanging $H$ and $K$ if necessary, the factors $H$ and $K$ have the form $2^4{:}\A_5$ and $2^8.\PSL_2(7)$ respectively ($K=2^8{}^{\boldsymbol{\cdot}}\PSL_2(7)$ if $H\cong\ASL_2(4)$, and $K=2^8{:}\PSL_2(7)$ if $H\ncong\ASL_2(4)$). For the second $G$, the factors have the form $2^8{:}\A_5$ and $2^4{:}\PSL_2(7)$ respectively.
In any case, $|G|=41287680<244823040=|\M_{24}|$.
\end{example}

In this section, we prove that the order $|G|=41287680$ in Example~\ref{Ex:Hopf} is the smallest possible $|G|$ for a biperfect triple $(G,H,K)$. The proof will be largely based on the classification of factorizations of finite simple groups and computation in \magma~\cite{BCP1997}.

\begin{theorem}\label{Thm:Hopf}
Let $(G,H,K)$ be a biperfect triple with $|G|\leqslant41287680$. Then $(G,H,K)$ is as described in Example~$\ref{Ex:Hopf}$.
\end{theorem}

Theorem~\ref{Thm:Hopf} will be proved at the end of this section, after preparation in this subsection and discussions in the next two subsections.

For a group factorization $G=HK$ with both $H$ and $K$ perfect, since $G=HK=H'K'\subseteq G'$, the group $G$ is perfect. If $|G|\leqslant41287680$, then by~\cite[Corollary~1.4]{HP1989}, $G/\Rad(G)=T_1\times\cdots\times T_k$ for some nonabelian simple groups $T_1,\ldots,T_k$. In particular, there exist epimorphisms $\varphi_1,\ldots,\varphi_k$ from $G$ to $T_1,\ldots,T_k$ respectively.

\begin{hypothesis}\label{Hyp:Hopf}
Let $G$ be a group with $|G|\leqslant41287680$ such that $G$ is not as described in Example~\ref{Ex:Hopf}, let $H$ and $K$ be subgroups of $G$ that are both perfect and self-normalizing such that $G=HK$ and $H\cap K=1$, and let $G/\Rad(G)=T_1\times\cdots\times T_k$ with nonabelian simple groups $T_1,\ldots,T_k$. For each $i\in\{1,\ldots,k\}$, let $\varphi_i$ be the epimorphism from $G$ to $T_i$.
\end{hypothesis}

\begin{lemma}\label{Lem:Hopf1}
Under Hypothesis~$\ref{Hyp:Hopf}$, the following statements hold:
\begin{enumerate}[\rm(a)]
\item $\Z(G)=1$;
\item for each normal subgroup $N$ of $G$, if $\,\overline{\phantom{x}}$ denotes the natural homomorphism $G\to G/N$, then $|N|=|H\cap N||K\cap N||\overline{H}\cap\overline{K}|$;
\item if $\,\overline{\phantom{x}}$ denotes the natural homomorphism $G\to G/\Rad(G)$, then for each $p\in\pi(\overline{H}\cap\overline{K})\Setminus\pi(\mathrm{H}_2(\overline{G},\mathbb{Z}))$, the integer $d:=\log_p|\Rad(G)|_p$ is such that $\GL_d(p)$ has a section isomorphic to $T_i$ for some $i\in\{1,\ldots,k\}$.
\end{enumerate}
\end{lemma}

\begin{proof}
Since $\Z(G)\leqslant\Nor_G(H)=H$ and $\Z(G)\leqslant\Nor_G(K)=K$, we deduce from $H\cap K=1$ that $\Z(G)=1$, proving~(a).

For~(b), the factorization $G=HK$ implies that $\overline{G}=\overline{H}\,\overline{K}$, and so $|\overline{G}||\overline{H}\cap\overline{K}|=|\overline{H}||\overline{K}|$. Substituting $|\overline{H}|=|H|/|H\cap N|$, $|\overline{K}|=|K|/|K\cap N|$ and $|\overline{G}|=|G|/|N|=|H||K|/|N|$, we obtain $|N|=|H\cap N||K\cap N||\overline{H}\cap\overline{K}|$.

For~(c), suppose for a contradiction that $\GL_d(p)$ has no section isomorphic to $T_i$ for any $i\in\{1,\ldots,k\}$. Since $|\overline{H}\cap\overline{K}|$ is divisible by $p$, we see from~(b) that $|\Rad(G)|$ is divisible by $p$.
In particular, there exist normal subgroups $N$ and $M$ of $G$ such that $N<M\leqslant\Rad(G)$, $|\Rad(G)/M|_p=1$ and $M/N=\C_p^\ell$ for some positive integer $\ell$.
Consider the action $\varphi$ of $G/M$ on $M/N$. As $\ell=\log_p|M/N|\leqslant\log_p|\Rad(G)|_p=d$, there is no section of $\GL_\ell(p)$ isomorphic to $T_i$ for any $i\in\{1,\ldots,k\}$. Hence $\mathrm{Im}(\varphi)$ is solvable, and so $\mathrm{Ker}(\varphi)\geqslant(G/M)^{(\infty)}$. Since $G$ is perfect and so is $G/M$, it follows that $\mathrm{Ker}(\varphi)\geqslant G/M$, or equivalently, $\mathrm{Im}(\varphi)=1$. This implies that $M/N\leqslant\Z(G/N)$. However, we derive from $|\Rad(G)/M|_p=1$ and $p\notin\pi(\mathrm{H}_2(\overline{G},\mathbb{Z}))$ that
\[
|\mathrm{H}_2(G/M,\mathbb{Z})|_p=|\mathrm{H}_2(G/\Rad(G),\mathbb{Z})|_p=1,
\]
a contradiction. This completes the proof.
\end{proof}

\begin{lemma}\label{Lem:Hopf10}
Let $G$ be a finite group such that $G/\Rad(G)=T_1\times T_2$ for some nonabelian simple groups $T_1$ and $T_2$, let $\Rad(G)=R_0>R_1>\cdots>R_n=1$ be a normal series of $G$ such that $R_{i-1}/R_i$ is a minimal normal subgroup of $G/R_i$ for each $i\in\{1,\ldots,n\}$, and let $H$ be a perfect subgroup of $G$ such that $HR_0/R_0\leqslant T_1$.
Suppose $\pi(T_2)\nsubseteq\pi(R_0)$ and that, for each $i\in\{1,\ldots,n\}$, either $(R_0/R_{i-1}).T_1$ or $(R_0/R_{i-1}).T_2$ acts trivially on $R_{i-1}/R_i$.
Then $\Nor_G(H)>H$.
\end{lemma}

\begin{proof}
We claim that there exists no $i\in\{0,1,\ldots,n-1\}$ such that $HR_{i+1}/R_{i+1}\geqslant R_i/R_{i+1}$ while $(R_0/R_i).T_2$ acts nontrivially on $R_i/R_{i+1}$.
Suppose for a contradiction that such an $i$ exists, and take $i$ to be the smallest.
Let $N=R_i/R_{i+1}$, $X=HR_i/R_i$ and $Y=((R_0/R_i).T_2)^{(\infty)}$.
In the group $N.(XY)$, the subgroup $N.X$ acts on the normal subgroup $N.Y$ such that the induced actions on $N$ and on $Y$ are both trivial.
Note that $N.X$ is perfect as a quotient of $H$.
It follows that $N.X$ centralizes $N.Y$, and in particular, $N$ centralizes $N.Y$, which is impossible as $(R_0/R_i).T_2$ acts nontrivially on $N$.
This proves our claim.

Now take any $p\in\pi(T_2)\Setminus\pi(R_0)$.
We prove by induction on $i$ that
\begin{equation}\label{Eqn:Hopf3}
\Nor_{G/R_i}(HR_i/R_i)\geqslant(HR_i/R_i)\times p
\end{equation}
for all $i\in\{0,1,\ldots,n\}$, which will imply $\Nor_G(H)=\Nor_{G/R_n}(HR_n/R_n)\geqslant(HR_n/R_n)\times p=H\times p$ and thus complete the proof.
The conclusion~\eqref{Eqn:Hopf3} is obvious for $i=0$.
Suppose that $\Nor_{G/R_i}(HR_i/R_i)\geqslant(HR_i/R_i)\times p$ for some $i\in\{0,1,\ldots,n-1\}$.
Let $\,\overline{\phantom{x}}$ be the natural homomorphism $G\to G/R_i$, let $\,\widetilde{\phantom{x}}$ be the natural homomorphism $G\to G/R_{i+1}$, and denote $M=\widetilde{R_i}$.

First assume that $\overline{R_0}.T_2$ acts trivially on $M$.
Consider the preimage $M.(\overline{H}\times p)$ of $\overline{H}\times p$ in $M.\overline{G}$.
Since $p$ is coprime to $|M|$, this preimage has the form $(M.\overline{H})\times p$.
Then as $\widetilde{H}\leqslant M.\overline{H}$, it follows that $\Nor_{\widetilde{G}}(\widetilde{H})\geqslant\widetilde{H}\times p$.

Next assume that $\overline{R_0}.T_2$ acts nontrivially on $M$.
Then by the condition of the lemma, $\overline{R_0}.T_1$ acts trivially on $M$.
Since $M$ is a minimal normal subgroup of $\widetilde{G}$ and since $\widetilde{H}\ngeqslant M$ by the above claim, it follows that $\widetilde{H}\cap M=1$.
Then as $\overline{H}\leqslant\overline{R_0}.T_1$ acts trivially on $M$, the subgroup $M.\overline{H}$ of $M.(\overline{H}\times p)$ has the form $M\times\overline{H}$ and contains a unique subgroup isomorphic to $\overline{H}$.
This together with $|M|_p=1$ implies that $\Nor_{\widetilde{G}}(\widetilde{H})\geqslant\Cen_{M.(\overline{H}\times p)}(\widetilde{H})\geqslant\widetilde{H}\times p$, completing the proof.
\end{proof}

\subsection{Unique nonsolvable composition factor}

\begin{lemma}\label{Lem:Hopf2}
Under Hypothesis~$\ref{Hyp:Hopf}$, $G/\Rad(G)\neq\PSL_2(q)$ for $q\geqslant5$.
\end{lemma}

\begin{proof}
Suppose for a contradiction that $G/\Rad(G)=\PSL_2(q)$ with $q\geqslant5$. Let $\,\overline{\phantom{x}}$ be the natural homomorphism $G\to G/\Rad(G)$. Then $\overline{H}$ and $\overline{K}$ are perfect subgroups of $\overline{G}=\PSL_2(q)$, and checking the candidates for the factorization $\overline{G}=\overline{H}\,\overline{K}$ in~\cite[Theorem~3.3]{LX}, we conclude that either $q=9$, or at least one of $\overline{H}$ or $\overline{K}$ equals $\overline{G}$.

\textsf{Case 1:} $q=5$ and $\overline{G}=\A_5$.
Since $\overline{H}$ and $\overline{K}$ are perfect, $\overline{H}=\overline{K}=\overline{G}=\A_5$, and so we derive from Lemma~\ref{Lem:Hopf1}(b) that $|\Rad(G)|$ is divisible by $2^2\cdot3\cdot5$.
Note that $G$ is perfect with $\Z(G)=1$ (as Lemma~\ref{Lem:Hopf1}(a) asserts) and $|\Rad(G)||\A_5|=|G|\leqslant41287680$.
It follows that $\pi(\Rad(G))=\{2,3,5\}$, $|\Rad(G)|_2\geqslant2^4$, $|\Rad(G)|_3\geqslant3^4$ and $|\Rad(G)|_5\geqslant5^2$.

Suppose $|\Rad(G)|_2=2^4$.
Then $G$ has the form $R.(V.\overline{G})$ with $\pi(R)=\{3,5\}$ and $V=2^4$ such that $V$ is an irreducible $\overline{G}$-module.
Let $\,\widetilde{\phantom{x}}$ be the natural homomorphism $G\to G/R$.
It follows that $|\widetilde{H}\cap V|$ and $|\widetilde{K}\cap V|$ are both in $\{1,2^4\}$.
However, this implies that $|H\cap\Rad(G)|_2|K\cap\Rad(G)|_2\in\{1,2^4,2^8\}$, contradicting Lemma~\ref{Lem:Hopf1}(b) as $|\overline{H}\cap\overline{K}|_2=2^2$ and $|\Rad(G)|_2=2^4$.

Thus, $|\Rad(G)|_2\geqslant2^5$.
If $|\Rad(G)|_3=3^4$, then $G$ has the form $R.(V.\A_5)$ or $R.(V.\SL_2(5))$, where $\pi(R)=\{2,5\}$ and $V=3^4$ is an irreducible $\A_5$- or $\SL_2(5)$-module, and we obtain a similar contradiction as in the previous paragraph.
Consequently, $|\Rad(G)|_3\geqslant3^5$.

Suppose $|\Rad(G)|_5\leqslant5^3$.
In the same vein as above, we obtain $G=R.(V.\SL_2(5))$ with $\pi(R)=\{2,3,5\}$ such that $V=5^2$ is an irreducible $\SL_2(5)$-module and $|R|_5=5$.
However, this implies that $\Z(G)\geqslant5$, which is impossible by Lemma~\ref{Lem:Hopf1}(a).

Now $|\Rad(G)|_2\geqslant2^5$, $|\Rad(G)|_3\geqslant3^5$ and $|\Rad(G)|_5\geqslant5^4$.
It follows that $|G|\geqslant2^5\cdot3^5\cdot5^4\cdot|\A_5|>41287680$, a contradiction.

\textsf{Case 2:} $q=7$.
This yields a contradiction $|G|\geqslant2^4\cdot3^7\cdot7^4\cdot|\PSL_2(7)|>41287680$ in the same vein as in Case~1.

\textsf{Case 3:} $q=9$ and $\overline{G}=\A_6$.
Since $\overline{G}=\overline{H}\,\overline{K}$ with $\overline{H}$ and $\overline{K}$ perfect, interchanging $H$ and $K$ if necessary, one of the following occurs:
\begin{enumerate}[(i)]
\item $\overline{H}=\overline{K}=\overline{G}$;
\item $\overline{H}=\overline{G}$, and $\overline{K}\cong\A_5$;
\item $\overline{H}\cong\overline{K}\cong\A_5$.
\end{enumerate}
In particular, $|\overline{H}\cap\overline{K}|_2\geqslant2$ and $|\overline{H}\cap\overline{K}|_5=5$.
We first suppose for a contradiction that $|\Rad(G)|_2\leqslant2^3$.
Note that $G$ is perfect and note from Lemma~\ref{Lem:Hopf1}(b) that $|\Rad(G)|$ is divisible by $|\overline{H}\cap\overline{K}|$.
We then conclude from $|\mathrm{H}_2(\overline{G},\mathbb{Z})|_2=2$ that $|\Rad(G)|_2=|\overline{H}\cap\overline{K}|_2=2$ as in~(iii).
It follows that the Hall $2'$-subgroup $N$ of $\Rad(G)$ is normal in $G$ such that $G/N=\SL_2(9)$.
However, this implies that $HN/N\cong KN/N=2.\A_5$, and so $|(HN/N)\cap(KN/N)|_2=2^2$, contradicting Lemma~\ref{Lem:Hopf1}(b).

Consequently, $|\Rad(G)|_2\geqslant2^4$. Moreover, since $|\mathrm{H}_2(\overline{G},\mathbb{Z})|_5=1$, there exist normal subgroups $N$ and $M$ of $G$ with $N<M\leqslant\Rad(G)$ and $|\Rad(G)/M|_5=1$ such that $M/N=5^\ell$ is an irreducible $(G/M)$-module for some integer $\ell>1$.
Since $\overline{G}=\A_6$, we have $\ell\geqslant4$, and so $|\Rad(G)|_5\geqslant5^4$.
Now that $|\Rad(G)|_2\geqslant2^4$ and $|\Rad(G)|_5\geqslant5^4$, we deduce from $|\Rad(G)||\A_6|\leqslant41287680$ that $\pi(\Rad(G))\subseteq\{2,3,5\}$.

We next claim that $|\Rad(G)|_5\geqslant5^5$.
Suppose otherwise that $|\Rad(G)|_5=5^4$.
Then the action of $G/M$ on $M/N$ induces the subgroup $\SL_2(9)$ of $\GL_4(5)$, and the only possible $(HM/M)$- or $(KM/M)$-invariant subspaces of $5^4$ are $1$, $5^2$ and $5^4$.
This implies that $|H\cap\Rad(G)|_5$ and $|K\cap\Rad(G)|_5$ are both in $\{1,5^2,5^4\}$.
It follows that $|\Rad(G)|_5\neq|H\cap\Rad(G)|_5|K\cap\Rad(G)|_5|\overline{H}\cap\overline{K}|_5$, contradicting Lemma~\ref{Lem:Hopf1}(b).
This proves the claim.

Now $|\Rad(G)|_3=1$ as $|\Rad(G)||\A_6|\leqslant41287680$ with $|\Rad(G)|_2\geqslant2^4$ and $|\Rad(G)|_5\geqslant5^5$.
For (i) or (ii), since $|\overline{H}\cap\overline{K}|_3\geqslant3$, Lemma~\ref{Lem:Hopf1}(b) would imply that $|\Rad(G)|_3\geqslant3$, a contradiction. Consequently, (iii) occurs.
Moreover, it follows from $|\Rad(G)||\A_6|\leqslant41287680$ and $|\Rad(G)|_5\geqslant5^5$ that $|\Rad(G)|_2\leqslant2^7$.
Note that the irreducible modules of $\A_6$ or $\SL_2(9)$ over $\mathbb{F}_2$ have dimension $4$ or $16$, while the irreducible modules of $\A_6$ or $\SL_2(9)$ over $\mathbb{F}_5$ have dimension $4$, $5$, $8$, $10$ or $20$.
We deduce from $G'=G$ and $\Z(G)=1$ that $G=R.(V.\widetilde{G})$ with $R=5^5$, $V=2^4$ and $\widetilde{G}=\A_6$ or $\SL_2(9)$, where $\,\widetilde{\phantom{x}}$ is the natural homomorphism $G\to G/(R.V)$.
It follows that the only $\widetilde{H}$- or $\widetilde{K}$-invariant subspaces of $V$ are $1$ and $2^4$, whence $|(HR/R)\cap V||(KR/R)\cap V|$ is either $1$ or at least $2^4$.
Accordingly, $|H\cap\Rad(G)|_2|K\cap\Rad(G)|_2$ is either at most $2^2$ or at least $2^4|\Z(\widetilde{G})|^2$.
However, as $|\overline{H}\cap\overline{K}|_2=2$, this contradicts Lemma~\ref{Lem:Hopf1}(b).

\textsf{Case 4:} $q=8$ or $q\geqslant11$.
Then at least one of $\overline{H}$ or $\overline{K}$ equals $\overline{G}$.
It follows that $\overline{H}\cap\overline{K}$ is a nontrivial perfect subgroup of $\PSL_2(q)$, whence $|\overline{H}\cap\overline{K}|$ has at least three prime divisors.
For each $p\in\pi(\overline{H}\cap\overline{K})\Setminus\pi(\mathrm{H}_2(\overline{G},\mathbb{Z}))$, denoting $d=\log_p|\Rad(G)|_p$, Lemma~\ref{Lem:Hopf1}(c) asserts that $\GL_d(p)$ has a section isomorphic to $\PSL_2(q)$. Thus, if $p$ is odd, then either $p^d=5^4$ or $p^d\geqslant3^6$. Note that $|\overline{H}\cap\overline{K}|$ has at least three prime divisors, and so does $\Rad(G)$ by Lemma~\ref{Lem:Hopf1}(b). As a consequence, since $|\mathrm{H}_2(\overline{G},\mathbb{Z})|=(2,q-1)$, there are at least two odd primes in $\pi(\overline{H}\cap\overline{K})\Setminus\pi(\mathrm{H}_2(\overline{G},\mathbb{Z}))$.
Therefore, $|\Rad(G)|\geqslant2\cdot5^4\cdot3^6$. This together with $|\Rad(G)||\PSL_2(q)|=|G|\leqslant41287680$ yields $2\cdot5^4\cdot3^6\cdot|\PSL_2(q)|\leqslant41287680$, which is not possible as $q\geqslant8$.
\end{proof}

\begin{lemma}\label{Lem:Hopf13}
Under Hypothesis~$\ref{Hyp:Hopf}$, if $G/\Rad(G)=\A_7$ or $\A_8$, then $(G,H,K)$ is as described in Example~$\ref{Ex:Hopf}$.
\end{lemma}

\begin{proof}
First suppose that $G/\Rad(G)=\A_7$.
Let $\,\overline{\phantom{x}}$ be the natural homomorphism $G\to G/\Rad(G)$.
Since $\overline{G}=\overline{H}\,\overline{K}$ with $\overline{H}$ and $\overline{K}$ perfect, $|\overline{H}\cap\overline{K}|$ is divisible by $2^2$.
By Lemma~\ref{Lem:Hopf1}(b), it follows that $|\Rad(G)|$ is divisible by $2^2$.
Then as $|\mathrm{H}_2(\A_7,\mathbb{Z})|_2=2$, we conclude that $|\Rad(G)|_2\geqslant2^4$.
If $|\Rad(G)|$ is divisible by some odd prime, then the conditions $|\Rad(G)||\A_7|\leqslant41287680$ and $\Z(G)=1$ force $G=(2^4\times3^6).\A_7$.
However, computation in \magma~\cite{BCP1997} shows that this $G$ does not admit any biperfect triple $(G,H,K)$.
Hence $\Rad(G)$ is a $2$-group.
This together with $|\Rad(G)||\A_7|\leqslant41287680$ and $\Z(G)=1$ yields that $\Rad(G)$ is one of the groups:
\begin{align*}
&2^4,\ \ 2^6,\ \ 2^{14},\ \ 2^{4+4},\ \ 2^{4+1+4},\ \ 2^{4+6},\ \ 2^{4+1+6},\ \ 2^{4+1+1+6},\\
&2^{6+4},\ \ 2^{6+1+4},\ \ 2^{6+6},\ \ 2^{6+1+6},\ \ 2^{6+1+1+6},\ \ 2^{4+4+4},\ \ 2^{4+4+1+4},\\
&2^{4+1+4+4},\ \ 2^{4+1+1+4+4},\ \ 2^{4+1+4+1+4},\ \ 2^{4+4+6},\ \ 2^{4+6+4},\ \ 2^{6+4+4}.
\end{align*}
For these candidates $G=\Rad(G).\A_7$, searching in \magma~\cite{BCP1997} for biperfect triples $(G,H,K)$ shows that $(G,H,K)$ is as described in Example~\ref{Ex:Hopf}.

Now suppose that $G/\Rad(G)=\A_8$.
As above, the conditions $|\Rad(G)||\A_8|\leqslant41287680$ and $\Z(G)=1$ imply that $\Rad(G)$ is one of the groups:
\[
2^4,\ \ 2^6,\ \ 2^{4+4},\ \ 2^{4+1+4},\ \ 2^{4+6},\ \ 2^{4+1+6},\ \ 2^{6+4},\ \ 2^{6+1+4}.
\]
Then computation in \magma~\cite{BCP1997} shows that none of these gives rise to a biperfect triple $(G,H,K)$. The proof is thus complete.
\end{proof}

\begin{lemma}\label{Lem:Hopf12}
Under Hypothesis~$\ref{Hyp:Hopf}$, either $k\geqslant2$, or $(G,H,K)$ is as described in Example~$\ref{Ex:Hopf}$.
\end{lemma}

\begin{proof}
Suppose for a contradiction that $k=1$ whereas $(G,H,K)$ is not as described in Example~$\ref{Ex:Hopf}$.
Let $\,\overline{\phantom{x}}$ be the natural homomorphism $G\to G/\Rad(G)$.
From Lemmas~\ref{Lem:Hopf2},~\ref{Lem:Hopf13} and the inequality $|T_1|\leqslant41287680$, we deduce that $T_1$ is one of the groups:
\begin{align*}
&\SL_3(3),\ \ \SU_3(3),\ \ \M_{11},\ \ \PSL_3(4),\ \ \SU_4(2),\ \ \Sz(8),\ \ \SU_3(4),\ \ \M_{12},\ \PSU_3(5),\ \ \J_1,\\
&\A_9,\ \ \SL_3(5),\ \ \M_{22},\ \ \J_2,\ \ \Sp_4(4),\ \ \Sp_6(2),\ \ \A_{10},\ \ \PSL_3(7),\ \ \PSU_4(3),\ \ \G_2(3),\ \ \PSp_4(5),\\
&\PSU_3(8),\ \ \SU_3(7),\ \ \PSL_4(3),\ \ \SL_5(2),\ \ \M_{23},\ \ \SU_5(2),\ \ \SL_3(8),\ \ ^2{}\F_4(2)',\ \ \A_{11},\ \ \Sz(32).
\end{align*}

Suppose that $T_1\in\{\SL_3(3),\SU_3(3),\M_{11},\PSL_3(4),\SU_4(2),\SU_3(4)\}$.
Since $\overline{G}=\overline{H}\,\overline{K}$ with $\overline{H}$ and $\overline{K}$ perfect, $\pi(\overline{H}\cap\overline{K})\supseteq\{2,3\}$.
Then Lemma~\ref{Lem:Hopf1}(b) implies that $\pi(\Rad(G))\supseteq\{2,3\}$.
If $T_1=\SL_3(3)$, then since $\Z(G)=1$, it follows that $|\Rad(G)|$ is divisible by $2^{12}3^3$, which leads to $|G|=|\Rad(G)||T_1|>41287680$, not possible.
Similarly, none of $\SU_3(3)$, $\M_{11}$, $\PSL_3(4)$, $\SU_4(2)$ or $\SU_3(4)$ is a possible candidate for $T_1$.

Now $T_1\notin\{\SL_3(3),\SU_3(3),\M_{11},\PSL_3(4),\SU_4(2),\SU_3(4)\}$.
By the classification of factorizations of $T_1$, we have $G\neq T_1$ and hence $\Rad(G)>1$.
If $T_1=\Sz(8)$, then as $\Z(G)=1$, it follows that $|\Rad(G)|$ is divisible by $2^{12}$, which leads to a contradiction that $|G|=|\Rad(G)||T_1|>41287680$.
The other candidates for $T_1$ are excluded similarly.
\end{proof}

\subsection{Two nonsolvable composition factors}

\begin{lemma}\label{Lem:Hopf3}
Under Hypothesis~$\ref{Hyp:Hopf}$, $G/\Rad(G)\neq\PSL_2(q)\times\PSL_2(q')$ for $q'\geqslant q\geqslant11$.
\end{lemma}

\begin{proof}
Suppose for a contradiction that $G/\Rad(G)=T_1\times T_2$ with $T_1=\PSL_2(q)$ and $T_2=\PSL_2(q')$ and $q'\geqslant q\geqslant11$. For $i\in\{1,2\}$, applying~\cite[Theorem~3.3]{LX} to the factorization $T_i=\varphi_i(H)\varphi_i(K)$ into perfect groups $\varphi_i(H)$ and $\varphi_i(K)$, we conclude that at least one of $\varphi_i(H)$ or $\varphi_i(K)$ is $T_i$. Without loss of generality, assume $\varphi_2(K)=T_2$. Let $\,\overline{\phantom{x}}$ be the natural homomorphism $G\to G/\Rad(G)$. It follows that either $\overline{K}\cong\PSL_2(q')$, or $|\overline{H}\cap\overline{K}|$ has at least three prime divisors. For the latter case, the same argument as in the proof of Lemma~\ref{Lem:Hopf2} gives a contradiction. Consequently, $\overline{K}\cong\PSL_2(q')$. Moreover, since at least one of $\varphi_1(H)$ or $\varphi_1(K)$ is $T_1$ while $|\overline{H}\cap\overline{K}|$ cannot have at least three prime divisors, we deduce that $\overline{H}\cong\PSL_2(q)$.
Accordingly, both $H$ and $K$ have order at most $41287680/|\PSL_2(q)|$ and hence at most $10^6$. Then, by~\cite{HP1989}, each of $H$ and $K$ is one of the following:
\[
\PSL_2(q),\ \ \SL_2(q),\ \ \PSL_2(q'),\ \ \SL_2(q'),\ \ \ASL_2(11)\text{ with }q=11.
\]
For the last candidate, we would have $|G|\geqslant|\ASL_2(11)||\PSL_2(11)|>41287680$, a contradiction.
Therefore, neither $H$ nor $K$ is $\ASL_2(11)$. Since $|G|=|H||K|$, we derive that $|\Rad(G)|=|G|/(|\PSL_2(q)||\PSL_2(q')|)$ divides $4$, and so $\Rad(G)\leqslant\Z(G)$. Thus, by Lemma~\ref{Lem:Hopf1}(a), $\Rad(G)=1$, which implies that $G=T_1\times T_2=\PSL_2(q)\times\PSL_2(q')$. Now $H\cong\PSL_2(q)$ and $K\cong\PSL_2(q')$ while both $H$ and $K$ are self-normalizing in $G$. It follows that $q=q'$, and
\[
H=\{(t,t^\alpha)\mid t\in\PSL_2(q)\}\ \text{ and }\ K=\{(t,t^\beta)\mid t\in\PSL_2(q)\}
\]
for some automorphisms $\alpha$ and $\beta$ of $\PSL_2(q)$. However, since the centralizer of the automorphism $\alpha\beta^{-1}$ in $\PSL_2(q)$ is nontrivial, it follows that $H\cap K\neq1$, a contradiction.
\end{proof}

In the following discussion, we will single out groups in the set
\begin{equation}\label{Eqn:Hopf1}
\{\A_5,\ \PSL_2(7),\ \A_6,\ \PSL_2(8),\ \A_7,\ \A_8,\ \M_{12},\ \A_9\}
\end{equation}
for the possibilities of $T_i$ in Hypothesis~\ref{Hyp:Hopf}.

\begin{lemma}\label{Lem:Hopf4}
Under Hypothesis~$\ref{Hyp:Hopf}$, if $k=2$, then either $T_1$ or $T_2$ lies in~\eqref{Eqn:Hopf1}.
\end{lemma}

\begin{proof}
Suppose for a contradiction that $G/\Rad(G)=T_1\times T_2$ such that neither $T_1$ nor $T_2$ lies in~\eqref{Eqn:Hopf1}. Then both $T_1$ and $T_2$ have order at least $|\PSL_2(11)|$. Let $\,\overline{\phantom{x}}$ be the natural homomorphism $G\to G/\Rad(G)$. For each $p\in\pi(\overline{H}\cap\overline{K})\Setminus\pi(\mathrm{H}_2(\overline{G},\mathbb{Z}))$, denoting $d=\log_p|\Rad(G)|_p$, Lemma~\ref{Lem:Hopf1}(c) asserts that $\GL_d(p)$ has a section isomorphic to $T_1$ or $T_2$. Thus, if $p$ is odd, then either $p^d=3^4$ or $p^d\geqslant11^2$.

Suppose that $|\overline{H}\cap\overline{K}|$ has at least three prime divisors. If both $\mathrm{H}_2(T_1,\mathbb{Z})$ and $\mathrm{H}_2(T_2,\mathbb{Z})$ are $2$-groups, then $\pi(\overline{H}\cap\overline{K})\Setminus\pi(\mathrm{H}_2(\overline{G},\mathbb{Z}))$ contains at least two odd primes, and so
\[
|G|=|\Rad(G)||T_1||T_2|\geqslant3^4\cdot11^2\cdot|\PSL_2(11)|\cdot|\PSL_2(11)|>41287680,
\]
a contradiction. If at least one of $\mathrm{H}_2(T_1,\mathbb{Z})$ or $\mathrm{H}_2(T_2,\mathbb{Z})$ is not a $2$-group, then by checking the nonabelian simple groups $T$ with $|\PSL_2(11)||T|\leqslant41287680$, we conclude that
\[
|T_1||T_2|\geqslant|\PSL_2(11)||\PSL_3(4)|\ \text{ and }\ \pi(\mathrm{H}_2(\overline{G},\mathbb{Z}))=\{2,3\},
\]
which leads to $|G|\geqslant11^2|\PSL_2(11)||\PSL_3(4)|>41287680$, again a contradiction.

Hence $|\overline{H}\cap\overline{K}|$ has at most two prime divisors. For $i\in\{1,2\}$, by the factorization $T_i=\varphi_i(H)\varphi_i(K)$ into perfect groups $\varphi_i(H)$ and $\varphi_i(K)$, we conclude from the classification of factorizations of $T_i$ that at least one of $\varphi_i(H)$ or $\varphi_i(K)$ is $T_i$.
Then we derive from $\overline{G}=\overline{H}\,\overline{K}$ that, interchanging $H$ and $K$ if necessary, $\overline{H}\cong T_1$ and $\overline{K}\cong T_2$. Since $|H||K|=|G|\leqslant41287680$ and both $T_1$ and $T_2$ have order at least $|\PSL_2(11)|$, it follows that both $H$ and $K$ have order at most $41287680/|\PSL_2(11)|$.
Then by~\cite{HP1989}, both $H$ and $K$ are quasisimple, and as $|\Rad(G)||T_1||T_2|=|H||K|$, it follows that $\pi(\Rad(G))\subseteq\{2,3\}$.
Therefore, $\pi(T_2)\nsubseteq\pi(\Rad(G))$, and so Lemma~\ref{Lem:Hopf10} gives $\Nor_G(H)>H$, a contradiction.
\end{proof}

\begin{lemma}\label{Lem:HopfTwo}
Under Hypothesis~$\ref{Hyp:Hopf}$, we have $k\neq2$.
\end{lemma}

\begin{proof}
Suppose for a contradiction that $k=2$. By Lemmas~\ref{Lem:Hopf3} and~\ref{Lem:Hopf4}, as $|T_1||T_2|\leqslant41287680$, we may assume that, interchanging $T_1$ and $T_2$ if necessary,
\[
T_1\in\{\A_5,\PSL_2(7),\A_6,\PSL_2(8),\A_7,\A_8,\M_{12},\A_9\}.
\]
We use the same reduction throughout as in the previous lemmas: Lemma~\ref{Lem:Hopf1}, the Schur multipliers and the minimum module degrees restrict the chief factors of $\Rad(G)$, while Lemma~\ref{Lem:Hopf10} excludes the cases in which one simple direct factor acts trivially on every chief factor. The remaining extensions, subject to the bound $|G|\leqslant41287680$, are excluded by computation in \magma~\cite{BCP1997}, apart from the case singled out below, where the candidate groups are too large for direct computation.

In the case when $T_1\cong T_2=\PSL_2(7)$ and $\Rad(G)=2^9$ or $2^{9+1}$, Lemma~\ref{Lem:Hopf1}(b) implies that $\pi(\overline{H}\cap\overline{K})\subseteq\{2\}$. Since $\overline{G}=\overline{H}\,\overline{K}$ with $\overline{H}$ and $\overline{K}$ perfect, interchanging $H$ and $K$ if necessary, $\overline{H}=T_1$ or $T_2$.
Without loss of generality, assume $\overline{H}=T_1$. Consider the irreducible $\overline{G}$-module $M=2^9$ (which is the tensor of the $\overline{H}$-module $U=2^3$ and the $T_2$-module $V=2^3$).
The $\overline{H}$-invariant subspace $H\cap M$ of $M$ has the form $U\otimes W$ for some subspace $W$ of $V$.
Since the stabilizer of $W$ in $T_2$ has order divisible by $3$, there exists a subgroup $P=3$ of $(\Rad(G)/M).T_2$ stabilizing $H\cap M$, and so the subgroup $(H\cap M).P=(H\cap M){:}3$ of $G$ normalizes $H$.
This implies that $\Nor_G(H)\geqslant H{:}3$, contradicting the condition that $H$ is self-normalizing.
\end{proof}

\subsection{Proof of Theorem~\ref{Thm:Hopf}}

Suppose for a contradiction that $(G,H,K)$ is not as described in Example~$\ref{Ex:Hopf}$. Lemmas~\ref{Lem:Hopf12} and~\ref{Lem:HopfTwo} imply that $k\geqslant3$.
Let $\,\overline{\phantom{x}}$ be the natural homomorphism $G\to G/\Rad(G)$.
Assume without loss of generality that $|T_1|\leqslant\cdots\leqslant|T_k|$.
Since $|T_1|\cdots|T_k|\leqslant|G|\leqslant41287680$, we only need to deal with the following cases.

\textsf{Case 1:} $k=3$ and $\A_5\notin\{T_1,T_2,T_3\}$.
We deduce from $|T_1||T_2||T_3|\leqslant41287680$ and $\Z(G)=1$ that one of the following occurs:
\begin{itemize}
\item $G=2^3.(\PSL_2(7)\times\PSL_2(7)\times\PSL_2(7))$;
\item $G=\PSL_2(7)\times\PSL_2(7)\times T_3$ with $T_3=\PSL_2(7)$, $\A_6$, $\PSL_2(8)$, $\PSL_2(11)$ or $\PSL_2(13)$;
\item $G=\PSL_2(7)\times\A_6\times T_3$ with $T_3=\A_6$, $\PSL_2(8)$ or $\PSL_2(11)$.
\end{itemize}
However, computation in \magma~\cite{BCP1997} shows that there exists no exact factorization $G=HK$ with both $H$ and $K$ perfect and self-normalizing, a contradiction.

\textsf{Case 2:} $k=3$ and $T_1=\A_5\notin\{T_2,T_3\}$.
In this case, as $|T_1||T_2||T_3|\leqslant41287680$ and $\Z(G)=1$, one of the following occurs:
\begin{enumerate}[(i)]
\item $\overline{G}=\A_5\times\PSL_2(7)\times\PSL_2(7)$, and $\Rad(G)=1$, $2^3$, $2^{3+1}$ or $2^4$;
\item $\overline{G}=\A_5\times\PSL_2(7)\times T_3$ with $T_3\in\{\A_6,\PSL_2(8)\}$, and $\Rad(G)=1$ or $2^3$;
\item $G=\A_5\times\PSL_2(7)\times T_3$, where $T_3=\A_7$ or $\PSL_2(q)$ with $11\leqslant q\leqslant19$;
\item $G=\A_5\times T_2\times T_3$, where $T_2\in\{\A_6,\PSL_2(8)\}$ and $T_3=\PSL_2(q)$ with $8\leqslant q\leqslant13$;
\item $G=\A_5\times\PSL_2(11)\times\PSL_2(11)$.
\end{enumerate}
Here we only show that (i) is not possible, while the others can be excluded similarly.
Suppose for a contradiction that (i) occurs.
Since $\overline{G}=\overline{H}\,\overline{K}$ with $\overline{H}$ and $\overline{K}$ perfect, interchanging $H$ and $K$ if necessary, either $\overline{H}\leqslant T_2\times T_3$, or $|\overline{H}\cap\overline{K}|$ is divisible by $2\cdot3$.
The latter is impossible as Lemma~\ref{Lem:Hopf1}(b) requires $\pi(\overline{H}\cap\overline{K})\subseteq\pi(\Rad(G))\subseteq\{2\}$.
Thus, $\overline{H}\leqslant T_2\times T_3$, and a similar argument as in the proof of Lemma~\ref{Lem:Hopf10} gives $\Nor_G(H)>H$, a contradiction.

\textsf{Case 3:} $k=3$ and $T_1=T_2=\A_5$.
It follows from $|T_1||T_2||T_3|\leqslant41287680$ and $\Z(G)=1$ that one of the following occurs:
\begin{itemize}
\item $\overline{G}=\A_5\times\A_5\times\A_5$, and either $\Rad(G)$ is a $2$-group or $\Rad(G)=3^4$;
\item $\overline{G}=\A_5\times\A_5\times T_3$, where $T_3\in\{\PSL_2(7),\A_6,\PSL_2(8),\PSL_2(11)\}$ and  $\Rad(G)$ is a $2$-group;
\item $G=\A_5\times\A_5\times T_3$, where $T_3=\A_7$, $\PSL_3(3)$, $\PSU_3(3)$, $\M_{11}$ or $\PSL_2(q)$ with $13\leqslant q\leqslant27$.
\end{itemize}
Then the same argument as in Case~2 yields a contradiction.

\textsf{Case 4:} $k=4$, $T_1=T_2=T_3=\A_5$ and $T_4\in\{\A_5,\PSL_2(7)\}$.
Since $|G|\leqslant41287680$ and $\Z(G)=1$, we obtain $G=\A_5\times\A_5\times\A_5\times T_4$.
However, computation in \magma~\cite{BCP1997} shows that there exists no exact factorization $G=HK$ with $H$ and $K$ perfect and self-normalizing, a contradiction. This completes the proof.

%%%%%%%%%%%%%%%%%%%%%%%%%%%%%%%%%%%%%%%%%%%%%%%%%%%%%%%%%%%%%%%%%%%
%%%%%%%%%%%%%%%%%%%%%%%%%%%%%%%%%%%%%%%%%%%%%%%%%%%%%%%%%%%%%%%%%%%

\section{Tables}

%%%%%%%%%%%%%%%%%%%%%%%%%%%%%%%%%%%%%%%%%%%%%%%%%%%%%%%%%%%%%%%%%%%
%%%%%%%%%%%%%%%%%%%%%%%%%%%%%%%%%%%%%%%%%%%%%%%%%%%%%%%%%%%%%%%%%%%

\begin{table}[htbp]
\captionsetup{justification=centering}
\caption{$(G_0,H_0,K_0)$ for linear groups ($c$ given in~\eqref{EqnLinear01})}\label{TabLinear}
\begin{tabular}{|c|l|l|l|l|l|}
\hline
  & $G_0$ & $H_0$ & $K_0$ & $H_0\cap K_0$ & Ref.\\
\hline
1 & $\SL_{ab}(q)$ & $\SL_a(q^b)$ & $q^{ab-1}{:}\SL_{ab-1}(q)$ & $q^{ab-b}{:}\SL_{a-1}(q^b)$ & \ref{ex:Linear} \\
  &  & $\Sp_a(q^b)'$ &  & $[q^c]{:}\Sp_{a-2}(q^b)$ & \ref{ex:Linear}\\
  &  & $\G_2(q^b)'$ ($a=6$, $q$ even) &  & $[(q^{5b},q^{6b}/4)]{:}\SL_2(q^b)$ & \ref{ex:Linear}\\
\hline
2 & $\SL_{2m}(q)$ & $\Sp_{2m}(q)'$ & $\SL_{2m-1}(q)$ &$\Sp_{2m-2}(q)'$ & \ref{ex:Sp<SL}\\
3 & $\SL_{2m}(2)$ & $\SiL_m(4)$ & $\SL_{2m-1}(2)$ & $\SL_{m-1}(4)$ & \ref{ex:SL_m<SL_2m}\\
  &  & $\GaSp_m(4)$ &  & $\Sp_{m-2}(4)$ & \ref{ex:Sp<SL-02}\\
4 & $\SL_{2m}(2){:}2$ & $\SL_m(4){:}2$ & $\SL_{2m-1}(2){:}2$ & $\SL_{m-1}(4)$ & \ref{ex:SL_m<SL_2m}\\
  &  & $\Sp_m(4){:}2$ &  & $\Sp_{m-2}(4)$& \ref{ex:Sp<SL-02}\\
5 & $\SL_{2m}(4){:}2$ & $\SL_m(16){:}4$ & $\SL_{2m-1}(4){:}2$ & $\SL_{m-1}(16)$ & \ref{ex:SL_m<SL_2m}\\
  &  & $\Sp_m(16){:}4$ &  & $\Sp_{m-2}(16)$ & \ref{ex:Sp<SL-02}\\
6 & $\SL_6(2^f)$ & $\G_2(2^f)$ & $\SL_5(2^f)$ & $\SL_2(2^f)$ & \ref{ex:G_2<Sp}\\
7 & $\SL_{12}(2)$ & $\GaG_2(4)$ & $\SL_{11}(2)$ &$\SL_2(4)$ & \ref{ex:Sp<SL-02}\\
8 & $\SL_{12}(2){:}2$ & $\G_2(4){:}2$ & $\SL_{11}(2){:}2$ &$\SL_2(4)$ & \ref{ex:Sp<SL-02}\\
9 & $\SL_{12}(4){:}2$ & $\G_2(16){:}4$ & $\SL_{11}(4){:}2$ &$\SL_2(16)$ & \ref{ex:Sp<SL-02}\\
\hline
10 & $\SL_4(2)$  & $\A_7$ & $2^3{:}\SL_3(2)$, $\SL_3(2)$ & $\PSL_2(7)$, $7{:}3$ & \ref{ex:Linear}, \ref{ex:Sp<SL}\\
11 & $\SL_4(3)$ & $2^{\boldsymbol{\cdot}}\Sy_5$, $8^{\boldsymbol{\cdot}}\A_5$, $2^{1+4}{}^{\boldsymbol{\cdot}}\A_5$ & $3^3{:}\SL_3(3)$ & $3$, $\Sy_3$, $\SL_2(3)$ & \ref{ex:Linear}\\
12 & $\SL_6(3)$ & $\SL_2(13)$ & $3^5{:}\SL_5(3)$ & $3$ & \ref{ex:Linear}\\
\hline
13 & $\SL_2(9)$ & $2^{\boldsymbol{\cdot}}\A_5$ & $\SL_2(5)$ & $2^{\boldsymbol{\cdot}}\D_{10}$ & \ref{ex:SL-exception1}\\
14 & $\SL_3(4){:}2$ & $3^{\boldsymbol{\cdot}}\M_{10}$ & $\SL_3(2){:}2$ & $\Sy_3$ & \ref{ex:SL-exception2}\\
\hline
\end{tabular}
\vspace{3mm}
\end{table}

\begin{table}[htbp]
\captionsetup{justification=centering}
\caption{$(G_0,H_0,K_0)$ for unitary groups with $K^{(\infty)}=(\N_1[G])^{(\infty)}$\\($c$ given in~\eqref{EqnUnitary7})}\label{TabUnitary}
\begin{tabular}{|c|l|l|l|l|l|l|}
\hline
  & $G_0$ & $H_0$ & $K_0$ & $H_0\cap K_0$ & Ref.\\
\hline
1 & $\SU_{2ab}(q)$ & $q^c{:}\SL_a(q^{2b})$ & $\SU_{2ab-1}(q)$ & $[q^{c-2b+1}].\SL_{a-1}(q^{2b})$ & \ref{LemUnitaryPm3} \\
  &  & $q^c{:}\Sp_a(q^{2b})$ &  & $[q^{c-2b+1}].\Sp_{a-2}(q^{2b})$ & \ref{LemUnitaryPm3} \\
  &  & $q^c{:}\G_2(q^{2b})$ ($a=6$, $q$ even) &  & $[q^{c-2b+1}].\SL_2(q^{2b})$ &  \ref{LemUnitaryPm3} \\
\hline
2 & $\SU_{2m}(q)$ & $\Sp_{2m}(q)$ & $\SU_{2m-1}(q)$ & $\Sp_{2m-2}(q)$ & \ref{LemUnitary09} \\
3 & $\SU_{2m}(2)$ & $\SL_m(4){:}2$ & $\SU_{2m-1}(2)$ &$\SL_{m-1}(4)$ & \ref{ex:Unitary-03} \\
  &  & $\Sp_m(4){:}2$ &  & $\Sp_{m-2}(4)$ & \ref{ex:Unitary-03} \\
4 & $\SU_{2m}(2){:}2$ & $\SL_m(4){:}2$ & $\SU_{2m-1}(2){:}2$ & $\SL_{m-1}(4)$ & \ref{ex:Unitary-03} \\
  &  & $\Sp_m(4){:}2$ & & $\Sp_{m-2}(4)$ & \ref{ex:Unitary-03} \\
5 & $\SU_{2m}(4){:}4$ & $\SL_m(16){:}4$ & $\SU_{2m-1}(4){:}4$ & $\SL_{m-1}(16)$ & \ref{ex:Unitary-03} \\
  &  & $\Sp_m(16){:}4$ & & $\Sp_{m-2}(16)$ &  \ref{ex:Unitary-03} \\
6 & $\SU_6(2^f)$ & $\G_2(2^f)$ & $\SU_5(2^f)$ & $\SL_2(2^f)$ & \ref{LemUnitary11} \\
\hline
7 & $\SU_6(2)$ & $3^{\boldsymbol{\cdot}}\PSU_4(3)$, $3^{\boldsymbol{\cdot}}\M_{22}$ & $\SU_5(2)$ & $3^5{:}\A_5$, $\PSL_2(11)$ & \ref{LemUnitary17} \\
8 & $\SU_{12}(2)$ & $\G_2(4){:}2$, $3^{\boldsymbol{\cdot}}\Suz$ & $\SU_{11}(2)$ & $\SL_2(4)$, $3^6{:}\PSL_2(11)$ & \ref{ex:Unitary-03}, \ref{LemUnitary19} \\
9 & $\SU_{12}(2){:}2$ & $\G_2(4){:}2$ & $\SU_{11}(2){:}2$ & $\SL_2(4)$ & \ref{ex:Unitary-03} \\
10 & $\SU_{12}(4){:}4$ & $\G_2(16){:}4$ & $\SU_{11}(4){:}4$ & $\SL_2(16)$ & \ref{ex:Unitary-03} \\
\hline
\end{tabular}
\vspace{3mm}
\end{table}

\begin{table}[htbp]
\captionsetup{justification=centering}
\caption{Exceptional $(G_0,H_0,K_0)$ for unitary groups\\(see Lemma~\ref{Lem:UnitaryProof} for remarks)}\label{TabUnitary-2}
\begin{tabular}{|l|l|l|l|l|l|}
\hline
 & $G_0$ & $H_0$ & $K_0$ & $H_0\cap K_0$\\
\hline
1 & $\PSU_4(3)$ & $\PSp_4(3)$ & $\PSL_3(4)$ & $2^4{:}\D_{10}$\\
2 & $\PSU_4(3){:}2$ & $3^4{:}(\A_5\times2)$ & $\PSL_3(4){:}2$ & $\A_5$\\
3 & $\PSU_4(3){:}2$ & $3^4{:}\A_6.2$ & $\PSL_3(4){:}2$ & $\A_6$\\
4 & $\PSU_4(3){:}\calO$ with $|\calO|=4$ & $3^4{:}\A_6.\calO$ & $\PSL_2(7).\calO$ & $\Sy_3$\\
5 & $\PSU_4(5){:}2$ & $5^4{:}(\PSL_2(25)\times2)^{\boldsymbol{\cdot}}4$ & $3^{\boldsymbol{\cdot}}\A_7{:}2$ & $\AGL_1(5)$\\
6 & $\PSU_9(2)$ & $\J_3$ & $2^{1+14}{:}\SU_7(2)$ & $2^{2+4}{:}(\Sy_3\times3)$\\
\hline
\end{tabular}
\vspace{3mm}
\end{table}

\begin{table}[htbp]
\captionsetup{justification=centering}
\caption{$(G_0,H_0,K_0)$ for orthogonal groups in odd dimension}\label{TabOmega}
\begin{tabular}{|c|l|l|l|l|l|}
\hline
  & $G_0$ & $H_0$ & $K_0$ & $H_0\cap K_0$ & Ref.\\
\hline
1 & $\Omega_{2ab+1}(q)$ & $(q^{ab(ab-1)/2}.q^{ab}){:}\SL_a(q^b)$ & $\Omega_{2ab}^-(q)$ & $[q^{(a^2b^2+ab-2b)/2}].\SL_{a-1}(q^b)$ & \ref{Ex:OmegaOdd-1} \\
  &  & $(q^{ab(ab-1)/2}.q^{ab}){:}\Sp_a(q^b)$ & &$[q^{(a^2b^2+ab-2b)/2}].\Sp_{a-2}(q^b)$ & \ref{Ex:OmegaOdd-1} \\
2 & $\Omega_7(q)$ & $\G_2(q)$, $\SL_3(q)$ ($q=3^f$) & $\Omega_6^-(q)$ & $\SU_3(q)$, $q^2-1$ & \ref{LemOmegaRow1,3,4} \\
3 & $\Omega_7(q)$ & $\G_2(q)$, $\SU_3(q)$ ($q=3^f$) & $\Omega_6^+(q)$ & $\SL_3(q)$, $q^2-1$ & \ref{LemOmega07} \\
  &  & ${^2}\G_2(q)$ ($q=3^f$, $f$ odd) &  & $\tfrac{q-1}{2}.2\phantom{{\tfrac{q^2-1}{2}}_f}$ & \ref{LemOmega07}\\
4 & $\Omega_7(q)$ & $q^5{:}\Omega_5(q)$, $q^4{:}\Omega_4^-(q)$ & $\G_2(q)$ & $[q^5]{:}\SL_2(q)$, $[q^3]$ & \ref{Lem:Omega7-01} \\
  &  & $\Omega_5(q)$ &  & $\SL_2(q)$ & \ref{LemOmega07} \\
5 & $\Omega_{13}(3^f)$ & $\PSp_6(3^f)$ & $\Omega_{12}^-(3^f)$ &$(\SL_2(3^f)\times\SL_2(3^{2f}))/2$ & \ref{LemOmega09} \\
6 & $\Omega_{25}(3^f)$ & $\F_4(3^f)$ & $\Omega_{24}^-(3^f)$ & $2.\Omega_8^-(3^f)$ & \ref{LemOmega10} \\
\hline
7 & $\Omega_7(3)$ & $3^4{:}\Sy_5$, $3^5{:}2^4{:}\A_5$ & $\G_2(3)$ &$3^2$, $\ASL_2(3)$ & \ref{LemXia28} \\
8 & $\Omega_7(3)$ & $3^3{:}\SL_3(3)$, $\Omega_6^+(3)$, $\G_2(3)$ & $\A_9$ & $\Sy_3$, $\Sy_5\times2$, $\PSL_2(7)$ & \ref{LemXia28} \\
9 & $\Omega_7(3)$ & $3^3{:}\SL_3(3)$, $\Omega_6^+(3)$, $\G_2(3)$ & $\Sp_6(2)$ & $\GL_2(3)$, $2^4{:}\Sy_5$, $2^3{}^{\boldsymbol{\cdot}}\PSL_2(7)$ & \ref{LemXia28} \\
10 & $\Omega_7(3)$ & $2^6{:}\A_7$, $\Sy_8$, $\A_9$ & $3^{3+3}{:}\SL_3(3)$ & $6^{\boldsymbol{\cdot}}\Sy_4$, $\Sy_3\times\Sy_3$, $3^3{:}\Sy_3$ & \ref{LemXia28} \\
  &  & $2^{\boldsymbol{\cdot}}\PSL_3(4)$, $\Sp_6(2)$ & & $3^2{:}4$, $\SU_3(2){:}\Sy_3$ & \ref{LemXia28} \\
\hline
11 & $\Omega_9(3)$ & $3^{6+4}{:}2.\Sy_5$, $3^{6+4}{:}8.\A_5$ & $\Omega_8^-(3)$ &$3^{3+3}.3$, $3^{3+3}.\Sy_3$ & \ref{Ex:OmegaOdd-1} \\
  &  & $3^{6+4}{:}2^{1+4}.\A_5$ &  & $3^{3+3}.\SL_2(3)$ & \ref{Ex:OmegaOdd-1} \\
12 & $\Omega_{13}(3)$ & $3^{15+6}{:}\SL_2(13)$ & $\Omega_{12}^-(3)$ &$3^{10+5}.3$ & \ref{Ex:OmegaOdd-1} \\
\hline
\end{tabular}
\end{table}

\begin{table}[htbp]
\caption{$(G_0,H_0,K_0)$ for orthogonal groups of minus type}\label{TabOmegaMinus}
\begin{tabular}{|c|l|l|l|l|l|l|}
\hline
  & $G_0$ & $H_0$ & $K_0$ & $H_0\cap K_0$ & Ref.\\
\hline
1 & $\Omega_{2m}^-(q)$ ($m$ odd) & $\SU_m(q)$ &$q^{2m-2}{:}\Omega_{2m-2}^-(q)$ & $q^{1+(2m-4)}{:}\SU_{m-2}(q)$ & \ref{ex:OmegaMinus-1}\\
2 & $\Omega_{10}^-(2)$ & $\A_{12}$, $\M_{12}$&$2^8{:}\Omega_8^-(2)$ & $(\A_8\times\A_4){:}2$, $2_+^{1+4}{:}\Sy_3$ & \ref{prop:OmegaMinus-1}\\
3 & $\Omega_{18}^-(2)$ & $3^{\boldsymbol{\cdot}}\J_3$ & $2^{16}{:}\Omega_{16}^-(2)$ & $2^{2+4}{:}(\Sy_3\times3)$ & \ref{prop:OmegaMinus-1}\\
\hline
4 & $\Omega_{2m}^-(q)$ ($m$ odd) & $\SU_m(q)$ & $\Omega_{2m-1}(q)$ & $\SU_{m-1}(q)$ & \ref{ex:OmegaMinus01}\\
5 & $\GO_{2m}^-(2)$ ($m/2$ odd) & $\SU_{m/2}(4){:}4$ & $\Sp_{2m-2}(2)\times2$ & $\SU_{m/2-1}(4){:}2$ & \ref{ex:OmegaMinus05}\\
6 & $\GaO_{2m}^-(4)$ ($m/2$ odd) & $\SU_{m/2}(16){:}8$ & $\Sp_{2m-2}(4){:}4$ & $\SU_{m/2-1}(16){:}2$ & \ref{ex:OmegaMinus05}\\
7 & $\GO_{2m}^-(2)$ ($m$ even) & $\GaO_m^-(4)$ & $\Sp_{2m-2}(2)\times2$ & $\Sp_{m-2}(4){:}2$ & \ref{ex:OmegaMinus02}\\
8 & $\GaO_{2m}^-(4)$ ($m$ even) & $\GaO_m^-(16)$ & $\Sp_{2m-2}(4){:}4$ & $\Sp_{m-2}(16){:}2$ & \ref{ex:OmegaMinus02}\\
\hline
9 & $\Omega_{2m}^-(2)$ ($m$ odd) & $\SU_m(2)$ & $\Omega^-_{2m-2}(2){:}2$ & $\SU_{m-2}(2)$ & \ref{LemOmegaMinus03}\\
10 & $\GO_{2m}^-(2)$ ($m$ odd) & $\SU_m(2){:}2$ & $\Omega^-_{2m-2}(2){:}2$ & $\SU_{m-2}(2)$ & \ref{LemOmegaMinus04}\\
11 & $\GaO_{2m}^-(4)$ ($m$ odd) & $\SU_m(4){:}4$ & $\Omega^-_{2m-2}(4){:}4$ & $\SU_{m-2}(4)$ & \ref{LemOmegaMinus04}\\
\hline
\end{tabular}
\end{table}

\begin{table}[htbp]
\captionsetup{justification=centering}
\caption{$(G_0,H_0,K_0)$ for orthogonal groups of plus type\\($c$ given in~\eqref{EqnOmegaPlus7}, $d=(2,m-1)$)}\label{TabOmegaPlus}
\small
\begin{tabular}{|l|l|l|l|l|l|}
\hline
 & $G_0$ & $H_0$ & $K_0$ & $H_0\cap K_0$ & Ref.\\
\hline
1 & $\Omega_{2ab}^+(q)$ & $q^c{:}\SL_a(q^b)$ & $\Omega_{2ab-1}(q)$ & $[q^{c-b+1}].\SL_{a-1}(q^b)$ & \ref{LemOmegaPlusPm3} \\
  &  & $q^c{:}\Sp_a(q^b)$ &  & $[q^{c-b+1}].\Sp_{a-2}(q^b)$ & \ref{LemOmegaPlusPm3} \\
  &  & $q^c{:}\G_2(q^b)'$ ($a=6$, $q$ even) &  & $[(q^{c-b+1},q^{c+1}/4)].\SL_2(q^b)$ & \ref{LemOmegaPlusPm3} \\
2 & $\Omega_{2ab}^+(2)$ & $2^c{:}\SL_a(2^b){:}b_2$ & $\Sp_{2ab-2}(2)$ & $[2^{c-b+2}].\SL_{a-1}(2^b).(b_2/2)$ & \ref{LemOmegaPlusPm4} \\
  & ($b$ even) & $2^c{:}\Sp_a(2^b){:}b_2$ &  & $[2^{c-b+2}].\Sp_{a-2}(2^b).(b_2/2)$ & \ref{LemOmegaPlusPm4} \\
  &  & $2^c{:}\G_2(2^b){:}b_2$ ($a=6$) &  & $[2^{c-b+2}].\SL_2(2^b).(b_2/2)$ & \ref{LemOmegaPlusPm4} \\
3 & $\Omega_{2ab}^+(4){:}2$ & $2^{2c}{:}\SL_a(4^b){:}2b_2$ & $\GaSp_{2ab-2}(4)$ & $[4^{c-b+2}].\SL_{a-1}(4^b).(b_2/2)$ & \ref{LemOmegaPlusPm4} \\
  & ($b$ even) & $2^{2c}{:}\Sp_a(4^b){:}2b_2$ &  & $[4^{c-b+2}].\Sp_{a-2}(4^b).(b_2/2)$ & \ref{LemOmegaPlusPm4} \\
  &  & $2^{2c}{:}\G_2(4^b){:}2b_2$ ($a=6$) &  & $[4^{c-b+2}].\SL_2(4^b).(b_2/2)$ & \ref{LemOmegaPlusPm4} \\
4 & $\Omega_{12}^+(3)$ & $3^{14}{:}\SL_2(13)$ & $\Omega_{11}(3)$ & $3^{9+1}$ & \ref{PropOmegaPlusPm} \\
\hline
5 & $\Omega_{2m}^+(q)$ & $\SL_m(q)$, $\Sp_m(q)$ & $\Omega_{2m-1}(q)$ & $\SL_{m-1}(q)$, $\Sp_{m-2}(q)$ & \ref{Prop:O^+=(SL,N_1)}, \ref{Prop:O^+=(tensor,N_1)} \\
 & & $\G_2(q)$ ($m=6$, $q$ even) & & $\SL_2(q)$ & \ref{Prop:O^+=(SL,N_1)}, \ref{Prop:O^+=(tensor,N_1)}\\
 & & $\SU_m(q)$ ($m$ even) & & $\SU_{m-1}(q)$ & \ref{Prop:O^+=(SU,N_1)}\\
 & & $\mathrm{Spin}_9(q)$ ($m=8$) & & $\mathrm{Spin}_7(q)$ & \ref{C_9-subgroups}\\
6 &$\Omega^+_{2m}(2)$ & $\Omega_m^+(4){:}2$, $\SL_{m/2}(4){:}2$ & $\Sp_{2m-2}(2)$ & $\Sp_{m-2}(4)$, $\SL_{m/2-1}(4)$ & \ref{LemOmegaPlus12}, \ref{Prop:O^+=(SL,N_1)}\\
 & & $\GaSp_{m/2}(4)$, $\Sp_{m/2}(4).4$ & & $\Sp_{m/2-2}(4)$, $\Sp_{m/2-2}(4).2$ & \ref{Prop:O^+=(SL,N_1)}, \ref{LemOmegaPlus35}\\
 & & $\GaG_2(4)$ ($m=12$) &  & $\SL_2(4)$ & \ref{Prop:O^+=(SL,N_1)}\\
 & & $\SU_{m/2}(4){:}4$ ($m/2$ even) & & $\SU_{m/2-1}(4).2$ & \ref{LemOmegaPlus35}\\
7 & $\Omega_{2m}^+(4){:}2$ & $\Omega_m^+(16){:}4$, $\SL_{m/2}(16){:}4$ & $\GaSp_{2m-2}(4)$ & $\Sp_{m-2}(16)$, $\SL_{m/2-1}(16)$ & \ref{LemOmegaPlus12}, \ref{Prop:O^+=(SL,N_1)}\\
 & & $\GaSp_{m/2}(16)$ & & $\Sp_{m/2-2}(16)$ & \ref{Prop:O^+=(SL,N_1)}, \ref{Prop:O^+=(tensor,N_1)}\\
 & & $\Sp_{m/2}(16).8$ & & $\Sp_{m/2-2}(16).2$ & \ref{LemOmegaPlus35}\\
 & & $\GaG_2(16)$ ($m=12$) & & $\SL_2(16)$ & \ref{Prop:O^+=(SL,N_1)}, \ref{Prop:O^+=(tensor,N_1)}\\
 & & $\SU_{m/2}(16){:}8$ ($m/2$ even) & & $\SU_{m/2-1}(16).2$ & \ref{LemOmegaPlus35}\\
\hline
8 & $\Omega_{12}^+(2)$ & $3^{\boldsymbol{\cdot}}\PSU_4(3)$, $3^{\boldsymbol{\cdot}}\M_{22}$ & $\Sp_{10}(2)$ &$3^5{:}\A_5$, $\PSL_2(11)$ & \ref{Prop:O^+=(SU,N_1)}\\
9 & $\Omega_{16}^+(2)$ & $\Omega_8^-(2).2$, $\GaSp_6(4)$ & $\Sp_{14}(2)$ & $\G_2(2)$, $\G_2(4)$, & \ref{LemOmegaPlus35}\\
 & & $\Omega_6^+(4){:}2$, $\Omega_6^-(4){:}2$ & & $\SL_3(4)$, $\SU_3(4)$, & \ref{LemOmegaPlus35}\\
 & & $\GaSp_4(4)$ & & $\SL_2(4)$ & \ref{C_9-subgroups}, \ref{LemOmegaPlus35}\\
10 & $\Omega_{16}^+(4){:}2$ & $\Omega_8^-(4).4$, $\GaSp_6(16)$ & $\GaSp_{14}(4)$ & $\G_2(4)$, $\G_2(16)$, & \ref{LemOmegaPlus35}\\
 & & $\Omega_6^+(16){:}4$, $\Omega_6^-(16){:}4$ & & $\SL_3(16)$, $\SU_3(16)$ & \ref{LemOmegaPlus35}\\
 & & $\GaSp_4(16)$ & & $\SL_2(16)$ & \ref{C_9-subgroups}, \ref{LemOmegaPlus35}\\
11 & $\Omega_{24}^+(2)$ & $\G_2(4){:}2$, $\G_2(4).4$, & $\Sp_{22}(2)$ & $\SL_2(4)$, $\SL_2(4).2$, & \ref{C_9-subgroups}, \ref{LemOmegaPlus35}\\
 &  & $3^{\boldsymbol{\cdot}}\Suz$, $\Co_1$ &  & $3^5{:}\PSL_2(11)$, $\Co_3$ & \ref{Prop:O^+=(SU,N_1)}, \ref{C_9-subgroups}\\
12 & $\Omega_{24}^+(4){:}2$ & $\GaG_2(16)$, $\G_2(16).8$ & $\GaSp_{22}(4)$ & $\SL_2(16)$, $\SL_2(16).2$ & \ref{Prop:O^+=(tensor,N_1)}, \ref{LemOmegaPlus35}\\
13 & $\Omega_{32}^+(2)$ & $\GaSp_8(4)$ & $\Sp_{30}(2)$ & $\Sp_6(4)$ & \ref{LemOmegaPlus35}\\
14 & $\Omega_{32}^+(4){:}2$ & $\GaSp_8(16)$ & $\GaSp_{30}(4)$ & $\Sp_6(16)$ & \ref{LemOmegaPlus35}\\
\hline
15 & $\Omega_{2m}^+(q)$ & $\SU_m(q)$ ($m$ even) & $q^{2m-2}{:}\Omega_{2m-2}^+(q)$ & $q^{1+(2m-4)}{:}\SU_{m-2}(q)$ & \ref{Ex:OmegaPlus17} \\
16 & $\Omega_{2m}^+(2)$ & $\SU_m(2){:}2$ ($m$ even) & $\Omega_{2m-2}^+(2)$ & $\SU_{m-2}(2)$ & \ref{LemOmegaPlus18} \\
17 & $\Omega_{2m}^+(2)$ & $\SU_m(2)$ ($m$ even) & $\Omega_{2m-2}^+(2){:}2$ & $\SU_{m-2}(2)$ & \ref{LemOmegaPlus19} \\
18 & $\Omega_{2m}^+(4){:}2$ & $\SU_m(4).4$ ($m$ even) & $\Omega_{2m-2}^+(4){:}2$ & $\SU_{m-2}(4)$ & \ref{LemOmegaPlus18}, \ref{LemOmegaPlus22} \\
\hline
19 & $\Omega_{2m}^+(q)$ & $q^{m(m-1)/2}{:}\SL_m(q)$ & $\Omega_{2m-2}^-(q)$ & $[q^{(m+1)(m-2)/2}]{:}\SL_{m-2}(q)$ & \ref{LemOmegaPlus13}\\
20 & $\Omega_{2m}^+(2)$ & $\SL_m(2)$ & $\Omega_{2m-2}^-(2).2$ & $\SL_{m-2}(2)$ & \ref{LemOmegaPlus15}\\
21 & $\Omega_{2m}^+(2){:}d$ & $\SL_m(2){:}2$ & $\Omega_{2m-2}^-(2).d$ & $\SL_{m-2}(2)$ & \ref{LemOmegaPlus14} \\
22 & $\Omega_{2m}^+(4){:}2$ & $\SL_m(4){:}2$ & $\Omega_{2m-2}^-(4).4$ & $\SL_{m-2}(4)$ & \ref{LemOmegaPlus16} \\
23 & $\Omega_8^+(4){:}2$ & $\Omega_8^-(2).2$ & $\Omega_6^-(4).4$ & $\SL_2(2)\times2$ & \ref{LemOmegaPlus23} \\
24 & $\Omega_8^+(16){:}4$ & $\Omega_8^-(4).4$ & $\Omega_6^-(16).8$ & $\SL_2(4)\times2$ & \ref{LemOmegaPlus23} \\
\hline
\end{tabular}
\vspace{3mm}
\end{table}

\begin{table}[htbp]
\captionsetup{justification=centering}
\caption{Exceptional $(G_0,H_0,K_0)$ for orthogonal groups of plus type\\(see Lemmas~\ref{PropOmegaPlusO+8} and~\ref{LemXia29} for remarks)}\label{TabOmegaPlus2}
%\small
\begin{tabular}{|l|l|l|l|l|l|}
\hline
 & $G_0$ & $H_0$ & $K_0$ & $H_0\cap K_0$\\
\hline
1 & $\POm_8^+(q)$ & $\Omega_7(q)$, $\Omega_6^+(q)$, $\Omega_6^-(q)$ & $\Omega_7(q)$ & $\G_2(q)$, $\SL_3(q)$, $\SU_3(q)$ \\
 & & $q^5{:}\Omega_5(q)$, $\Omega_5(q)$ & & $[q^5]{:}\SL_2(q)$, $\SL_2(q)$ \\
 & & $q^4{:}\Omega^-_4(q)$, $\Omega_8^-(q^{1/2})$ & & $[q^3]$, $\G_2(q^{1/2})$ \\
\hline
2 & $\Omega_8^+(2)$ & $\Sy_5$, $\A_5{:}4$, $2^4{:}\A_5$, $\A_6$, $2^5{:}\A_6$ & $\Sp_6(2)$ & $1$, $2$, $\Q_8$, $3$, $4^2{:}3{:}2$ \\
 & & $\A_7$, $2^6{:}\A_7$, $\A_8$ & & $7{:}3$, $2^3{}^{\boldsymbol{\cdot}}\SL_3(2)$, $\SL_3(2)$ \\
 & & $\A_8$, $\A_9$ & & $\AGaL_1(8)$, $\PGaL_2(8)$ \\
3 & $\Omega_8^+(2)$ & $2^6{:}\A_7$, $\A_8$, $\Sy_8$, $\A_9$ & $\SU_4(2)$ & $\SL_2(3)$, $3$, $\Sy_3$, $9{:}3$  \\
4 & $\Omega_8^+(2)$ & $\A_8$ & $\SU_4(2).2$ & $\Sy_3$ \\
5 & $\Omega_8^+(2)$ & $2^4{:}\A_5$, $2^5{:}\A_6$, $2^6{:}\A_7$, $\A_8$, $2^6{:}\A_8$ & $\A_9$ &$1$, $\A_4$,
$\SL_3(2)$, $7{:}3$, $\AGL_3(2)$ \\
\hline
6 & $\POm_8^+(3)$ & $3^4{:}\Sy_5$, $3^4{:}(\A_5\times4)$ & $\Omega_7(3)$ & $3^2$, $\Sy_3\times3$ \\
 & & $(3^5{:}2^4){:}\A_5$, $\A_9$ & & $\AGL_3(2)$, $\SL_3(2)$ \\
 & & $\SU_4(2)$, $\Sp_6(2)$, $\Omega_8^+(2)$ & & $\SL_2(3)$, $2^3{}^{\boldsymbol{\cdot}}\SL_3(2)$, $2^6{:}\A_7$ \\
7 & $\POm_8^+(3)$ & $2^6{:}\A_7$, $\A_8$ & $3^6{:}\PSL_4(3)$ & $(\SL_2(3)\times3){:}2$, $\Sy_3\times3$ \\
 & & $2^6{:}\A_8$, $\A_9$ & & $(2^3{:}\Sy_4){:}\Sy_3$, $3^3{:}\Sy_3$ \\
 & & $2^{\boldsymbol{\cdot}}\PSL_3(4)$, $\Sp_6(2)$, $\Omega_6^-(3)$ & & $3^2{:}4$, $3^{\boldsymbol{\cdot}}\AGL_2(3)$, $3^{3+2}{:}\SL_2(3)$ \\
8 & $\POm_8^+(3)$ & $3^6{:}\SL_3(3)$, $3^{3+3}{:}\SL_3(3)$ & $\Omega_8^+(2)$ & $(\SL_2(3)\times3){:}2$, $(\SL_2(3)\times3){:}2$ \\
 & & $3^{6+3}{:}\SL_3(3)$, $3^6{:}\PSL_4(3)$ & & $3^2{}^{\boldsymbol{\cdot}}\AGL_2(3)$, $(\PSp_4(3)\times3){:}2$ \\
\hline
9 & $\Omega_8^+(4){:}2$ & $\SiL_2(16)$ & $\GaSp_6(4)$ & $1$ \\
\hline
\end{tabular}
\vspace{3mm}
\end{table}

\begin{table}[htbp]
\captionsetup{justification=centering}
\caption{Infinite families of $(G_0,H_0,K_0)$ for symplectic groups\\($c$ given in~\eqref{EqnSymplectic2}, $d$ given in~\eqref{EqnSymplectic13})}\label{TabSymplectic1}
\small
\begin{tabular}{|l|l|l|l|l|l|l|}
\hline
 & $G_0$ & $H_0$ & $K_0$ &$H_0\cap K_0$& Ref.\\
\hline
1 & $\Sp_{2ab}(q)$ & $\Sp_{2a}(q^b)$ & $\Pa_1[G_0]^{(\infty)}$ & $[q^d]{:}\Sp_{2a-2}(q^b)$ & \ref{ex:K<P1<Sp} \\
 & & $\G_2(q^b)'$ ($a=3$, $q$ even) &  & $[(q^{5b},q^{6b}/4)]{:}\SL_2(q^b)'$ & \ref{ex:K<P1<Sp}\\
\hline
2 & $\Sp_{2ab}(q)$  & $\Sp_{2a}(q^b)$ & $\Omega_{2ab}^+(q){:}(2,b)$ &$\Omega^+_{2a}(q^b){:}(2,b)$ & \ref{LemSymplectic04}\\
 & ($q$ even) & $\G_2(q^b)$ ($a=3$) & & $\SL_3(q^b).(2,b)$ & \ref{LemSymplectic49}\\
3 & $\Sp_{2m}(2^f)$ & $\Sz(2^{fm/2})$ ($f$ odd, $m/2$ odd) & $\mathrm{O}_{2m}^+(2^f)$ & $\D_{2(2^{fm/2}-1)}$ & \ref{LemSymplectic21}\\
4 & $\Sp_{2m}(2)$ & $\mathrm{O}_{2m}^-(2)$, $\SU_m(2){:}2$ ($m$ odd) & $\Omega_{2m}^+(2)$ & $\Sp_{2m-2}(2)$, $\SU_{m-1}(2)$ & \ref{prop:Sp(2)=O^-O^+}\\
5 & $\Sp_{2m}(2)$ & $\Omega_{2m}^-(2)$, $\SU_m(2)$ ($m$ odd) & $\GO_{2m}^+(2)$ & $\Sp_{2m-2}(2)$, $\SU_{m-1}(2)$ & \ref{prop:Sp(2)=O^-O^+}\\
 & & $\GaO_m^-(4)$, $\SU_{m/2}(4).4$ ($m/2$ odd) & & $\Sp_{m-2}(4){:}2$, $\SU_{m/2-1}(4){:}2$ & \ref{prop:Sp(2)=O^-O^+}\\
6 & $\GaSp_{2m}(4)$ & $\GaO_{2m}^-(4)$, $\SU_m(4).4$ ($m$ odd) & $\Omega_{2m}^+(4){:}2$ & $\Sp_{2m-2}(4)$, $\SU_{m-1}(4)$ &  \ref{prop:Sp(2)=O^-O^+}\\
\hline
7 & $\Sp_{2ab}(q)$ & $q^c{:}\SL_a(q^b){:}b_2$ ($b$ even) & $\Omega_{2ab}^-(q)$ & $[q^{c-b}].\SL_{a-1}(q^b).b_2$ & \ref{LemSymplecticPm4} \\
 & ($q$ even) & $q^c{:}\Sp_a(q^b){:}b_2$ ($b$ even) &  & $[q^{c-b}].\Sp_{a-2}(q^b).b_2$ & \ref{LemSymplecticPm4} \\
 & & $q^c{:}\G_2(q^b){:}b_2$ ($a=6$, $b$ even) &  & $[q^{c-b}].\SL_2(q^b).b_2$ & \ref{LemSymplecticPm4} \\
8 & $\Sp_{2ab}(q)$ & $q^c{:}\SL_a(q^b)$ & $\Omega_{2ab}^-(q){:}(2,b)$ & $[q^{c-b}].\SL_{a-1}(q^b).(2,b)$ & \ref{LemSymplecticPm5} \\
 & ($q$ even) & $q^c{:}\Sp_a(q^b)$ &  & $[q^{c-b}].\Sp_{a-2}(q^b).(2,b)$ & \ref{LemSymplecticPm5} \\
 & & $q^c{:}\G_2(q^b)'$ ($a=6$) &  & $[(q^{c-b},q^c/4)].\SL_2(q^b).(2,b)$ & \ref{LemSymplecticPm5} \\
9 & $\Sp_{2ab}(2)$ & $2^c{:}\SL_a(2^b){:}b_2$ ($b$ even) & $\mathrm{O}_{2ab}^-(2)$ & $[2^{c-b+2}].\SL_{a-1}(2^b).(b_2/2)$&\ref{LemSymplecticPm7}\\
 & & $2^c{:}\Sp_a(2^b){:}b_2$ ($b$ even) &  & $[2^{c-b+2}].\Sp_{a-2}(2^b).(b_2/2)$ & \ref{LemSymplecticPm7} \\
 & & $2^c{:}\G_2(2^b){:}b_2$ ($a=6$, $b$ even) &  & $[2^{c-b+2}].\SL_2(2^b).(b_2/2)$ & \ref{LemSymplecticPm7} \\
10 & $\Sp_{2ab}(2)$ & $2^c{:}\SL_a(2^b)$ & $\mathrm{O}_{2ab}^-(2)$ & $[2^{c-b+1}].\SL_{a-1}(2^b)$ & \ref{LemSymplecticPm6} \\
 & & $2^c{:}\Sp_a(2^b)$ &  & $[2^{c-b+1}].\Sp_{a-2}(2^b)$ & \ref{LemSymplecticPm6} \\
 & & $2^c{:}\G_2(2^b)'$ ($a=6$) &  & $[(2^{c-b+1},2^{c-1})].\SL_2(2^b)$ & \ref{LemSymplecticPm6} \\
11 & $\GaSp_{2ab}(4)$ & $2^{2c}{:}\SL_a(4^b){:}2b_2$ ($b$ even) & $\GaO_{2ab}^-(4)$ & $[4^{c-b+1}].\SL_{a-1}(4^b).b_2$ & \ref{LemSymplecticPm7} \\
 & & $2^{2c}{:}\Sp_a(4^b){:}2b_2$ ($b$ even) &  & $[4^{c-b+1}].\Sp_{a-2}(4^b).b_2$ & \ref{LemSymplecticPm7} \\
 & & $2^{2c}{:}\G_2(4^b){:}2b_2$ ($a=6$, $b$ even) &  & $[4^{c-b+1}].\SL_2(4^b).b_2$ & \ref{LemSymplecticPm7} \\
\hline
12 & $\Sp_{2ab}(q)$ & $\Sp_{2a}(q^b){:}b_2$ ($b$ even) & $\Omega^-_{2ab}(q)$ & $\mathrm{O}^-_{2a}(q^b).(b_2/2)$ & \ref{LemSymplectic05}\\
 & ($q$ even) & $\G_2(q^b){:}b_2$ ($a=3$, $b$ even) & & $\SU_3(q^b).b_2$ & \ref{LemSymplectic50}\\
13 & $\Sp_{2m}(2^f)$ & $\Sp_4(2^{fm/4}){:}2$ ($f$ odd, $m/4$ odd) & $\Omega_{2m}^-(2^f)$ & $\D_{2(2^{fm/2}-1)}$ & \ref{LemSymplectic42}\\
14 & $\Sp_{2ab}(q)$ & $\Sp_{2a}(q^b)$ & $\Omega_{2ab}^-(q){:}(2,b)$ & $\Omega^-_{2a}(q^b){:}(2,b)$ & \ref{LemSymplectic04}\\
 & ($q$ even) &$\G_2(q^b)$ ($a=3$) &  & $\SU_3(q^b).(2,b)$ & \ref{LemSymplectic49}\\
15 & $\Sp_{2m}(2)$ & $\SL_m(2){:}2$ ($m$ odd) & $\Omega_{2m}^-(2)$ & $\SL_{m-1}(2)$ & \ref{prop:Sp(2)=O^+O^-}\\
16 & $\Sp_{2m}(2)$ & $\SL_m(2)$, $\Sp_m(2)$ & $\GO_{2m}^-(2)$ & $\SL_{m-1}(2)$, $\Sp_{m-2}(2)$ & \ref{prop:Sp(2)=O^+O^-} \\
 & & $\SU_m(2)$ ($m$ even), $\Omega_m^+(4){:}2$ &  & $\SU_{m-1}(2)$, $\Sp_{m-2}(4)$ & \ref{prop:Sp(2)=O^+O^-}\\
 & & $\Sp_m(2)\times2$ &  & $\Sp_{m-2}(2)\times2$ & \ref{LemSymplectic18}\\
 & & $\SL_{m/2}(4){:}2$, $\GaSp_{m/2}(4)$ &  & $\SL_{m/2-1}(4)$, $\Sp_{m/2-2}(4)$ & \ref{prop:Sp(2)=O^+O^-} \\
 & & $\Sp_{m/2}(4).4$ &  & $\Sp_{m/2-2}(4){:}2$ & \ref{prop:Sp(2)=O^+O^-}\\
 & & $\Sp_{m/2}(4){:}2^2$ &  & $\Sp_{m/2-2}(4)\times2$ & \ref{LemSymplectic15}\\
 & & $\SU_{m/2}(4){:}4$ ($m/2$ even) &  & $\SU_{m/2-1}(4){:}2$ & \ref{prop:Sp(2)=O^+O^-} \\
17 & $\GaSp_{2m}(4)$ & $\Sp_{2m}(2)\times2$ & $\GaO_{2m}^-(4)$ & $\Sp_{2m-2}(2)\times2$ & \ref{LemSymplectic20} \\
 & & $\SL_m(4){:}2$, $\Sp_m(4){:}2$ &  & $\SL_{m-1}(4)$, $\Sp_{m-2}(4)$ & \ref{prop:Sp(2)=O^+O^-}\\
 & & $\Sp_m(4){:}2^2$ & & $\Sp_{m-2}(4)\times2$ & \ref{LemSymplectic18} \\
 & & $\SU_m(4){:}2$ ($m$ even) &  & $\SU_{m-1}(4)$ & \ref{prop:Sp(2)=O^+O^-} \\
18 & $\GaSp_{2m}(16)$ & $\Sp_{2m}(4){:}4$ & $\GaO_{2m}^-(16)$ & $\Sp_{2m-2}(4)\times2$ & \ref{LemSymplectic20} \\
\hline
19 & $\Sp_{2m}(2)$ & $\GaSp_m(4)$ & $\Sp_{2m-2}(2)$ & $\Sp_{m-2}(4)$ & \ref{LemSymplectic13}\\
 & & & $\Omega_{2m-1}(2)\times2$ & $\Sp_{m-2}(4)\times2$ & \ref{LemSymplectic19}\\
20 & $\GaSp_{2m}(4)$ & $\GaSp_m(16)$ & $\Sp_{2m-2}(4){:}2$ & $\Sp_{m-2}(16)$ & \ref{LemSymplectic13}\\
 & & & $\Omega_{2m-1}(4).2^2$ & $\Sp_{m-2}(16)\times2$ & \ref{LemSymplectic19}\\
 & & & $\Omega_{2m-1}(4).4$ & $\Sp_{m-2}(16){:}2$ & \ref{LemSymplectic19}\\
\hline
21 & $\Sp_6(2^f)$ & $\Sp_4(2^f)$, $2^{4f}{:}\Omega_4^-(2^f)$ & $\G_2(2^f)$ & $\SL_2(2^f)$, $[2^{3f}]$ & \ref{LemSymplectic12}\\
 & ($f\geqslant2$) & $2^{5f}{:}\Sp_4(2^f)$ & & $[2^{5f}]{:}\SL_2(2^f)$ & \ref{LemSymplectic12}\\
\hline
\end{tabular}
\vspace{3mm}
\end{table}

\begin{table}[htbp]
\captionsetup{justification=centering}
\caption{Sporadic $(G_0,H_0,K_0)$ for symplectic groups}\label{TabSymplectic2}
\begin{tabular}{|l|l|l|l|l|l|}
\hline
 & $G_0$ & $H_0$ & $K_0$ & $H_0\cap K_0$ & Ref.\\
\hline
1 & $\GaSp_4(4)$ & $\Sy_6$, $\A_6\times2$ & $\GaO_4^-(4)$ & $\Sy_3$, $\Sy_3$ & \ref{LemSymplectic52} \\
2 & $\Sp_4(9)$ & $3^{2+4}{:}(\SL_2(5)\times8)$ & $\SL_2(81)$ & $3^3{:}4$ & \ref{lem:K<P1-Sp(4,q)} \\
3 & $\Sp_4(9){:}2$ & $3^{2+4}{:}\SL_2(5){:}2$ & $\SiL_2(81)$ & $3^3{:}2$ & \ref{ex:S5<S6<Sp(4,9).2} \\
4 & $\Sp_4(11)$ & $11^{1+2}{:}\SL_2(5)$ & $\SL_2(11^2)$ & $11$ & \ref{lem:K<P1-Sp(4,q)} \\
5 & $\Sp_4(19)$ & $19^{1+2}{:}(\SL_2(5)\times9)$ & $\SL_2(19^2)$ & $19{:}3$ & \ref{lem:K<P1-Sp(4,q)} \\
6 & $\Sp_4(29)$ & $29^{1+2}{:}(\SL_2(5)\times7)$ & $\SL_2(29^2)$ & $29$ & \ref{lem:K<P1-Sp(4,q)} \\
7 &$ \Sp_4(59)$ & $59^{1+2}{:}(\SL_2(5)\times29)$ & $\SL_2(59^2)$ & $59$ & \ref{lem:K<P1-Sp(4,q)} \\
\hline
8 & $\Sp_6(2)$ & $2^4{:}\A_5$ & $\PGaL_2(8)$ & $1$ & \ref{LemSymplectic52} \\
9 & $\Sp_6(2)$ & $2^4{:}\A_5$, $\Sy_5\times2$, $\Sy_6$, $\A_6\times2$ & $\SU_3(3)$ & $4$, $1$, $3$, $3$ & \ref{LemSymplectic52} \\
 & & $\Sy_7$, $\Omega_6^-(2)$, $\mathrm{O}_6^+(2)$ & & $7{:}3$, $3^{1+2}{:}4$, $\SL_3(2)$ & \ref{LemSymplectic52} \\
10 & $\Sp_6(2)$ & $\Sy_5$, $\A_6$, $\A_7$ & $\G_2(2)$ & $1$, $3$, $7{:}3$ & \ref{LemSymplectic52} \\
11 & $\Sp_6(3)$ & $\SL_2(13)$ & $3^{1+4}{:}\Sp_4(3)$ & $3$ & \ref{LemSymplectic52} \\
12 & $\Sp_6(3)$ & $3^{1+4}{:}2^{1+4}{}^{\boldsymbol{\cdot}}\A_5$ & $\SiL_2(27)$ & $3$ & \ref{LemSymplectic52} \\
13 & $\Sp_6(4)$ & $\J_2$ & $\Omega_6^-(4)$ & $5^2{:}\Sy_3$ & \ref{LemSymplectic10} \\
14 & $\GaSp_6(4)$ & $\SU_3(3)\times2$, $\G_2(2)$ & $\GaO_6^-(4)$ & $\Sy_3$ & \ref{LemSymplectic10} \\
15 & $\GaSp_6(4)$ & $\SiL_2(16)$ & $\GaG_2(4)$ & $1$ & \ref{LemSymplectic10} \\
16 & $\GaSp_6(16)$ & $\G_2(4){:}4$ & $\GaO_6^-(16)$ & $\SL_2(4)\times2$ & \ref{LemSymplectic32} \\
\hline
17 & $\Sp_8(2)$ & $\PSL_2(17)$ & $\mathrm{O}_8^+(2)$ & $\D_{18}$ & \ref{LemSymplectic52} \\
18 & $\Sp_8(2)$ & $\PGaL_2(9)$, $\Sy_{10}$ & $\Omega_8^-(2)$ & $\Sy_3$, $(\A_7\times3){:}2$ & \ref{LemSymplectic52} \\
 && $2^{10}{:}\A_6$, $2^{10}{:}\A_7$ && $2^5{}^{\boldsymbol{\cdot}}(\Sy_4\times2)$, $2^6{}^{\boldsymbol{\cdot}}\SL_3(2)$ &\ref{LemSymplectic52}\\
19 & $\Sp_8(2)$ & $\Sy_5$, $\Sy_5\times2$, $\A_5{:}4$, $\A_6$ & $\mathrm{O}_8^-(2)$ & $1$, $2$, $2$, $3$ & \ref{LemSymplectic52} \\
 & & $\A_6\times2$, $\Sy_6$, $\M_{10}$, $\PGL_2(9)$ & & $\Sy_3$, $\Sy_3$, $\Sy_3$, $\Sy_3$ & \ref{LemSymplectic52}\\
 & & $\A_7$, $\A_8$, $\A_8$ & & $7{:}3$, $\PSL_2(7)$, $\AGaL_1(8)$ & \ref{LemSymplectic52}\\
 & & $\A_9$, $\A_{10}$ & & $\PGaL_2(8)$, $(\A_7\times3){:}2$ & \ref{LemSymplectic52}\\
 & & $\Sp_6(2)$, $2^5{:}\A_6$, $2^6{:}\A_7$ & & $\G_2(2)$, $4^2{:}\Sy_3$, $2^3{}^{\boldsymbol{\cdot}}\SL_3(2)$ & \ref{LemSymplectic52}\\
20 & $\GaSp_8(4)$ & $\Omega_8^-(2){:}2$, $\Sp_6(4){:}2$, $\Omega_6^+(4){:}2$ & $\GaO_8^-(4)$ & $\G_2(2)$, $\G_2(4)$, $\SL_3(4)$ & \ref{prop:Sp(2)=O^+O^-} \\
 & & $\Omega_6^-(4).4$, $\Sp_4(4){:}2$ & & $\SU_3(4){:}2$, $\SL_2(4)$ & \ref{prop:Sp(2)=O^+O^-}\\
\hline
21 & $\Sp_{12}(2)$ & $\J_2{:}2$ & $\Omega_{12}^-(2)$ & $5^2{:}\Sy_3.2$ $\phantom{\overline{l}}$ & \ref{prop:SpaO-<Sp} \\
22 & $\Sp_{12}(2)$ & $\G_2(2)$, $3^{\boldsymbol{\cdot}}\PSU_4(3)$, $3^{\boldsymbol{\cdot}}\M_{22}$ & $\mathrm{O}_{12}^-(2)$ & $\SL_2(2)$, $3^5{:}\A_5$, $\PSL_2(11)$ & \ref{prop:Sp(2)=O^+O^-} \\
 & & $\SU_3(3)\times2$ & & $\SL_2(2)$ & \ref{LemSymplectic51}\\
 & & $\J_2$ & & $5^2{:}\Sy_3.2$ & \ref{prop:SpaO-<Sp}\\
23 & $\Sp_{12}(2)$ & $\GaG_2(4)$ & $\Sp_{10}(2)$ & $\SL_2(4)$ & \ref{LemSymplectic36} \\
 & & & $\Omega_{11}(2)\times2$ & $\SL_2(4)\times2$ & \ref{LemSymplectic36}\\
24 & $\GaSp_{12}(4)$ & $\G_2(4){:}2$ & $\GaO_{12}^-(4)$ & $\SL_2(4)$ & \ref{LemSymplectic30} \\
 & & $\G_2(4){:}2^2$ & & $\SL_2(4)\times2$ & \ref{LemSymplectic51} \\
25 & $\GaSp_{12}(4)$ & $\GaG_2(16)$ & $\Sp_{10}(4){:}2$ & $\SL_2(16)$ & \ref{LemSymplectic36} \\
 & & & $\Omega_{11}(4).2^2$ & $\SL_2(16)\times2$ & \ref{LemSymplectic36} \\
 & & & $\Omega_{11}(4).4$ & $\SL_2(16){:}2$ & \ref{LemSymplectic36} \\
\hline
26 & $\Sp_{16}(2)$ & $\Omega_9(2)$, $\Omega_8^-(2){:}2$, $\GaSp_6(4)$ & $\mathrm{O}_{16}^-(2)$ & $\Omega_7(2)$, $\G_2(2)$, $\G_2(4)$ & \ref{prop:Sp(2)=O^+O^-} \\
 & & $\Omega_6^+(4){:}2$, $\Omega_6^-(4){:}2$, $\GaSp_4(4)$ & & $\SL_3(4)$, $\SU_3(4)$, $\SL_2(4)$ & \ref{prop:Sp(2)=O^+O^-}\\
27 & $\GaSp_{16}(4)$ & $\Omega_9(4){:}2$ & $\GaO_{16}^-(4)$ & $\Omega_7(4)$ & \ref{LemSymplectic30} \\
\hline
28 & $\Sp_{24}(2)$ & $3^{\boldsymbol{\cdot}}\Suz$, $\Co_1$ & $\mathrm{O}_{24}^-(2)$ & $3^5{:}\PSL_2(11)$, $\Co_3$ & \ref{prop:Sp(2)=O^+O^-} \\
 & & $\G_2(4){:}2$, $\G_2(4).4$ & & $\SL_2(4)$, $\SL_2(4){:}2$ & \ref{prop:Sp(2)=O^+O^-} \\
 & & $\G_2(4){:}2^2$ & & $\SL_2(4)\times2$ & \ref{LemSymplectic47} \\
\hline
29 & $\Sp_{32}(2)$ & $\GaSp_8(4)$ & $\mathrm{O}_{32}^-(2)$ & $\Sp_6(4)$ & \ref{prop:Sp(2)=O^+O^-} \\
\hline
\end{tabular}
\vspace{3mm}
\end{table}

%%%%%%%%%%%%%%%%%%%%%%%%%%%%%%%%%%%%%%%%%%%%%%%%%%%%%%%%%%%%%%%%%%%
%%%%%%%%%%%%%%%%%%%%%%%%%%%%%%%%%%%%%%%%%%%%%%%%%%%%%%%%%%%%%%%%%%%

\section*{Acknowledgments}

Guralnick acknowledges support from Simons Foundation Fellowship 00019819.
Li and Xia acknowledge support of NNSFC 11931005.
Wang acknowledges support of NNSFC~12571022 and~12061083.

\end{document}